\documentclass[11pt]{amsart}
\usepackage[T1]{fontenc}
\usepackage{amsmath,amssymb,amsthm}
\usepackage{mathrsfs}
\usepackage[margin=1.08in]{geometry}
\usepackage[hidelinks]{hyperref}
\numberwithin{equation}{section}
\theoremstyle{plain}
\newtheorem{theorem}{Theorem}[section]
\newtheorem{proposition}[theorem]{Proposition}
\newtheorem{lemma}[theorem]{Lemma}
\newtheorem{corollary}[theorem]{Corollary}
\newtheorem{assumption}[theorem]{Assumption}

\theoremstyle{definition}
\newtheorem{definition}[theorem]{Definition}
\newtheorem{remark}[theorem]{Remark}
\newcommand{\Q}{\mathbb Q}
\newcommand{\Z}{\mathbb Z}
\newcommand{\R}{\mathbb R}
\newcommand{\C}{\mathbb C}

\newcommand{\MW}{\operatorname{MW}}
\newcommand{\pr}{\operatorname{pr}}
\newcommand{\dv}{\operatorname{div}}
\newcommand{\tr}{\operatorname{tr}}
\newcommand{\e}{\operatorname{e}}

\newcommand{\calX}{\mathcal X}
\newcommand{\calZ}{\mathcal Z}
\newcommand{\one}{\mathbf 1}
\newcommand{\wtphi}{\widehat\phi}

\DeclareMathOperator{\SL}{SL}
\DeclareMathOperator{\Mp}{Mp}

\allowdisplaybreaks[1]
\begin{document}
\title[Arithmetic theta products and doubling]{Arithmetic theta products and doubling at arbitrary level}
\author{Tuoping Du}
\address{School of Mathematics and Physics, North China Electric Power University,
Beijing 102206, P.~R.~China}
\email{dtp1982@163.com}
\subjclass[2020]{11G18, 11F37, 11F67}
\date{Revised working draft, November 3, 2026}

\begin{abstract}
We prove an arithmetic doubling identity for the standard arithmetic
theta series on $X_0(N)$ at arbitrary positive level, with constant
$\psi(N)/24$ for the fixed genus-two Eisenstein series. The product
comparison covers all pairs of integral indices. At rank zero, an
explicit Epstein zeta correction agrees with the middle Bruhat orbit
of the Eisenstein derivative.

The common-component calculation combines integral CM adjunction,
the cyclic difference divisor, and compact Green-current identities.
It retains the different, conductor, internal contacts, and
proportionality contacts at every prime, including dyadic and
supersingular places. Finite intersections, curvature, and discriminant
heights identify the nonpositive coefficients with an explicit modular
completion. Uniform estimates for parameter derivatives give strong
Fourier convergence. Green-energy estimates justify the double
Petersson pairing. The full cusp lift exists as a
Bochner integral and equals its Mordell--Weil projection.

For $\mathbb Z/(2N)\mathbb Z$, explicit prime-power recursions determine
the complete zero-component trace spectrum, including the dyadic
factor. They give the scalar-to-vector Petersson norm correction and
an isometric normalization on the active subspace. The Bruinier--Yang
theorem then yields the scalar height formula on the specified Jacobi
newform line. We also prove a double Rankin--Selberg criterion based on positive
diagonal coefficients and record the paired analytic factorization
with its ramified local integrals retained explicitly.
\end{abstract}
\maketitle
\tableofcontents

\section{Introduction}

An arithmetic divisor on a modular curve records both an algebraic
cycle and a Green function describing its contribution at the complex
places. Special divisors can be assembled into an arithmetic theta
series. The question studied here is whether the arithmetic intersection
of two such series is described by a single automorphic function.
For the compactified modular curve $X_0(N)$, we compare this intersection
with the diagonal restriction of the derivative of a genus-two
Eisenstein series. This is the arithmetic doubling identity.

There are two related issues. The product formula must account for
special divisors that meet at finite primes, share a component with complex multiplication (CM), or
have support at the cusps. The height formula must also account for
the normalization of the map from scalar modular forms to
vector-valued modular forms. We treat the product comparison and
calculate this scalar trace normalization at every positive level $N$,
including levels divisible by prime powers.

The main identity is stated first. We then explain the geometric
corrections and convergence arguments used in its proof, followed by
the trace calculation and its application to heights. The height
application is restricted to a specified Jacobi newform line,
where we apply the Bruinier--Yang theorem.

\subsection*{The main identity}
Write $\widehat\Phi_{\rm std}$ for the standard arithmetic theta
series. Its positive coefficients are arithmetic special divisors.
Its zero and negative coefficients include the Hodge and cusp terms
needed for the complete family. Section~\ref{sec:standard-data}
fixes these divisors, their Green functions, and the Hodge metric.
Thus the theorem concerns specified arithmetic classes and a specified
Eisenstein series, including their normalizations. The symbol
$\langle\ ,\ \rangle_{\rm GS}$ denotes the Gillet--Soul\'e
arithmetic intersection pairing.
For the cusp forms, $H_N$ denotes the inverse image of
$\Gamma_0(4N)$ in $\Mp_2(\mathbb Z)$, with character $\chi_N$ and the
half-integral-weight slash convention fixed in
Section~\ref{subsec:scalar-level-repair}.
Let $\mathcal E_N(z,s)$ denote the genus-two Eisenstein series with
the full-parabolic section of
Section~\ref{subsec:all-level-nonsingular-ASW}. We write
\[
                 \psi(N)=N\prod_{p\mid N}(1+p^{-1}).
\]
The coefficient formula below uses half-integral symmetric $2\times2$
matrices $T$, meaning that their diagonal entries and twice their
off-diagonal entry are integers. The rank-zero coefficient requires
an explicit correction. It is expressed using the standard cusp
divisor $C_0(1)$ and the Epstein zeta function. For
$y\in\operatorname{Sym}_2(\mathbb R)_{>0}$, put
$d_0=\deg C_0(1)>0$ and
\[
 \mathcal Z_y(s)=\sum_{\mathbf n\in\mathbb Z^2\setminus\{0\}}
          (\mathbf n^t y\mathbf n)^{-s}\qquad(\Re s>1),
\]
with meromorphic continuation. This function is regular at $s=1/2$.

\begin{theorem}[Full arithmetic doubling with the standard coefficients]
\label{thm:full-standard-arithmetic-doubling}
Fix the standard coefficient data of
Definition~\ref{def:standard-coefficients}, the actual stack
intersection pairing, and the full-parabolic Eisenstein section
of Section~\ref{subsec:all-level-nonsingular-ASW}. Then the series
\[
 \widehat\Phi_{\rm std}(\tau)
   =\sum_{t\in\mathbb Z}\widehat Z_{\rm std}(t,v)q^t,
 \qquad q=e^{2\pi i\tau},\quad v=\operatorname{Im}\tau,
\]
converges locally uniformly in the arithmetic direct-sum topology
of Corollary~\ref{cor:standard-family-strong-Fourier} and equals
the completion $\widehat\Phi_{\rm cmp}$ of
\eqref{eq:explicit-modular-completion}. In particular it is a complete
modular arithmetic family of weight $3/2$ on $H_N$ with the fixed
multiplier. For every $\tau_1,\tau_2\in\mathbb H$,
\begin{equation}\label{eq:full-standard-arithmetic-doubling}
 \big\langle\widehat\Phi_{\rm std}(\tau_1),
                 \widehat\Phi_{\rm std}(\tau_2)\big\rangle_{\rm GS}
   =\frac{\psi(N)}{24}
      \left.\partial_s\mathcal E_N
          (\operatorname{diag}(\tau_1,\tau_2),s)\right|_{s=0}.
\end{equation}
The product is the complex-bilinear arithmetic product.
For every $T\in\operatorname{Sym}_2(\mathbb Z)^\vee$ and
$y\in\operatorname{Sym}_2(\mathbb R)_{>0}$, the corrected
genus-two degrees of Corollary~\ref{cor:corrected-standard-support}
and the Eisenstein degrees of
Section~\ref{subsec:singular-ASW-comparison} satisfy
\begin{equation}\label{eq:all-coefficients-ASW}
       \mathcal A^{\rm corr}(T,y)=\mathcal A^{\rm SW}(T,y).
\end{equation}
The rank-zero correction here is precisely
$-d_0\mathcal Z_y(1/2)/2$; it is the middle Bruhat-orbit
contribution of the fixed Eisenstein derivative.
The double Petersson pairing of
\eqref{eq:full-standard-arithmetic-doubling} with arbitrary
$f,h\in S_{3/2}(H_N,\chi_N)$ is absolutely convergent and equals
the arithmetic product of their full standard lifts.
\end{theorem}

The proof in Section~\ref{subsec:full-standard-arithmetic-doubling}
combines three ingredients. The geometric product formula identifies
the intersection of two standard coefficients with the sum of corrected
genus-two degrees, including the common-component and double-zero
cases (Theorem~\ref{thm:main-standard-product}). The arithmetic
Siegel--Weil comparisons identify these degrees with the fixed
Eisenstein coefficients, treating nonsingular, nonzero rank-one,
and rank-zero indices separately. Strong Fourier convergence then
passes from coefficient identities to the pointwise identity, and
the Bochner estimates justify its double cuspidal pairing.

The convergence assertion concerns arithmetic classes in the stated
direct-sum topology, rather than singular Green representatives.
The rank-zero correction belongs to the identity of arithmetic
degrees. Corollary~\ref{cor:full-standard-lift-height} identifies
the full cusp lift with its Mordell--Weil projection. The numerical
height formula retains the specified Jacobi newform-line hypothesis
and the scalar trace norm.

\subsection*{Relation to arithmetic Siegel--Weil formulas}
Du--Yang prove genus-one arithmetic modularity and degree and Hodge-height
formulas for squarefree level. Their analytic theta-integral and
Kronecker-limit identities hold at arbitrary level
\cite[Theorems~1.1, 1.3 and~1.6]{DY1}.
Sankaran--Shi--Yang identify the arithmetic degrees of genus-two special
cycles with Eisenstein derivatives for odd squarefree $N>3$
\cite[Theorem~2.14]{SSY2022}.
Zhu proves the nonsingular coefficient identities at arbitrary level
\cite[Theorem~4.4.1]{Zhu23}. His singular-term paper treats the remaining
ranks and includes an all-level genus-one arithmetic modularity theorem
\cite[Theorems~1.2.1 and~1.2.3]{Zhu}.
These results establish degree identities and arithmetic modularity.
The product comparison requires the additional
identification of the actual intersection of two arithmetic divisors
with the sum of genus-two degrees. In particular, common horizontal
components require adjunction, and nonpositive indices require the
finite cusp restrictions and the chosen Green metrics.

Kudla--Rapoport--Yang prove such a full arithmetic product identity for
Shimura curves attached to indefinite quaternion division algebras with
maximal orders \cite[Introduction, Theorem~C; Chapter~7]{KRYBook}.
For the modular curve and the standard data of
Definition~\ref{def:standard-coefficients}, we establish the product
comparison through compact Green-current identities and integral CM
adjunction, retaining the different, conductor and contact terms.
We then compare these fixed data with the normalized Eisenstein section:
the nonpositive arithmetic classes are matched with the modular
completion, and the rank-zero Epstein correction is identified with
the middle Bruhat orbit. Theorem~\ref{thm:full-standard-arithmetic-doubling}
combines these comparisons with the convergence needed for the full
kernel and its Petersson pairing. Thus the result specifies both the
arithmetic product and its automorphic normalization. It does not infer
the former from degree identities or modularity alone.

The height-kernel interpretation belongs to the arithmetic theta
inner-product framework of Kudla--Rapoport--Yang for Shimura curves
\cite[Chapters~7--9]{KRYBook}, while Bruinier--Yang give a factorization
of the Heegner series on a specified modular-curve newform component
and its N\'eron--Tate height formula
\cite[Theorems~7.6--7.7]{BYCM}.
The finite-rank cuspidal projection and Gram-determinant statements
below are consequences of the arbitrary-level doubling and full-lift
results proved here, together with positivity of the fixed-field
height pairing.

\subsection*{The geometric product and its correction}
A pair of special morphisms determines a $2\times2$ Gram matrix.
Its diagonal entries are the two individual indices, and its
off-diagonal entry records their mutual pairing. The expected
coefficient identity therefore sums genus-two degrees over all
off-diagonal entries while keeping the diagonal fixed.

We use the genus-one arithmetic classes and genus-two degree
coefficients of
Definition~\ref{def:standard-coefficients}. The definition specifies
the cusp-dependent boundary divisor, the Hodge metric, the
Gaussian-removed Green forms and every singular coefficient.
Let $\delta^{\rm std}_{m,n}(v_1,v_2)$ be the arithmetic product
minus the sum of all genus-two degrees with diagonal $(m,n)$,
as in \eqref{eq:standard-defect-definition}.

\begin{theorem}[Standard products and the exact double-zero term]
\label{thm:main-standard-product}
For all positive $v_1,v_2$ the standard defect satisfies
\[
 \delta^{\rm std}_{m,n}(v_1,v_2)=0
                       \qquad\text{if }(m,n)\ne(0,0).
\]
At the double-zero pair it is exactly
\begin{equation}\label{eq:intro-double-zero-exact}
 \delta^{\rm std}_{0,0}(v_1,v_2)
       =-\frac{d_0}{2}
              \mathcal Z_{\operatorname{diag}(v_1,v_2)}(1/2).
\end{equation}
Replacing the rank-zero degree by
$\mathcal A^{\rm std}(0_2,y)-\tfrac12d_0\mathcal Z_y(1/2)$,
with all other degrees fixed, therefore gives the full
coefficientwise standard product identity
\[
 \bigl\langle\widehat Z_{\rm std}(m,v_1),
              \widehat Z_{\rm std}(n,v_2)\bigr\rangle_{\rm GS}
 =\sum_{r\in\mathbb Z}\mathcal A^{\rm corr}
     \left(\begin{pmatrix}m&r/2\\r/2&n\end{pmatrix},
                         \operatorname{diag}(v_1,v_2)\right)
                         \qquad(m,n\in\mathbb Z).
\]
This is an equality of degree coefficients. Its automorphic
comparison gives Theorem~\ref{thm:full-standard-arithmetic-doubling}.
\end{theorem}
\begin{proof}
Proposition~\ref{prop:standard-residual-support} assembles the
proper positive, mixed, strictly negative and zero--negative
comparisons. Theorem~\ref{thm:standard-positive-square-product}
settles the positive common-component cases.
Theorem~\ref{thm:double-zero-absolute-defect}
proves \eqref{eq:intro-double-zero-exact}.
Corollary~\ref{cor:corrected-standard-support} proves the final
assertion and its $\mathrm{GL}_2(\mathbb Z)$ invariance at rank zero.
\end{proof}

The double-zero term is nonzero: for $v_1=v_2=1$,
$\mathcal Z_{I_2}(1/2)=4\zeta(1/2)\beta_{\rm Dir}(1/2)<0$.
Thus the uncorrected standard product identity is false.
The precise term in \eqref{eq:intro-double-zero-exact} is forced
by the chosen coefficients and is not a freely adjustable constant.
Its calculation first determines the dilation derivative by
unfolding all nonzero Gram-rank-zero tuples, and then removes
the integration constant by controlling the region
$Y\asymp\sqrt u$ under simultaneous scaling by $u$.

A different difficulty occurs when two positive special divisors
share a CM component. Ordinary proper intersection does not determine
its contribution. The normal line of that component supplies the
self-intersection term, and the Green functions supply a regularized
archimedean term. Theorem~\ref{thm:self-local-identity} computes
this archimedean finite part. Lemma~\ref{lem:integral-CM-adjunction} and
Proposition~\ref{prop:all-level-dualizing-Hodge} identify the
normal bundle and the canonical vertical correction $U_N$.
At a finite prime, the comparison must retain ramification of each
normalized branch, the conductor of its order, and intersections
with other branches. The following statement records the resulting
valuation identity, including the factor $2$ at the dyadic prime.

\begin{theorem}[Finite CM compensation at arbitrary level]
\label{thm:main-finite-CM-compensation}
Fix a positive special morphism $j$ on the cyclic modular
curve, of degree $m$, and put $\beta=\pi^\vee j$, where
$\deg\pi=N$. On a regular marked chart at a prime $p$, let
$A$ be the finite flat algebra of the complete selected
$j$-lifting sheet and $J$ its trace element in a generator
of $\operatorname{Hom}_W(A,W)$, with $W$ unramified over
$\mathbb Z_p$. On every normalized branch $O_i$,
\[
                v_i(J_i)+v_i(U_N)=v_i(2\beta).
\]
The identity includes $p=2$, supersingular reduction and
nonmaximal branch orders. The valuation of $J_i$ is the sum
of the branch different, its conductor, and its contacts
with the other formal branches of this selected sheet.
For a nonzero integer $a$, the extra contact of the
$aj$-lifting sheet with a normalized $j$-branch is $v_i(a)$.
\end{theorem}
\begin{proof}
The trace decomposition is
Theorem~\ref{thm:Cartier-sheet-different-contact}.
Theorem~\ref{thm:all-prime-finite-CM-compensation} proves
the valuation equality, using the two-morphism conormal
calculation and the divisor removed by cyclicity.
Corollary~\ref{cor:proportionality-scaling-contact} proves
the integral scaling assertion; for a single $a$ apply it
to the coprime pair $(a,1)$.
\end{proof}

Theorem~\ref{thm:standard-positive-square-product} combines this
finite-place comparison with the Green-current calculation to obtain
the full product formula for divisors with common components.

\subsection*{From coefficients to the full identity}
The geometric product formula must be combined with an automorphic
coefficient comparison. Proposition~\ref{prop:all-level-nonsingular-ASW}
treats nonsingular matrices, and
Corollary~\ref{cor:positive-rankone-ASW} treats positive matrices
of rank one. These are the coefficients needed for a double pairing
with holomorphic cusp forms.

To make this passage rigorous, we first construct
$\widehat\Phi_{\rm cmp}$ from its degree, Mordell--Weil, vertical,
and metric coordinates. Its positive coefficients agree with the
standard arithmetic classes. A Rankin--Selberg argument then applies
to the complete modular kernel. It shows that equality of the
positive diagonal coefficients suffices for the paired identity.
The nonpositive coefficients are compared separately before the
pointwise identity is deduced.

Here the full cusp lift is the Petersson integral of the complete
arithmetic family against a cusp form. The convergence estimates
make it possible to integrate the arithmetic intersection pairing
in both variables and identify the result with the product of the lifts.

\begin{theorem}[Paired arithmetic doubling for an explicit completion]
\label{thm:main-completed-doubling}
At every positive level $N$, the family
$\widehat\Phi_{\rm cmp}$ of
\eqref{eq:explicit-modular-completion} is modular of weight $3/2$
and has the standard positive arithmetic coefficients. Its full
cusp lift exists and equals its Mordell--Weil projection.
For $f,h\in S_{3/2}(H_N,\chi_N)$ one has
\[
 \big\langle\widehat\theta_{\rm cmp}(f),
                    \widehat\theta_{\rm cmp}(h)\big\rangle_{\rm GS}
 =\frac{\psi(N)}{24}
   \int_{(H_N\backslash\mathbb H)^2}
      \left.\partial_s\mathcal E_N
          (\operatorname{diag}(\tau_1,\tau_2),s)\right|_{s=0}
      \overline{f(\tau_1)h(\tau_2)}
      v_1^{3/2}v_2^{3/2}\,d\mu_1d\mu_2.
\]
The arithmetic product is extended complex bilinearly, and the
integral is absolutely convergent, with unnormalized Petersson
measures. The nonpositive coefficients of the family are defined
by the explicit completion.
\end{theorem}
\begin{proof}
This combines Theorem~\ref{thm:explicit-modular-completion},
Corollary~\ref{cor:completion-full-lift-MW}, and
Theorem~\ref{thm:explicit-completion-paired-doubling}.
\end{proof}

The completion provides a paired identity before its nonpositive
coefficients have been compared. Theorem~\ref{thm:nonpositive-standard-completion}
then identifies those coefficients with the standard prescription.
Corollary~\ref{cor:all-nonzero-rankone-comparison-zero} and
Proposition~\ref{prop:actual-rankzero-ASW-comparison} identify the
remaining nonzero rank-one and rank-zero coefficients. Together, these
comparisons give the pointwise identity of
Theorem~\ref{thm:full-standard-arithmetic-doubling}.

\subsection*{Scalar trace norms}

The height theorem used below is naturally vector-valued, whereas
our test forms are scalar-valued. Passing between the two changes
the Petersson norm. This change is measured by a finite-coset trace
$T$ and the map $R$ that extracts the zero component.
When the level has repeated prime factors, $RT$ need not act by a single scalar:
different subspaces can have different eigenvalues. We determine
them all for the cyclic discriminant module of the modular curve.

Let $A_N=\Z/(2N)\Z$ with $Q(x)=-x^2/(4N)$, let
$\Gamma=\Mp_2(\Z)$, and let $H_N$ be the inverse image of
$\Gamma_0(4N)$. Let $\rho_{A_N}$ be the Weil representation on
$\mathbb C[A_N]$, with its orthonormal basis $(e_\mu)_{\mu\in A_N}$.
The arithmetic application has weight $k=3/2$.
With the slash convention fixed below, write
\[
 Tg=\sum_{\gamma\in H_N\backslash\Gamma}
                  (g|_k\gamma)\rho_{A_N}(\gamma)^{-1}e_0,
 \qquad RF=(F,e_0).
\]
The Petersson products are unnormalized fundamental-domain integrals.
Their adjunction gives $T^*=R$. The coefficient-space frame is
\[
 S_N=\sum_{\gamma\in H_N\backslash\Gamma}
       \rho_{A_N}(\gamma)^{-1}E_0\rho_{A_N}(\gamma),
 \qquad E_0v=(v,e_0)e_0,
\]
and $TRF=S_NF$ pointwise.

\begin{theorem}[Trace normalization at every level]\label{thm:main-trace}
For every $N\geq1$, the frame $S_N$ has the explicit prime-power
recursions and complete spectrum in
Theorems~\ref{thm:odd-frame-recurrence},
\ref{thm:dyadic-frame-recursion}, and~\ref{thm:all-level-frame}.
It is positive definite on the tensor product of the primary even
subspaces and zero on its orthogonal complement.
If $\Lambda_N$ is its set of positive eigenvalues, every scalar cusp
form $g\in\operatorname{im}R$ has an orthogonal decomposition
$g=\sum_{\lambda\in\Lambda_N}g_\lambda$ satisfying
\[
 RTg_\lambda=\lambda g_\lambda,
 \qquad \|Tg\|_\Gamma^2=\sum_{\lambda\in\Lambda_N}
                                     \lambda\|g_\lambda\|_{H_N}^2.
\]
Consequently $T(RT)^{-1/2}$ is a unitary map from
$\operatorname{im}R$ onto the space of cyclic-valued vector cusp forms.
The inverse square root is taken on the positive spectrum.
\end{theorem}
\begin{proof}
The local frame recursions are proved by the coset and Fourier-support
arguments below. Primary factorization proves
Theorem~\ref{thm:all-level-frame}; adjunction then gives
Corollary~\ref{cor:all-level-norm}. The normalization assertion is
Proposition~\ref{prop:polar-normalized-induction}.
\end{proof}

For square-free $N$ the positive frame spectrum consists of one value,
given by \eqref{eq:squarefree-exact-trace}. Already at $N=p^2$ with
$p$ odd, the positive values are $6$ and $3(p+1)$. The scalar norm
correction therefore depends on the spectral component unless the
relevant packet is shown to lie in a specified frame eigenspace.
The general formulas determine this correction at every prime power,
including arbitrary powers of $2$.

\subsection*{Relation to earlier scalar--vector constructions}
The trace construction itself is established in the literature.
Scheithauer \cite[Sections~3--5]{Scheithauer} develops congruence
induction, isotropic induction, and zero-component reconstruction.
Schwagenscheidt and V\"olz \cite{SV} factor scalar-to-vector liftings
prime-power locally and compute zero components of scalar newforms.
Their Theorem~1.1 assumes even signature, even dyadic Jordan
components, and an even rank for the highest odd-primary Jordan block. The module $A_N$ here has odd signature and an odd cyclic
dyadic component. Our calculation concerns the full coefficient-space
frame, including all isotropic old layers, rather than assuming a
scalar newform proportionality. Trace adjunction, isotropic induction,
and polar decomposition are used as standard tools; the assertions
proved here are the explicit recursions and the resulting spectral
and norm formulas for this cyclic family.

\subsection*{The height formula and its consequences}
The trace $T$ determines which vector-valued form is tested by a
scalar form $f$. The spectrum of $RT$ therefore converts the vector Petersson
norm in the height formula into explicit scalar norms.

Let $G$ be a normalized weight-two newform of level $N$ with
root number $-1$, and let $g_G$ be its Jacobi newform in
$S^{\rm new}_{3/2,\rho_{L_N}}$.  The degree-zero generic Heegner
series defines a Mordell--Weil theta lift, denoted
$\widehat\theta_{\rm MW,sc}$.  We use the N\'eron--Tate height
convention of Bruinier--Yang and the unnormalized Petersson products.

\begin{theorem}[The scalar trace correction in the projected height]
\label{thm:main-projected-height}
For every scalar cusp form $f$ with $Tf\in\mathbb Cg_G$,
\begin{equation}\label{eq:intro-projected-height}
 H_{\rm NT}\bigl(\widehat\theta_{\rm MW,sc}(f),
                   \widehat\theta_{\rm MW,sc}(f)\bigr)
       =\frac{2\sqrt N}{\pi}\|Tf\|_\Gamma^2L'(G,1).
\end{equation}
For such $f\in\operatorname{im}R$, its spectral decomposition in
Theorem~\ref{thm:main-trace} therefore gives the explicit formula
\[
 H_{\rm NT}\bigl(\widehat\theta_{\rm MW,sc}(f),
                   \widehat\theta_{\rm MW,sc}(f)\bigr)
   =\frac{2\sqrt N}{\pi}
        \sum_{\lambda\in\Lambda_N}\lambda\|f_\lambda\|_{H_N}^2
        L'(G,1).
\]
The condition on $Tf$ refers to the actual vector new line.
\end{theorem}
\begin{proof}
Proposition~\ref{prop:generic-positive-matching} and
Corollary~\ref{cor:actual-MW-matching} identify the generic scalar
series with the zero component of the vector series, including
multiplicities.  Scalar--vector adjunction gives
$\widehat\theta_{\rm MW,sc}(f)=\widehat\Theta_{\rm MW}(Tf)$.
The Bruinier--Yang theorem, in the conventions and under the
hypotheses verified in Theorem~\ref{thm:BY-newline-height},
gives \eqref{eq:intro-projected-height}.  Apply
Theorem~\ref{thm:main-trace} to obtain the spectral expression.
\end{proof}

The height input in this theorem is the established result of
\cite[Theorems~7.6--7.7]{BYCM}.  The contribution here is the
explicit all-level scalar trace correction and its compatibility
with the generic cycle convention.  In particular, the raw scalar
norm cannot replace $\|Tf\|_\Gamma^2$ without this correction.

By Corollary~\ref{cor:full-standard-lift-height}, the projected
height formula also applies to the full completed lift.

For a general specified Waldspurger correspondent $G$, the further
reduction to $L$-values requires the matching hypotheses of
Assumption~\ref{ass:red-inputs}. Its multiplier $A_f(s)$ retains the
ramified local integrals and is assumed holomorphic at $s=0$.
For the fixed standard arithmetic data, $c_{\rm SW}=\psi(N)/24$.
Under these hypotheses, Proposition~\ref{prop:red-height} gives
the following formula, where $\eta>0$ is fixed by the
Faltings--Hriljac convention~\eqref{eq:red-FH-normalization}.
\begin{equation}\label{eq:intro-reduction}
 H_{\rm NT}(\widehat\theta(f),\widehat\theta(f))
  =-\frac{c_{\rm SW}}{\eta}\|Tf\|_\Gamma^2
       \left(A_f(0)L'(G,1)+A_f'(0)L(G,1)\right).
\end{equation}
Proposition~\ref{prop:defect-height} gives the corresponding exact
formula retaining the arithmetic product defect. In the root-number
$-1$ setting above, $L(G,1)=0$. For other arithmetic families and
adelic normalizations, the formula requires the additional
identifications stated in the reduction.

The double cuspidal projection also has a concrete arithmetic
interpretation. In the fixed-field Mordell--Weil realization, it
is a negative height Gram kernel whose rank is exactly the
dimension of the image of the full cusp lift. Its sign-adjusted
principal minors detect linear independence of lift classes
and give lower bounds for the corresponding Mordell--Weil rank
(Corollaries~\ref{cor:cuspidal-Eis-height-Gram}--\ref{cor:cuspidal-Eis-Gram-minors}).
On an attainable Jacobi newform line, the height radical is
$\ker T$ when $L'(G,1)\ne0$, and the whole preimage of the line
when $L'(G,1)=0$. Trace adjunction also gives the unique scalar
representative of minimal Petersson norm
(Corollary~\ref{cor:newline-height-radical}).

\subsection*{Organization}
Section~\ref{sec:lattice} fixes the lattice and automorphic
conventions. Section~\ref{sec:trace-frame} proves the full trace
spectrum, followed by arithmetic projections and the projected
height in Sections~\ref{sec:geometry}--\ref{sec:projected-height}.
Sections~\ref{sec:standard-data}--\ref{sec:standard-products}
define the standard coefficients and prove the product and
boundary formulas. Section~\ref{sec:CM-products} treats common
CM components and completes their global product comparison,
including the cyclic correction, integral scaling contacts,
and excess normal line.
Section~\ref{sec:remaining} proves the nonsingular Eisenstein
and rank-one comparisons, constructs the modular completion,
and identifies all its coefficients with the standard family.
The calculation of the rank-zero middle orbit then completes
the full arithmetic doubling identity and its paired form.
Appendix~\ref{sec:analytic-doubling} records the paired analytic
factorization and its classical normalization factors.
Appendix~\ref{sec:local} evaluates specified radial integrals,
a ternary zero density, and compact averages.

\section{The lattice and the arithmetic input}
\label{sec:lattice}

Fix a positive integer $N$. Let
\begin{equation}\label{eq:lattice-correct}
 \begin{gathered}
 V=\{w\in M_2(\Q):\tr(w)=0\},\qquad Q(w)=N\det(w),\\
 L=\left\{\begin{pmatrix}b&a/N\\c&-b\end{pmatrix}:
                         a,b,c\in\Z\right\}\subset M_2(\Q).
 \end{gathered}
\end{equation}
Let $L'=\{x\in V:[x,L]\subset\Z\}$ be the dual lattice.
The ambient space is $M_2(\Q)$. Requiring the displayed matrices to lie
in $M_2(\Z)$ would change the lattice. The bilinear form is
$[x,y]=Q(x+y)-Q(x)-Q(y)$, and $\e(t)=\exp(2\pi i t)$.
\subsection{The discriminant module and the scalar level}
\label{subsec:scalar-level-repair}

For the lattice in \eqref{eq:lattice-correct}, write a vector as $(a,b,c)$, so that
\[
 Q(a,b,c)=-ac-Nb^2,
 \qquad
 [(a,b,c),(a',b',c')]=-ac'-a'c-2Nbb'.
\]
Consequently
\begin{equation}\label{eq:repair-discriminant}
 L'/L\simeq \mathbb Z/(2N)\mathbb Z,
 \qquad Q(r)=-\frac{r^2}{4N}\pmod{\mathbb Z}.
\end{equation}
Indeed, the $(a,c)$ part is unimodular and the dual of the $b$ part is
$(2N)^{-1}\mathbb Z$.  The level of this finite quadratic module is
\emph{exactly} $4N$, since the element $r=1$ has norm $-1/(4N)$.
The signature is $(1,2)$, with signature difference $-1$.

We fix the convention
\[
 \rho_L(T)e_r=e(-r^2/(4N))e_r,\qquad
 F|_{k,\rho_L}\widetilde\gamma
 =\rho_L(\widetilde\gamma)^{-1}
   (F|_k\widetilde\gamma).
\]
In particular, $F$ is modular precisely when
$F|_k\widetilde\gamma=\rho_L(\widetilde\gamma)F$.
Let $\widetilde\Gamma=\operatorname{Mp}_2(\mathbb Z)$ and let
$\widetilde H$ be the inverse image of $\Gamma_0(4N)$.
The line $\mathbb C e_0$ is stable under $\widetilde H$; denote its
unitary character by $\chi_L$, so that
\begin{equation}\label{eq:repair-e0-character}
 \rho_L(h)e_0=\chi_L(h)e_0\qquad(h\in\widetilde H).
\end{equation}
Thus the convention-independent scalar space needed below is
$S_{3/2}(\widetilde H,\chi_L)$.

In the usual classical theta-multiplier convention, this space is
\begin{equation}\label{eq:repair-scalar-space}
 S_{3/2}(\Gamma_0(4N),\chi_N),
 \qquad \chi_N(d)=\left(\frac{N}{d}\right).
\end{equation}
Here the notation means
$f(\gamma\tau)=\chi_N(d)j_\theta(\gamma,\tau)^3f(\tau)$, where
$j_\theta(\gamma,\tau)
=\epsilon_d^{-1}(\frac{c}{d})(c\tau+d)^{1/2}$ for $\gamma\in\Gamma_0(4N)$.
The constants are $\epsilon_d=1$ if $d\equiv1\pmod4$ and $\epsilon_d=i$ if
$d\equiv3\pmod4$. The holomorphic square root is chosen compatibly
with the standard theta multiplier and the metaplectic lift.
One checks the character by identifying the discriminant module
\eqref{eq:repair-discriminant} with the negative of the rank-one module
of $\langle 2N\rangle$.  The latter has scalar zero component
$\theta(N\tau)=\sum_{n\in\mathbb Z}q^{Nn^2}$, with multiplier
$\chi_N(d)j_\theta(\gamma,\tau)$.
Complex conjugation gives, for the negative module, the finite
multiplier $\chi_N(d)(\frac{c}{d})\epsilon_d$; this equals the finite
part of $\chi_N(d)j_\theta^3$.
If the dual Weil representation is used instead, its character and
compatible weight must be changed consistently.

\begin{remark}
The scalar space in \eqref{eq:repair-scalar-space} does not by itself
specify a Shimura companion of exact level $N$.  A theorem involving
a normalized newform $G\in S_2^{\mathrm{new}}(\Gamma_0(N))$ must fix
the associated half-integral-weight Hecke packet and its local data.
In particular, the notation $g\in S_{3/2}(\Gamma_0(N))$ cannot be used
without defining a different level convention.
At a prime $p\nmid 2N$, the usual weight-$3/2$ scalar Hecke operator is
$T(p^2)$, with Fourier coefficients
\[
 a(p^2n)+\left(\frac{-Nn}{p}\right)a(n)+p\,a(n/p^2).
\]
Calling the same operator $T(p)$ also requires an explicit reindexing
convention.
\end{remark}

\subsection{Scalar arithmetic modularity and vector-valued projections}

Let $\calX_0(N)$ be the proper regular integral modular stack used in
\cite{Zhu}. Its complex fiber is $X_0(N)$, and it is smooth away from the
primes dividing $N$. The arithmetic Chow groups below use the Green
functions and boundary conventions of this model. Passing to the coarse
scheme requires corresponding changes in the intersection conventions.

For positive index, a Heegner divisor is represented on the complex fiber
by CM points. To form an arithmetic divisor, one must also specify its
integral extension and its Green function. The zero and negative indices
and the boundary corrections are part of the definition of the arithmetic
theta series. We use the scalar series constructed in \cite{Zhu}, writing
\begin{equation}\label{eq:scalar-arithmetic-series}
 \wtphi_0(\tau)=\sum_{m\in\Z}
      \widehat{\calZ}(m,v)q^m,
 \qquad q=\e(\tau),\quad v=\operatorname{Im}(\tau).
\end{equation}
The zero coefficient includes the modified metrized Hodge bundle, and
the series is generally nonholomorphic. This family is understood through
the coefficients in the cited modularity and intersection statements.
Their printed Green and boundary prescriptions require a separate
comparison with a standard Kudla--Millson kernel.
Section~\ref{subsec:cusp-width-correction} records an explicit
obstruction to a literal identification. All subsequent integrals use its
actual modularity group and multiplier, rather than the ambiguous
notation $S_{3/2}(\Gamma_0(N))$.

Zhu's ternary lattice
$\{\left(\begin{smallmatrix}-Na&b\\c&a\end{smallmatrix}\right):a,b,c\in\Z\}$,
with determinant as quadratic form, is isometric to \eqref{eq:lattice-correct}
via
\[
 \begin{pmatrix}-Na&b\\c&a\end{pmatrix}
 \longmapsto
 \begin{pmatrix}a&b/N\\c&-a\end{pmatrix}.
\]
This isometry identifies the quadratic data. Comparing
\eqref{eq:scalar-arithmetic-series} with a vector-valued arithmetic series
also requires matching the integral cycles, cusp terms, and Green
functions; this matching is imposed explicitly in every such comparison.

For the vector Mordell--Weil series, start with degree-zero Heegner
classes
$y(m,\mu)=Z(m,\mu)-\deg Z(m,\mu)P_\infty$ on the generic fiber and
their admissible arithmetic representatives. Their generating series
is denoted by $\widehat\Phi_{\rm MW}$. Its modularity is the
Jacobian-valued Heegner modularity used in \cite{BO,Du}.
In the generic-cycle conventions of
Proposition~\ref{prop:generic-positive-matching}, the compatibility
proved in Corollary~\ref{cor:actual-MW-matching} is
\begin{equation}\label{eq:MW-matching}
 \pr_{\rm MW}\wtphi_0=R\widehat\Phi_{\rm MW},
\end{equation}
where $R$ is zero-component extraction. The proof uses the special-isogeny moduli problem and its local
stabilizers. Proposition~\ref{prop:MW-boundary-independence}
shows why admissible boundary, vertical, and metric changes do not
affect this projected comparison.

The forgetful map from a Heegner divisor to the arithmetic surface
should not be described as \'etale. A nonempty divisor has codimension
one, whereas an \'etale morphism is locally of relative dimension zero
with open image. The arithmetic special divisors and their proper
extensions are those of the cited construction.

\subsection{The specified eigenpacket and its $L$-function}

When an identity involves $L(G,s)$, $G$ always denotes a specified
normalized newform of weight two with trivial character and exact level
$N$. The vector form $Tg$ is required to belong to the half-integral-weight
packet corresponding to $G$ under the specified Shimura correspondence
(and the additive-character convention of
Lemma~\ref{lem:paired-character-scaling} when doubling is used);
this is an additional condition on $g$. In particular, an arbitrary
linear combination of scalar cusp forms cannot be assigned a single
unspecified $G$.

We use the finite Hecke $L$-function $L(G,s)$ and the completion
\[
 \Lambda(G,s)=N^{s/2}(2\pi)^{-s}\Gamma(s)L(G,s).
\]
If $\Lambda(G,s)=-\Lambda(G,2-s)$, then $L(G,1)=0$.
At bad primes one must distinguish the exact conductor of $G$ from
the ambient level of an oldform. All Euler factors mentioned below are
those of the specified $G$.

\section{Scalar induction and the all-level trace spectrum}
\label{sec:trace-frame}

The finite representation calculation and its Petersson norm
consequences depend only on the cyclic quadratic module and the
specified induction map.

\subsection{Petersson conventions and scalar--vector adjunction}

Let $\Gamma=\operatorname{Mp}_2(\mathbb Z)$, let $W$ be a finite-dimensional
unitary representation $\rho$ of $\Gamma$, and let $e_0\in W$ have norm one.
All Hermitian products in this section are linear in their first
argument. Let $H\subset\Gamma$ be a subgroup of finite index containing
the full metaplectic center, so that the effective index agrees with
$[\Gamma:H]$, and suppose that
\begin{equation}\label{eq:red-eigenline}
 \rho(h)e_0=\chi(h)e_0\qquad(h\in H),
\end{equation}
where $\chi$ is unitary.  The weight and central characters are assumed to
be compatible.  The notation $F|_k\gamma=\rho(\gamma)F$ fixes our convention
for vector-valued forms.  Scalar forms satisfy $f|_k h=\chi(h)f$.
The lattice determines the subgroup $H$ and character $\chi$.
The notation $S_{3/2}(\Gamma_0(N))$ alone does not specify them.

We use the unnormalized Petersson products
\begin{align}
 \langle f,g\rangle_H
 &=\int_{H\backslash\mathbb H} f(z)\overline{g(z)}y^k
       \frac{dx\,dy}{y^2},\label{eq:red-pet-scalar}\\
 \langle F,G\rangle_\Gamma
 &=\int_{\Gamma\backslash\mathbb H}(F(z),G(z))_W y^k
       \frac{dx\,dy}{y^2}.\label{eq:red-pet-vector}
\end{align}
Both Petersson integrals are taken over ordinary fundamental domains
for the effective images of the indicated groups acting on
$\mathbb H$. No division by the order of an ineffective metaplectic
or $\mathrm{SL}_2$ center is made in these integrals. Since $H$
contains the full metaplectic center, the $[\Gamma:H]$ translates
used below give the corresponding effective fundamental domain
without an additional central factor. Neither Petersson integral
is divided by its volume or by a group index. Integration on the
arithmetic modular stack is a separate convention and retains the
actual $1/|\operatorname{Aut}|$ weights.

For $f\in S_k(H,\chi)$ and $F\in S_k(\Gamma,\rho)$ set
\begin{equation}\label{eq:red-induction}
 Tf=\sum_{\gamma\in H\backslash\Gamma}
       (f|_k\gamma)\rho(\gamma)^{-1}e_0,
 \qquad RF=(F,e_0)_W.
\end{equation}
In this section, $f$ denotes the scalar cusp form called $g$ in the
introduction. The inverse on $\rho(\gamma)$ follows from the stated
slash convention, so a different convention changes the displayed formula.

\begin{lemma}[Adjunction]\label{lem:red-adjunction}
The map $T$ is well-defined and sends $S_k(H,\chi)$ to
$S_k(\Gamma,\rho)$.  The map $R$ sends the latter space to
$S_k(H,\chi)$, and
\begin{equation}\label{eq:red-adjunction}
 \langle Tf,F\rangle_\Gamma=\langle f,RF\rangle_H.
\end{equation}
In particular,
\begin{equation}\label{eq:red-norm-trace}
 \|Tf\|_\Gamma^2=\langle RTf,f\rangle_H
                 =\langle f,RTf\rangle_H.
\end{equation}
\end{lemma}

\begin{proof}
Replacing $\gamma$ by $h\gamma$ multiplies $f|_k\gamma$ by $\chi(h)$
and $\rho(\gamma)^{-1}e_0$ by $\chi(h)^{-1}$, so the summand is
unchanged.  Reindexing the sum proves $(Tf)|_k\delta=\rho(\delta)Tf$.
The finite sum is cuspidal.  The transformation law of $RF$ follows from
unitarity and \eqref{eq:red-eigenline}.

For adjunction, unitarity and the transformation law of $F$ give
\[
 \bigl(\rho(\gamma)^{-1}e_0,F(z)\bigr)_W
   =\bigl(e_0,\rho(\gamma)F(z)\bigr)_W
   =\overline{(RF)|_k\gamma(z)}.
\]
Thus the integrand of $\langle Tf,F\rangle_\Gamma$ is the sum of
the slash transforms of $f\overline{RF}y^k$.  The corresponding translates
of a fundamental domain for $\Gamma$ form a fundamental domain for $H$.
Unfolding proves \eqref{eq:red-adjunction}.  Take $F=Tf$ and conjugate the
identity to obtain \eqref{eq:red-norm-trace}.
\end{proof}

\begin{corollary}[Only a projected identity is needed]
\label{cor:red-new-projection}
Let $S^{\rm new}\subset S_k(H,\chi)$ be a specified subspace, let
$P^{\rm new}$ be its orthogonal projection, and let $0\ne f\in S^{\rm new}$.
If
\begin{equation}\label{eq:red-projection-hypothesis}
 P^{\rm new}RTf=c_f f,
\end{equation}
then
\begin{equation}\label{eq:red-cf}
 c_f=\frac{\|Tf\|_\Gamma^2}{\|f\|_H^2}\geq0.
\end{equation}
Consequently, $RTf-cf$ lying in a subspace orthogonal to $S^{\rm new}$
is sufficient to obtain $\|Tf\|^2=c\|f\|^2$.  The stronger identity
$RTf=cf$ is unnecessary.
\end{corollary}

\begin{proof}
Pair \eqref{eq:red-projection-hypothesis} with $f$ and apply
\eqref{eq:red-norm-trace}.
\end{proof}

\begin{remark}\label{rem:red-oldform-gap}
The assertion that Fourier coefficients vanish when $(n,N)=1$ does not,
by itself in the present argument, prove membership in a specified
half-integral-weight oldspace.  An applicable oldform criterion and its
hypotheses are needed to establish
\eqref{eq:red-projection-hypothesis}.  In particular, the numerical
assertion $c_f=\nu/2$ is an additional input, not a consequence of
matching Hecke eigenvalues.  The argument below can use the norm
$\|Tf\|_\Gamma$ directly and does not require this assertion.
\end{remark}

\subsection{A finite representation calculation for the trace constant}

Retain the unitary coefficient space $W$, $e_0$, $R$, and $T$ of
\eqref{eq:red-induction}. Let $E_0(v)=(v,e_0)_We_0$, let
$I=[\Gamma:H]$, and put
\begin{equation}\label{eq:frame-operator}
 S=\sum_{\gamma\in H\backslash\Gamma}
             \rho(\gamma)^{-1}E_0\rho(\gamma),
 \qquad W_{\rm cyc}=\operatorname{span}\{\rho(\gamma)e_0:\gamma\in\Gamma\}.
\end{equation}

\begin{lemma}[Finite-frame normalization]\label{lem:frame-operator}
For every vector-valued modular form $F$, $TRF=SF$ pointwise.
The operator $S$ commutes with $\rho(\Gamma)$, has trace $I$, is zero on
$W_{\rm cyc}^{\perp}$, and is positive definite on $W_{\rm cyc}$.
If $W_{\rm cyc}$ is irreducible, then
\[
 S|_{W_{\rm cyc}}=\frac{I}{\dim W_{\rm cyc}}\operatorname{id}.
\]
In this case, if $F$ is a $W_{\rm cyc}$-valued cusp form and $g=RF\ne0$, then
\begin{equation}\label{eq:frame-norm-ratio}
 \frac{\|Tg\|_\Gamma^2}{\|g\|_H^2}=
 \frac{I}{\dim W_{\rm cyc}}.
\end{equation}
\end{lemma}
\begin{proof}
The transformation law of $F$ gives
$(RF)|_k\gamma=(\rho(\gamma)F,e_0)_W$, hence $TRF=SF$.
Conjugation of $E_0$ by $\rho(h)$ leaves $E_0$ unchanged for $h\in H$.
The sum is therefore independent of coset representatives, and right
translation of the left cosets proves that $S$ commutes with $\rho$.
Each summand is a rank-one orthogonal projection. Thus
\[
 (Sv,v)_W=\sum_{\gamma\in H\backslash\Gamma}
                        |(\rho(\gamma)v,e_0)_W|^2.
\]
The vectors $\rho(\gamma)^{-1}e_0$ in this finite sum span
$W_{\rm cyc}$, because every group element is a left coset representative
times an element of $H$ after inversion. It follows that
$\ker S=W_{\rm cyc}^{\perp}$ and $\tr S=I$. Schur's lemma proves the
scalar formula in the irreducible case. Now $RTg=I g/\dim W_{\rm cyc}$;
pair this with $g$ and use \eqref{eq:red-norm-trace}.
\end{proof}

The scalar formula requires irreducibility. In the reducible case,
diagonalizing $S$ on its cyclic subspace retains its distinct eigenvalues.
This finite calculation does not evaluate the ramified doubling operator.
For $H$ above $\Gamma_0(4N)$, the effective index is
$I=4N\prod_{p\mid4N}(1+p^{-1})$ with the Petersson conventions of
\eqref{eq:red-pet-scalar}--\eqref{eq:red-pet-vector}.



\subsection{The odd-primary frame at every exponent}
\label{subsec:odd-frame-general}

Fix an odd prime $p$ and a $p$-adic unit $a$. For $n\geq1$, put
\[
 D_n=\mathbb Z/p^n\mathbb Z,\qquad Q_n(x)=a x^2/p^n,
 \qquad V_n=\mathbb C[D_n],
\]
and put $D_0=\{0\}$, $V_0=\mathbb C$.
Let $\iota_n e_x=e_{-x}$ and $P_n^+=(1+\iota_n)/2$.
Write $V_n^+=P_n^+V_n$ and $V_n^-=(1-P_n^+)V_n$.
Use the odd finite Weil representation $\rho_n$ and its unitary
conjugation action. In coordinates its generators are
\begin{align*}
 \rho_n(T)e_x&=e(a x^2/p^n)e_x,\\
 \rho_n(S)e_x&=\gamma_n^{-1}p^{-n/2}
                    \sum_{y\bmod p^n}e(-2axy/p^n)e_y,
 \qquad
 \gamma_n=p^{-n/2}\sum_{x\bmod p^n}e(a x^2/p^n).
\end{align*}
Set $\gamma_0=1$. Common scalar phases in a realization of the Weil
action disappear from all conjugated rank-one projectors below.
The upper triangular subgroup modulo $p^n$ preserves $\mathbb Ce_0$:
upper unipotents fix $e_0$, and diagonal units act on it by a unit
scalar. Thus, writing $H_n=\Gamma_0(p^n)\subset\operatorname{SL}_2(\mathbb Z)$,
we may define
\[
 S_n=\sum_{\gamma\in H_n\backslash\operatorname{SL}_2(\mathbb Z)}
             \rho_n(\gamma)^{-1}E_0\rho_n(\gamma),
 \qquad E_0(v)=(v,e_0)e_0.
\]
For $n=0$ set $H_0=\operatorname{SL}_2(\mathbb Z)$ and $S_0=1$.
The same operators are obtained from the finite congruence groups.
As before,
\[
 I_0=1,\qquad I_n=[\operatorname{SL}_2(\mathbb Z):H_n]
                  =p^{n-1}(p+1)\quad(n\geq1).
\]

For $n\geq2$ define the isometry
\begin{equation}\label{eq:odd-frame-embedding}
 J_n:V_{n-2}\longrightarrow V_n,\qquad
 J_ne_z=p^{-1/2}\sum_{t\bmod p}e_{pz+t p^{n-1}}.
\end{equation}
Its image consists of functions supported on $pD_n$ and invariant
under translation by $p^{n-1}$. This is the isotropic lift from the
quotient of $D_n$ by its order-$p$ isotropic subgroup.

\begin{lemma}[Equivariance of the isotropic lift]
\label{lem:odd-frame-J}
For $n\geq2$, the map $J_n$ intertwines $\rho_{n-2}$ and $\rho_n$,
and $J_nP_{n-2}^+=P_n^+J_n$. Moreover
$J_n^*e_0=p^{-1/2}e_0$.
\end{lemma}
\begin{proof}
The parity and adjoint assertions follow directly from
\eqref{eq:odd-frame-embedding}. For the $T$ action, if
$x=pz+t p^{n-1}$, then
$a x^2/p^n\equiv a z^2/p^{n-2}\pmod{\mathbb Z}$.
For the Fourier action, sum first over $t$ in
\eqref{eq:odd-frame-embedding}. The sum
$\sum_t e(-2aty/p)$ vanishes unless $p\mid y$ and equals $p$
otherwise. The remaining Fourier sum is exactly the one on $D_{n-2}$,
with the factor $p^{-1/2}$ required by $J_n$.
The phases agree because $\gamma_n=\gamma_{n-2}$.
Indeed, splitting a quadratic Gauss sum modulo $p^n$ into classes
$x=r+t p^{n-1}$ first forces $p\mid r$, and gives
\[
 \sum_{x\bmod p^n}e(a x^2/p^n)
 =p\sum_{z\bmod p^{n-2}}e(a z^2/p^{n-2}).
\]
This also holds when $n=2$, with the final sum interpreted as $1$.
Since $S$ and $T$ generate the modular group, the intertwining follows.
\end{proof}

\begin{theorem}[Exact odd-primary frame recurrence]
\label{thm:odd-frame-recurrence}
One has $S_1=2P_1^+$. For every $n\geq2$,
\begin{equation}\label{eq:odd-frame-recurrence}
 S_n=2\bigl(P_n^+-J_nP_{n-2}^+J_n^*\bigr)
              +\kappa_nJ_nS_{n-2}J_n^*,
 \qquad
 \kappa_n=\begin{cases}p+1,&n=2,\\p,&n\geq3.\end{cases}
\end{equation}
Thus the primitive even subspace
$V_n^+\ominus J_nV_{n-2}^+$ has eigenvalue $2$, while the old
restriction is $\kappa_n S_{n-2}$ under $J_n$.
The operator is zero on $V_n^-$ and positive definite on $V_n^+$.
In particular, the cyclic subspace generated by $e_0$ is exactly $V_n^+$.
\end{theorem}

\begin{proof}
The frame is self-adjoint and commutes with $\rho_n$ by the general
finite-frame lemma. Parity commutes with the Weil action and fixes
$e_0$. Every projector in the frame therefore has image in $V_n^+$,
and $S_n$ vanishes on $V_n^-$.

\emph{The unit-coordinate calculation.}
Represent the projective bottom rows modulo $p^n$ by
\[
 (1,b)\ (b\bmod p^n),\qquad
 (pc,1)\ (c\bmod p^{n-1}).
\]
Use matrices $ST^b$ and $n^-(pc)$, respectively.
The vectors $\rho_n(ST^b)^{-1}e_0$ have coordinates, up to a unit
phase, $p^{-n/2}e(-ba x^2/p^n)$. Their projection sum consequently
has matrix
\begin{equation}\label{eq:odd-frame-squares}
 A_n(x,y)=\mathbf1_{x^2\equiv y^2\ (p^n)}.
\end{equation}
For a lower unipotent with parameter $t=pc$, its inverse image of
$e_0$ is supported on $pD_n$. To see this, when $t\ne0$ write
$t=p^ku$, with $k\geq1$ and $u$ a unit. Its coordinate at $x$ is,
up to a scalar, the quadratic Fourier sum
\[
 \sum_{r\bmod p^n}e\bigl((at r^2-2axr)/p^n\bigr).
\]
Translation $r\mapsto r+p^{n-k}$ leaves the quadratic factor
unchanged and multiplies the sum by $e(-2ax/p^k)$.
It therefore vanishes unless $p^k\mid x$.
The parameter $t=0$ gives $e_0$ and has the same support property.
Thus these projectors kill every unit-coordinate vector $e_y$.
For $y$ a unit, equality of squares in
\eqref{eq:odd-frame-squares} forces $x=y$ or $x=-y$ modulo $p^n$,
because $p$ is odd. Consequently
\begin{equation}\label{eq:odd-frame-unit-eigen}
 S_ne_y=e_y+e_{-y}\qquad(p\nmid y).
\end{equation}

\emph{What the unit vectors generate.}
Let $U_n$ be the invariant subspace generated by
$e_y+e_{-y}$ for unit $y$. For $n\geq2$ we claim
\begin{equation}\label{eq:odd-frame-unit-span}
 U_n=V_n^+\ominus J_nV_{n-2}^+.
\end{equation}
The old space $J_nV_{n-2}^+$ is invariant by
Lemma~\ref{lem:odd-frame-J}, and all of its vectors are supported on
$pD_n$. It is therefore orthogonal to $U_n$.
Conversely, let $v\in V_n^+$ be orthogonal to $U_n$.
Its unit coordinates vanish. Since $U_n$ is invariant, the same is
true of its Fourier transform. Fourier support in $pD_n$ is equivalent
to translation invariance of $v$ by $p^{n-1}$: the two subgroups are
annihilators for the pairing $2axy/p^n$.
Thus $v$ is supported on $pD_n$ and is invariant under that
translation, which is exactly the image of $J_n$.
Its evenness places it in $J_nV_{n-2}^+$, proving the claim.
For $n=1$, the same reasoning would force both $v$ and its Fourier
transform to be supported at $0$, which is impossible for a nonzero
vector. Hence $U_1=V_1^+$.
By \eqref{eq:odd-frame-unit-eigen} and equivariance, $S_n=2$ on $U_n$.
This proves $S_1=2P_1^+$ and the primitive-even assertion for $n\geq2$.
Self-adjointness now implies preservation of its orthogonal complement
$J_nV_{n-2}^+$.

\emph{The old restriction.}
Using the intertwining and $J_n^*e_0=p^{-1/2}e_0$ gives
\[
 J_n^*S_nJ_n
 =\frac1p\sum_{\gamma\in H_n\backslash\operatorname{SL}_2(\mathbb Z)}
      \rho_{n-2}(\gamma)^{-1}E_0\rho_{n-2}(\gamma).
\]
Because $H_n\subset H_{n-2}$ and the lower-level $e_0$ line is
$H_{n-2}$-stable, each lower-level projector occurs
$I_n/I_{n-2}$ times. Therefore
\[
 J_n^*S_nJ_n=\frac{I_n}{pI_{n-2}}S_{n-2}
           =\kappa_nS_{n-2}.
\]
Combining the primitive, old, and odd subspaces proves
\eqref{eq:odd-frame-recurrence}.
The positivity assertion follows inductively from $S_0=1$ and
$S_1=2P_1^+$. The image of a frame is the cyclic subspace generated
by its seed, proving the last assertion.
\end{proof}

\begin{corollary}[Complete spectrum and multiplicities]
\label{cor:odd-frame-spectrum}
The kernel of $S_n$ is $V_n^-$, of dimension $(p^n-1)/2$.
The remaining spectrum is as follows.
\begin{enumerate}
\item If $n=2m+1$, for each $0\leq r\leq m$ there is an eigenvalue
      $2p^r$. Its multiplicity is
      \[
      \begin{cases}
       (p^{n-2r}-p^{n-2r-2})/2,&r<m,\\
       (p+1)/2,&r=m.
      \end{cases}
      \]
\item If $n=2m\geq2$, for each $0\leq r<m$ the eigenvalue $2p^r$
      has multiplicity $(p^{n-2r}-p^{n-2r-2})/2$.
      There is in addition the eigenvalue $(p+1)p^{m-1}$ of
      multiplicity one.
\end{enumerate}
For $n=0$ the only eigenvalue is $1$, with multiplicity one.
All these lists have trace $I_n$.
\end{corollary}
\begin{proof}
The even space has dimension $(p^n+1)/2$.
For $j\geq2$ its primitive even part has dimension
$(p^j-p^{j-2})/2$, while the terminal spaces at $j=1$ and $j=0$
have dimensions $(p+1)/2$ and $1$.
Iterate \eqref{eq:odd-frame-recurrence}; every step above exponent two
multiplies old eigenvalues by $p$, while the final step from exponent
two to zero multiplies by $p+1$. The trace follows either by summing
these multiplicities or from the defining sum of $I_n$ rank-one
projectors.
\end{proof}

\begin{remark}[An exact matrix algorithm]
The recurrence gives an exact algorithm with rational entries and
requires no floating-point Gauss sums. Start from $S_0=(1)$ and
$S_1(x,y)=\delta_{x,y}+\delta_{x,-y}$.
For $n\geq2$ one may equivalently compute
\[
 S_n(x,y)=\delta_{x,y}+\delta_{x,-y}
 +\frac{\mathbf1_{p\mid x,\ p\mid y}}p
  \left(\kappa_nS_{n-2}(\bar x/p,\bar y/p)
       -\delta_{\bar x/p,\bar y/p}
       -\delta_{\bar x/p,-\bar y/p}\right),
\]
where barred quotients are reduced modulo $p^{n-2}$.
The formula shows explicitly that the frame is independent of the unit
$a$ in these coordinates. Computing only its spectrum requires a
number of recurrence steps linear in $n$.
\end{remark}


\subsection{The dyadic finite frame at every exponent}
\label{subsec:dyadic-frame-general}

Let $n\geq1$, let $a$ be odd, and put
\[
 D_n=\mathbb Z/2^n\mathbb Z,
 \qquad q_n(x)=\frac{a x^2}{2^{n+1}}\pmod{\mathbb Z},
 \qquad W_n=\mathbb C[D_n].
\]
Thus $D_n$ has level $2^{n+1}$. The dyadic primary component of our
discriminant module corresponds to $n=v_2(N)+1$.
Use the orthonormal basis $e_x$, the standard Weil representation
$\rho_n$, and the inverse image $H_n$ of $\Gamma_0(2^{n+1})$ in
$\Gamma=\operatorname{Mp}_2(\mathbb Z)$.
The line $\mathbb C e_0$ is stable under $H_n$.
The finite frame is
\[
 \mathcal S_n=
 \sum_{\gamma\in H_n\backslash\Gamma}
       \rho_n(\gamma)^{-1}E_0\rho_n(\gamma),
 \qquad E_0v=(v,e_0)e_0.
\]
All central lifts give the same conjugated projections.
As for the general finite-frame lemma, $\mathcal S_n$ is self-adjoint,
positive, commutes with $\rho_n(\Gamma)$, and has trace $3\cdot2^n$.
Let $W_n^\pm$ be the eigenspaces of $\iota e_x=e_{-x}$, and let
$P_n^+=(1+\iota)/2$.

For $n\geq3$ define an isometric map
\begin{equation}\label{eq:dyadic-old-embedding}
 U_n:W_{n-2}\longrightarrow W_n,
 \qquad
 U_ne_z=\frac{e_{2z}+e_{2z+2^{n-1}}}{\sqrt2}.
\end{equation}
It is well-defined for $z$ modulo $2^{n-2}$.
Its image is the space of functions supported on $2D_n$ and invariant
under translation by $2^{n-1}$.
Equivalently it is the isotropic induction space for
$H=\{0,2^{n-1}\}$, since $H^\perp=2D_n$ and
$H^\perp/H\simeq D_{n-2}$ as quadratic modules.

\begin{theorem}[Exact dyadic recursion]\label{thm:dyadic-frame-recursion}
The initial frames are
\[
 \mathcal S_1=3\operatorname{id}_{W_1},
 \qquad \mathcal S_2=4P_2^+.
\]
For $n\geq3$, put
\[
 W_{n,\mathrm{prim}}^+
   =W_n^+\ominus U_nW_{n-2}^+.
\]
The following orthogonal decomposition is preserved by the frame:
\[
 W_n=W_n^-\oplus W_{n,\mathrm{prim}}^+
                   \oplus U_nW_{n-2}^+.
\]
On its three summands, respectively, the frame acts by
\begin{equation}\label{eq:dyadic-frame-recursion}
 0,\qquad 4\operatorname{id},\qquad
 2U_n\mathcal S_{n-2}U_n^*.
\end{equation}
In particular the cyclic space generated by $e_0$ is exactly $W_n^+$
at every exponent. This assertion does not assume irreducibility of
the primitive space.
\end{theorem}

\begin{proof}
We give the coset computation and the generation argument explicitly.
Write $\mathrm e(t)=\exp(2\pi i t)$.
Up to their standard scalar phases, the Weil generators act by
\[
 \rho_n(T)e_x=\mathrm e(q_n(x))e_x,
 \qquad
 \rho_n(S)e_x=
  \frac{\kappa_n}{\sqrt{2^n}}
       \sum_{y\bmod2^n}\mathrm e(-a xy/2^n)e_y,
 \qquad |\kappa_n|=1.
\]

\emph{The isotropic embedding.}
The maps $U_n$ intertwine these Weil representations.
For $T$, the two coordinates in \eqref{eq:dyadic-old-embedding} have
quadratic value $q_{n-2}(z)$ modulo integers.
For $S$, summing those two coordinates gives the factor
$1+(-1)^x$, and hence the required smaller Fourier transform on even
coordinates. The scalar phases agree: the normalized Gauss sums
\[
 g_n=2^{-n/2}\sum_{x\bmod2^n}\mathrm e(a x^2/2^{n+1})
\]
satisfy $g_n=g_{n-2}$ for $n\geq3$.
Indeed, pairing $x$ with $x+2^{n-1}$ cancels the odd terms;
writing $x=2z$ in the even terms leaves twice the lower Gauss sum.
Their normalizing factors give the stated equality, and therefore
$\kappa_n=\kappa_{n-2}$ in the standard Weil convention.
The embedding also intertwines $\iota$, so both its even image and its
orthogonal complement inside $W_n^+$ are Weil invariant.

\emph{The coset vectors.}
Represent the projective bottom rows modulo $2^{n+1}$ by
\[
 (1,b),\quad b\bmod2^{n+1},
 \qquad
 (2c,1),\quad c\bmod2^n.
\]
Representing matrices are $ST^b$ and
$n^-(2c)=\left(\begin{smallmatrix}1&0\\2c&1\end{smallmatrix}\right)$.
The first family's vectors $\rho_n(ST^b)^{-1}e_0$ have coordinates,
up to a common phase, $2^{-n/2}\mathrm e(-ba x^2/2^{n+1})$.
Their projection sum has matrix
\begin{equation}\label{eq:dyadic-first-family}
 A_n(x,y)=2\,\mathbf1_{x^2\equiv y^2\ (2^{n+1})}.
\end{equation}
The second family's vectors have coordinates, up to harmless scalar
phases,
\begin{equation}\label{eq:dyadic-second-family}
 v_c(x)=2^{-n}\sum_{z\bmod2^n}
                 \mathrm e\bigl(a(cz^2-xz)/2^n\bigr).
\end{equation}
For $n\geq2$, shifting $z$ by $2^{n-1}$ multiplies the summand by
$(-1)^x$. Thus $v_c(x)=0$ whenever $x$ is odd, for every $c$.
The support assertion also covers the degenerate cases $c=0$ and
$v_2(c)=n-1$, without a nonsingular Gauss-sum hypothesis.

\emph{The initial frames.}
For $n=1$, \eqref{eq:dyadic-first-family} gives $A_1=2\operatorname{id}$,
and the two vectors in \eqref{eq:dyadic-second-family} are $e_0$ and
$e_1$. Hence $\mathcal S_1=3\operatorname{id}$.
For $n=2$, the vectors with $c=0,2$ are $e_0,e_2$.
The vectors with $c=1,3$ form an orthonormal basis of
$\operatorname{span}\{e_0,e_2\}$: after adjusting their scalar phases
they are $(e_0\pm i^a e_2)/\sqrt2$.
The second projection sum is consequently $2E_0+2E_2$.
The first sum has entries $2$ on the separate coordinates $0,2$ and
the all-$2$ matrix on $\{1,3\}$.
Adding gives $\mathcal S_2=4P_2^+$.

\emph{The primitive even space.}
Fix $n\geq3$.
For odd $y$, the congruence $x^2\equiv y^2\pmod{2^{n+1}}$ is equivalent
to $x\equiv y$ or $x\equiv-y\pmod{2^n}$.
To see this, one of $v_2(x-y),v_2(x+y)$ is exactly $1$, and their sum
must be at least $n+1$.
The second family has no odd coordinates. Consequently
\[
 \mathcal S_n(e_y+e_{-y})=4(e_y+e_{-y})\qquad(y\text{ odd}).
\]
These unit symmetric vectors lie in $W_{n,\mathrm{prim}}^+$, and their
Weil orbit spans it. For a proof, let $M$ be their invariant orbit span,
and suppose $f\in W_{n,\mathrm{prim}}^+\cap M^\perp$.
Since $f$ is even and orthogonal to all the unit symmetric vectors,
it is supported on even coordinates. The same is true of its Fourier
transform, because $M^\perp\cap W_{n,\mathrm{prim}}^+$ is invariant.
Fourier support on even coordinates means that $f$ is invariant under
translation by $2^{n-1}$: this follows immediately by Fourier inversion,
or by multiplying its Fourier coefficients by $(-1)^x$.
Thus $f$ is supported on even coordinates and has that translation
invariance, so $f\in U_nW_{n-2}$.
Its defining orthogonality to the old even space forces $f=0$.
This proves $M=W_{n,\mathrm{prim}}^+$.
Since the frame commutes with every Weil operator, it acts by $4$ on
this orbit span.

\emph{The old space and the odd kernel.}
The vector $e_0$ and all its Weil translates are even, so the frame
annihilates $W_n^-$.
Self-adjointness and the primitive calculation show that the remaining
old even space is also preserved.
The embedding identity and $U_n^*e_0=e_0/\sqrt2$ give
\[
 U_n^*\mathcal S_nU_n
  =\frac12\sum_{\gamma\in H_n\backslash\Gamma}
       \rho_{n-2}(\gamma)^{-1}E_0\rho_{n-2}(\gamma).
\]
The index $[H_{n-2}:H_n]$ is $4$.
Conjugation by $H_{n-2}$ fixes the smaller $E_0$, so this sum is four
copies of the smaller frame, and its value is $2\mathcal S_{n-2}$.
This proves \eqref{eq:dyadic-frame-recursion}.
Finally the two initial frames and the recursion are positive definite
on the entire even space. The kernel/image characterization of a frame
therefore identifies its cyclic space with $W_n^+$.
\end{proof}

\begin{corollary}[All dyadic eigenvalues and multiplicities]
\label{cor:dyadic-frame-spectrum}
Put $m=n-1=v_2(N)$. The zero eigenvalue has multiplicity $2^m-1$.
The nonzero spectrum is as follows.

If $m=2r$ is even, the primitive layers have eigenvalues and multiplicities
\[
 2^{j+2}\quad\text{with multiplicity}\quad
       3\cdot4^{r-j-1},\qquad 0\leq j<r,
\]
and the final layer has eigenvalue $3\cdot2^r$ with multiplicity $2$.
For $r=0$ the primitive list is empty, giving the initial eigenvalue
$3$ with multiplicity $2$.

If $m=2r+1$ is odd, the primitive layers have
\[
 2^{j+2}\quad\text{with multiplicity}\quad
       6\cdot4^{r-j-1},\qquad 0\leq j<r,
\]
and the final layer has eigenvalue $2^{r+2}$ with multiplicity $3$.
For $r=0$ this gives just the initial eigenvalue $4$ with multiplicity $3$.
Every displayed spectrum is independent of the odd unit $a$.
\end{corollary}

\begin{proof}
One has $\dim W_n^+=2^{n-1}+1$ and
$\dim W_n^-=2^{n-1}-1$.
For $n\geq3$,
\[
 \dim W_{n,\mathrm{prim}}^+
    =(2^{n-1}+1)-(2^{n-3}+1)=3\cdot2^{n-3}.
\]
At each induction step the old eigenvalues are multiplied by $2$.
Iterating the theorem to $n=1$ or $n=2$ gives the stated lists.
\end{proof}

\begin{remark}[Scope]
These formulas determine the dyadic coefficient-space frame and,
through scalar--vector adjunction, the trace norm on each frame
eigenspace. Evaluating a dyadic local doubling operator, proving the
arithmetic product identity, and establishing cuspidal orthogonality
of the full arithmetic theta series require separate arguments.
\end{remark}

\subsection{The spectrum at every global level}
\label{subsec:all-level-frame}

Write $N=2^m\prod_{p\mid N,\ p\text{ odd}}p^{n_p}$.
Let $A_N=\Z/(2N)\Z$ with $Q(x)=-x^2/(4N)$, and let $S_N$
be its coefficient-space frame for the subgroup above $\Gamma_0(4N)$.
For every odd $p$, let $\mathcal L_{p,n_p}$ denote the multiset of
positive eigenvalues in Corollary~\ref{cor:odd-frame-spectrum}, with the
stated multiplicities. Let $\mathcal L_{2,m}$ denote the dyadic multiset
of Corollary~\ref{cor:dyadic-frame-spectrum}.

\begin{theorem}[All-level frame spectrum]\label{thm:all-level-frame}
The cyclic subspace generated by $e_0$ is the tensor product of the
local even subspaces. Its dimension is
\begin{equation}\label{eq:all-level-cyclic-dimension}
 d_N=(2^m+1)\prod_{p\mid N,\ p\text{ odd}}\frac{p^{n_p}+1}{2}.
\end{equation}
The frame is zero on its orthogonal complement and positive definite
on this cyclic subspace. Its positive spectrum is the multiset
\begin{equation}\label{eq:all-level-product-spectrum}
 \left\{\lambda_2\prod_{p\mid N,\ p\text{ odd}}\lambda_p:
       \lambda_2\in\mathcal L_{2,m},\quad
       \lambda_p\in\mathcal L_{p,n_p}\right\}.
\end{equation}
Multiplicities multiply in this formula; equal numerical products are
combined by adding their multiplicities. In particular,
\[
 \tr S_N=[\SL_2(\Z):\Gamma_0(4N)]
 =3\cdot2^{m+1}\prod_{p\mid N,\ p\text{ odd}}p^{n_p-1}(p+1).
\]
All eigenvalues and their multiplicities are independent of the unit
coefficients in the local primary quadratic forms.
\end{theorem}
\begin{proof}
The orthogonal primary decomposition of $A_N$ has a dyadic cyclic
factor of order $2^{m+1}$ and, for each odd $p$, a cyclic factor of
order $p^{n_p}$. Their quadratic forms have exactly the shapes in
Theorems~\ref{thm:odd-frame-recurrence} and
\ref{thm:dyadic-frame-recursion}, with unit coefficients.
The Weil representation is the tensor product of these primary
representations. On conjugated rank-one projectors all scalar
metaplectic phases cancel. The resulting action factors through the
product of the local congruence groups.

The map from $\SL_2(\Z)$ to this product is onto, and the projective
bottom-row cosets for $\Gamma_0(4N)$ decompose into the product of
the local projective cosets by the Chinese remainder theorem.
Thus the defining projector sum is the tensor product of the local
frames. The local kernel and positivity results identify its cyclic
space, and the spectral theorem gives the product list and
multiplicities. The local even dimensions are $2^m+1$ and
$(p^{n_p}+1)/2$, proving \eqref{eq:all-level-cyclic-dimension}.
The trace also follows directly from the sum of one rank-one projector
for each coset. Unit independence was proved in both local calculations.
\end{proof}

\begin{corollary}[Scalar norm decomposition]\label{cor:all-level-norm}
Let $\Lambda_N$ be the set of distinct positive values in
\eqref{eq:all-level-product-spectrum}. In a compatible fixed weight,
put $\mathcal E_{\rm act}=\operatorname{im}R$ inside the scalar cusp
space and $K=RT$. There is an orthogonal decomposition
\[
 \mathcal E_{\rm act}=\bigoplus_{\lambda\in\Lambda_N}\mathcal E_\lambda,
 \qquad K|_{\mathcal E_\lambda}=\lambda\operatorname{id}.
\]
Some of these scalar spaces may be zero. If $g=\sum_\lambda g_\lambda$,
then
\begin{equation}\label{eq:all-level-scalar-norm}
 \|Tg\|_\Gamma^2=\sum_{\lambda\in\Lambda_N}
                         \lambda\|g_\lambda\|_H^2.
\end{equation}
On the full scalar cusp space $K$ is zero on
$(\operatorname{im}R)^\perp$. In particular, when $N$ is square-free,
\begin{equation}\label{eq:squarefree-exact-trace}
 RT=c_NP_{\operatorname{im}R},\qquad
 c_N=\begin{cases}
 3\,2^{\omega(N)},&N\text{ odd},\\
 4\,2^{\omega(N)-1},&N\text{ even},
 \end{cases}
\end{equation}
where $\omega(N)$ is the number of distinct prime divisors of $N$.
\end{corollary}
\begin{proof}
Adjunction gives $K=T^*T$ and $TR=S_N$ on vector cusp forms.
Every value of $Tf$ belongs to the cyclic coefficient space.
Let $S_N^\dagger$ be the inverse of $S_N$ on that space,
extended by zero on its orthogonal complement. It commutes
with $\rho$ and preserves vector-valued modularity and
cuspidality. For a cyclic-valued vector cusp form $F$,
\[
       T\bigl(R(S_N^\dagger F)\bigr)=S_NS_N^\dagger F=F.
\]
Thus the image of $T$ is exactly the cyclic-valued vector cusp
space. The positive spectra of $T^*T$ and $TT^*$ agree on
their actual nonzero spectral subspaces, which proves the decomposition.
Adjunction now gives \eqref{eq:all-level-scalar-norm}. At square-free
level every local positive eigenvalue is constant: $2$ at odd primes,
$3$ for the order-two dyadic module, and $4$ when $v_2(N)=1$.
Their products give \eqref{eq:squarefree-exact-trace}.
\end{proof}

\begin{corollary}[Explicit norm bounds]\label{cor:all-level-bounds}
For every $g\in\operatorname{im}R$,
\[
 a_N\|g\|_H^2\leq\|Tg\|_\Gamma^2\leq b_N\|g\|_H^2,
\]
where
\[
 a_N=\begin{cases}3,&m=0,\\4,&m\geq1\end{cases}
        \prod_{p\mid N,\ p\text{ odd}}2,
 \qquad
 b_N=b_2(m)\prod_{p\mid N,\ p\text{ odd}}b_p(n_p),
\]
and
\[
 b_2(m)=\begin{cases}
   3\cdot2^{m/2},&m\text{ even},\\
   2^{(m+3)/2},&m\text{ odd},
 \end{cases}
 \qquad
 b_p(n)=\begin{cases}
   2p^{(n-1)/2},&n\text{ odd},\\
   (p+1)p^{n/2-1},&n\text{ even}.
 \end{cases}
\]
These are the extreme eigenvalues on the finite cyclic coefficient
space. They need not both be attained by cusp forms in every weight.
\end{corollary}
\begin{proof}
Read the smallest and largest positive values from the local lists,
take their products, and apply \eqref{eq:all-level-scalar-norm}.
\end{proof}

\begin{remark}[Finite primitive spaces and automorphic conductor]
The primitive subspaces in the local recursions are the specified
orthogonal complements of isotropic lifts. Identifying them with an
exact-conductor Shimura packet requires compatibility of the chosen
Shimura correspondence with level lowering; no such identification
is made here. The existence of nonzero cusp forms in a finite
coefficient summand remains a separate question.
\end{remark}


\subsection{Polar normalization of the scalar induction}
\label{subsec:polar-normalized-induction}

Let $\mathcal E=S_k(H,\chi)$ and $\mathcal V=S_k(\Gamma,\rho)$,
with the Petersson products already fixed, and let
$T:\mathcal E\to\mathcal V$ be the unnormalized trace.
By adjunction $T^*=R$. Put
\[
 A=RT=T^*T,\qquad
 \mathcal E_{\rm act}=(\ker T)^\perp=\operatorname{im}R,
 \qquad \mathcal V_{\rm act}=\operatorname{im}T.
\]
These are finite-dimensional Hilbert spaces.
Define $A^{-1/2}$ by the positive inverse square root on
$\mathcal E_{\rm act}$ and by zero on $\ker T$.

\begin{proposition}[Canonical partial isometry]
\label{prop:polar-normalized-induction}
The map
\[
 U=TA^{-1/2}
\]
is a unitary map from $\mathcal E_{\rm act}$ onto
$\mathcal V_{\rm act}$. On the full spaces it satisfies
\[
 U^*U=P_{\mathcal E_{\rm act}},\qquad
 UU^*=P_{\mathcal V_{\rm act}}.
\]
If $\mathcal S=TR$ is the coefficient-space frame acting on vector
cusp forms, then equivalently $U=\mathcal S^{-1/2}T$, with the same
zero-on-kernel convention.
\end{proposition}
\begin{proof}
One has $U^*=A^{-1/2}R$, hence
$U^*U=A^{-1/2}AA^{-1/2}=P_{\mathcal E_{\rm act}}$.
Its image is $\operatorname{im}T$, so $UU^*$ is the corresponding
orthogonal projection.
Finally $\mathcal ST=TA$; applying functional calculus, or checking
on each positive singular-value subspace, gives the alternative formula.
\end{proof}

Suppose the scalar/vector arithmetic matching has been established, so
that $\widehat\theta_0(f)=\widehat\Theta(Tf)$.
Define the \emph{polar-normalized scalar lift} by
\[
 \widehat\theta_{\rm iso}(g)
 :=\widehat\theta_0(A^{-1/2}g)
 =\widehat\Theta(Ug),\qquad g\in\mathcal E_{\rm act}.
\]
This definition changes the normalization of the test-form lift.
On an $A$-eigenspace of eigenvalue $\lambda>0$ it gives
$\lambda^{-1/2}\widehat\theta_0(g)$, since the square root is positive
real. Substituting it for the original scalar lift therefore requires
recording the change of definition. The map $U$ preserves Petersson
norms, but this does not prove that the arithmetic lift preserves heights.

\begin{proposition}[Packet preservation requires compatibility]
\label{prop:polar-packet-compatibility}
Suppose a collection of corresponding scalar and vector Hecke
operators $H_\ell^{\rm sc}$ and $H_\ell^{\rm vec}$ satisfies
\[
 TH_\ell^{\rm sc}=H_\ell^{\rm vec}T,
 \qquad
 RH_\ell^{\rm vec}=H_\ell^{\rm sc}R.
\]
Then $A$ and $A^{-1/2}$ commute with every $H_\ell^{\rm sc}$, and
$U$ intertwines these scalar and vector operators. In particular,
the positive inverse square root preserves their joint eigenspaces.
\end{proposition}
\begin{proof}
Compose the two intertwining identities to obtain
$AH_\ell^{\rm sc}=H_\ell^{\rm sc}A$.
An operator commuting with $A$ preserves its eigenspaces and commutes
with its spectral functions, including the zero-on-kernel inverse
square root. The formula for $U$ now proves the assertion.
\end{proof}

For self-adjoint matching Hecke operators, the second intertwining
identity follows by taking adjoints of the first.
Neither identity follows from polar decomposition alone.
For an unramified Hecke packet, it is enough to establish them for the
specified operators away from $2N$. Additional local packet conditions
must also be preserved if the cited vector formula requires them.
Under this compatibility, a vector theorem
\[
 H_{\rm NT}(\widehat\Theta(F),\widehat\Theta(F))
      =C_G\|F\|_\Gamma^2L'(G,1)
\]
on the chosen packet gives
\[
 H_{\rm NT}(\widehat\theta_{\rm iso}(g),
                 \widehat\theta_{\rm iso}(g))
      =C_G\|g\|_H^2L'(G,1)
\]
for active scalar forms in that packet. Polar normalization removes
the induction Gram factor. The constant $C_G$ and the vector theorem
remain independent inputs.

\begin{remark}[A packet-specific construction without a commutation claim]
Fix instead a specified vector packet $\mathcal V_G$ and its Petersson
orthogonal projection $P_G$. Set
\[
 T_G=P_GT,\qquad A_G=T_G^*T_G=RP_GT,
 \qquad U_G=T_GA_G^{-1/2},
\]
where the inverse is taken on $(\ker T_G)^\perp$ and is zero on the
kernel. The preceding partial-isometry proof applies verbatim.
The image of $U_G$ lies in $\mathcal V_G$ by construction, without
assuming that $RT$ preserves that packet.
Define the projected raw scalar functional by
$\widehat\theta_{0,G}(f)=\widehat\Theta(P_GTf)$, and its normalized
version by
$\widehat\theta_{G,\rm iso}(g)=\widehat\theta_{0,G}(A_G^{-1/2}g)$.
A vector height theorem on $\mathcal V_G$ then transfers to this
functional with $\|g\|_H^2$ as its Petersson factor.
Identifying $\widehat\theta_{0,G}$ with a particular projection of the
original arithmetic scalar lift still requires arithmetic Hecke
compatibility. The construction canonically normalizes the specified
projected lift and proves no corresponding statement for the
unprojected raw lift.
\end{remark}

\section{Arithmetic projections at arbitrary level}
\label{sec:geometry}

\subsection{Vertical components and their intersection matrices}

For $p\mid N$, put $n=v_p(N)$ and $N_p=Np^{-n}$. In the notation of
\cite[Theorem~5.4.4]{Zhu}, the special fiber is
\begin{equation}\label{eq:special-fiber}
 F_p=\dv(p)=
 \sum_{\substack{-n\leq a\leq n\\a\equiv n\ (2)}}
 \varphi\bigl(p^{(n-|a|)/2}\bigr)\calX_p^a(N).
\end{equation}
Here $\varphi$ is Euler's totient function, with $\varphi(1)=1$.
For $n=1$ there are two components. For $n=2$ there are three,
with middle multiplicity $p-1$. The two-component description at
square-free level therefore needs modification at higher level.

For a finite vertical divisor $D$, write $\widehat D=(D,0)$.
The principal arithmetic divisor relation gives
\begin{equation}\label{eq:fiber-constant}
 (\dv(p),-2\log p)=0,
 \qquad \widehat F_p=2\log p\,\one,
 \qquad \one=(0,1).
\end{equation}
Consequently the full-fiber directions at different primes determine the
same constant line in the real arithmetic Chow group. They cannot be
treated as independent direct summands. Quotient each vertical component
space by its full-fiber direction and choose compatible representatives
$Y_{p,i}$ of a basis. Put
\[
 M_p=(\langle Y_{p,i},Y_{p,j}\rangle_{\rm GS})_{i,j}.
\]
The intersection pairing on this quotient is negative definite in the
standard regular arithmetic-surface setting; the corresponding stack
intersection statement is understood with the conventions of \cite{Zhu}.
In particular $M_p$ is invertible.

\begin{lemma}[The vertical dual basis]\label{lem:vertical-dual}
For any such basis, its intersection-dual basis is
\begin{equation}\label{eq:vertical-dual}
 Y_{p,i}^{\vee}=\sum_j(M_p^{-1})_{ji}Y_{p,j}.
\end{equation}
The component in this span of an arithmetic class $z$ is
\begin{equation}\label{eq:vertical-projection}
 z_{p}^{\rm corr}=\sum_i
             \langle z,Y_{p,i}\rangle_{\rm GS}Y_{p,i}^{\vee}.
\end{equation}
The simpler formula
$Y_{p,i}^{\vee}=Y_{p,i}/\langle Y_{p,i},Y_{p,i}\rangle_{\rm GS}$
requires that the chosen basis be orthogonal.
\end{lemma}
\begin{proof}
For every $k$, matrix multiplication gives
$\langle Y_{p,k},Y_{p,i}^{\vee}\rangle_{\rm GS}=\delta_{ki}$.
Pairing the right side of \eqref{eq:vertical-projection} with each
$Y_{p,k}$ therefore recovers $\langle z,Y_{p,k}\rangle_{\rm GS}$.
Nondegeneracy proves uniqueness. If the basis is orthogonal, $M_p$ is
diagonal and its inverse has the asserted diagonal entries.
\end{proof}

\subsection{Admissible decomposition}

Let $\nu$ be a smooth conjugation-invariant positive $(1,1)$-form on
$X_0(N)$ of total mass one. Write $A^0$ for the real smooth invariant
functions of $\nu$-mean zero, and $a(f)=(0,f)$. Choose an admissible
degree-one arithmetic class $\widehat\omega$ whose curvature is $\nu$.
Fix lifts $Y_{p,i}$ satisfying
$\langle\widehat\omega,Y_{p,i}\rangle_{\rm GS}=0$.
In the admissible decomposition of arithmetic Chow space
\cite[Chapter~4]{KRYBook}, use the convention
\begin{equation}\label{eq:arith-decomposition}
 \mathcal A_{\R}=
 \widetilde{\MW}_{\R}\oplus\R\widehat\omega
 \oplus\left(\R\one\oplus\bigoplus_{p\mid N}
                   \operatorname{span}_{\R}\{Y_{p,i}\}\right)
 \oplus a(A^0).
\end{equation}
The Mordell--Weil summand is represented by admissible degree-zero
classes with the appropriate vertical corrections. The precise extension
at cusps is inherited from the arithmetic Chow theory being used.

Under the standard normalization
$\langle\widehat\omega,\one\rangle_{\rm GS}=1/2$, an expansion
$z=z_{\rm MW}+d\widehat\omega+z_{\rm corr}+c\one+a(f)$
has $d=\deg z$ and
\begin{equation}\label{eq:constant-projection}
 c=2\left(\langle z,\widehat\omega\rangle_{\rm GS}
      -d\langle\widehat\omega,\widehat\omega\rangle_{\rm GS}\right).
\end{equation}
Indeed, the other summands pair trivially with $\widehat\omega$:
for $a(f)$ this follows from $\int f\nu=0$. Formulas
\eqref{eq:vertical-projection} and \eqref{eq:constant-projection}
retain the full intersection matrix and the constant arithmetic class.

\begin{proposition}[The orthogonality required of a full lift]
\label{prop:MW-orthogonality}
Suppose coefficientwise decomposition and integration commute for
\eqref{eq:scalar-arithmetic-series}. If its degree, corrected vertical,
constant, and smooth components all have zero Petersson pairing with
$g$, then
\[
 \widehat\theta_0(g)=
 \langle\pr_{\rm MW}\wtphi_0,g\rangle_H
 \in\widetilde{\MW}_{\C}.
\]
It is then orthogonal to the complementary summand for the arithmetic
intersection pairing. Conversely, membership in the Mordell--Weil
summand forces all other projections of this integral to vanish.
\end{proposition}
\begin{proof}
Integrate the direct-sum decomposition
\eqref{eq:arith-decomposition} and use the stated vanishings. The
converse follows from uniqueness of the direct-sum components.
Orthogonality is part of the admissible Mordell--Weil decomposition.
\end{proof}

Eisenstein descriptions can establish some of these vanishings once their
Petersson orthogonality is justified. The assertion that each component is
modular does not establish them: modular forms can have nonzero cuspidal
projections. In particular, the smooth component and the constant
arithmetic component must be treated explicitly.

\subsection{Boundary changes and the Mordell--Weil projection}
\label{subsec:MW-boundary-independence}

Comparison after Mordell--Weil projection requires less information
than comparison of complete arithmetic classes. The Manin--Drinfeld
theorem removes Green functions and cuspidal divisors from the projected
comparison, without a vanishing assertion for an isolated boundary
Fourier series.

\begin{proposition}[Independence of the projected class from boundary data]
\label{prop:MW-boundary-independence}
Use the admissible decomposition \eqref{eq:arith-decomposition}, based
on a degree-one cuspidal class $P_\infty$. Let $\widehat D_1$ and
$\widehat D_2$ be arithmetic divisor classes in the same arithmetic
Chow space. If their generic-fiber divisor classes differ by a complex
linear combination of cuspidal divisor classes, then
\[
 \pr_{\rm MW}\widehat D_1=\pr_{\rm MW}\widehat D_2.
\]
In particular, admissible changes of Green functions, vertical
corrections, and cuspidal divisors do not change this projection.
The statement applies only when the changed data define classes in
the stated arithmetic Chow space.
\end{proposition}
\begin{proof}
Generic-fiber restriction identifies the Mordell--Weil summand with
the degree-zero Picard group tensored with $\C$, and sends its
projection to
\begin{equation}\label{eq:MW-generic-projection}
 \bigl[D_\Q-(\deg D_\Q)P_\infty\bigr].
\end{equation}
Indeed, the degree-one summand restricts to $P_\infty$, whereas
vertical and metric classes restrict to zero. This also proves that
\eqref{eq:MW-generic-projection} characterizes the projection.

Pass to a finite extension over which all cusps in question are
rational. By the Manin--Drinfeld theorem \cite{Drinfeld,Elkik},
every degree-zero cuspidal divisor has torsion class in the
Jacobian, hence zero class after tensoring with $\C$.
Consequently, for a cuspidal divisor $C$,
\[
       [C-(\deg C)P_\infty]=0.
\]
The same assertion holds over the original field: restriction
followed by norm is multiplication by the extension degree,
which is invertible after tensoring with $\C$. For stack classes, use the rational Picard group of the generic coarse
curve with a consistent degree convention. Normalizing the base cusp
by that degree, including the common generic stabilizer factor, gives
the same rational divisor relation. Apply
\eqref{eq:MW-generic-projection} to
$\widehat D_1-\widehat D_2$. Changes of metrics and vertical
classes already have zero generic restriction, proving the last
assertion.
\end{proof}

\begin{corollary}[Matching via positive generic cycles]
\label{cor:MW-generic-matching}
Let $\widehat\phi_{\rm sc}=\sum_m\widehat D(m,v)q^m$ be an
arithmetic family in this Chow space. Write the zero component
of the Heegner series with values in the Jacobian as
\[
 R\widehat\Phi_{\rm MW}=\sum_{m>0}\widehat y(m,0)q^m,
\]
where
\[
 y(m,0)=[Z(m,0)-\deg Z(m,0)P_\infty].
\]
Suppose that, for a fixed scalar $a_N$, the following generic-fiber
conditions hold:
\begin{enumerate}
\item for every $m>0$, the class of $D(m,v)-a_N Z(m,0)$ is
      cuspidal;
\item for $m\leq0$, the class of $D(m,v)$ is cuspidal.
\end{enumerate}
Then, coefficientwise and hence for every convergent pairing,
\begin{equation}\label{eq:MW-generic-matching}
 \pr_{\rm MW}\widehat\phi_{\rm sc}
       =a_N R\widehat\Phi_{\rm MW}.
\end{equation}
Thus this projected family is independent of $v$ in its Fourier
coefficients and has no nonpositive Fourier terms. For the modified
zero coefficient, a Hodge class satisfies the second condition
whenever its generic rational Picard class is cuspidal, as in
\cite[Theorem~6.4.3]{Zhu}.
\end{corollary}
\begin{proof}
For $m>0$, apply Proposition~\ref{prop:MW-boundary-independence}
to $\widehat D(m,v)$ and any arithmetic extension of
$a_N Z(m,0)$. Formula~\eqref{eq:MW-generic-projection} gives
$a_N\widehat y(m,0)$. For $m\leq0$ the same proposition gives
zero. These coefficientwise equalities prove the series identity;
linearity proves the assertion for convergent pairings.
\end{proof}

Equality after projection does not identify the full lifts or replace
the compactified arithmetic product identity. The constant metric class
$\one=(0,1)$, for example, has zero Mordell--Weil projection but pairs
nontrivially with a degree-one class. Determining $a_N$ still requires
a separate check of the positive generic-cycle multiplicities, including
the stack convention and the treatment of $x$ and $-x$.

\begin{proposition}[The positive scalar cycle and the zero vector component]
\label{prop:generic-positive-matching}
Let $\mathcal X_0(N)_\mathbb Q$ be the usual modular stack of cyclic
degree-$N$ isogenies, and let $X_0(N)_\mathbb Q$ be its coarse curve.
For $m>0$, let $\mathcal Z(m)$ be the special-isogeny moduli stack
of \cite[Definition~7.1.1]{Zhu}, whose objects are $(E\xrightarrow\pi E',j)$ with
\[
 \deg\pi=N,\qquad \deg j=m,\qquad
       j^\vee\pi+\pi^\vee j=0.
\]
Over $\mathbb Q$ the displayed orthogonality condition is equivalent
to $j\pi^\vee+\pi j^\vee=0$, by composing with $\pi^\vee$
on the left and $\pi$ on the right and using $\pi^\vee\pi=[N]$.
Its generalized Cartier divisor on $\mathcal X_0(N)_\mathbb Q$
has rational Picard descent to the coarse curve equal to
\[
 Z_{\rm vec}(m,0)=
 \sum_{x\in\Gamma_0(N)\backslash L_N[m]}
          \frac{2}{|\Gamma_x|}[z_x],
 \qquad
 L_N=\left\{\begin{pmatrix}b&A/N\\C&-b\end{pmatrix}:
                    b,A,C\in\mathbb Z\right\},\quad Q=N\det.
\]
Here $\Gamma_x$ is the stabilizer in $\Gamma_0(N)\subset
\mathrm{SL}_2(\mathbb Z)$; the displayed notation means that
each orbit contributes the indicated rational point divisor.
Both $x$ and $-x$ are retained.  This is the divisor normalization
of the untwisted zero component in \cite[(1.16)]{Du}.
In particular the scalar $a_N$ in the generic-cycle comparison is
$1$, with these rational Picard and Heegner-divisor conventions.
\end{proposition}
\begin{proof}
To identify the moduli problem as well as the quadratic lattice,
choose integral Betti bases on $E,E'$ in which the cyclic isogeny
$\pi$ has matrix $P=\operatorname{diag}(N,1)$.
Such bases exist by the Smith normal form.  If $J$ is the integral
matrix of $j$, its dual is represented by $\operatorname{adj}(J)$.
The orthogonality equation is therefore
\[
 \operatorname{adj}(J)P+\operatorname{adj}(P)J=0.
\]
Writing $J=\left(\begin{smallmatrix}r&s\\t&u\end{smallmatrix}\right)$,
this equation is equivalent to $r=-Nu$.  Set
\[
 x=\frac1N\operatorname{adj}(J)P
   =-P^{-1}J
   =\begin{pmatrix}u&-s/N\\-t&-u\end{pmatrix}.
\]
Then $x\in L_N$ and $N\det x=\det J=m$.
Conversely $x\in L_N[m]$ gives the integral matrix $J=-Px$,
which satisfies the two required equations.  The rational
endomorphism represented by $x$ is
$\alpha/N$ with $\alpha=j^\vee\pi$.  Thus $x$ is holomorphic
precisely at its associated CM point $z_x$ in the orthogonal
uniformization.  This is also the special-endomorphism bridge in
\cite[Lemma~4.3.15]{Zhu23}.

Changes of the compatible Betti bases have the form
$(g,PgP^{-1})$ with $g\in\Gamma_0(N)$; on special matrices they
act by
$J\mapsto PgP^{-1}Jg^{-1}$ and $x\mapsto gxg^{-1}$.
Consequently the above construction is a bijection of the
analytic special-cycle groupoids, including their stabilizers.
It does not quotient by the additional relation $x\sim-x$.

Each special homomorphism contributes one reduced Cartier branch
on an orbifold chart in characteristic zero.  One way to see the
reducedness is the usual CM equation
$Cz^2-2bz-A/N=0$ for
$x=\left(\begin{smallmatrix}b&A/N\\C&-b\end{smallmatrix}\right)$.
Its discriminant is $-4m/N<0$, so its roots are simple.
At a CM point $z$, write $w=|\Gamma_z|$.
The stabilizer of a positive-norm $x$ fixing $z$ is exactly
$\Gamma_z$: rotations fixing $z$ commute with the associated
complex structure and hence with $x$.

It remains to distinguish divisor descent from stack degree.
The central $\{\pm1\}$ acts trivially on the local spatial
disk.  The effective ramification degree of the coarse map is
therefore $e=w/2$.  If $D_z$ denotes one reduced Cartier stack
branch and $[z]$ its coarse point, then locally
\[
       p^*[z]=eD_z,\qquad
       D_z\longleftrightarrow \frac{2}{w}[z]
       \quad\text{in rational Picard groups}.
\]
Summing this relation over all special-homomorphism branches gives
the asserted formula.  The homomorphisms $j$ and $-j$ are distinct branches
over the same base point; exactly the same multiplicity is
present in the two distinct vector orbits $x,-x$.
The stack integration weight $1/w$ is a different operation and
is not used in this rational Picard descent.
\end{proof}

\begin{remark}
The lattice map used in this proof is the moduli-equivariant map
$J\mapsto-P^{-1}J$.  An abstract determinant isometry differing
by conjugation by a determinant-minus-one matrix should not be
substituted without checking its action on the connected complex
domain.  The comparison also requires the vector divisor to use
the displayed $2/|\Gamma_x|$ convention.  Replacing it by an
arithmetic stack-degree pushforward changes the convention and
requires a new scalar comparison.
\end{remark}

\begin{corollary}[The projected scalar family at every level]
\label{cor:actual-MW-matching}
For the generic special divisors of \cite[Definition~7.1.1]{Zhu},
take their degree-zero classes based on $P_\infty$, and use the
modified Hodge and cuspidal generic coefficients at nonpositive
indices. Their admissible Mordell--Weil generating series is
\[
                  \widehat\phi_{\rm sc,MW}=R\widehat\Phi_{\rm MW}.
\]
In particular \eqref{eq:MW-matching} holds for every arithmetic
realization of these generic coefficients in the admissible Chow
space. No Green-function identification is required for this
projected statement.
\end{corollary}
\begin{proof}
Proposition~\ref{prop:generic-positive-matching} verifies the
positive-index hypothesis of Corollary~\ref{cor:MW-generic-matching}
with $a_N=1$. At negative indices the generic corrections are
cuspidal. The generic Hodge class at zero is rationally cuspidal
by \cite[Theorem~6.4.3]{Zhu}; hence it too has zero degree-zero
Jacobian class. The cited corollary proves the result.
\end{proof}

\subsection{Degree, vertical, and Hodge cusp pairings at every level}
\label{subsec:geometric-eisenstein-vanishing}

Explicit Eisenstein formulas for the degree, vertical intersection
functionals, and Hodge intersection functional of the scalar arithmetic
series imply their cusp orthogonality at arbitrary level.

\begin{lemma}[Eisenstein derivatives are orthogonal to cusps]
\label{lem:eisenstein-cusp-orthogonality}
Let $E_H(\tau,s)$ be a meromorphic weight-$k$ Eisenstein family
induced from constant Fourier terms at the cusps of $H$, with the
character used for the Petersson pairing. For a holomorphic cusp form
$g$ of that weight and character,
\[
 \langle E_H(\cdot,s),g\rangle_H=0
\]
as a meromorphic identity. Hence all regular values, Laurent
coefficients, and parameter derivatives have zero cusp pairing.
The same holds for a lower-level family pulled back to $H$.
\end{lemma}
\begin{proof}
In a right half-plane, unfold the Eisenstein sum at each inducing
cusp. The integral along its horizontal period is the complex
conjugate of the constant Fourier term of $g$ at that cusp, hence
zero. Cuspidal decay justifies this unfolding. Continue the zero
identity meromorphically. On compact parameter sets away from poles,
the Eisenstein series and its parameter derivatives have polynomial
growth, with possible logarithmic factors, at the cusps. The
exponential decay of $g$ permits differentiation under the integral
and extraction of Laurent coefficients. For a lower-level group
$H'\supset H$, unfold the finite covering first: its pairing with
$g$ is the pairing at $H'$ with the finite trace of $g$, which is
again cuspidal. The last assertion therefore uses the same
normalization of the measure.
\end{proof}

For the remainder of this subsection write $\widehat\phi_M$ for
Zhu's scalar arithmetic series at level $M$, and put
\[
 d_M(\tau)=\deg\widehat\phi_M(\tau),\qquad
 V_{p,a}^{M}(\tau)=
       [\widehat\phi_M(\tau),\widehat{\calX}_p^a(M)]_{\rm GS}.
\]
All comparisons use his actual boundary and zero coefficients. In
particular, $d_M$ includes the degree of the cusp correction; it is
not just the series of positive CM degrees.

\begin{proposition}[Three vanishing functionals]
\label{prop:degree-vertical-hodge-vanishing}
Let $g\in S_{3/2}(\Gamma_0(4N),\chi_N)$. Use an arithmetic
family realizing the normalized degree, vertical, and Hodge identities
stated in \cite{Zhu}. Assuming the arithmetic coefficient pairings
commute with this Petersson integral, one has
\begin{align}
 \langle d_N,g\rangle_{H_N}&=0,\label{eq:degree-cusp-zero}\\
 \langle V_{p,a}^N,g\rangle_{H_N}&=0
      &&(p\mid N),\label{eq:vertical-cusp-zero}\\
 \left\langle[\widehat\phi_N,\widehat\omega_N]_{\rm GS},
                    g\right\rangle_{H_N}&=0.
                    \label{eq:hodge-cusp-zero}
\end{align}
Here $a$ ranges over the components in \eqref{eq:special-fiber},
and $\widehat\omega_N$ is Zhu's metrized Hodge class. Consequently
the degree and every corrected vertical projection of the lift vanish.
For a family with a changed Green function or cusp correction, these
three scalar identities are hypotheses to be rechecked.
\end{proposition}
\begin{proof}
Use the notation $\mathcal E_M(\tau,s)$ and $E_M(\tau,s)$ for the
two genus-one Eisenstein normalizations in \cite{Zhu}; their
parameters need not agree with the genus-two parameter.
Proposition~8.1.1 there gives
\[
 d_M(\tau)=\frac{2}{\varphi(M)}\mathcal E_M(\tau,1).
\]
Lemma~\ref{lem:eisenstein-cusp-orthogonality} proves
\eqref{eq:degree-cusp-zero}.

For $n=v_p(N)\geq1$, Proposition~8.2.4 there gives
\begin{equation}\label{eq:outer-vertical-degree-series}
 V_{p,n}^{N}=\frac{\log p}{2}
     \left(d_N-(p-1)\sum_{i=1}^{\lfloor n/2\rfloor}
                                  d_{Np^{-2i}}\right).
\end{equation}
Atkin--Lehner invariance gives $V_{p,-n}^N=V_{p,n}^N$.
For an interior component, Corollary~8.2.3 reduces the pairing
from level $N$ to level $Np^{-2}$, with the same component index.
Repeat until it is an outer component, or until the prime no
longer divides the level. In the latter case its fiber is a
principal full fiber, and its pairing is a constant multiple of
$d_M$ by \eqref{eq:fiber-constant}. Thus every vertical pairing
is a linear combination of degree Eisenstein series at levels
$M=Np^{-2i}$. Since $N/M$ is a square, their classical characters
agree: $\chi_M=\chi_N$ on the common group. Their lower-level
pullbacks are orthogonal to $g$ by the last assertion of the lemma.
This proves \eqref{eq:vertical-cusp-zero}; multiplying by the inverse
intersection matrix in \eqref{eq:vertical-projection} proves the
vanishing of the corrected vertical coordinates as well.

For the Hodge pairing, use the full-series identity in
\cite[Proposition~8.3.2]{Zhu} and the divisor of $\Delta_N$:
\[
 \begin{split}
 12\varphi(N)[\widehat\phi_N,\widehat\omega_N]_{\rm GS}
  &=[\widehat\phi_N,\widehat\Delta_N]_{\rm GS}
       +\sum_{p\mid N}[\widehat\phi_N,(f_p(N),0)]_{\rm GS}\\
  &=12\mathcal E_N'(\tau,1)
       +\sum_{p\mid N}[\widehat\phi_N,(f_p(N),0)]_{\rm GS}.
 \end{split}
\]
Each $f_p(N)$ is a fixed linear combination of vertical components.
The derivative term is orthogonal to cusp forms by
Lemma~\ref{lem:eisenstein-cusp-orthogonality}; the vertical terms
are orthogonal by the preceding calculation. This proves
\eqref{eq:hodge-cusp-zero}. The local coefficient
$\beta'_{p,t}(0)$ in Zhu's Corollary~8.3.3 depends on $t$;
it cannot be replaced by a single scalar multiplying the full
Eisenstein series.
\end{proof}

\begin{corollary}[The remaining projection is the smooth one]
\label{cor:MW-from-smooth}
Use the admissible decomposition \eqref{eq:arith-decomposition}
with the degree-one class based on $P_\infty$, scaled by its
stack degree when necessary. Under the interchanges
above, if the smooth component of the lift vanishes, then the full
lift belongs to $\widetilde{\MW}_{\C}$.
\end{corollary}
\begin{proof}
The preceding proposition removes its degree and corrected vertical
coordinates. If the smooth coordinate also vanishes, write the lift
as $x+c\one$, with $x\in\widetilde{\MW}_{\C}$. The generic-fiber
divisor relation for $\Delta_N$ in \cite[Theorem~6.4.3]{Zhu} expresses
the Hodge line as a nonzero multiple of $P_\infty$ in the rational
Picard group; all other divisor terms are vertical. Hence
$\widehat\omega_N$ has zero Mordell--Weil coordinate in this
decomposition. Admissibility gives
$H_{\rm GS}(x,\widehat\omega_N)=0$: the vertical coordinates
are orthogonal to $x$, and any residual metric function pairs
trivially with $x$ because its admissible curvature is zero. This
uses the same extended Green-current pairing for the Hodge metric. Its pairing with the lift is
therefore $c\deg(\omega_N)/2$. The Hodge pairing vanishes by
\eqref{eq:hodge-cusp-zero}, and $\deg(\omega_N)>0$, so $c=0$.
\end{proof}

\subsection{A direct Kudla--Millson cusp orthogonality lemma}

The argument uses the following quadratic and differential conventions.
Let $(V,Q)$ have signature $(1,2)$,
put $B(x,y)=Q(x+y)-Q(x)-Q(y)$, and choose the positive vector
$h_z\in V(\mathbb R)$ associated with $z\in\mathfrak H$, normalized
by $Q(h_z)=1$.  Write
\[
 \Omega(z)=\frac{i}{2}\frac{dz\wedge d\bar z}{(\operatorname{Im}z)^2},
 \qquad
 \mathcal R(x,z)=\frac12 B(x,h_z)^2-B(x,x).
\]
Thus $\mathcal R(x,z)\geq0$, and
$\mathcal R(x,z)+2Q(x)=\tfrac12B(x,h_z)^2$.
For an even lattice $L\subset V$ the two theta kernels in these
conventions are
\begin{align}
 \Theta_S(\tau,z)
  &=v\sum_{\mu\in L'/L}\sum_{x\in L+\mu}
       e^{-2\pi v\mathcal R(x,z)}q^{Q(x)}e_\mu,
       \label{eq:KM-proof-Siegel}\\
 \Theta_{\rm KM}(\tau,z)
  &=\sum_{\mu\in L'/L}\sum_{x\in L+\mu}
       \left(vB(x,h_z)^2-\frac1{2\pi}\right)
       e^{-2\pi v\mathcal R(x,z)}q^{Q(x)}e_\mu\,\Omega(z),
       \label{eq:KM-proof-kernel}
\end{align}
where $\tau=u+iv$ and $q=e^{2\pi i\tau}$.  These are the untwisted
conventions in \cite[Sections~3.1--3.2]{Alfes}.  Their modular weights
are $-1/2$ and $3/2$, respectively, with the same unitary Weil
representation.  The arguments below also apply after a fixed finite
Schwartz combination or a compatible scalar component has been taken.

\begin{lemma}[Exact raising identity and pointwise cusp orthogonality]
\label{lem:KM-direct-cusp}
With $R_k=2i\partial_\tau+k/v$ one has
\begin{equation}\label{eq:KM-exact-raising}
 \Theta_{\rm KM}(\tau,z)
   =-\frac1\pi R_{-1/2,\tau}\Theta_S(\tau,z)\,\Omega(z).
\end{equation}
Consequently, for every holomorphic cusp form $g$ of weight $3/2$
for the matching representation,
\begin{equation}\label{eq:KM-pointwise-cusp-zero}
 \int_{\Gamma\backslash\mathfrak H}
      (\Theta_{\rm KM}(\tau,z),g(\tau))
      v^{3/2}\,\frac{du\,dv}{v^2}=0
 \qquad(z\in\mathfrak H).
\end{equation}
The integral converges absolutely and locally uniformly in $z$,
with all derivatives in $z$ on compact subsets.  There is no
square-freeness condition on $L$ or its level.
\end{lemma}

\begin{proof}
For fixed $x,z$ put $A_x=e^{-2\pi v\mathcal R(x,z)}q^{Q(x)}$.
Since $2i\partial_\tau v=1$ and
$2i\partial_\tau q^{Q(x)}=-4\pi Q(x)q^{Q(x)}$, one obtains
\[
 R_{-1/2}(vA_x)
 =\left\{\frac12-2\pi v\bigl(\mathcal R(x,z)+2Q(x)\bigr)\right\}A_x
 =\left\{\frac12-\pi vB(x,h_z)^2\right\}A_x.
\]
This gives \eqref{eq:KM-exact-raising}.  Termwise differentiation is
valid because
$\mathcal R(x,z)+Q(x)$ is a positive definite quadratic form in
$x$, uniformly comparable with a fixed positive definite form when
$z$ ranges over a compact set.

For integration by parts, put $\kappa=-1/2$ and
$h(\tau)=\Theta_S(\tau,z)$, and use a Hermitian product linear
in the first variable. The boundary term is included below.  Holomorphy of $g$
gives, component by component,
\begin{equation}\label{eq:KM-Stokes-exact}
 (R_\kappa h,g)v^\kappa\,du\wedge dv
   =-d\bigl((h,g)v^\kappa d\bar\tau\bigr).
\end{equation}
The one-form on the right is invariant: the weights of $h$ and $g$
are $\kappa$ and $\kappa+2$, and the representation is unitary.
Integrate over a truncated fundamental domain.  Paired sides cancel.
For fixed $z$ the Gaussian estimate gives $h=O_z(v)$ as $v\to\infty$,
whereas $g=O(e^{-cv})$ for some $c>0$.  Thus the remaining horizontal
boundary integral tends to zero.  On a finite-index scalar modular
group the same estimates hold at every cusp after applying its
scaling matrix; there are only finitely many cusps.  They also hold
locally uniformly with all $z$ derivatives.  Equations
\eqref{eq:KM-Stokes-exact} and \eqref{eq:KM-exact-raising} prove the
assertion.  The explicit kernels establish the pointwise assertion.
The conclusion of \cite[Theorem~5.1]{Alfes} concerns orthogonality
of certain integrated lifts and is not itself a pointwise assertion.
\end{proof}


\subsection{Absolute joint integrability of the standard kernel}
\label{subsec:KM-joint-integrability}

The absolute value below is taken \emph{after} summing the
Kudla--Millson kernel in the lattice variable. We do not claim
integrability of the sum of the absolute values of the individual
original lattice summands. Poisson cancellation in one isotropic
direction is essential.

Put
\[
 V=\{\lambda\in M_2(\mathbb Q):\operatorname{tr}\lambda=0\},
 \qquad Q(\lambda)=N\det\lambda,
 \qquad
 L_N=\left\{\begin{pmatrix}b&-a/N\\c&-b\end{pmatrix}:
                         a,b,c\in\mathbb Z\right\}.
\]
Let $B$ be the polarization of $Q$, let $z=x+iy$ be the
geometric variable, and let $\tau=u+iv$ be the modular variable.
Use the positive vector and standard majorant
\[
 h_z=\frac1{\sqrt N y}
       \begin{pmatrix}-x&|z|^2\\-1&x\end{pmatrix},
 \qquad
 \mathcal R(\lambda,z)=\tfrac12B(\lambda,h_z)^2-2Q(\lambda),
 \qquad \Omega=\frac{dx\wedge dy}{y^2}.
\]
For a rational lattice coset $\Lambda+\mu\subset V$, define
the scalar coefficient of its standard KM kernel by
\begin{equation}\label{eq:KM-bound-component}
 \vartheta_{\Lambda+\mu}(\tau,z)=
 \sum_{\lambda\in\Lambda+\mu}
   \left(vB(\lambda,h_z)^2-\frac1{2\pi}\right)
   e^{-2\pi v\mathcal R(\lambda,z)}e^{2\pi i\tau Q(\lambda)}.
\end{equation}
For the Weil kernel associated with $L_N$, the components are
$\vartheta_{L_N+\mu}\Omega$, with $\mu\in L_N'/L_N$.

\begin{lemma}[Uniform geometric $L^1$ bound]
\label{lem:KM-uniform-geometric-L1}
Fix $N$ and $v_0>0$. For the complete standard KM kernel
associated with $L_N$, including the zero lattice vector,
there is a constant $C_{N,v_0}$ such that
\begin{equation}\label{eq:KM-uniform-geometric-L1}
 \int_{\Gamma_0(N)\backslash\mathbb H}
       \|\Theta_{\rm KM}(u+iv,z)\|
       \leq C_{N,v_0}
       \qquad(u\in\mathbb R, v\geq v_0).
\end{equation}
The norm of a form means the total variation density, using
any fixed norm on the finite-dimensional coefficient space.
The same conclusion holds for a fixed finite rational
Schwartz combination and after passing to a finite-index
geometric subgroup. For fixed $\tau$, the kernel extends
smoothly and without atoms to the geometric cusps.
\end{lemma}

\begin{proof}
First consider the geometric cusp at infinity and one
rectangular rational coset in the coordinates $(a,b,c)$:
\[
       a\in\alpha+D\mathbb Z,\quad
       b\in\beta+D\mathbb Z,\quad
       c\in\gamma+D\mathbb Z,
 \qquad \alpha,\beta,\gamma\in D^{-1}\mathbb Z,
\]
where $D$ is a positive integer. Put
\[
       b'=b-cx,\qquad a_0=2Nbx-Ncx^2,\qquad t=a-a_0.
\]
Direct expansion gives
\begin{align}\label{eq:KM-completed-square}
 Q(\lambda)&=ct-N(b')^2,\\
 B(\lambda,h_z)&=-\frac{t+Ncy^2}{\sqrt N y},\notag\\
 \mathcal R(\lambda,z)+Q(\lambda)
       &=\frac{t^2}{2Ny^2}+N(b')^2+\frac{Nc^2y^2}{2}.
       \notag
\end{align}
The original summand is thus a polynomial times a positive
definite Gaussian at each fixed $\tau,z$, so its lattice sum
converges absolutely there.

Poisson summation in $a$ has frequencies $m\in D^{-1}\mathbb Z$.
To compute it, put $A=v/(Ny^2)$ and $K=m-cu$. The Fourier
transform of the $t$-dependent factor is obtained from
\[
 \int_{\mathbb R}e^{-\pi At^2}e^{-2\pi iKt}\,dt
       =A^{-1/2}e^{-\pi K^2/A}.
\]
Applying the polynomial $A(t+Ncy^2)^2-1/(2\pi)$ under
this integral gives
\begin{align*}
 &\left(-\frac{K^2}{A}-2iNcy^2K+AN^2c^2y^4\right)
       A^{-1/2}e^{-\pi K^2/A}\\
 &\hspace{12mm}=-\frac{N^{3/2}y^3}{v^{3/2}}
       (m-c\overline\tau)^2 e^{-\pi Ny^2(m-cu)^2/v}.
\end{align*}
The $1/(2\pi)$ term cancels exactly. Combining the remaining
$c$-Gaussian from \eqref{eq:KM-completed-square}, the absolute
value of the Poisson expansion is bounded by
\begin{equation}\label{eq:KM-Poisson-majorant}
 \frac{N^{3/2}y^3}{Dv^{3/2}}
 \sum_{\substack{c\in\gamma+D\mathbb Z\,,\ m\in D^{-1}\mathbb Z\\
                  (m,c)\ne(0,0)}}
  |m-c\tau|^2 e^{-\pi Ny^2|m-c\tau|^2/v}
       \sum_{b\in\beta+D\mathbb Z}
            e^{-2\pi Nv(b-cx)^2}.
\end{equation}
The omitted factors are complex numbers of modulus one.
This transformed series also converges absolutely at fixed
$\tau,z$. The zero pair $(m,c)=(0,0)$ is absent because its
polynomial is zero, even when the original rectangular coset
contains the zero lattice vector.

The elementary shifted Gaussian estimate
\[
 \sum_{b\in\beta+D\mathbb Z}e^{-2\pi Nv(b-cx)^2}
       \leq C_{N,D,v_0}\qquad(v\geq v_0)
\]
is uniform in $x,c$. Integrate
\eqref{eq:KM-Poisson-majorant} over a cusp strip
$0\leq x\leq w$, $y\geq Y>0$, against $dx\,dy/y^2$.
Tonelli applies to this nonnegative majorant. For every
$\omega=m-c\tau\ne0$,
\[
 \int_Y^\infty y|\omega|^2
          e^{-\pi Ny^2|\omega|^2/v}\,dy
       =\frac{v}{2\pi N}
          e^{-\pi NY^2|\omega|^2/v}.
\]
The resulting integral is at most
\begin{equation}\label{eq:KM-integrated-majorant}
 C\,v^{-1/2}
   \sum_{m,c\in D^{-1}\mathbb Z}
       \exp\left(-\pi NY^2
           \left(\frac{(m-cu)^2}{v}+c^2v\right)\right).
\end{equation}
Adding the zero pair at this stage only enlarges the bound.
Write $m=j/D$, $c=k/D$ and $a=\pi NY^2/D^2>0$.
Uniformly in real shifts,
\[
 \sum_{j\in\mathbb Z}e^{-a(j-ku)^2/v}
       \leq C_a(1+\sqrt v),
 \qquad
 \sum_{k\in\mathbb Z}e^{-ak^2v}\leq C_{a,v_0}.
\]
Hence \eqref{eq:KM-integrated-majorant} is bounded by
$C v^{-1/2}(1+\sqrt v)$, which is uniformly bounded for
$v\geq v_0$.

An arbitrary rational lattice coset is a finite disjoint
union of such rectangular cosets: choose $D$ so that
$DL_N\subset\Lambda$ and
$\Lambda+\mu\subset D^{-1}L_N$, then use its cosets modulo
$DL_N$. Therefore the cusp-strip estimate holds for every
fixed rational lattice coset.

At a different geometric cusp, choose
$\sigma\in\mathrm{SL}_2(\mathbb Z)$ taking infinity to that
cusp. The identities for $h_z$ and invariance of $Q$ show
that pulling the kernel to the cusp coordinate conjugates
the lattice to $\sigma^{-1}(L_N+\mu)\sigma$.
It is a rational lattice coset, so the preceding argument
applies. There are finitely many cusps, each with finite
width. On the remaining compact geometric region,
$\mathcal R(\lambda,z)+Q(\lambda)$ is uniformly positive
definite. The polynomial-Gaussian bound there is uniform
for $v\geq v_0$. Summing the finitely many regions proves
\eqref{eq:KM-uniform-geometric-L1}.

For fixed $\tau$, the nonzero pairs in
\eqref{eq:KM-Poisson-majorant} have
$|m-c\tau|\geq\min(1,v)/D$. Thus at a geometric cusp
the coefficient and all its geometric derivatives are
bounded by a polynomial in $y$ times $e^{-c_\tau y^2}$.
In a local parameter $q_z=e^{2\pi iz/w}$ this extends
smoothly by zero, even after the factors in $\Omega$ are
included. In particular the extension is an atom-free
locally integrable density, proving the last assertion.
\end{proof}

\begin{corollary}[The compact weak Fubini condition]
\label{cor:KM-compact-Fubini-proved}
Let $H$ be a finite-index scalar modular group for a
compatible component of the standard kernel, and let
$g\in S_{3/2}(H,\chi)$ be holomorphic and cuspidal at
every cusp of $H$. For every smooth function $a$ on the
compact geometric curve, one has
\begin{equation}\label{eq:KM-joint-absolute-integral}
 \int_{H\backslash\mathbb H}\int_X
     \left|a(z)\bigl(\Theta_{{\rm KM},0}(\tau,z),g(\tau)\bigr)\right|
           v^{3/2}\,d\mu(\tau)<\infty.
\end{equation}
Consequently its two integrations commute. Together with
Lemma~\ref{lem:KM-direct-cusp}, this proves that the
cuspidally paired standard curvature kernel is zero
as a current on the entire compact geometric curve.
\end{corollary}

\begin{proof}
The function $a$ is bounded. A fundamental domain for
$H$ is a finite union of translates of the standard
$\mathrm{SL}_2(\mathbb Z)$ fundamental domain. Transform
each integral back to that domain. Modular covariance
replaces the kernel by a fixed finite combination of
its Weil components and replaces $g$ by the slash
transform at the corresponding cusp. The absolute
integrand with the Petersson measure is unchanged by
this transformation. Each such slash transform of
$g$ is $O(e^{-c v})$ as $v\to\infty$, with some $c>0$.
Lemma~\ref{lem:KM-uniform-geometric-L1} therefore bounds
the tail of the double integral by a constant times
\[
       \int_{v_0}^\infty e^{-cv}v^{-1/2}\,dv<\infty.
\]
The remaining modular region is compact, and the same
lemma applies there. This proves
\eqref{eq:KM-joint-absolute-integral}; ordinary Fubini
now applies. The pointwise inner integral vanishes
off the finite cusp set by the exact raising identity.
Since the compact extension is a density, its integral
against $a$ is zero. Thus the paired current vanishes.
\end{proof}

\paragraph{Scope of the conclusion.}
This proves the atom-free extension and compact weak
Fubini requirements for the standard KM theta kernel
at arbitrary $N$. If the actual arithmetic curvature
series is identified with this standard kernel and
curvature commutes with its arithmetic Petersson
integral, the compact-zero-curvature proposition then
applies without assuming its Fubini condition.
The arithmetic cusp divisor and zero Green conventions are
compared separately in
Sections~\ref{subsec:cusp-width-correction} and
\ref{subsec:negative-cusp-correction}.
Corollary~\ref{cor:completion-full-lift-MW} justifies passage
through the arithmetic Petersson integral.


\subsection{The smooth projection and the compact Green operator}
\label{subsec:KM-Green-operator}

Let $X$ be the compact complex modular curve, and choose a
smooth positive form $\nu$ with $\int_X\nu=1$. The statements
below can also be made in its finite orbifold charts. Fix the
normalized compact Green kernel $G_\nu(z,w)$ with
\[
 dd_z^cG_\nu(z,w)=\delta_w-\nu(z),
 \qquad \int_XG_\nu(z,w)\nu(z)=0.
\]
Its logarithmic diagonal singularity and smooth regular part
on compact charts imply
\begin{equation}\label{eq:KM-Green-L1-bound}
 C_\nu:=\sup_{w\in X}\int_X|G_\nu(z,w)|\nu(z)<\infty.
\end{equation}
Indeed a finite family of local charts gives a uniform bound
for the integral of $|\log|z-w||$ near the diagonal, and the
kernel is bounded on the compact complement. The same proof
works on finite quotient charts.

For a complex finite measure $\lambda$ on $X$ put
\[
 \mathcal G_\nu\lambda(z)=\int_XG_\nu(z,w)\,d\lambda(w).
\]
Tonelli and \eqref{eq:KM-Green-L1-bound} give
\begin{equation}\label{eq:KM-Green-measure-L1}
 \|\mathcal G_\nu\lambda\|_{L^1(X,\nu)}
       \leq C_\nu\|\lambda\|_{\rm TV},
 \qquad
 dd^c\mathcal G_\nu\lambda
       =\lambda-\lambda(X)\nu.
\end{equation}
The potential has mean zero. It is the unique mean-zero
distribution satisfying this last equation, since the
difference of two such potentials is a distributionally
harmonic function on compact connected $X$, hence is smooth
and constant.

\begin{corollary}[Weak interchange for the smooth coordinate]
\label{cor:KM-Green-interchange}
Let $\Theta(\tau)$ be the compact extension of a compatible
component of the standard KM kernel. Put
$d(\tau)=\int_X\Theta(\tau)$. Suppose the actual mean-zero
smooth coordinate $f(\tau,\cdot)$ of an arithmetic family
satisfies
\begin{equation}\label{eq:KM-actual-smooth-curvature}
       dd^cf(\tau,\cdot)=\Theta(\tau)-d(\tau)\nu
\end{equation}
on the entire compact curve, with the same zero and boundary
conventions. Then
$f(\tau)=\mathcal G_\nu\Theta(\tau)$ as distributions.
For every compatible holomorphic cusp form $g$,
\[
 \int_{H\backslash\mathbb H}\int_X
       |f(\tau,z)\overline{g(\tau)}|\nu(z)
       v^{3/2}\,d\mu(\tau)<\infty.
\]
Its Petersson integral therefore exists in $L^1(X,\nu)$,
has mean zero, and commutes with $dd^c$ weakly. Moreover
\[
       \langle d,g\rangle_H=0,
       \qquad \langle f,g\rangle_H=0
                  \quad\text{as a distribution on }X.
\]
This establishes the degree and smooth-coordinate
vanishings under the actual curvature comparison. It does
not assert existence of the Petersson integral of the
other arithmetic coordinates.
\end{corollary}

\begin{proof}
Uniqueness of the normalized potential gives
$f=\mathcal G_\nu\Theta$. The uniform total variation
bound in Lemma~\ref{lem:KM-uniform-geometric-L1}, followed
by \eqref{eq:KM-Green-measure-L1}, gives a uniform
$L^1(X,\nu)$ bound for $f(u+iv)$ for $v\geq v_0$.
It also gives $|d(u+iv)|\leq\|\Theta(u+iv)\|_{\rm TV}$.
After the finitely many modular slash transforms, the
same bounds apply at every cusp of $H$. Cuspidal decay
of $g$ proves the displayed absolute integral, just as
in Corollary~\ref{cor:KM-compact-Fubini-proved}.

Let $a$ be a smooth test function on $X$. The absolute
integral permits interchange of the Petersson integral
with testing $f$ against $dd^ca$. Equations
\eqref{eq:KM-actual-smooth-curvature} and
\eqref{eq:KM-Green-measure-L1} then give
\[
 dd^c\langle f,g\rangle_H
       =\langle\Theta,g\rangle_H
                       -\langle d,g\rangle_H\nu.
\]
Corollary~\ref{cor:KM-compact-Fubini-proved} makes the
first term zero as a compact current. Testing that
same corollary against $1$ shows
$\langle d,g\rangle_H=0$. Hence the potential integral
is distributionally harmonic. It has mean zero by
absolute Fubini, so elliptic regularity and compactness
make it identically zero.
\end{proof}

\subsection{The compact-current consequence, with the zero term included}

Applying the analytic lemma to an arithmetic family requires the
normalization and convergence assumptions below, but no twice
differentiable extension of the Green functions to the cusps.
For an arithmetic divisor write
\[
 \operatorname{curv}(D,G)=dd^cG+\delta_{D(\mathbb C)}.
\]
All such expressions in this paragraph are currents on the compact
complex curve $X$, including its cusps.

\begin{proposition}[Zero curvature of the cuspidally paired family]
\label{prop:KM-compact-zero-curvature}
Let $\widehat\phi(\tau)$ be an arithmetic theta family.  Assume that
its actual compact Green identities, including its prescribed cusp
divisors and its zero coefficient, identify
\begin{equation}\label{eq:KM-compact-curvature-match}
 \operatorname{curv}\widehat\phi(\tau)
    =\Theta_{{\rm KM},0}(\tau,\cdot)
\end{equation}
as currents on $X$, where the right side is the matching scalar
component of \eqref{eq:KM-proof-kernel}.  Assume also that the
arithmetic Petersson integral against a cusp form $g$ exists and
that curvature commutes with this integral in the weak sense after
testing against smooth functions on $X$. The compact Fubini
interchange for this standard kernel is supplied by
Corollary~\ref{cor:KM-compact-Fubini-proved}. Then
\[
 \operatorname{curv}\langle\widehat\phi,g\rangle=0
 \quad\text{as a current on }X.
\]
In an admissible arithmetic Chow decomposition with curvature-zero
Mordell--Weil and vertical terms, this simultaneously removes the
degree component and the mean-zero smooth component of the lift.
\end{proposition}

\begin{proof}
Test the weak curvature identity against a smooth function $a$.
The assumed curvature interchange and the proved compact Fubini
identity reduce its value to the pointwise cusp pairing in
Lemma~\ref{lem:KM-direct-cusp}, integrated against $a$ on $X$.
That pairing is zero. Lemma~\ref{lem:KM-uniform-geometric-L1}
also proves that the compact standard kernel is atom-free, so no
cusp Dirac term is omitted. This proves zero curvature as a compact
current.

Write the degree and smooth curvature contributions of the lift as
$d\nu+dd^cf$, where $\int_X\nu=1$ and $f$ has mean zero.
Integration of the zero current gives $d=0$.  It follows that
$dd^cf=0$ distributionally on compact $X$.  Weyl's lemma makes $f$
a smooth harmonic function.  Compactness and the mean-zero
normalization then give $f=0$.
\end{proof}

\paragraph{Checking the missing lattice vector at zero.}
In \eqref{eq:KM-proof-kernel} the vector $x=0$ contributes exactly
$-\Omega/(2\pi)$ in the zero component.  In Zhu's arithmetic
normalization the zero coefficient contains
$H=-\widehat\omega_N-W_N^*\widehat\omega_N$ and the
term $-(0,\log v)$.  Since $\log v$ is constant in the geometric
variable and the Petersson metric satisfies
\[
 c_1(\widehat\omega_N)=\frac{\Omega}{4\pi},
 \qquad W_N^*\Omega=\Omega,
\]
the extra zero-coefficient curvature is
$c_1(H)=-\Omega/(2\pi)$, precisely the missing
$x=0$ term.  This Hodge contribution completes the matching of the nonzero
lattice vectors. Omitting it would invalidate the raising identity
for the arithmetic curvature series.
Zhu's Theorem~7.1.2 states a compact Green equation with the
prescribed cusp divisor. Its application to the present standard
family requires the normalization checks in the next paragraph
and the cusp expansions in
Sections~\ref{subsec:cusp-width-correction} and
\ref{subsec:negative-cusp-correction}. The log--log terms are
locally integrable currents; their cusp flux is
\(O(1/|\log\varepsilon|)\), so they contribute no further Dirac atom.

\paragraph{Green-function normalization.}
The displayed notation in \cite[Section~7.1.2]{Zhu} is not by itself
a consistent substitute for the conventions above.  Section~2.1
defines $(x,y)$ by polarization of $Q$, while the later coefficient
condition is printed $(x,x)=t$; the Fourier label in the present
paper is $Q(x)=t$.  Equation~(50) also prints an exponential without
the modular imaginary part, despite the scaled Green function in
the preceding display.  With that same polarized bilinear form,
the printed expression
$-|B(x,c(z))|^2/B(c(z),\overline{c(z)})$ equals
$\mathcal R(x,z)/2$, whereas the standard normalization used in
\eqref{eq:KM-proof-kernel} is $\mathcal R$.
Matching the actual Green normalization to
\eqref{eq:KM-proof-kernel} therefore requires both the polarization
factor and the factor of $v$ in its Gaussian.  For instance, at $N=1$, $z=i$, and
$x=\operatorname{diag}(a,-a)$, the standard and printed expressions
give $2a^2$ and $a^2$, respectively. A
correct statement of \eqref{eq:KM-compact-curvature-match} must
record this convention; a literal uncorrected use of those printed
symbols does not prove it.
The standard normalization can instead be read without this ambiguity
in \cite[(2.3), (2.5)--(2.8), (5.3)--(5.4)]{DY1}.
Section~5 of that paper explicitly allows an arbitrary even lattice,
and Corollary~5.5 is the compact Green current identity; the
square-free restriction on its arithmetic integral model is imposed
only subsequently.  This provides an all-level analytic reference for
the standard Kudla Green construction.  Its coefficients of geometric
cusp divisors, which depend on cusp widths and the lattice at each
cusp, are computed for the standard arithmetic family in
Sections~\ref{subsec:cusp-width-correction} and
\ref{subsec:negative-cusp-correction}.

\section{The projected arithmetic height}
\label{sec:projected-height}

The trace calculation transfers the corrected Bruinier--Yang height
theorem to the degree-zero generic theta lift, under the stated
Jacobi newform hypothesis.

\subsection{Transfer from an existing projected vector-valued formula}

A proved vector-valued formula gives a second route to the scalar
Mordell--Weil formula. Its formal transfer preserves the norm factor
from scalar--vector induction and requires no recomputation of the
doubling integrals.

Let $\widehat\Phi_{\rm MW}(z)$ be a $W$-valued modular form with
coefficients in $\widetilde{\rm MW}_\mathbb C$, using the same $\rho$
and Petersson conventions as above.  Write
$\widehat\phi_{0,\rm MW}=R\widehat\Phi_{\rm MW}$ and define
\begin{align}
 \widehat\Theta_{\rm MW}(F)
   &=\int_{\Gamma\backslash\mathbb H}
          (\widehat\Phi_{\rm MW}(z),F(z))_W y_z^k\,d\mu(z),
       \label{eq:red-vector-MW-lift}\\
 \widehat\theta_{0,\rm MW}(g)
   &=\int_{H\backslash\mathbb H}
          \widehat\phi_{0,\rm MW}(z)\overline{g(z)}y_z^k\,d\mu(z).
       \label{eq:red-scalar-MW-lift}
\end{align}
Here $(\widehat\Phi_{\rm MW},F)_W$ denotes contraction in the $W$
factor; it leaves the arithmetic coefficient unchanged.

\begin{proposition}[Transfer of a projected vector formula]
\label{prop:red-existing-vector-transfer}
Assume that a vector-valued arithmetic inner product theorem, in the
normalizations fixed above, gives
\begin{equation}\label{eq:red-vector-theorem-input}
 H_{\rm NT}\bigl(\widehat\Theta_{\rm MW}(F),
                    \widehat\Theta_{\rm MW}(F)\bigr)
   =C_N^{\rm vec}\|F\|_\Gamma^2L'(G,1)
\end{equation}
for the specified nonzero form $F=Tg$, with its specified Shimura
correspondent $G$.  Then
\begin{equation}\label{eq:red-scalar-transfer}
 H_{\rm NT}\bigl(\widehat\theta_{0,\rm MW}(g),
                    \widehat\theta_{0,\rm MW}(g)\bigr)
   =C_N^{\rm vec}\|Tg\|_\Gamma^2L'(G,1).
\end{equation}
Equivalently, for $g\ne0$ the scalar formula is
\begin{equation}\label{eq:red-scalar-ratio-formula}
 H_{\rm NT}\bigl(\widehat\theta_{0,\rm MW}(g),
                    \widehat\theta_{0,\rm MW}(g)\bigr)
   =C_N^{\rm vec}r(g)\|g\|_H^2L'(G,1),
 \qquad r(g)=\frac{\|Tg\|_\Gamma^2}{\|g\|_H^2}.
\end{equation}
\end{proposition}

\begin{proof}
Apply the conjugate of \eqref{eq:red-adjunction} coefficientwise to
$\widehat\Phi_{\rm MW}$. This can be done by applying a linear
functional to the arithmetic coefficients, using adjunction, and
varying the functional. Cuspidal decay justifies the unfolding.  This proves
\begin{equation}\label{eq:red-MW-adjunction}
 \widehat\theta_{0,\rm MW}(g)=\widehat\Theta_{\rm MW}(Tg).
\end{equation}
Substitute $F=Tg$ into \eqref{eq:red-vector-theorem-input}.
\end{proof}

\begin{corollary}[Arithmetic coefficient rescaling]
\label{cor:arithmetic-cycle-rescaling}
Suppose a scalar arithmetic family satisfies
\[
 \operatorname{pr}_{\rm MW}\widehat\phi_{\rm sc}
       =a_N R\widehat\Phi_{\rm MW}
\]
with a specified $a_N\in\mathbb C$. If its full cusp lift equals
its Mordell--Weil projection and the vector formula
\eqref{eq:red-vector-theorem-input} applies to $F=Tg$, then
\[
 H_{\rm NT}(\widehat\theta_{\rm sc}(g),\widehat\theta_{\rm sc}(g))
   =|a_N|^2 C_N^{\rm vec}\|Tg\|_\Gamma^2 L'(G,1).
\]
\end{corollary}
\begin{proof}
Adjunction gives
$\widehat\theta_{\rm sc}(g)=a_N\widehat\Theta_{\rm MW}(Tg)$.
Apply the Hermitian height pairing and the cited vector formula.
\end{proof}

For the generic special-isogeny cycles and the vector divisor
convention used here, Proposition~\ref{prop:generic-positive-matching}
and Corollary~\ref{cor:actual-MW-matching} prove $a_N=1$.
The weights $2/|\operatorname{Aut}|$ in the coarse Heegner divisor
and $1/|\operatorname{Aut}|$ in the stack arithmetic degree refer
to different operations. Rational Picard descent uses the effective
stabilizer modulo $\{\pm1\}$. The numerical arithmetic pairing
still requires its own stack and height normalization.

\begin{remark}[Precise use of a cited theorem]
The projected vector-valued formula stated in
\cite[Theorem~0.11 and formula~(0.12)]{Du} is a possible input to
\eqref{eq:red-vector-theorem-input}.  Its proof in Section~9.2 is formulated for a newform in a fixed
Shimura packet, so that is the input scope used here. To apply it,
its hypotheses on $F$,
its pairing convention, its Petersson measure, its definition of the
arithmetic generating series, and its constant $C_N^{\rm vec}$ must
agree with those here.  The transfer assumes the cited theorem in these conventions.
Its general-level assertion concerns a Mordell--Weil projection,
and does not by itself extend to the full arithmetic theta lift.  The latter requires the additional
identity
\[
 \widehat\theta_0(g)=\widehat\theta_{0,\rm MW}(g),
\]
which amounts to the vanishing of the non-Mordell--Weil components after
cuspidal pairing.  Their modularity alone does not imply that vanishing.
Nor does the transfer prove $r(g)=1$ or $r(g)=\nu/2$.
\end{remark}

\begin{proposition}[A packet height determined by one Heegner coefficient]
\label{prop:packet-one-coefficient-height}
Fix a weight-two newform $G$ and assume the actual
Gross--Kohnen--Zagier factorization in its chosen vector new line,
\[
 \widehat\Phi_{\rm MW}^{G}(\tau)=F_G(\tau)\otimes P_G,
 \qquad F_G(\tau)=\sum_{m,\mu}b_G(m,\mu)q^m e_\mu.
\]
Here $P_G$ is an admissible Mordell--Weil class; this factorization
is the theorem of \cite[Theorem~7.7]{BO} under its stated newform
hypotheses, with the untwisted formulation provided by
Gross--Kohnen--Zagier.  Suppose the scalar arithmetic family's
Mordell--Weil projection equals $R\widehat\Phi_{\rm MW}$ with
the generic divisor normalization above.  For every scalar cusp
form $f$ for which the lift converges, one then has
\[
 \operatorname{pr}_G\widehat\theta_{\rm MW,sc}(f)
       =\langle RF_G,f\rangle_H P_G.
\]
Consequently, whenever $b_G(m_0,\mu_0)\ne0$,
\begin{equation}\label{eq:packet-height-single-Heegner}
 H_{\rm NT}\bigl(\operatorname{pr}_G\widehat\theta_{\rm MW,sc}(f),
                \operatorname{pr}_G\widehat\theta_{\rm MW,sc}(f)\bigr)
 =\frac{|\langle RF_G,f\rangle_H|^2}
        {|b_G(m_0,\mu_0)|^2}
     H_{\rm NT}\bigl(y^G(m_0,\mu_0),y^G(m_0,\mu_0)\bigr).
\end{equation}
If in addition $Tf\in\mathbb C F_G$, the fraction in front of
the Heegner height can be written
\[
       \frac{\|Tf\|_\Gamma^2\,\|F_G\|_\Gamma^2}
            {|b_G(m_0,\mu_0)|^2}.
\]
\end{proposition}
\begin{proof}
Pair the factorized generating series with $Tf$ and use scalar--vector
adjunction.  This gives the first identity.  Its coefficient at
$(m_0,\mu_0)$ is
$y^G(m_0,\mu_0)=b_G(m_0,\mu_0)P_G$.
The Hermitian scaling law for the height proves
\eqref{eq:packet-height-single-Heegner}.
Finally, if $Tf=cF_G$, adjunction gives
$|\langle RF_G,f\rangle_H|^2
 =|\langle F_G,Tf\rangle_\Gamma|^2
 =\|Tf\|_\Gamma^2\|F_G\|_\Gamma^2$.
\end{proof}

The identity combines the cited factorization with the generic-cycle
comparison. The corrected Bruinier--Yang theorem evaluates the height
constant in this vector newform normalization.

\subsection{The Bruinier--Yang height formula at arbitrary level}
\label{subsec:BY-height}

Fix the N\'eron--Tate pairing on $J_0(N)$ in the convention of
\cite[Theorem~7.7]{BYCM}.  In particular its comparison with the
arithmetic intersection pairing remains the separately specified
constant $\eta$ in \eqref{eq:red-FH-normalization}.
The Petersson norm in the following result is the ordinary
unnormalized integral over the standard
$\mathrm{SL}_2(\mathbb Z)$ fundamental domain, with the
orthonormal basis of $\mathbb C[L_N'/L_N]$. This is
\cite[(3.2)]{BYCM} and is exactly $\|\cdot\|_\Gamma$ under
the effective-domain convention fixed above. Thus the transfer
of the Bruinier--Yang constant introduces no additional central factor.

\begin{theorem}[Bruinier--Yang and its theta-height consequence]
\label{thm:BY-newline-height}
Let $N$ be any positive integer and let $G$ be a normalized
weight-two newform of level $N$ and root number $-1$.
Let $g_G\in S^{\rm new}_{3/2,\rho_{L_N}}$ be its Jacobi--Shimura
newform, with real coefficients in its Hecke field.
Use the generic Heegner divisor convention of
\cite[(7.6)]{BYCM}: $Z(m,\mu_r)=P_{-4Nm,r}+P_{-4Nm,-r}$,
retaining both signs.  Write $\mathcal A^G$ for the $G$-projection
of the associated degree-zero Jacobian-valued generating series.
There is $P_G\in J_0(N)(\mathbb Q)\otimes\mathbb R$ such that
\begin{equation}\label{eq:BY-factor-and-height}
 \mathcal A^G=g_G\otimes P_G,\qquad
 H_{\rm NT}(P_G,P_G)=
    \frac{2\sqrt N}{\pi\|g_G\|_\Gamma^2}L'(G,1).
\end{equation}
Consequently, for every vector cusp form $F\in\mathbb Cg_G$,
\begin{equation}\label{eq:BY-vector-height}
 H_{\rm NT}\bigl(\widehat\Theta_{\rm MW}(F),
                       \widehat\Theta_{\rm MW}(F)\bigr)
       =\frac{2\sqrt N}{\pi}\|F\|_\Gamma^2L'(G,1).
\end{equation}
Here the Mordell--Weil generating series is identified by its
generic Heegner coefficients, as proved above; changes by cusps,
vertical divisors and Green data have no effect on that projection.
\end{theorem}
\begin{proof}
The lattice in \cite[(7.3)]{BYCM} is exactly $L_N$, with
$Q=N\det$; changing the name of its upper-right integer does not
change the lattice.  Section~7.1 allows every positive integer
$N$.  Lemma~7.4 supplies a harmonic Maass form $h_G$ with
$\xi h_G=g_G/\|g_G\|_\Gamma^2$ and $c^+_{h_G}(0,0)=0$.
Theorems~7.6--7.7 give \eqref{eq:BY-factor-and-height} with
$P_G=y(h_G)$.  There is no squarefree-level restriction in
these theorems.  The choice of cusp divisor defining degree zero
does not change its Jacobian class, by Manin--Drinfeld.

If $F=a g_G$, coefficientwise Petersson pairing of the factorization
gives
\[
 \widehat\Theta_{\rm MW}(F)
       =\overline a\,\|g_G\|_\Gamma^2 P_G.
\]
The other Jacobian isotypic components make no contribution by
Hecke equivariance and Petersson orthogonality of distinct newform
eigensystems.
Taking the Hermitian height and substituting
\eqref{eq:BY-factor-and-height} proves
\eqref{eq:BY-vector-height}.
\end{proof}

The result also gives, for any scalar cusp form $f$,
\begin{equation}\label{eq:BY-scalar-projected-height}
 H_{\rm NT}\bigl(\operatorname{pr}_G\widehat\theta_{\rm MW,sc}(f),
                 \operatorname{pr}_G\widehat\theta_{\rm MW,sc}(f)\bigr)
 =\frac{2\sqrt N}{\pi}
    \frac{|\langle Rg_G,f\rangle_H|^2}{\|g_G\|_\Gamma^2}
       L'(G,1).
\end{equation}
Indeed the projected lift is $\langle Rg_G,f\rangle_H P_G$.
If $Tf\in\mathbb Cg_G$, adjunction makes the fraction in
\eqref{eq:BY-scalar-projected-height} equal to $\|Tf\|_\Gamma^2$.
Thus, once the full lift is known to equal this Mordell--Weil
projection, its exact scalar height is
\[
 \frac{2\sqrt N}{\pi}\|Tf\|_\Gamma^2L'(G,1)
 =\frac{2\sqrt N}{\pi}r(f)\|f\|_H^2L'(G,1).
\]
Membership of the test form's image under $T$ in the indicated new line
must be verified; packet preservation is not assumed implicitly.

\begin{corollary}[Raw and isometrically normalized scalar lifts]
\label{cor:BY-scalar-normalizations}
Suppose $f\ne0$, $0\ne Tf\in\mathbb Cg_G$, and $L'(G,1)\ne0$.
Then the formula
\[
 H_{\rm NT}(\widehat\theta_{\rm MW,sc}(f),
             \widehat\theta_{\rm MW,sc}(f))
       =\frac{2\sqrt N}{\pi}\|f\|_H^2L'(G,1)
\]
holds if and only if $r(f)=1$.  In particular it fails for the raw
scalar lift of every such nonzero $f\in\operatorname{im}R$ at
squarefree level, since the trace eigenvalue in
\eqref{eq:squarefree-exact-trace} is greater than one.

There is a canonical normalization for which that coefficient is
correct.  Put $K=RT$ and $U=TK^{-1/2}$ on
$\operatorname{im}R$, and define
\[
 \widehat\theta^{\rm iso}_{\rm MW,sc}(f)
       =\widehat\theta_{\rm MW,sc}(K^{-1/2}f).
\]
For $f\in\operatorname{im}R$ with $Uf\in\mathbb Cg_G$, one has
\begin{equation}\label{eq:BY-isometric-height}
 H_{\rm NT}(\widehat\theta^{\rm iso}_{\rm MW,sc}(f),
             \widehat\theta^{\rm iso}_{\rm MW,sc}(f))
       =\frac{2\sqrt N}{\pi}\|f\|_H^2L'(G,1).
\end{equation}
\end{corollary}
\begin{proof}
Compare Theorem~\ref{thm:BY-newline-height} with the proposed raw
formula and cancel its nonzero common factors.  At squarefree
level, Corollary~\ref{cor:all-level-norm} evaluates $r(f)$ by
\eqref{eq:squarefree-exact-trace}.  For the second assertion,
Proposition~\ref{prop:polar-normalized-induction} gives
$\|TK^{-1/2}f\|_\Gamma^2=\|f\|_H^2$.  Apply
\eqref{eq:BY-vector-height} and scalar--vector adjunction.
\end{proof}

\paragraph{The product defect.}
In the notation and under the convergence hypotheses of
Proposition~\ref{prop:defect-height}, suppose that the full lift
equals the Mordell--Weil lift above and $0\ne Tf\in\mathbb Cg_G$.
If $A_f$ is regular at zero, the existing residual identity now
evaluates the defect exactly:
\begin{equation}\label{eq:BY-defect-evaluation}
 d_N(f,f)=-\|Tf\|_\Gamma^2
       \left(\eta\frac{2\sqrt N}{\pi}+C_NA_f(0)\right)L'(G,1).
\end{equation}
This follows by substitution, using $L(G,1)=0$.
In particular the paired defect vanishes when $L'(G,1)=0$.
When $L'(G,1)\ne0$ and $Tf\ne0$, its vanishing is precisely the
analytic normalization equality
\begin{equation}\label{eq:BY-normalization-test}
             C_NA_f(0)=-\eta\frac{2\sqrt N}{\pi}.
\end{equation}
The Bruinier--Yang theorem determines the height side of this
comparison at every level on the specified new line. It neither
identifies the ramified section of the specified genus-two arithmetic
Eisenstein series nor evaluates its doubling multiplier. It therefore
does not by itself prove the general paired arithmetic product identity
for that particular kernel.

\subsection{The Gross--Zagier comparison}

For the trivial class-group character in Gross--Zagier, the Rankin
$L$-function still cannot be replaced by the untwisted $L$-function
over $\Q$.
Let $K=\Q(\sqrt D)$ be an imaginary quadratic field of fundamental
discriminant $D<0$, satisfying the Heegner hypotheses, and let $\chi_D$
be its quadratic Dirichlet character. For trivial class-group character,
the base-change $L$-function factors as
\begin{equation}\label{eq:basechange-factor}
 L_K(G,s)=L(G,s)L(G\otimes\chi_D,s).
\end{equation}
Thus, when $L(G,1)=0$,
\begin{equation}\label{eq:basechange-derivative}
 L_K'(G,1)=L'(G,1)L(G\otimes\chi_D,1).
\end{equation}
The quadratic-twist central value must be retained unless a precisely
stated normalization or a Fourier-coefficient identity absorbs it.

Under the hypotheses and height conventions of
\cite[Theorem~6.3]{GZ}, let $H_K$ be the Hilbert class field, $h_K$ its
class number, $2u_K$ the number of roots of unity in $K$, and
$c_{\xi,G}$ the $G$-projection of
$\sum_{\sigma\in\operatorname{Gal}(H_K/K)}\xi(\sigma)^{-1}c^\sigma$.
Then
\begin{equation}\label{eq:GZ-correct-context}
 L_K'(G,\xi,1)=
 \frac{8\pi^2\|G\|^2}{h_Ku_K^2\sqrt{|D|}}
                   \widehat h_{H_K}(c_{\xi,G}).
\end{equation}
Here $c$ is the Heegner divisor class used there, and
$\widehat h_{H_K}$ is their height over $H_K$. Passing to a height
normalized over $K$ or $\Q$, or replacing the class-group sum by its
average, changes the constants. Equation~\eqref{eq:GZ-correct-context}
is recalled as a cited theorem, not proved by the formal reductions in
this paper. For $\xi=1$, \eqref{eq:basechange-derivative} explains the
factor missing from a formula involving only $L'(G,1)$ and an
unrenormalized individual Heegner class.

For the exact Jacobi new line, the factorization and coefficient
normalization are supplied by the corrected Bruinier--Yang theorem
in Section~\ref{subsec:BY-height}.  Its height formula is the
projected arithmetic input used here.  The quadratic-twist factor
in an individual Heegner coefficient is thereby accounted for;
it does not disappear from the Gross--Zagier identity itself.

\section{Standard coefficients and boundary conventions}
\label{sec:standard-data}

The coefficient calculations use the standard Kudla Green functions
and their cusp expansions. The singular genus-two degrees are
specified explicitly, without assuming modularity of a modified
family.

\subsection{The coefficient model used in the product theorems}
\label{subsec:coefficient-model}

We fix the normalizations for the coefficientwise product
comparisons. Let $\mathcal X=\mathcal X_0(N)$ be the regular
compactified stack, and let $X=\mathcal X(\mathbb C)$ carry its
orbifold measure. Finite lengths, complex integrals and cusp sums
all use the same stabilizer weights; equivalently, they may be
computed on an auxiliary fine level and divided by its degree.
We use
\[
 dd^c=\frac{i}{2\pi}\partial\bar\partial,\qquad
 \widehat{\deg}(0,T)=\frac12\int_XT,\qquad
 \mathfrak a_1(a)=(0,a).
\]
Let $B$ be the polarization of $Q$. For $z=x+iY$, put
\[
 h_z=\frac1{\sqrt N\,Y}
     \begin{pmatrix}-x&x^2+Y^2\\-1&x\end{pmatrix},\qquad
 R(\lambda,z)=\frac12B(\lambda,h_z)^2-2Q(\lambda),
\]
and define
\[
 F_m(v,z)=\sum_{\substack{\lambda\in L\setminus\{0\}\\Q(\lambda)=m}}
                  \beta_1(2\pi vR(\lambda,z)),\qquad
 \beta_1(a)=\int_1^\infty e^{-at}\frac{dt}{t}.
\]
These functions are interpreted as Green currents at positive CM
points. At a cusp with representative $a/c$, $c\mid N$, set
$q_c=\gcd(c,N/c)$ and $M_c=N/q_c^2$; infinity is represented
by denominator $N$. Define
\[
 g_c(m,v)=
 \begin{cases}
 \dfrac{\sqrt N}{2\pi q_c\sqrt v},&m=0,\\[4pt]
 \dfrac{\sqrt N}{2\pi q_c\sqrt v}
       \beta_{3/2}(4\pi|m|v),&m=-M_cj^2,\ j\in\mathbb Z_{>0},\\[4pt]
 0,&\text{otherwise},
 \end{cases}
 \qquad
 \beta_{3/2}(a)=\int_1^\infty e^{-at}t^{-3/2}\,dt.
\]
Write $C_m(v)=\sum_c g_c(m,v)[c]$, interpreted on closed
arithmetic cusp orbits and then with the corresponding complex
weights. The next two subsections prove that these are exactly
the boundary coefficients of $F_m$.

Let $\widehat\lambda$ be the source Hodge bundle with Petersson
squared norm $c_HY$ in a Tate frame, where
$c_H=2\sqrt\pi e^{-\gamma/2}$, and put
\[
 H=-\widehat\lambda-W_N^*\widehat\lambda,\qquad
 \operatorname{curv}H=\varphi_0
                  =-\frac{dx\wedge dY}{2\pi Y^2}.
\]
For $m>0$, $Z(m)$ denotes the actual finite flat special-morphism
Cartier divisor with the multiplicities used in the generic cycle
comparison. Define the standard arithmetic coefficients by
\[
 \widehat Z_{\rm std}(m,v)=
 \begin{cases}
       (Z(m),F_m(v)),&m>0,\\
       (C_m(v),F_m(v)),&m<0,\\
       H+(C_0(v),F_0(v))-\mathfrak a_1(\log v),&m=0.
 \end{cases}
\]
We use the extended arithmetic pairing with the integral cusp
frames and regular parts verified in
Lemma~\ref{lem:all-level-normal-cusp-parameter} and
Lemma~\ref{lem:complete-zero-cusp-regular-part}.
In particular the boundary divisor is part of the arithmetic class.

\begin{definition}[Standard genus-two degree coefficients]
\label{def:standard-coefficients}
The indices are half-integral matrices
\[
 T\in\operatorname{Sym}_2(\mathbb Z)^\vee
   =\left\{\begin{pmatrix}m&r/2\\r/2&n\end{pmatrix}:
                                      m,n,r\in\mathbb Z\right\},
 \qquad y\in\operatorname{Sym}_2(\mathbb R)_{>0}.
\]
For nonsingular $T$, let $\mathcal A^{\rm std}(T,y)$ be the
finite special-cycle degree plus one half of the integral of the
Gaussian-removed Green form $g_2^o(T,y)$ used in
Section~\ref{subsec:proper-compact-heat-star}. The finite term is
zero unless $T$ is positive definite. It uses the scheme
lengths of the special-morphism stack, with ordered morphisms,
residue degrees and stabilizer factors, rather than reduced
point counts. Lemma~\ref{lem:standard-proper-finite-coefficients}
identifies these multiplicities with the proper divisor product.
Both terms carry the special-moduli multiplicities fixed above.

For $\operatorname{rank}T=1$, let $t\ne0$ be its signed integral content
under $\mathrm{GL}_2(\mathbb Z)$ equivalence, so that $T$ is
equivalent to $\operatorname{diag}(0,t)$. Put
\[
 s=t^{-1}\operatorname{tr}(Ty)>0,\qquad \nu=\det(y)/s>0,
\]
and prescribe
\[
 \mathcal A^{\rm std}(T,y)
     =\langle\widehat Z_{\rm std}(t,s),H\rangle_{\rm GS}
                   -\frac12\log\nu\,\deg Z^*_{\rm std}(t,s),
\]
where $Z^*_{\rm std}(t,s)$ is the divisor component of
$\widehat Z_{\rm std}(t,s)$. For $T=0_2$, prescribe
\[
 \mathcal A^{\rm std}(0_2,y)
      =\widehat{\deg}(H^2)-\frac12\log\det(y)\deg H.
\]
These are degree prescriptions, with no genus-two modularity
assumption. Their diagonal product defect is
\begin{align}\label{eq:standard-defect-definition}
 \delta^{\rm std}_{m,n}(v_1,v_2)
    ={}&\langle\widehat Z_{\rm std}(m,v_1),
                       \widehat Z_{\rm std}(n,v_2)\rangle_{\rm GS}\notag\\
       &-\sum_{r\in\mathbb Z}
          \mathcal A^{\rm std}\!\left(
              \begin{pmatrix}m&r/2\\r/2&n\end{pmatrix},
                       \operatorname{diag}(v_1,v_2)\right).
\end{align}
The sum includes both singular endpoints when they occur.
Convergence and the permitted regroupings are proved in the
coefficient calculations below.
\end{definition}

\begin{lemma}[Intrinsic parameters of a rank-one coefficient]
\label{lem:standard-rank-one-intrinsic-parameters}
Every nonzero half-integral rank-one matrix has a unique
expression
\[
                 T=t\mathbf a\mathbf a^t,
 \qquad t\in\mathbb Z\setminus\{0\},\quad
 \mathbf a\in\mathbb Z^2\text{ primitive},
\]
where $\mathbf a$ is unique up to sign.  It is integrally
equivalent to $\operatorname{diag}(0,t)$.  For $y>0$ its
rank-one parameters are intrinsically
\[
 s=\mathbf a^ty\mathbf a=t^{-1}\operatorname{tr}(Ty)>0,
                  \qquad\nu=\det(y)/s>0.
\]
In particular the standard rank-one degree prescription does
not depend on the choice of diagonalizing matrix.
\end{lemma}
\begin{proof}
If $T_{12}=r/2$, rank one gives $r^2=4T_{11}T_{22}$,
so $r$ is even and $T$ is integral.  Its rational image
is generated by a primitive integer column $\mathbf a$;
symmetry and rank one give $T=t\mathbf a\mathbf a^t$
with $t\in\mathbb Q^\times$.  Choose an integer column
$\mathbf b$ with $\mathbf b^t\mathbf a=1$.  Then
$t=\mathbf b^tT\mathbf b\in\mathbb Z$.
To prove uniqueness, suppose that another primitive integer column
$\mathbf c$ gives the same nonzero rank-one matrix. Equality of the
rational image lines gives $\mathbf c=q\mathbf a$ for a nonzero
rational number $q$. Pairing with an integral Bezout vector for
$\mathbf a$ shows that $q$ is an integer. Pairing in the opposite
direction with a Bezout vector for $\mathbf c$ shows that $q^{-1}$
is an integer. Hence $q=\pm1$, and the scalar coefficient is unchanged.
Complete $\mathbf a^t$
to the second row of $\gamma\in\mathrm{GL}_2(\mathbb Z)$.
Then
\[
 \gamma^{-t}T\gamma^{-1}=\operatorname{diag}(0,t),
 \qquad (\gamma y\gamma^t)_{22}=\mathbf a^ty\mathbf a,
 \qquad \det(\gamma y\gamma^t)=\det y.
\]
The second diagonal entry and its Schur complement are thus
$s$ and $\nu$.  This proves the claim in either sign of $t$.
\end{proof}

\subsection{The zero-index cusp correction at arbitrary level}
\label{subsec:cusp-width-correction}

\begin{lemma}[Zero-index cusp coefficients for the standard Green kernel]
\label{lem:zero-cusp-width}
Let
\[
 L=\left\{\begin{pmatrix}b&a/N\\c&-b\end{pmatrix}:a,b,c\in\mathbb Z\right\},
 \qquad Q(w)=N\det w,
\]
and let $a/c$ be a cusp, with $\gcd(a,c)=1$.  For the standard
Kudla Green function $\Xi(0,0,v)$ of \cite[Section~5]{DY1}, its
coefficient of $-\log|q_{a/c}|^2$ is
\[
 g_{a/c}^{\rm std}(0,v)
 =\frac{\sqrt N}{2\pi\sqrt v}
       \frac{\gcd(N,c)}{\gcd(N,c^2)}.
\]
Here $q_{a/c}$ is the local parameter with the classical cusp width,
and the formula at $c=0$ uses $\gcd(N,0)=N$.
\end{lemma}
\begin{proof}
Let $\ell_\infty=\left(\begin{smallmatrix}0&1\\0&0\end{smallmatrix}\right)$
and $\sigma=\left(\begin{smallmatrix}a&b\\c&d\end{smallmatrix}\right)$
with determinant one.  Direct conjugation gives
\[
 \sigma\ell_\infty\sigma^{-1}
 =\begin{pmatrix}-ac&a^2\\-c^2&ac\end{pmatrix}.
\]
The smallest positive rational $\beta$ such that this matrix times
$\beta$ belongs to $L$ is
\[
 \beta=\frac1{\gcd(Na^2,ac,c^2)}
       =\frac1{\gcd(N,c)}.
\]
The classical width is $\kappa=N/\gcd(N,c^2)$, as follows by
conjugating $\left(\begin{smallmatrix}1&m\\0&1\end{smallmatrix}\right)$
by $\sigma$.  \cite[Theorem~5.1(3)]{DY1} gives
$g_{a/c}^{\rm std}(0,v)=\kappa/(2\pi\beta\sqrt{Nv})$.
\end{proof}

For example, when $N=p^2$, the standard coefficients at $\infty$
and $1/p$ are respectively
\[
 \frac{p}{2\pi\sqrt v},\qquad\frac1{2\pi\sqrt v}.
\]
Their ratio is $p$.  By contrast, the Cartier divisor called
$\operatorname{Cusp}(X_0(N))$ in \cite[Section~5.3]{Zhu} is defined by
pulling back the sum of the two boundary divisors on
$X_0(1)\times X_0(1)$.  Its coefficients are the sum of the
ramification indices of $p_1$ and $p_2$ at each cusp.  For $N=p^2$
they are $p^2+1$ at $\infty$ and $2$ at $1/p$; indeed the first
projection has cusp widths $1,1$ at these two cusps, and the second
has widths $p^2,1$, since $p_2=p_1\circ W_N$.
Their ratio is $(p^2+1)/2$, which differs from $p$ for every prime.
A uniform scalar multiple of that pullback cusp divisor therefore
fails to reproduce the zero-index standard Green cusp coefficients,
and a common global stack-normalization factor cannot change the
ratio. This obstructs identification of the boundary prescription
in Zhu's equations~(48)--(49) with the standard all-level Du--Yang
kernel by corrections to the Fourier label, Gaussian factor and
definition of $R$ alone. A separately specified geometric correction
remains possible, but both the correction and its effect on the
arithmetic identities require proof.

\begin{corollary}[A necessary standard zero-index boundary divisor]
\label{cor:standard-zero-boundary}
With reduced cusp points as basis, put
\[
 C_N^{\rm std}=\sum_{[a/c]}
    \frac{\gcd(N,c)}{\gcd(N,c^2)}[a/c].
\]
The standard zero-index Green function has logarithmic cusp divisor
\[
       \frac{\sqrt N}{2\pi\sqrt v}\,C_N^{\rm std}.
\]
For $N=p^n$ and a cusp of denominator $p^k$, its coefficient in
$C_N^{\rm std}$ is $p^{-\min(k,n-k)}$.
\end{corollary}
\begin{proof}
Sum Lemma~\ref{lem:zero-cusp-width} over the cusp classes.
The last assertion is $p^{k-\min(n,2k)}=p^{-\min(k,n-k)}$.
\end{proof}

This determines the boundary term forced by the standard analytic
Green function. Replacing a cusp divisor in an arithmetic family
changes that family. The degree, Hodge, and genus-two identities
for the changed family require verification; the corollary alone
is not an arithmetic Siegel--Weil theorem.

\subsection{The negative-index correction at each cusp}
\label{subsec:negative-cusp-correction}

\begin{proposition}[Cusp-dependent negative-index Green coefficients]
\label{prop:standard-negative-boundary}
Let
\[
 L_N=\left\{\begin{pmatrix}b&A/N\\C&-b\end{pmatrix}:
                    b,A,C\in\mathbb Z\right\},\qquad Q=N\det.
\]
Represent a cusp by $a/c$ with $c\mid N$ and $\gcd(a,c)=1$;
the cusp $\infty$ is treated separately.  Put
\[
 q_c=\gcd(c,N/c),\quad M_c=N/q_c^2,
 \quad \kappa_c=\frac N{\gcd(N,c^2)},\quad\beta_c=\frac1c.
\]
At $\infty$ put $q_\infty=1$, $M_\infty=N$,
$\kappa_\infty=1$ and $\beta_\infty=1/N$.
Then the negative rank-one lattice at this cusp is
\[
  (L_N\cap I_c^\perp)/(L_N\cap I_c)
       \simeq (\mathbb Z,j\mapsto-M_cj^2).
\]
For the standard zero-component Kudla Green function
$\Xi(n,0,v)$, the coefficient of $-\log|q_{a/c}|^2$ at this cusp is
\[
 g_c(n,v)=\frac{\sqrt N}{2\pi q_c\sqrt v}
 \begin{cases}
  \beta_{3/2}(4\pi |n|v),& n=-M_cj^2\text{ for }j\in\mathbb Z_{>0},\\
  0,&n<0\text{ otherwise},\\
  1,&n=0.
 \end{cases}
\]
The assertion is independent of the numerator $a$.  For negative
indices outside the stated support, the standard Green function has
no logarithmic divisor coefficient at this cusp.
\end{proposition}
\begin{proof}
Write $\sigma=\left(\begin{smallmatrix}a&b_0\\c&d\end{smallmatrix}\right)
\in\mathrm{SL}_2(\mathbb Z)$ and $I_c=\mathbb Q\sigma\ell_\infty\sigma^{-1}$.
The lattice $\sigma^{-1}L_N\sigma\cap I_\infty^\perp$ consists of
upper-triangular matrices
$\left(\begin{smallmatrix}\xi&\eta\\0&-\xi\end{smallmatrix}\right)$.
If $w=\left(\begin{smallmatrix}b&A/N\\C&-b\end{smallmatrix}\right)$
lies in $I_c^\perp$, its action on $(a,c)^t$ has eigenvalue $\xi$.
Thus
\[
 aC=c(\xi+b),\qquad cA=Na(\xi-b).
\]
The first equation gives $c\xi=k\in\mathbb Z$; the second, with
$N=c(N/c)$ and $\gcd(a,c)=1$, gives
$c\mid (N/c)k$.  Consequently $\xi\in q_c^{-1}\mathbb Z$.
Conversely, for $\xi=j/q_c$, choose $b\in\mathbb Z$ such that
$c\xi+cb$ is divisible by $a$ (possible because $a,c$ are coprime).
The two displayed equations then give integers $A,C$.
The case $a=0,c=1$ follows directly, and gives the same conclusion.
The quadratic form is $Q=-N\xi^2$, proving the quotient formula.
The isotropic fibers in $\eta$ have period $\beta_c$; this follows
from $\sigma\ell_\infty\sigma^{-1}$ exactly as in the zero-index
cusp-width calculation.

For completeness we compute the logarithmic term directly.  In
the transformed spatial coordinate $z=x+iY$, for a fixed nonzero
diagonal entry $\xi$ the fiber is $\eta\in\eta_0+\beta_c\mathbb Z$,
and the standard majorant is
\[
 \mathcal R\!\left(\begin{pmatrix}\xi&\eta\\0&-\xi\end{pmatrix},z\right)
 =\frac{N(2\xi x+\eta)^2}{2Y^2}+2N\xi^2.
\]
Poisson summation in $\eta$ shows that its zero frequency in
$\sum_\eta\beta_1(2\pi v\mathcal R)$ is
\[
 \frac{Y}{\beta_c\sqrt{Nv}}
      \int_1^\infty e^{-4\pi Nv\xi^2s}s^{-3/2}\,ds.
\]
The two diagonal entries $\xi=\pm j/q_c$ give twice this term.
Since $-\log|q_{a/c}|^2=4\pi Y/\kappa_c$, the required
coefficient is
$\kappa_c\beta_{3/2}(4\pi|n|v)/(2\pi\beta_c\sqrt{Nv})$.
The zero-index coefficient is
$\kappa_c/(2\pi\beta_c\sqrt{Nv})$ by the zero-index calculation.
Finally $\kappa_c/\beta_c=N/q_c$, giving the displayed formula.
Nonzero Poisson frequencies and vectors outside $I_c^\perp$
have no linear term in $Y$; this is the cusp expansion in
\cite[Theorem~5.1]{DY1}.  At $\infty$ the same calculation has
$\xi\in\mathbb Z$, as asserted.
\end{proof}

Thus the complete standard boundary coefficient series at this cusp is
\[
 b_{N,c}(\tau)=\frac{\sqrt N}{2\pi q_c\sqrt v}
 \left(1+\sum_{j\ge1}
    \beta_{3/2}(4\pi M_cj^2v)q^{-M_cj^2}\right).
\]
Its shadow is
$\xi_{3/2}b_{N,c}=-\sqrt N\,\vartheta_{M_c}/(4\pi q_c)$.
This is a statement about the standard analytic Green kernel;
the series alone is not asserted to be modular.

At $N=p^2$, Proposition~\ref{prop:standard-negative-boundary}
gives support $n=-p^2j^2$ at $\infty$, but support $n=-j^2$
at a middle cusp of denominator $p$. In particular the standard
coefficient for $n=-1$ at $\infty$ is zero. This also follows
immediately from $B(x,\ell_\infty)=-Nc$: perpendicularity forces
$c=0$ and then $Q(x)=-p^2b^2$. In contrast, the literal boundary
prescription in \cite[equations~(48)--(49)]{Zhu}, at the same
quadratic index $n=-1$, gives at $\infty$ the positive coefficient
\[
       (p^2+1)\frac{p}{2\pi\sqrt v}\beta_{3/2}(4\pi v).
\]
The discrepancy thus extends to negative-index support and cannot
be removed by rescaling the Green Gaussian or applying a common
stack factor. The comparison uses the same Fourier index
$Q(x)=n$ on both sides. It tests the proposed standard-kernel
identification but does not determine the cuspidal projection of
the complete arithmetic product defect.

\subsection{The exact lowering term of the standard arithmetic family}
\label{subsec:standard-arithmetic-lowering}

Use the standard Kudla Green function and the cusp-dependent
boundary coefficients proved above.  At a cusp of denominator
$c\mid N$, put $q_c=\gcd(c,N/c)$, $M_c=N/q_c^2$ and
$\kappa_c^{\rm bdry}=\sqrt N/(2\pi q_c)$; use
$q_\infty=1$ and $M_\infty=N$ at infinity.  Write
$\widehat\phi_{\rm std}(\tau)$ for this arithmetic family,
including its usual zero Hodge coefficient and $-(0,\log v)$,
where $\tau=u+iv$.  Its lowering operator is
$L_{3/2}=-2iv^2\partial_{\bar\tau}$.

\begin{lemma}[Current convergence before lowering]
\label{lem:standard-Green-current-convergence}
On every compact subset of the modular upper half plane, the
complete Green series
\[
 \sum_{0\ne x\in L_N}
       \beta_1(2\pi v\mathcal R(x,z))q^{Q(x)}
\]
converges absolutely in $L^1(X,\nu)$, where $\nu$ is any fixed
smooth positive form on the compactified geometric curve.
The series can be differentiated in $u,v$ as currents, locally
uniformly in $\tau$.  The boundary-divisor coefficient series
and all their modular-variable derivatives also converge
locally uniformly.  These assertions concern divisor and Green
data; they do not assert modularity or realization in a specific
arithmetic Chow space.
\end{lemma}
\begin{proof}
Fix $0<v_0\leq v\leq v_1$.  At a geometric cusp, write its
scaling coordinate as $z=x+iY$, with $x$ bounded and $Y\geq1$.
For a fixed rational coordinate grid containing the transformed
lattice, the positive majorant $M_z=\mathcal R+Q$ satisfies
\[
       M_z(x)\geq cY^{-2}\|x\|^2.
\]
For $Q(x)\leq0$ and heat time $s\geq1$, one has
$s\mathcal R(x,z)+Q(x)=sM_z(x)+(1-s)Q(x)\geq sM_z(x)$.
Consequently the sum of the absolute values of the nonpositive
terms is at most $\sum_{x\ne0}\beta_1(2\pi v_0M_z(x))$.
An excluded-zero Gaussian estimate on a fixed three-dimensional
grid gives
\[
 \sum_{x\ne0}e^{-c s\|x\|^2/Y^2}
   \leq C Y^3s^{-3/2}e^{-c's/Y^2}.
\]
Integration against $ds/s$ bounds this Green sum by $CY^3$.
This is integrable against every smooth compactified density,
whose expression in the cusp coordinate decays exponentially
in $Y$.

For positive norm $m$, regroup the lattice vectors in their
finitely many CM $\Gamma_0(N)$-orbits and unfold their positive
Green sums.  The number $h_N(m)$ of these orbits, with all signs
retained, is $O_N(1+m)$.  Indeed they give integral positive
binary quadratic forms of discriminant $-4Nm$ with fixed
congruence conditions.  A reduced form has leading coefficient
$A\leq\sqrt{4Nm/3}$ and middle coefficient $|B|\leq A$, giving
$O_N(1+m)$ possibilities; passing to the fixed finite-index
subgroup changes only the constant.
If $r$ is hyperbolic distance to a CM point of norm $m$, then
$\mathcal R=2m\sinh^2r$.  Hence, by polar coordinates,
\[
 \int_{\mathbb H}\beta_1(2\pi v\mathcal R)\,
                      \frac{dx\,dy}{y^2}
 \leq2\pi\int_0^\infty s\beta_1(4\pi vm s^2)\,ds
 =\frac1{4vm}.
\]
Here $\int_0^\infty\beta_1(t)\,dt=1$.  The ratio of $\nu$
to hyperbolic area is bounded, including at the cusps.
Unfolding therefore bounds the positive part in $L^1(X,\nu)$
by $C\sum_{m\geq1}e^{-2\pi v_0m}h_N(m)/m<\infty$.
Orbifold stabilizer weights only reduce this upper bound.
The same orbit count proves total-variation convergence of the
positive horizontal divisor-current series and all its
modular-variable derivatives.

Modular-variable derivatives add polynomial factors in $Q$,
$\mathcal R$ and heat time.  On the nonpositive part,
$|Q|\leq M_z$ and $\mathcal R\leq2M_z$, so a fixed fraction
of the Gaussian absorbs every such factor, retaining a
polynomial cusp bound.  On the positive part, derivatives of
$\beta_1(2\pi v\mathcal R)$ of positive order are bounded by
constants times a Gaussian in $v_0\mathcal R$; its same radial
integral is $O((v_0m)^{-1})$.  The powers of $m$ from
differentiating $q^m$ are absorbed by its exponential decay.
Thus these bounds hold for each specified derivative order and
justify current differentiation of the full sum.
Finally the boundary series and their derivatives are
Gaussian-convergent unary series, which proves the remaining
assertion.
\end{proof}

\begin{proposition}[Lowering with its cusp divisor retained]
\label{prop:standard-arithmetic-lowering}
As differentiable standard divisor and Green-current data one has
\begin{equation}\label{eq:standard-arithmetic-lowering}
 L_{3/2}\widehat\phi_{\rm std}(\tau)
     =\bigl(C_{\rm low}(\tau),-\Theta_S(\tau,\cdot)\bigr),
\end{equation}
where $\Theta_S=v\sum_{x\in L_N}e^{-2\pi v\mathcal R(x,z)}q^{Q(x)}$
is the full scalar Siegel theta function, including $x=0$, and
\[
 C_{\rm low}(\tau)=
   -\sum_{\mathfrak c}
       \frac{\sqrt N}{4\pi q_{\mathfrak c}}
       v^{1/2}\overline{\vartheta_{M_{\mathfrak c}}(\tau)}
           [\mathfrak c],\qquad
 \vartheta_M(\tau)=\sum_{j\in\mathbb Z}q^{Mj^2}.
\]
If these data realize differentiable classes in the chosen
extended arithmetic Chow space, the same identity holds there.
Under that realization, the lowering class is not identically
zero and cannot be a pure metric class: its generic-fiber degree
is not identically zero.  The displayed cusp divisor is therefore
necessary.  No modularity assertion is included in this statement.
\end{proposition}
\begin{proof}
For every nonzero lattice vector,
\[
 \frac{\partial}{\partial v}
     \beta_1(2\pi v\mathcal R(x,z))
    =-v^{-1}e^{-2\pi v\mathcal R(x,z)}.
\]
Since $q^{Q(x)}$ is holomorphic, applying $L_{3/2}$ to a
coefficient times $q^{Q(x)}$ acts as $v^2\partial_v$ on the
coefficient.  Thus the Green part is minus the Siegel kernel
with $x=0$ omitted.  At index zero the derivative of
$-(0,\log v)$ is $(0,-v)$, supplying exactly the missing
zero vector.  The fixed Hodge bundle has derivative zero.

The positive-index horizontal divisors are independent of $v$.
At a cusp the negative boundary coefficient is
$\kappa_c^{\rm bdry}v^{-1/2}
\beta_{3/2}(4\pi M_cj^2v)$ for $j\ge1$, and the zero
coefficient is $\kappa_c^{\rm bdry}v^{-1/2}$.  Therefore
\[
 v^2\partial_v g_c(-M_cj^2,v)
     =-\kappa_c^{\rm bdry}v^{1/2}e^{-4\pi M_cj^2v},
 \qquad
 v^2\partial_v g_c(0,v)
     =-\tfrac12\kappa_c^{\rm bdry}v^{1/2}.
\]
Here one uses the integration-by-parts identity
\[
 \beta_{3/2}(t)=2e^{-t}-2t\beta_{1/2}(t),
 \qquad \beta_s(t)=\int_1^\infty e^{-tu}u^{-s}\,du,
 \quad t>0.
\]
Multiplication by $q^{-M_cj^2}$ converts
$e^{-4\pi M_cj^2v}q^{-M_cj^2}$ into
$\overline{q^{M_cj^2}}$.  Summing gives precisely the stated
$C_{\rm low}$.  Lemma~\ref{lem:standard-Green-current-convergence}
justifies the operations after summing the full series, as compact
currents.  The Siegel kernel and its modular-variable derivatives
also have a Gaussian majorant $O(1+Y)$ in a cusp coordinate,
by successive summation in the three lattice coordinates, so its
compact current interpretation is well-defined.

For the class assertion, let $d_{\mathfrak c}>0$ be the generic
degree of the reduced cusp class in the fixed degree convention.
The constant Fourier term in $u$ of the generic degree of
$C_{\rm low}$ is
\[
       -\frac{\sqrt N}{4\pi}\sqrt v
          \sum_{\mathfrak c}\frac{d_{\mathfrak c}}
                                     {q_{\mathfrak c}}<0.
\]
This proves that its generic degree is not identically zero.
Principal divisor classes and pure metric classes have generic
degree zero, so neither can remove this conclusion.
\end{proof}

The formula identifies the complete geometric input for a
lowering calculation of the arithmetic product kernel.  No
identity with the lowering of an Eisenstein derivative, and no
harmonicity of their defect, is asserted here.

\subsection{Exact algebra for the singular diagonal terms}
\label{subsec:singular-diagonal-algebra}

Write $\mathfrak a_1(f)=(0,f)$ for a metric class of
codimension one and $\mathfrak a_2(T)=(0,T)$ for a
codimension-two current class.  We use Zhu's arithmetic degree
normalization
\[
 \widehat{\deg}\,\mathfrak a_2(T)=\frac12\int_X T.
\]
All products below are products of arithmetic classes in a fixed
intersection theory.  In particular, no product of two Dirac
currents, or of a singular Green function with its own Dirac
current, is used.  For an arithmetic divisor $A$ with compact
curvature current $\omega_A$, the constant-metric rule is
\begin{equation}\label{eq:singular-constant-metric-rule}
 \mathfrak a_1(\ell)\cdot A=\mathfrak a_2(\ell\omega_A),
 \qquad
 \mathfrak a_1(\ell)\cdot\mathfrak a_1(\ell')=0,
 \qquad
 \widehat{\deg}\bigl(\mathfrak a_1(\ell)\cdot A\bigr)
       =\frac{\ell}{2}\deg A.
\end{equation}
Here $\ell,\ell'$ are constants, and $\deg A=\int_X\omega_A$.
These formulas also follow by rescaling one metric by a constant,
so they do not require a singular-current multiplication.

Put
\begin{equation}\label{eq:singular-isotropic-piece}
 \begin{aligned}
 H&=-\widehat\omega_N-W_N^*\widehat\omega_N,\\
 B(v)&=\bigl(Z^*(0,v),g(0,v)\bigr),\qquad \ell(v)=\log v.
 \end{aligned}
\end{equation}
Zhu's Definition~7.1.3 \cite{Zhu} gives the equality of divisor classes
\begin{equation}\label{eq:singular-zero-splitting}
       \widehat Z(0,v)=H+B(v)-\mathfrak a_1(\ell(v)).
\end{equation}
The term $B(v)$ contains the nonzero isotropic-vector contribution
and its boundary divisor; it must not be silently discarded.

\begin{proposition}[A zero diagonal entry and a rank-one term]
\label{prop:singular-zero-positive-algebra}
Let $t\ne0$ and $v_1,v_2>0$.  Suppose that
$A=\widehat Z(t,v_2)$ and $B(v_1)$ are classes in the stated
intersection theory, with globally defined, locally integrable
Green current $g_t=g(t,v_2)$ and compact Green equation
\[
       \omega_A=\delta_{Z^*(t,v_2)(\C)}+dd^c g_t.
\]
Then Zhu's rank-one definition yields
\begin{equation}\label{eq:singular-rank-one-class}
 \widehat Z(0,v_1)\cdot\widehat Z(t,v_2)
 -\widehat Z\!\left(
     \begin{pmatrix}0&0\\0&t\end{pmatrix},
       \operatorname{diag}(v_1,v_2)\right)
       =B(v_1)\cdot\widehat Z(t,v_2)
\end{equation}
in codimension-two classes modulo exact currents.  In particular,
\begin{equation}\label{eq:singular-rank-one-degree}
 \widehat{\deg}\bigl(\widehat Z(0,v_1)\cdot\widehat Z(t,v_2)\bigr)
 -\widehat{\deg}\,\widehat Z\!\left(
     \begin{pmatrix}0&0\\0&t\end{pmatrix},
       \operatorname{diag}(v_1,v_2)\right)
       =\langle B(v_1),\widehat Z(t,v_2)\rangle_{\rm GS}.
\end{equation}
\end{proposition}
\begin{proof}
Definition~7.2.2 gives the rank-one class
\[
       A\cdot H-
       \mathfrak a_2\bigl(\ell(v_1)
                    \delta_{Z^*(t,v_2)(\C)}\bigr).
\]
Expand the first product using
\eqref{eq:singular-zero-splitting} and use
\eqref{eq:singular-constant-metric-rule}.  Subtracting the
rank-one class gives
\[
 B(v_1)\cdot A-
       \mathfrak a_2\bigl(\ell(v_1)dd^c g_t\bigr).
\]
The second summand is exact because $g_t$ is a global distribution
on the compact curve and $dd^c g_t$ belongs to
$\operatorname{im}\partial+\operatorname{im}\bar\partial$.
Distributional Stokes also gives zero integral. Both statements
follow, since the multiplier of $dd^c g_t$ is constant.
\end{proof}

\begin{proposition}[The rank-zero algebra]
\label{prop:singular-rank-zero-algebra}
Let $B_i=B(v_i)$, $\ell_i=\log v_i$, and let
$\Omega_Z$ denote the form in Zhu's Definition~7.2.4.
Assume the same well-defined products and compact curvatures as
above.  Without imposing a convention on $\Omega_Z$, one has
\begin{align}\label{eq:singular-rank-zero-general}
 &\widehat Z(0,v_1)\cdot\widehat Z(0,v_2)
       -\widehat Z(0_2,\operatorname{diag}(v_1,v_2))\notag\\
 &\quad=H\cdot(B_1+B_2)+B_1\cdot B_2
       -\mathfrak a_2(\ell_1\omega_{B_2}
                               +\ell_2\omega_{B_1})
       -\mathfrak a_2\bigl((\ell_1+\ell_2)
                                      (\omega_H+\Omega_Z)\bigr).
\end{align}
For the compatible Petersson convention
\[
       \omega_H=-\frac{dx\wedge dy}{2\pi y^2},
       \qquad \Omega_Z=\frac{dx\wedge dy}{2\pi y^2},
\]
the final summand is zero.  Thus, writing
$d_i=\deg B_i$, the exact numerical discrepancy at the rank-zero
matrix alone is
\begin{align}\label{eq:singular-rank-zero-degree}
 &\widehat{\deg}\bigl(\widehat Z(0,v_1)\cdot\widehat Z(0,v_2)\bigr)
       -\widehat{\deg}\,\widehat Z(0_2,
                                  \operatorname{diag}(v_1,v_2))\notag\\
 &\quad=\langle H,B_1+B_2\rangle_{\rm GS}
            +\langle B_1,B_2\rangle_{\rm GS}
            -\frac12(\ell_1d_2+\ell_2d_1).
\end{align}
In particular the explicit logarithmic Hodge correction in the
genus-two definition has the correct sign and coefficient; the
remaining terms involve the nonzero isotropic contribution.
\end{proposition}
\begin{proof}
Expand the two factors in
\eqref{eq:singular-zero-splitting}.  The product of the two
constant metric classes is zero.  The remaining terms are
\[
 H^2+H(B_1+B_2)+B_1B_2
 -\mathfrak a_2\bigl((\ell_1+\ell_2)\omega_H
                     +\ell_1\omega_{B_2}
                     +\ell_2\omega_{B_1}\bigr).
\]
Definition~7.2.4 is
$\widehat Z(0_2,y)=H^2+
 \mathfrak a_2(\log\det y\,\Omega_Z)$.
Subtract it and then apply the degree convention.
The curvature convention displayed above follows from
$H=-\widehat\omega_N-W_N^*\widehat\omega_N$,
$c_1(\widehat\omega_N)=dx\wedge dy/(4\pi y^2)$, and invariance
of the hyperbolic form under $W_N$.  Zhu's final computation in
\S9.3 uses precisely the displayed integral for $\Omega_Z$.
\end{proof}

\paragraph{Relation to the complete Fourier coefficient.}
The left sides of
\eqref{eq:singular-rank-one-degree} and
\eqref{eq:singular-rank-zero-degree} compare the product with
\emph{one} singular genus-two matrix.  They are not the complete
product defects $\delta_{0,t}$ and $\delta_{0,0}$.
With $y=\operatorname{diag}(v_1,v_2)$, those defects satisfy
\begin{align}\label{eq:singular-complete-defects}
 \delta_{0,t}
   &=\langle B(v_1),\widehat Z(t,v_2)\rangle_{\rm GS}
       -\sum_{r\ne0}\widehat{\deg}\,
          \widehat Z\!\left(
            \begin{pmatrix}0&r/2\\r/2&t\end{pmatrix},y\right),\notag\\
 \delta_{0,0}
   &=\langle H,B_1+B_2\rangle_{\rm GS}
          +\langle B_1,B_2\rangle_{\rm GS}
          -\tfrac12(\ell_1d_2+\ell_2d_1)
       -\sum_{r\ne0}\widehat{\deg}\,
          \widehat Z\!\left(
            \begin{pmatrix}0&r/2\\r/2&0\end{pmatrix},y\right).
\end{align}
These formulas apply whenever the indicated sums converge in
the fixed coefficient-extraction convention.  Every matrix in
the sums has determinant $-r^2/4<0$; its genus-two contribution
cannot be discarded simply because the finite special cycle is
empty. For $t>0$, Corollary~\ref{cor:mixed-standard-arithmetic}
identifies the two terms in the first line in the standard
normalization; Corollary~\ref{cor:zero-negative-standard-product}
does so for $t<0$. Thus $\delta^{\rm std}_{0,t}=0$ for every
$t\ne0$. For the second line,
Theorem~\ref{thm:double-zero-absolute-defect} evaluates the
entire defect as an Epstein zeta value. The unmodified standard
prescription therefore has a nonzero double-zero defect, and
Corollary~\ref{cor:corrected-standard-support} removes it by an
explicit rank-zero degree correction. Transferring these comparisons
to $\delta_{m,n}$ for a specified arithmetic Siegel--Weil series
requires the normalization comparison described above.

\paragraph{Necessary scope.}
The compact Green equation used above is essential.  If a printed
cusp divisor fails to match the Green asymptotics, its omitted
Dirac defect survives the rank-one calculation.  If instead
$\omega_A=\delta_{Z^*}+dd^c g_t+\varepsilon_A$, the right side of
\eqref{eq:singular-rank-one-class} acquires the additional term
$-\mathfrak a_2(\ell(v_1)\varepsilon_A)$.
The logarithmic singular corrections are therefore exactly compatible
under the stated compact current equation. This algebra does not
justify a boundary normalization that fails to satisfy that
equation.

\begin{lemma}[Both rank-one endpoints for positive indices]
\label{lem:singular-positive-endpoints}
Let $m,n>0$ be integers and suppose $mn$ is a square.  There
are unique positive integers $t,a,b$ with $\gcd(a,b)=1$ such that
$m=ta^2$, $n=tb^2$.  The two rank-one matrices in the
$(m,n)$ diagonal coefficient are
\[
 T_\pm=t\begin{pmatrix}a^2&\pm ab\\\pm ab&b^2\end{pmatrix}.
\]
Put $v_*=a^2v_1+b^2v_2$ and $v_\perp=v_1v_2/v_*$.
Zhu's Definition~7.2.2 gives the same class at the two endpoints,
namely
\begin{equation}\label{eq:singular-positive-endpoints}
 \widehat Z(T_\pm,\operatorname{diag}(v_1,v_2))
 =\widehat Z(t,v_*)\cdot H
     -\mathfrak a_2\bigl(\log v_\perp\,
                       \delta_{Z^*(t,v_*)(\C)}\bigr).
\end{equation}
Thus both copies must be counted.  In particular for $m=n=t$,
$v_*=v_1+v_2$ and $v_\perp=v_1v_2/(v_1+v_2)$.
\end{lemma}
\begin{proof}
Write $m=kA^2$, $n=kB^2$ with the common square-free part $k$,
and put $d=\gcd(A,B)$, $t=kd^2$, $a=A/d$, $b=B/d$.
A primitive vector $(a,\pm b)$ can be completed to an integral
unimodular basis, so $T_\pm$ has the rank-one content $t$ of
Lemma~7.2.1.  Moreover
$t^{-1}\operatorname{tr}(T_\pm\operatorname{diag}(v_1,v_2))
=a^2v_1+b^2v_2$.  Substitution into Definition~7.2.2 proves
the formula.  This statement evaluates the singular genus-two
terms; it makes no assertion about the separate arithmetic
self-intersection of the genus-one divisor.
\end{proof}

\section{Compact arithmetic product comparisons}
\label{sec:standard-products}

We compare individual diagonal Fourier coefficients, including all
matrices with that diagonal. The estimates in this section justify
the required spatial and heat integrations.

\subsection{The exact residual condition for arithmetic doubling}
\label{subsec:exact-product-defect}

Fix Zhu's arithmetic coefficients, including the modified Hodge
coefficient at zero, and put
\[
 \widehat\phi_1(\tau)=\sum_{m\in\mathbb Z}
          \widehat Z(m,v)q^m,
 \qquad \tau=u+iv.
\]
Use the same compactification, Green-current normalization and
arithmetic-degree map in the genus-two series
$\widehat\phi_2(Z)$. Write $C_N=\psi(N)/24$.
For $z,w\in\mathbb H$, define the scalar product defect by
\begin{equation}\label{eq:defect-kernel-exact}
 \mathcal D_N(z,w)=
 H_{\rm GS}(\widehat\phi_1(z),\widehat\phi_1(w))
 -C_N E'\!\left(\operatorname{diag}(z,-\overline w),
                    0,\Delta(N)^2\right).
\end{equation}
This definition requires no multiplication or summation of
uncontrolled current-valued Fourier series. Whenever those operations
are justified, the coefficient of $q_z^m\overline{q_w^n}$ is
\begin{equation}\label{eq:defect-coefficient-exact}
 \delta_{m,n}(v_z,v_w)
  =\widehat{\deg}\left(
       \widehat Z(m,v_z)\cdot\widehat Z(n,v_w)\right)
   -\sum_{r\in\mathbb Z}\widehat{\deg}\,
       \widehat Z\left(
          \begin{pmatrix}m&r/2\\r/2&n\end{pmatrix},
          \operatorname{diag}(v_z,v_w)\right).
\end{equation}
The sum over $r$ is essential. One fixed genus-two coefficient is
not the coefficient of a product of two genus-one series. Formula
\eqref{eq:defect-coefficient-exact} follows from Zhu's
Theorem~7.2.5 and the restriction of a genus-two Fourier expansion
to the diagonal. It is a definition of a residual quantity, not an
assertion that it vanishes.

For a finite-dimensional scalar cusp space $S$ with the required
weight and character, assume the following pairings converge and set
\begin{equation}\label{eq:defect-pairing}
 d_N(f,h)=\int_{(H\backslash\mathbb H)^2}
       \mathcal D_N(z,w)\overline{f(z)}h(w)
       v_z^{3/2}v_w^{3/2}\,d\mu(z)d\mu(w).
\end{equation}
This form is antilinear in $f$ and linear in $h$. If $(e_a)$ is an
orthonormal basis of $S$, the double cuspidal projection is
\begin{equation}\label{eq:defect-projection}
 (\Pi_S\otimes\overline{\Pi_S})\mathcal D_N(z,w)
       =\sum_{a,b}d_N(e_a,e_b)e_a(z)\overline{e_b(w)}.
\end{equation}
This finite sum can be used as the definition of the projection even
if the unpaired kernel is not square integrable.

\begin{proposition}[A residual version of the arithmetic reduction]
\label{prop:defect-height}
Suppose the arithmetic lift of $f$ exists in the admissible
Mordell--Weil subspace and
$H_{\rm GS}=-\eta H_{\rm NT}$ there. Suppose also that the
classical Eisenstein kernel in \eqref{eq:defect-kernel-exact}
satisfies the analytic doubling identity, with all differentiation
and integration interchanges justified. In the notation of
Proposition~\ref{prop:red-height}, and with $K_{00}$ normalized
to be this classical Eisenstein kernel, one then has the exact identity
\begin{align}\label{eq:defect-height-exact}
 H_{\rm NT}(\widehat\theta(f),\widehat\theta(f))
 =-\frac{1}{\eta}\bigg[&C_N\|Tf\|_\Gamma^2
       \bigl(A_f(0)L'(G,1)+A_f'(0)L(G,1)\bigr)
             +d_N(f,f)\bigg].
\end{align}
For a fixed $f$, the required arithmetic replacement is equivalent
to $d_N(f,f)=0$. For all $f\in S$ it is equivalent to
\begin{equation}\label{eq:defect-minimal-condition}
       (\Pi_S\otimes\overline{\Pi_S})\mathcal D_N=0.
\end{equation}
If a different normalization of $K_{00}$ is used, its explicitly
computed comparison scalar must also multiply $C_N$.
\end{proposition}

\begin{proof}
Insert \eqref{eq:defect-kernel-exact} into the double integral
defining the Hermitian pairing of the two lifts. The first summand
is $C_N$ times the derivative of the analytic doubling integral;
the second is \eqref{eq:defect-pairing}. The height convention
then gives \eqref{eq:defect-height-exact}. The equivalence for a
fixed $f$ follows immediately. For all $f\in S$, complex
polarization recovers every $d_N(f,h)$ from the diagonal values
$d_N(f,f)$, proving the last assertion.
\end{proof}

\paragraph{The coefficient comparisons.}
For positive $m,n$ with $mn$ not a square, the positive special
divisors have no common horizontal component. Zhu's fiber-product
description and proper finite-intersection calculation
\cite[Corollary~4.3.8 and Proposition~8.3.1]{Zhu23}
identify the finite-place contribution with the corresponding sum
of nonsingular genus-two cycles. Theorem~\ref{thm:compact-proper-green-star}
below proves the matching archimedean identity, including vanishing
of the cusp flux, for the explicitly normalized Gaussian-removed
genus-two kernel. Thus the positive nonsquare comparison reduces
to identifying that kernel with the genus-two coefficients in the
arithmetic Siegel--Weil theorem. This normalization comparison
must be made before concluding $\delta_{m,n}=0$ for that series.

When $m,n>0$ and $mn$ is a square, common-component
self-intersections can occur, and rank-one matrices occur at
$r=\pm2\sqrt{mn}$. Thus an improper arithmetic adjunction
calculation is required in addition to the proper calculation.
The local analytic finite part, with its precise tangent-metric
constant, is computed in Theorem~\ref{thm:self-local-identity};
its underlying normal-bundle identity is proved in
Lemma~\ref{lem:integral-CM-adjunction}.  Actual Cartier branches
and the all-level metric comparison are supplied in
Section~\ref{subsec:all-level-canonical-correction}. Their finite
conductors, vertical contacts and distinct-branch intersections
are included in the common-component comparison of
Theorem~\ref{thm:standard-positive-square-product}. Lemma~\ref{lem:singular-positive-endpoints} evaluates both
rank-one endpoint terms without assuming this self-intersection
identity. For $mn=0$, the zero genus-one coefficient and rank-one
genus-two terms must be compared; the rank-zero matrix occurs
only for $m=n=r=0$. Zhu defines these terms separately in
Definitions~7.1.3, 7.2.2 and 7.2.4, with logarithmic and Hodge
corrections. The exact cancellations and the remaining terms are
computed in Section~\ref{subsec:singular-diagonal-algebra}.
Corollary~\ref{cor:mixed-standard-arithmetic} closes every mixed
coefficient with one positive and one nonpositive index for the
standard Gaussian-removed coefficients.
Theorem~\ref{thm:strictly-negative-standard-product} settles all
strictly negative pairs, including the split rank-one contributions;
Corollary~\ref{cor:zero-negative-standard-product} settles the
zero--negative pairs. The only nonpositive pair at which the
standard defect can be supported is $(0,0)$, and
Theorem~\ref{thm:double-zero-absolute-defect} evaluates its
uncorrected defect exactly. The rank-zero correction of
Corollary~\ref{cor:corrected-standard-support} removes that term.
The positive square-product self-intersections are settled in
Theorem~\ref{thm:standard-positive-square-product}.
Section~\ref{sec:remaining} proves the arithmetic Siegel--Weil
comparison for the fixed sections. Other normalizations require
a separate comparison.
Nonpositive indices retain their cusp-divisor and Green-current
contributions.  A partition by Fourier indices need not give
separately modular pieces.

\paragraph{A precise sufficient criterion for discarding a residual.}
Let $\mathcal U$ be the compatible unary-theta subspace and take
$S\subset\mathcal U^\perp$. Suppose a residual kernel admits an
integrably convergent expansion in separable tensors such that in
each tensor at least one factor has zero Petersson pairing with
$S$. Its double pairing with $S\times S$ is then zero, by Fubini.
This applies in particular if that factor belongs to $\mathcal U$,
or is an Eisenstein series or its derivative whose cusp pairing
has separately been proved to vanish. It is the spectral
condition on the factor, rather than a formal label such as
``zero coefficient'' or ``square class'', that proves the
vanishing. A coefficientwise cut of a modular Fourier expansion
usually does not preserve modularity and cannot be substituted
for this condition.

\subsection{The finite part of a proper product}

\begin{lemma}[The proper finite coefficient identity]
\label{lem:standard-proper-finite-coefficients}
Let $m,n>0$ and $mn\notin\mathbb Q^{\times2}$.  Use the
actual positive special-morphism generalized Cartier divisors
and their multiplicities.  Their intersection is proper, and
\[
 (Z(m),Z(n))_{\rm fin}
   =\sum_{r\in\mathbb Z}
        \deg_{\rm fin}\mathcal Z\!\left(
                   \begin{pmatrix}m&r/2\\r/2&n\end{pmatrix}\right).
\]
Only positive-definite Gram matrices on the right can be
nonempty, and the sum of their finite degrees is finite.
\end{lemma}
\begin{proof}
Horizontality requires a separate argument; Cartierness alone
does not imply it. The generic point of every special-fibre
component of the modular curve has transcendental $j$-invariant
and therefore rational endomorphism algebra $\mathbb Q$.
A special isogeny $j$ orthogonal to $\pi$ would yield a nonzero
endomorphism $\alpha=\pi^\vee j$ satisfying
$\alpha^2=-Nm$, which is impossible in $\mathbb Q$.
Thus a positive special divisor contains no vertical
component.  On a regular surface its Cartier equation
defines a Cohen--Macaulay quotient without embedded
associated primes; absence of vertical minimal primes
makes every rational prime a nonzerodivisor on this
quotient.  This proves the required flatness locally
over $\mathbb Z$.  This is the geometric content of the
flatness argument in the proof of
\cite[Proposition~4.3.7]{Zhu23}.  In a ramified chart one
must exclude every vertical component, rather than merely
test divisibility of the equation by the whole element $p$.

The positive divisors also miss the cuspidal locus.  They
have no common horizontal component: on the
complex fibre a common point would place their two positive
special vectors on its one positive rational line, forcing
$n/m$ to be a rational square.  They have no common vertical
component since they are flat horizontal divisors.  Thus the
intersection is zero-dimensional and has no generic fibre.
Properness of the modular stack and finiteness of the
special-morphism maps make this intersection finite over
$\operatorname{Spec}\mathbb Z$.

After an etale scheme cover, generalized Cartierness splits
each special-morphism stack into closed Cartier sheets.  Fix
one ordered pair of sheets at an intersection point.  In the
regular local surface ring $R$ let their equations be $f,g$.
The preceding absence of common codimension-one components
makes $(f,g)$ a regular sequence.  For completeness, $R/(f)$
is a one-dimensional Cohen--Macaulay ring, with no embedded
associated primes; $g$ belongs to none of its minimal primes
and is a nonzerodivisor.  The same holds after interchanging
$f,g$.  Hence the Cartier resolutions give
\[
 \operatorname{Tor}_i^R(R/(f),R/(g))=0\quad(i>0),\qquad
 i_R(f,g)=\operatorname{length}_R R/(f,g).
\]
Thus the finite product multiplicity equals the scheme length
on the fibre product, including at nonreduced points and bad
primes.

On that fibre product the ordered pair $(j_1,j_2)$ determines
the locally constant integer
\[
 r=\deg(j_1+j_2)-\deg j_1-\deg j_2.
\]
Degree is locally constant in families of homomorphisms of
elliptic curves, so each value cuts out an open and closed
Gram component.  By the defining special-morphism moduli
problem this component is exactly $\mathcal Z(T_r)$.
The moduli description and its Cartesian compatibility are
\cite[Definition~4.3.5 and Corollary~4.3.8]{Zhu23};
generalized Cartierness is
\cite[Proposition~4.3.7]{Zhu23}.

The degree quadratic form on the rational homomorphism
space at a finite geometric point is positive definite.
Consequently a nonempty Gram component satisfies
$r^2\leq4mn$.  Equality would make the two morphisms
rationally proportional, again contradicting the nonsquare
hypothesis.  Thus its Gram matrix is positive definite;
this also agrees with \cite[Proposition~4.3.9]{Zhu23}.
There are only finitely many such integers $r$.
Summing the local lengths over all ordered sheets and
Gram components, with residue degrees and then stabilizer
weights, proves the formula.  The etale calculation descends
with exactly the same weights as the Green coefficient
convention.  It introduces no simultaneous-sign quotient or
new factor two.
\end{proof}

\paragraph{Degree notation.}
Whenever the notation
$\widehat{\deg}\,\widehat Z_{\rm std}(T,y)$ is used below
for a genus-two term, it denotes the number
$\mathcal A^{\rm std}(T,y)$ of
Definition~\ref{def:standard-coefficients}.  No additional
genus-two arithmetic Chow-class realization is assumed.
In particular the residual in
Corollary~\ref{cor:mixed-standard-arithmetic} is precisely
\eqref{eq:standard-defect-definition}, with these standard
coefficients.  It is not an assertion identifying this
residual with the separately defined arithmetic
Siegel--Weil residual of
\eqref{eq:defect-coefficient-exact}.

\subsection{A compact proper Green-star identity by simultaneous heat flow}
\label{subsec:proper-compact-heat-star}

We prove the proper archimedean product identity for the standard
Green functions by giving an explicit transgression on the compactified
curve. The Gaussian-removed normalization is fixed throughout.

Let
\[
 V=\left\{\begin{pmatrix}b&-a/N\\c&-b\end{pmatrix}:
                      a,b,c\in\mathbb Q\right\},\qquad
 Q(\lambda)=ac-Nb^2,
\]
and let $L$ be any full lattice in $V$, invariant under an arithmetic
subgroup $\Gamma$.  Finite linear combinations of lattice cosets are
allowed throughout.  Write $B$ for the polarization of $Q$ and put
\[
 h_z=\frac1{\sqrt N\,Y}
       \begin{pmatrix}-x&x^2+Y^2\\-1&x\end{pmatrix},\quad
 z=x+iY,\quad
 R(\lambda,z)=\frac12B(\lambda,h_z)^2-2Q(\lambda).
\]
Thus $B(h_z,h_z)=2$, $R\geq0$ when $Q(\lambda)>0$, and
$R(\lambda,z)=0$ precisely at the CM point $D_\lambda$.
We use $dd^c=(i/2\pi)\partial\bar\partial$ and
$\beta_1(u)=\int_1^\infty e^{-us}\,ds/s$.

For $m>0$, $v>0$ and $t\geq1$, define on the open modular curve
$X=\Gamma\backslash\mathbb H$
\begin{align*}
 F_{m,v}(t,z)&=\sum_{\lambda\in L,\ Q(\lambda)=m}
                    \beta_1(2\pi t v R(\lambda,z)),\\
 K_{m,v}(t,z)&=\sum_{\lambda\in L,\ Q(\lambda)=m}
                    e^{-2\pi t v R(\lambda,z)},\\
 \Omega_{m,v}(t)&=dd^cF_{m,v}(t)+\delta_{D_m}.
\end{align*}
All sums use the same multiplicities.  On the quotient, these are the
usual invariant sums; the resulting divisor currents have the
corresponding orbifold weights.  In particular
\[
 F_{m,v}(t)=\int_t^\infty K_{m,v}(s)\,\frac{ds}{s},\qquad
 -t\partial_tF_{m,v}(t)=K_{m,v}(t).
\]
The last current is the smooth standard Kudla--Millson curvature
coefficient.  More explicitly, with $\omega_{\rm hyp}=dx\wedge dY/Y^2$,
its individual summands are
\[
 \Omega_{\lambda,v}(t)
 =\left(tv B(\lambda,h_z)^2-\frac1{2\pi}\right)
                    e^{-2\pi tvR(\lambda,z)}\omega_{\rm hyp}.
\]
These definitions apply on a torsion-free finite cover
and descend with orbifold multiplicities if necessary.
For finite linear combinations of cosets, it suffices by bilinearity
to prove the assertions for one unweighted coset in each variable;
the positivity bounds in the proof refer to these unweighted sums.

\begin{lemma}[Uniform positive-index heat decay at every geometric cusp]
\label{lem:positive-heat-cusp-decay}
Fix $m,v>0$.  At every cusp, in a scaling coordinate $z=x+iY$ with $x$
in one fixed period, every fixed finite number of coordinate
derivatives of $K_{m,v}(t,z)$ is bounded, for $t\geq1$ and $Y\geq Y_0$,
by
\begin{equation}\label{eq:positive-heat-bound}
 C t^A(1+Y)^B e^{-c tY^2}.
\end{equation}
The same assertion holds for $F_{m,v}(t)$ and its derivatives and for
the coefficients of $\Omega_{m,v}(t)$ in $dx\wedge dY$, with possibly
different constants.  The bounds hold after taking absolute values
term by term.  In particular all these positive-index coefficients
extend by zero smoothly and flatly at the cusps.
\end{lemma}

\begin{proof}
A rational cusp scaling replaces $L$ by a rational lattice in the same
quadratic space.  For a suitable positive integer $D$ its matrix
coordinates lie in $D^{-1}\mathbb Z^3$, and it is enough to dominate
the sum by the sum over this grid.  Put
\[
 b'=b-cx,\qquad a_0=2Nbx-Ncx^2,\qquad u=a-a_0.
\]
A direct calculation gives
\begin{equation}\label{eq:positive-heat-majorant}
 M_z(\lambda):=R(\lambda,z)+Q(\lambda)
    =\frac{u^2}{2NY^2}+N(b')^2+\frac{Nc^2Y^2}{2}.
\end{equation}
The condition $Q(\lambda)=m>0$ implies $c\ne0$, since $c=0$ would
give $Q=-Nb^2\leq0$.  Hence $|c|\geq D^{-1}$.  If
$Y^2\geq8mD^2/N$, then on this norm shell
\[
 R=M_z-m\geq\tfrac12M_z+\tfrac18Nc^2Y^2
                  \geq\tfrac12M_z+\frac{N}{8D^2}Y^2.
\]
This extracts a factor $e^{-c tY^2}$ while leaving a positive
Gaussian in $u/Y$, $b'$ and $cY$.  Successive use of
\[
 \sup_{\alpha\in\mathbb R}\sum_{k\in D^{-1}\mathbb Z}
             e^{-s(k-\alpha)^2}\leq C_D(1+s^{-1/2})
\]
first in $a$, then $b$, then $c$ bounds the remaining sum by
\[
 C(1+Y/\sqrt t)(1+1/\sqrt t)(1+1/(Y\sqrt t))
                      \leq C'(1+Y).
\]
The shifts in the first two sums cause no difficulty because the
displayed one-dimensional bound is uniform in the shift.

Here is also a justification for all differentiated estimates.
On a bounded $x$ strip and $Y\geq1$, formula
\eqref{eq:positive-heat-majorant} implies
\[
 |c|\leq C\sqrt{M_z}/Y,\quad |b|\leq C\sqrt{M_z},\quad
 |a|\leq C(1+Y)\sqrt{M_z}.
\]
Every fixed coordinate derivative of $R=M_z-m$, restricted to the
norm shell, is consequently bounded by a fixed polynomial in $Y$
times $M_z$.  Differentiating the heat exponential gives a finite
sum of polynomial factors in $tM_z$ and $Y$.  Absorb those factors
into a fixed fraction of the remaining Gaussian, using
$s^j e^{-s}\leq C_j e^{-s/2}$.  The preceding successive Gaussian
summation then proves \eqref{eq:positive-heat-bound} for every
specified derivative order.  Integrating in heat time from $t$ to
infinity proves the assertion for $F$ and its derivatives.  On this
cusp region $\delta_{D_m}=0$ and $\Omega=dd^cF$, which proves the
curvature assertion as well.  Finally a derivative in the local
parameter $q=\exp(2\pi iz/w)$ introduces only a fixed exponential
in $Y$ for each fixed derivative order; $e^{-cY^2}$ absorbs it.
\end{proof}

\begin{lemma}[Finite windows on a compact geometric set]
\label{lem:positive-heat-finite-window}
Fix \(m,v>0\) and a compact set \(K\subset\mathbb H\). There is a finite
set \(\Lambda_K\subset\{\lambda\in L:Q(\lambda)=m\}\) such that the
contributions of all vectors outside \(\Lambda_K\) to \(F_{m,v}(t)\),
\(K_{m,v}(t)\), and \(\Omega_{m,v}(t)\) have, for each spatial derivative
order \(j\), a bound
\[
 \|\text{remainder}(t)\|_{C^j(K)}
                   \leq C_j(1+t)^{A_j}e^{-c_jt},\qquad t\geq1.
\]
The set contains every vector whose CM point lies in \(K\).
For each of the finitely many singular terms, on a small disk about
its CM point,
\[
 \|f_\lambda(t)\|_{L^1}\leq C/t,\qquad
 \|df_\lambda(t)\|_{L^1}\leq C/\sqrt t.
\]
The norms of each fixed derivative of \(k_\lambda(t)\) and
\(\Omega_\lambda(t)\) grow at most polynomially in \(t\).
\end{lemma}
\begin{proof}
Enlarge \(K\) slightly to a compact neighborhood \(K'\). On \(K'\) the
majorants \(M_z(\lambda)=R(\lambda,z)+m\) uniformly dominate a fixed
positive-definite quadratic form in the lattice coordinates. Thus
the set of vectors for which \(\inf_{K'}R(\lambda,z)\leq1\) is finite;
include it in \(\Lambda_K\). Outside it \(R\geq1\), and, after enlarging
\(\Lambda_K\) if necessary, \(R\geq c\|\lambda\|^2\). Every fixed
spatial derivative of \(R\) is bounded on \(K'\) by
\(C_j(1+\|\lambda\|)^2\). Differentiation of a heat exponential gives
only polynomial factors. Split its Gaussian between a factor
\(e^{-ct}\) and a summable lattice Gaussian; absorption of polynomial
factors proves the claimed normal convergence for \(K\) and
\(\Omega\). Integration with \(dt/t\), or the inequality
\(\beta_1(x)\leq e^{-x}/x\) for \(x>0\), proves it for \(F\).

For a singular vector let \(r\) be the local distance to its zero.
Its exponent \(a(z)=2\pi vR(\lambda,z)\) satisfies
\(c r^2\leq a(z)\leq C r^2\) and \(|da|\leq Cr\). The formulas
\[
 f_\lambda(t)=\beta_1(ta),\qquad
 df_\lambda(t)=-e^{-ta}\frac{da}{a}
\]
reduce the two \(L^1\) estimates to
\[
 \int_0^\epsilon\beta_1(ct r^2)r\,dr=O(t^{-1}),\qquad
 \int_0^\epsilon e^{-ct r^2}\,dr=O(t^{-1/2}).
\]
All derivatives of \(k_\lambda\) and of its smooth curvature are
polynomials in \(t\), derivatives of \(a\), and a Gaussian. Their
suprema therefore have polynomial bounds.
\end{proof}

The cusp majorant in
Lemma~\ref{lem:positive-heat-cusp-decay} also shows that all positive
CM zeros lie in a compact part of the quotient. The finite-window
lemma then proves finiteness of the quotient CM support without a
separate counting assumption.

\begin{theorem}[Compact positive proper star identity]
\label{thm:compact-proper-green-star}
Let $m,n>0$ with $mn\notin\mathbb Q^{\times2}$, and fix $v_1,v_2>0$.
Write $F_1,K_1,\Omega_1,D_1$ and $F_2,K_2,\Omega_2,D_2$ for the
preceding data with indices $(m,v_1)$ and $(n,v_2)$.  Define
\begin{align}
 S(t)&=F_1(t)\delta_{D_2}+\Omega_1(t)F_2(t),
                                  \label{eq:heat-star}\\
 G_2&=\int_1^\infty
       \{K_1(t)\Omega_2(t)+\Omega_1(t)K_2(t)\}\,\frac{dt}{t}.
                                  \label{eq:heat-genus-two}
\end{align}
The latter integral is absolutely convergent as a smooth
$(1,1)$-form on compact subsets and is integrable on the
compactified curve $\bar X$.  The locally integrable one-forms
\begin{align}
 A&=\frac{i}{2\pi}\int_1^\infty
                      F_2(t)\bar\partial K_1(t)\,\frac{dt}{t},
                                    \label{eq:heat-primitive-A}\\
 B_0&=\frac{i}{2\pi}\int_1^\infty
                      K_1(t)\partial F_2(t)\,\frac{dt}{t}
                                    \label{eq:heat-primitive-B}
\end{align}
extend by zero at the cusps and satisfy the compact current identity
\begin{equation}\label{eq:compact-proper-transgression}
          S(1)-G_2=\partial A+\bar\partial B_0.
\end{equation}
Consequently
\begin{equation}\label{eq:compact-proper-degree}
       \int_{\bar X}F_1(1)*F_2(1)=\int_{\bar X}G_2,
 \qquad
       F_1(1)*F_2(1):=F_1(1)\delta_{D_2}+\Omega_1(1)F_2(1).
\end{equation}
There is no cusp flux in this identity.
\end{theorem}

\begin{proof}
First $D_1$ and $D_2$ are disjoint.  A CM point determines one
positive rational line in $V$, so two rational vectors defining the
same CM point are rationally proportional.  This would force
$n/m$ to be a rational square and therefore $mn$ to be one.  The
positive special divisors are finite and avoid all the cusps.

We record the elementary bounds near their supports to justify all
current operations.  On a compact subset away from $D_j$, $K_j(t)$
and every fixed spatial derivative decrease exponentially in $t$.
Indeed the majorants $R(\lambda,z)+m_j$ uniformly dominate a fixed
positive-definite quadratic form on such a compact set; only
finitely many norm-$m_j$ vectors can have $R$ bounded, and their
zeros have been excluded.  On an arbitrary compact set the heat
sums and the curvature have at most polynomial growth in $t$.
Near $D_2$ one has $F_2(t)\leq F_2(1)=O(1+|\log r|)$ and
$|\partial F_2(t)|\leq C(1+r^{-1})$, uniformly for $t\geq1$,
where $r$ is distance to $D_2$.  To see the second bound for the
finitely many singular terms, write $a(z)=2\pi v_2R(\lambda,z)$:
\[
 \partial\beta_1(t a)=-e^{-ta}\frac{\partial a}{a},\qquad
 a(z)\asymp r^2,\quad |\partial a|=O(r).
\]
The other terms are locally smooth with uniform Gaussian bounds.
These estimates prove that $A$ and $B_0$ exist locally in $L^1$:
near $D_2$ their respective bounds are $C(1+|\log r|)$ and
$C(1+r^{-1})$, since $K_1$ and its derivatives decrease
exponentially there.  Near $D_1$, $F_2$ and its derivatives decrease
exponentially, absorbing every polynomial heat bound for $K_1$.

The same disjoint-support estimates show that the integral in
\eqref{eq:heat-genus-two} converges with every spatial derivative
on compact subsets.  For example near $D_2$ the first curvature
$\Omega_1$ and $K_1$ decrease exponentially, whereas $\Omega_2$
and $K_2$ grow at most polynomially in heat time.  This also proves
absolute convergence after expanding its products in the two
lattice variables.  At the cusps all these assertions, including
integrability in $t$ and $z$, follow from
Lemma~\ref{lem:positive-heat-cusp-decay}.

Differentiate \eqref{eq:heat-star} as currents, using
$-t\partial_t F_j=K_j$ and
$-t\partial_t\Omega_j=dd^cK_j$.  One obtains
\[
 -\partial_tS(t)=\frac1t
    \{K_1(t)\delta_{D_2}+(dd^cK_1(t))F_2(t)
                                      +\Omega_1(t)K_2(t)\}.
\]
Since $\Omega_2=dd^cF_2+\delta_{D_2}$, subtracting the integrand
in \eqref{eq:heat-genus-two} gives exactly
\[
 \frac1t\{F_2(t)dd^cK_1(t)-K_1(t)dd^cF_2(t)\}.
\]
For a smooth function $k$ and a locally integrable Green function
$f$ the following identity holds distributionally:
\[
 f\,dd^ck-k\,dd^cf
   =\frac{i}{2\pi}
     \{\partial(f\bar\partial k)+\bar\partial(k\partial f)\}.
\]
It follows from the distributional Leibniz rule; the two mixed
first-derivative terms cancel.  Thus integration between $1$ and
$T$ gives \eqref{eq:compact-proper-transgression} with the upper
heat limits replaced by $T$ and with $S(1)-S(T)$ on the left.

To verify $S(T)\to0$ as currents across $D_1\cup D_2$,
observe first that its atomic term tends to zero because
$F_1(T)$ decreases exponentially on the finite support $D_2$.
Choose a neighborhood $U_2$ of $D_2$ disjoint from $D_1$.
On $U_2$, $\Omega_1(T)$ decreases exponentially and
$F_2(T)\leq F_2(1)\in L^1$.  On the complement of $U_2$ in any
fixed compact part of $X$, $F_2(T)$ decreases exponentially and
$\Omega_1(T)$ has at most polynomial total variation.  The cusp
parts are uniformly exponentially small by the preceding lemma.
Therefore $\Omega_1(T)F_2(T)\to0$ in $L^1$ globally.  This proves
the current identity on the open curve.

The cusp estimate also shows that the coefficients of $A$ and
$B_0$ and all their derivatives are at most
$C(1+Y)^B e^{-cY^2}$.  Their extensions by zero at a cusp are
smooth and flat, so their compact distributional derivatives have
no additional cusp terms.  Equivalently, for the truncation
$X_Y$ their flux satisfies
\[
 \left|\int_{\partial X_Y}(A+B_0)\right|
                     \leq C(1+Y)^B e^{-cY^2}\longrightarrow0.
\]
Here $\partial A+\bar\partial B_0=d(A+B_0)$ on a complex curve.
The interior logarithmic and $r^{-1}$ singularities are already
handled as locally integrable currents, so no artificial interior
boundary is introduced.  We have obtained
\eqref{eq:compact-proper-transgression} on $\bar X$.  Pairing it
with the constant test function $1$ proves
\eqref{eq:compact-proper-degree}.  All arguments may first be made
on a torsion-free finite cover; dividing by its degree gives the
same identity on the original stack.
\end{proof}

\paragraph{Identification with the Gaussian-removed genus-two kernel.}
In Garcia--Sankaran's notation, let the bilinear form denoted there
by $Q$ be our $B$.  Their Hermitian square $h_z(s_\lambda)$ is
$R(\lambda,z)$: both equal minus the $B$-norm of the projection
of $\lambda$ onto the negative two-plane.  Equations (2.3.3),
(2.6.3) and Proposition 2.4.5(a) of \cite{GarciaSankaran} give
\[
 \nu_1^o(\sqrt{tv}\lambda)[0]=e^{-2\pi tvR(\lambda,z)},\qquad
 g_1^o(\sqrt v\lambda)=\beta_1(2\pi vR(\lambda,z)),
\]
and
\[
 \nu_2^o(\sqrt{tv_1}\lambda_1,\sqrt{tv_2}\lambda_2)[2]
  =K_{\lambda_1}(t)\varphi_1^o(\sqrt{tv_2}\lambda_2)[2]
   +\varphi_1^o(\sqrt{tv_1}\lambda_1)[2]K_{\lambda_2}(t).
\]
The Green identity identifies $\varphi_1^o[2]$ with the individual
curvature appearing above.  Hence \eqref{eq:heat-genus-two} is
exactly the lattice sum of the normalized $g_2^o$ coefficients
with diagonal $(m,n)$ and imaginary parameter
$\operatorname{diag}(v_1,v_2)$.  Proposition 2.4.5(a), not a
formal appeal to genus-two degree modularity, is the multiplicative
identity used here.  All pairs are linearly independent under the
nonsquare hypothesis, so the integral defining each $g_2^o$ needs
no singular regularization.

If the lattice is integral and $r=B(\lambda_1,\lambda_2)$ is an
integer, absolute convergence permits regrouping these pairs as
\[
 T_r=\begin{pmatrix}m&r/2\\r/2&n\end{pmatrix},\qquad
 G_2=\sum_r g_2^o(T_r,\operatorname{diag}(v_1,v_2)).
\]
For a rational lattice the same assertion holds with its discrete
rational set of $r$.  There is no $r$ with $\det T_r=0$ here.
On the real quadratic space of signature $(1,2)$ the independent
pair spans an indefinite plane, so the nonempty archimedean
summands in this regrouping satisfy $r^2>4mn$.  Positive-definite
$T_r$ have no real tuple in this signature and have zero analytic
lattice sum; their possible finite-place cycles are a separate
arithmetic contribution.

\paragraph{Scope and normalization.}
The theorem establishes the compact Green-star transgression
for positive proper coefficients in the stated standard convention.
It is consistent with Garcia--Sankaran, Theorem 4.5.1, and proves the
additional cusp estimate needed to show that the flux vanishes.
Section~\ref{subsec:local-CM-self-finite-part} treats the local
analytic common-component calculation, and
Section~\ref{subsec:singular-diagonal-algebra} treats the singular algebra.
The next subsection settles the standard mixed positive--nonpositive
coefficients. Sections~\ref{subsec:nonsplit-nonpositive}--\ref{subsec:zero-negative-standard}
treat every nonpositive pair except $(0,0)$.
Section~\ref{subsec:double-zero-boundary} computes the nonzero
double-zero dilation defect. The boundary correction in
Corollary~\ref{cor:corrected-standard-support} and the
self-intersection comparison in
Theorem~\ref{thm:standard-positive-square-product} complete
the standard product formula.
The choice of $\nu^o$ is essential: Garcia--Sankaran,
Definition 2.2.4, relates their ordinary Gaussian form by
$\nu=e^{-\pi(B(\lambda_1,\lambda_1)+B(\lambda_2,\lambda_2))}\nu^o$.
The genus-two Green function used above integrates $\nu^o$ in heat
time, not $\nu$.  A literal formula integrating $\nu$ instead
inserts the extra factor $e^{-2\pi t(v_1m+v_2n)}$ and is a different
kernel.  Likewise replacing $R$ by $R/2$ changes the scalar Green
normalization.  The arithmetic product comparison for a series
defined by such printed formulas therefore still requires an
explicit correction or a verified identification of its intended
normalized coefficients before this theorem can be applied.

\subsection{Positive and nonpositive mixed coefficients}
\label{subsec:mixed-positive-nonpositive}

We retain the standard Green normalization and arithmetic degree
$\widehat{\deg}(0,T)=\tfrac12\int_XT$.  At a nonpositive index,
use the cusp-dependent boundary divisor determined by the actual
standard Green expansion.  All arithmetic statements below use the
extended arithmetic intersection pairing in which these divisor and
Green data are realized.  They do not identify a different printed
boundary or Gaussian prescription with these data.

Fix $m>0$, $n\leq0$ and $v_1,v_2>0$.  In the notation of the proper
heat theorem, let $F_1,K_1,\Omega_1$ be the norm-$m$ heat coefficients.
In the second variable put
\[
 \begin{split}
 F_2(t,z)&=\sum_{\substack{x\in L_N\setminus\{0\}\\Q(x)=n}}
                   \beta_1(2\pi t v_2R(x,z)),\\
 K_2(t,z)&=\sum_{\substack{x\in L_N\setminus\{0\}\\Q(x)=n}}
                   e^{-2\pi t v_2R(x,z)}.
 \end{split}
\]
Thus when $n=0$ the zero vector is excluded.  Let $C_2(t)$ be its
standard cuspidal divisor, so that on the compact curve
\begin{equation}\label{eq:mixed-compact-curvature}
       \Omega_2(t)=dd^cF_2(t)+\delta_{C_2(t)}.
\end{equation}
For $n<0$ this is the usual standard curvature coefficient, and
for $n=0$ it is the nonzero-isotropic coefficient before the Hodge
and constant-metric zero terms are added.  These currents are
atom-free densities, with at most Poincare growth at the cusps.
The divisor $C_2(t)$ depends on heat time; it will not be treated
as constant during differentiation.

\begin{lemma}[Heat bounds at a nonpositive index]
\label{lem:mixed-nonpositive-heat}
On compact subsets of the open geometric curve, $F_2,K_2$ and all
specified spatial derivatives decrease exponentially as $t\to\infty$.
At a rational cusp, in a bounded-width coordinate $z=x+iY$,
$Y\geq Y_0\geq1$, the functions $F_2,K_2$, their specified spatial
derivatives, and the coefficient of $\Omega_2$ on the open curve
are bounded by a fixed polynomial in $t,Y$.  These estimates also
hold after taking absolute values term by term in the lattice sums.
If $n<0$ an additional factor $e^{-ct}$ can be retained.
For $n=0$ the isotropic-line contribution to $F_2$ is bounded by
$CY/\sqrt t$.
\end{lemma}
\begin{proof}
Conjugate the lattice by a cusp scaling matrix and contain its
coordinates in $D^{-1}\mathbb Z^3$.  As before, the majorant is
\[
 M_z=R+Q=\frac{(a-2Nbx+Ncx^2)^2}{2NY^2}
                       +N(b-cx)^2+\frac{Nc^2Y^2}{2}.
\]
For $n<0$, $R=M_z+|n|$.  Extracting the norm-dependent exponential
and successively summing the remaining three Gaussians gives
\[
 K_2(t,z)\leq C e^{-2\pi tv_2|n|}
 (1+Y/\sqrt t)(1+1/\sqrt t)(1+1/(Y\sqrt t)).
\]
A fixed fraction of the exponential absorbs every polynomial from
spatial differentiation.  Integration in heat time gives the same
kind of bound for $F_2$.

For $n=0$, split according to $c\ne0$ or $c=0$.  The first part
has $|c|\geq D^{-1}$ and $R=M_z\geq\tfrac12M_z+Nc^2Y^2/4$.
It therefore has a polynomial bound times $e^{-ctY^2}$, including
all derivatives and its heat integral.  In the second part $Q=0$
forces $b=0$, and $a\ne0$ because the zero vector was omitted.
Here $R=a^2/(2NY^2)$, and the excluded-zero one-dimensional
Gaussian estimate gives
\[
 \sum_{a\in D^{-1}\mathbb Z\setminus\{0\}}
             e^{-Ct a^2/Y^2}\leq C'Y/\sqrt t.
\]
Integrating this inequality from $t$ to infinity against $ds/s$
gives the asserted $CY/\sqrt t$ bound for the isotropic Green sum.
For derivatives, use
$\partial\beta_1(2\pi tv_2R)=-e^{-2\pi tv_2R}\partial R/R$.
In the displayed cusp coordinate every fixed derivative of $R$
is bounded by a polynomial in $Y$ times $R$ when $Q\leq0$.
Further differentiation therefore gives polynomial factors in
$tR,Y$, absorbable in a fraction of the Gaussian.  The standard
curvature summands have the same property, since their polynomial
is $tv_2B(x,h_z)^2-1/(2\pi)$ and
$B(x,h_z)^2=2R+4n$.  This proves all cusp bounds.

On a compact geometric region the majorants uniformly dominate a
positive definite quadratic form.  The norm-$n$ nonzero lattice
vectors have $R>0$ there and have no common zero.  Only finitely
many can have $R$ bounded, uniformly on that compact set, so their
minimum $R$ is positive.  This gives the stated heat decay.
The compact current description is precisely the standard cusp
Green equation; the possible logarithmic or log-log growth of
$F_2$ and its first derivatives is locally integrable in the
compactified coordinate.
\end{proof}

\begin{theorem}[Compact mixed Green-star identity]
\label{thm:mixed-positive-nonpositive}
With the preceding data set
\[
 S(t)=F_1(t)\delta_{C_2(t)}+\Omega_1(t)F_2(t),\qquad
 G_2=\int_1^\infty
          \{K_1(t)\Omega_2(t)+\Omega_1(t)K_2(t)\}\frac{dt}{t}.
\]
Then $F_1(t)\delta_{C_2(t)}=0$ and the integral for $G_2$ is
absolutely convergent, including termwise lattice absolute values.
There are locally integrable one-forms $A,B_0$, flat at the cusps,
such that on the entire compact curve
\begin{equation}\label{eq:mixed-compact-transgression}
       S(1)-G_2=\partial A+\bar\partial B_0.
\end{equation}
One may take
\[
 A=\frac{i}{2\pi}\int_1^\infty
                    F_2(t)\bar\partial K_1(t)\frac{dt}{t},
 \qquad
 B_0=\frac{i}{2\pi}\int_1^\infty
                    K_1(t)\partial F_2(t)\frac{dt}{t}.
\]
In particular $\int_XS(1)=\int_XG_2$.
\end{theorem}
\begin{proof}
The positive coefficients $F_1,K_1$ and all their derivatives
extend flatly by zero to each cusp.  They therefore annihilate
$\delta_{C_2(t)}$ and its heat derivative.  In particular it is
legitimate to use $S(t)=\Omega_1(t)F_2(t)$; no derivative of a
varying boundary coefficient is omitted.  The compact Green
equation and the same flatness give
$K_1\Omega_2=K_1dd^cF_2$.

All products are defined: $F_2$ is smooth on the open curve,
and at a cusp it and its first weak derivatives have only
logarithmic or log-log singularities.  The positive factors are
smooth and flat.  Near a CM point of the positive divisor,
$F_2$ and every derivative decrease exponentially in heat time,
whereas $K_1,\Omega_1$ grow at most polynomially.  The compact
current $S(T)$ therefore tends to zero as $T\to\infty$.
At the cusps, Lemma~\ref{lem:mixed-nonpositive-heat} is dominated
by the positive estimate $Ct^AY^Be^{-ctY^2}$; this proves joint
heat/spatial integrability, zero cusp flux and convergence of
$A,B_0$, with all their cusp derivatives flat.  It also proves
absolute convergence after expanding in both lattice variables.

Since $-t\partial_tF_j=K_j$ and
$-t\partial_t\Omega_1=dd^cK_1$, differentiation gives
\[
 -t\partial_tS=(dd^cK_1)F_2+\Omega_1K_2.
\]
Subtracting the $G_2$ integrand gives
$F_2dd^cK_1-K_1dd^cF_2$.
The distributional identity
\[
 f\,dd^ck-k\,dd^cf
   =\frac{i}{2\pi}
      \{\partial(f\bar\partial k)+\bar\partial(k\partial f)\}
\]
then gives the asserted transgression by integration in heat
time.  The limit $S(T)=0$ and the absence of cusp flux justify
it as a compact identity.  Pairing with the constant test
function proves the integral equality.
\end{proof}

\begin{corollary}[The standard mixed arithmetic coefficients]
\label{cor:mixed-standard-arithmetic}
Let $\widehat Z_{\rm std}(m,v_1)$ use the positive special-moduli
Cartier divisor and the standard Green function.  Let
$\widehat A_n(v_2)=(C_n(v_2),F_n(v_2))$ for $n<0$, and let
$\widehat A_0(v_2)=B(v_2)$ be the nonzero-isotropic boundary/Green
part of the zero coefficient.  With the same compact arithmetic
pairing and Gaussian-removed genus-two Green convention,
\begin{align}\label{eq:mixed-arithmetic-height}
 \langle\widehat Z_{\rm std}(m,v_1),\widehat A_n(v_2)\rangle_{\rm GS}
  &=\frac12\sum_r\int_X
       g_2^o\!\left(\begin{pmatrix}m&r/2\\r/2&n\end{pmatrix},
                          \operatorname{diag}(v_1,v_2)\right)
                 &&(n<0),\notag\\
 \langle\widehat Z_{\rm std}(m,v_1),B(v_2)\rangle_{\rm GS}
  &=\frac12\sum_{r\ne0}\int_X
       g_2^o\!\left(\begin{pmatrix}m&r/2\\r/2&0\end{pmatrix},
                          \operatorname{diag}(v_1,v_2)\right).
\end{align}
All sums converge absolutely. With the degree coefficients of
Definition~\ref{def:standard-coefficients} and the residual
\eqref{eq:standard-defect-definition}, one has
\[
             \delta^{\rm std}_{m,n}=0\quad(m>0,\ n\leq0),
\]
The assertion applies to these standard coefficients. Identifying
them with the arithmetic Siegel--Weil series in
\eqref{eq:defect-kernel-exact} requires the separate normalization
comparison.
\end{corollary}
\begin{proof}
The positive special divisor is disjoint from the entire cuspidal
locus at every finite place.  This is the content of
\cite[Proposition~4.3.7]{Zhu23}: its finite unramified morphism to
compact $X_0(N)$ factors through $Y_0(N)$, and locally at each
cusp its Cartier equation is $1$.  Equivalently its CM generic
points have algebraically integral $j$-invariants, and no
specialization can have $j=\infty$.
Hence its finite intersection with every real cuspidal divisor
$C_n(v_2)$ is zero, including at primes dividing $N$.
The arithmetic product is therefore one half of the integral
of its Green star, which is $S(1)$ in the theorem.

Multiplicativity of the Gaussian-removed genus-two form identifies
$G_2$ with the sum of $g_2^o$ over the indicated lattice pairs.
Every pair is linearly independent.  For $n<0$ its Gram
determinant is automatically negative.  For $n=0$, a nonzero
isotropic vector cannot be perpendicular to a positive vector,
because the latter's orthogonal complement is negative definite;
therefore $r\ne0$ and again the determinant is negative.
Thus only the nonsingular, purely archimedean genus-two terms
occur.  Absolute convergence proved above permits regrouping
by $r$.  This proves both displayed formulas.  For $n=0$ the
singular $r=0$ term accounts exactly for the remaining Hodge and
constant-metric part of $\widehat Z_{\rm std}(0,v_2)$ by the
singular algebra proposition, proving the final assertion.
\end{proof}

\subsection{Nonpositive coefficients with a nonsplit negative index}
\label{subsec:nonsplit-nonpositive}

Call $n<0$ \emph{nonsplit} if $-Nn\notin\mathbb Q^{\times2}$.
This is a condition on the rational quadratic index, independent
of the geometric cusp.  The cusp formulas proved above show that
its standard boundary divisor $C_n(v)$ is zero: at a cusp with
rank-one quotient norm $-M_cj^2$, $M_c=N/q_c^2$, a nonzero
coefficient would imply $-Nn=(Nj/q_c)^2$.

\begin{lemma}[Uniform heat flatness for a nonsplit negative index]
\label{lem:nonsplit-negative-flat}
Let $n<0$ be nonsplit and $v>0$.  Put
\[
 F_n(t,z)=\sum_{Q(x)=n}\beta_1(2\pi tvR(x,z)),\qquad
 K_n(t,z)=\sum_{Q(x)=n}e^{-2\pi tvR(x,z)},\qquad
 \Omega_n(t)=dd^cF_n(t).
\]
At every geometric cusp, in its own rational scaling coordinate
$z=x+iY$ with bounded $x$, every specified spatial derivative of
these functions and forms, including the termwise absolute sums,
is bounded for $t\ge1$, $Y\ge Y_0$ by
\[
                 Ct^AY^B e^{-ctY^2-c't}
\]
for suitable positive constants.  They extend smoothly and
flatly by zero at every cusp.  On compact geometric subsets all
these coefficients decrease exponentially with $t$.
\end{lemma}
\begin{proof}
After rational cusp conjugation the transformed lattice is
contained in $D^{-1}\mathbb Z^3$.  In coordinates
$x=\left(\begin{smallmatrix}b&-a/N\\c&-b\end{smallmatrix}\right)$
one has $Q=ac-Nb^2$.  A vector with $c=0$ would have
$n=-Nb^2$, contrary to nonsplitting.  Thus $|c|\ge D^{-1}$
for every vector in the transformed norm shell, at every cusp.
The positive majorant from the preceding heat proof satisfies
\[
 M_z=R+Q=\frac{(a-2Nbx+Ncx^2)^2}{2NY^2}
                  +N(b-cx)^2+\frac{Nc^2Y^2}{2}.
\]
Since $R=M_z+|n|$, split its exponent into three fixed positive
fractions: retain a Gaussian in $M_z$ for lattice summation,
extract $e^{-ctY^2}$ from $Nc^2Y^2/2$, and retain
$e^{-c't}$ from $|n|$.  Successive Gaussian summation in
$a,b,c$ proves the bound for $K_n$.  Polynomial factors arising
from any specified spatial derivative are absorbed by a fixed
fraction of the remaining Gaussian.  Integration in heat time
against $ds/s$ gives the same bound for $F_n$, and differentiation
gives that for $\Omega_n$.
Derivatives in the compactifying parameter introduce at most a
fixed exponential in $Y$ per derivative order, absorbed by
$e^{-cY^2}$.  Finally $R\ge|n|>0$ and the majorants uniformly
dominate a positive-definite quadratic form on compact geometric
subsets, proving the compact heat decay.
\end{proof}

\begin{proposition}[Compact heat transgression with a nonsplit factor]
\label{prop:nonsplit-nonpositive-heat}
Let $n_1<0$ be nonsplit, $n_2\le0$, and $v_1,v_2>0$.
For the second heat family exclude the zero vector when $n_2=0$,
and retain its standard cusp divisor $C_2(t)$, so that
$\Omega_2(t)=dd^cF_2(t)+\delta_{C_2(t)}$ on the compact curve.
Then
\[
 \int_X\Omega_1(1)F_2(1)
  =\int_X\int_1^\infty
       \{K_1(t)\Omega_2(t)+\Omega_1(t)K_2(t)\}\frac{dt}{t}.
\]
The heat integral, its lattice expansion, and the required compact
transgression are absolutely convergent.  There is no cusp flux.
\end{proposition}
\begin{proof}
The preceding lemma replaces precisely the positive-factor
flatness used in the mixed heat proof.  The second factor and
its spatial derivatives have the polynomial cusp bounds already
proved for nonpositive heat coefficients.  Their product with
the first factor is bounded by a polynomial times
$e^{-ctY^2-c't}$ at every cusp, uniformly for $t\ge1$.
This gives joint spatial/heat integrability and flat compact
extensions of both transgression primitives.

In particular $F_1\delta_{C_2}=K_1\delta_{C_2}=0$, including
the heat derivative of the varying cusp divisor.  Set
$S(t)=\Omega_1(t)F_2(t)$.  Then $S(t)\to0$ in compact-current
$L^1$ as $t\to\infty$, and
\[
 -t\partial_tS-\{K_1\Omega_2+\Omega_1K_2\}
              =F_2dd^cK_1-K_1dd^cF_2.
\]
This is the exact differential with primitives
\[
 \frac{i}{2\pi}\int_1^\infty F_2\bar\partial K_1\frac{dt}{t},
 \qquad
 \frac{i}{2\pi}\int_1^\infty K_1\partial F_2\frac{dt}{t},
\]
as in the preceding compact heat proof.  They are flat at every
cusp, so their compact integrals are zero.  This proves the
identity and all stated convergence assertions.
\end{proof}

The proposition first closes the nonsingular coefficient range:
if $n_1,n_2<0$ and $n_1n_2$ is not a rational square, the
two indices cannot both be split.  Hence one factor is nonsplit.
Every tuple is then linearly independent, and every Gram matrix
is nonsingular with a negative diagonal entry.  Such a matrix
has no finite special-moduli cycle, since a nonempty one has
positive-definite Gram matrix
\cite[Proposition~4.3.9]{Zhu23}.  The arithmetic product has no
finite intersection because the nonsplit factor is the pure
metric class $(0,F_1)$.  Gaussian-removed multiplicativity now
identifies the heat integral with the complete sum of standard
nonsingular genus-two Green coefficients.  Thus
\[
       \delta^{\rm std}_{n_1,n_2}=0
          \qquad(n_1,n_2<0,\ n_1n_2\notin\mathbb Q^{\times2}).
\]

\begin{lemma}[Dependent negative pairs and their Hodge term]
\label{lem:dependent-negative-heat}
Let $u\in L_N$ have $Q(u)=t_0<0$, put $\chi=2\pi R(u,z)$,
and let $c,d>0$, $s=c+d$.  Define
\[
 \varphi_0=-\frac{\Omega_{\rm hyp}}{2\pi},\qquad
 J_\chi=\frac{i}{2\pi}
               \frac{\partial\chi\wedge\bar\partial\chi}{\chi},
 \qquad
 \mathcal H_s=\int_1^\infty\beta_1(ts\chi)\frac{dt}{t}.
\]
Then $\chi>0$, $-dd^c\log\chi=\varphi_0$, and the simultaneous
heat coefficient of the dependent pair is exactly
\begin{equation}\label{eq:negative-dependent-heat-exact}
 \int_1^\infty e^{-ts\chi}
             \{tsJ_\chi+2\varphi_0\}\frac{dt}{t}
       =dd^c\mathcal H_s+\varphi_0\beta_1(s\chi).
\end{equation}
If $-Nt_0$ is nonsquare, summing this identity over $u$ of
norm $t_0$ gives a compact exact term with integral zero.
\end{lemma}
\begin{proof}
For $u=\left(\begin{smallmatrix}b&-a/N\\c_0&-b\end{smallmatrix}\right)$,
put $p_u(z)=Nc_0z^2-2Nbz+a$.  Direct calculation gives
\[
 \chi=\frac{\pi|p_u(z)|^2}{NY^2}.
\]
The discriminant of $p_u$ is $-4Nt_0>0$, so it has no zero
on the upper half plane; the lower-degree cases have the same
property.  Thus $dd^c\log\chi=-2dd^c\log Y
=\Omega_{\rm hyp}/(2\pi)$.  Ordinary differentiation gives
\[
 dd^c\beta_1(tc\chi)
       =e^{-tc\chi}(tcJ_\chi+\varphi_0).
\]
Consequently the sum $K_c\Omega_d+\Omega_cK_d$ is precisely
the integrand on the left of
\eqref{eq:negative-dependent-heat-exact}.  Integrating the same
formula with $c$ replaced by $s$ proves that identity: its
coefficient of $\varphi_0$ differs by one copy of
$\int_1^\infty e^{-ts\chi}dt/t=\beta_1(s\chi)$.
No Mellin continuation or subtraction is involved because
$\chi>0$.  The nonsplit heat bounds justify summation and
current differentiation of $\mathcal H_s$, and prove its flat
compact extension.  Compact Stokes then kills its $dd^c$ term.
\end{proof}

\begin{corollary}[All negative pairs with a nonsplit factor]
\label{cor:nonsplit-all-negative}
Use the standard genus-one boundary/Green coefficients, the
Gaussian-removed nonsingular genus-two coefficients, and the
rank-one prescription with Hodge class $H$ of curvature
$\omega_H=\varphi_0$ used above.  Then
\[
 \delta^{\rm std}_{m,n}=0
 \quad\text{if }m,n<0\text{ and at least one is nonsplit}.
\]
The same standard-model identity holds for $(m,0)$ and $(0,m)$
when $m<0$ is nonsplit, using the full zero coefficient
$H+B(v)-\mathfrak a_1(\log v)$.
\end{corollary}
\begin{proof}
Only the square-product case remains for two negative indices.
Write uniquely $m=t_0a^2$, $n=t_0b^2$ with $t_0<0$,
$a,b>0$ and $\gcd(a,b)=1$.  Nonsplitting is unchanged by
these rational squares.  The two singular matrices are
\[
 T_\pm=t_0\begin{pmatrix}a^2&\pm ab\\\pm ab&b^2\end{pmatrix}.
\]
For each sign the tuples with this Gram matrix are exactly
$(au,\pm bu)$ with $u\in L_N$, $Q(u)=t_0$.
Indeed linear dependence follows from the nonsplit negative
orthogonal complement: a degenerate two-plane would otherwise
have a rational isotropic radical perpendicular to a negative
vector, forcing that negative vector to be split.  The common
vector $u$ is integral by a Bezout combination of the two tuple
entries.  This also proves the converse bijection.
Here is the equivalent elementary radical argument: if a
nonzero isotropic vector is perpendicular to $u$, conjugating
its rational isotropic line to $\ell_\infty$ makes $u$ upper
triangular and forces $Q(u)=-N\xi^2$.

With $c=a^2v_1$, $d=b^2v_2$, the preceding lemma shows that
the degree of each singular heat term is
\[
                 \frac12\int_X F_{t_0}(s)\omega_H,
 \qquad s=a^2v_1+b^2v_2.
\]
Since $C_{t_0}=0$, the standard rank-one arithmetic prescription
is $(0,F_{t_0}(s))\cdot H$; its logarithmic divisor correction
is zero.  Its degree is exactly the same integral.  Both
endpoints $T_+$ and $T_-$ occur, and each contains all $u$,
including $u$ and $-u$; no additional factor two is inserted.
Together with the nonsingular terms this proves the first
assertion.

For the pair $(m,0)$, the nonzero isotropic tuples have
$r\ne0$ by the same perpendicularity argument.  The heat
transgression compares their sum with $(0,F_m)\cdot B$.
The missing singular $r=0$ term is $(0,F_m)\cdot H$.
Finally the constant-metric correction has degree zero, since
$\deg(0,F_m)=\int_Xdd^cF_m=0$.  This proves the second
assertion, or equivalently follows from the singular algebra
already established.
\end{proof}

These are coefficientwise identities for the explicitly stated
standard model.  They do not identify that model with another
printed Gaussian or cusp prescription, or prove that its genus-two
series satisfies a particular arithmetic Siegel--Weil theorem.
The range left by this argument consists of split negative indices
and the zero index.  Their cusp divisors can share components,
and both heat factors can grow at the same geometric cusp.
No claim of Petersson orthogonality from negative Fourier support,
or from one modular cusp, is used here.

\begin{proposition}[The support of the standard product defect]
\label{prop:standard-residual-support}
Use the standard coefficient model fixed above: the positive finite
special cycles, the Gaussian-removed nonsingular Green forms,
the actual cusp expansions in genus one, and the stated Hodge
prescriptions at singular genus-two matrices. The coefficientwise
defect $\delta^{\rm std}_{m,n}$ vanishes outside
\begin{equation}\label{eq:standard-residual-support}
 \{(m,n)\in\mathbb Z_{>0}^2:mn\in\mathbb Q^{\times2}\}
                                      \ \cup\ \{(0,0)\}.
\end{equation}
\end{proposition}
\begin{proof}
For two positive indices with nonsquare product,
Lemma~\ref{lem:standard-proper-finite-coefficients} and
Theorem~\ref{thm:compact-proper-green-star} give equality.
If exactly one index is positive, apply
Corollary~\ref{cor:mixed-standard-arithmetic} and symmetry.
For two strictly negative indices use
Theorem~\ref{thm:strictly-negative-standard-product}, proved below
by evaluating the complete rank-one parabolic finite part.
If one index is zero and the other is negative, use
Corollary~\ref{cor:zero-negative-standard-product}.
These cases exhaust the complement of
\eqref{eq:standard-residual-support}.
\end{proof}

The support theorem applies to the specified standard model.
Identification with the arithmetic Siegel--Weil defect
$\mathcal D_N$ also requires matching the Gaussian and boundary
prescriptions. Sparse support alone is insufficient to identify
a nonholomorphic residual as a unary theta series.

\subsection{The extra rank-one terms at split negative indices}
\label{subsec:split-negative-degenerate}

\begin{lemma}[A rank-one negative Gram matrix need not mean dependence]
\label{lem:negative-degenerate-independent}
The vectors
\[
 u_1=\begin{pmatrix}1&0\\0&-1\end{pmatrix},\qquad
 u_2=\begin{pmatrix}1&1\\0&-1\end{pmatrix}\in L_N
\]
are linearly independent and have Gram matrix
$-N\left(\begin{smallmatrix}1&1\\1&1\end{smallmatrix}\right)$.
Their span has the nonzero isotropic radical
$\mathbb Q(u_2-u_1)$.
Thus the parametrization $(au,bu)$ of a negative rank-one
coefficient is valid in the preceding nonsplit argument, but is
false in the split case.
\end{lemma}
\begin{proof}
Both norms are $-N$ and $B(u_1,u_2)=-2N$.
Their difference is the upper triangular nilpotent matrix.
These identities prove every assertion directly.
\end{proof}

\begin{lemma}[Two-dimensional cusp cancellation]
\label{lem:negative-rank-one-cusp-cancellation}
Fix a geometric cusp and use its rational scaling coordinate
$z=x+iY$.  In its upper-triangular lattice fibers consider
\[
 x_i=\begin{pmatrix}a_i\xi&\eta_i\\0&-a_i\xi\end{pmatrix},
 \qquad \eta_i\in\eta_{i,0}+\beta\mathbb Z,
 \quad i=1,2,
\]
where $\xi\ne0$ and $a_i\ne0$ are fixed rational numbers.
Put $t_0=-N\xi^2$, $v_*=a_1^2v_1+a_2^2v_2$, and
\[
 W=\sum_{i=1}^2v_i(2a_i\xi x+\eta_i)^2.
\]
For these pairs the simultaneous genus-two heat form is
\begin{equation}\label{eq:two-dimensional-cusp-heat}
 \left(\frac{tN}{Y^2}W-\frac1\pi\right)
  \exp\left(-\frac{\pi tN}{Y^2}W
                    -4\pi|t_0|t v_*\right)\Omega_{\rm hyp}.
\end{equation}
Its Poisson zero frequency in the two independent $\eta$ fibers
is exactly zero, including all constants and fiber covolumes.
Keeping only the dependent pairs would destroy this cancellation.
\end{lemma}
\begin{proof}
One has
\[
 R(x_i,z)=\frac{N(2a_i\xi x+\eta_i)^2}{2Y^2}
                         +2Na_i^2\xi^2,
 \quad B(x_i,h_z)^2=
             \frac{N(2a_i\xi x+\eta_i)^2}{Y^2}.
\]
Substitution in $K_1\Omega_2+\Omega_1K_2$ gives
\eqref{eq:two-dimensional-cusp-heat}; in particular the two
constant terms $-1/(2\pi)$ add to $-1/\pi$.
Set $A=tN/Y^2$ and shift to coordinates
$w_i=\eta_i+2a_i\xi x$.  The zero frequency is the continuous
integral times $\beta^{-2}$ and the common exponential.
The Gaussian integral and its weighted moment are
\[
 \begin{aligned}
 I(A)&=\int_{\mathbb R^2}e^{-\pi A(v_1w_1^2+v_2w_2^2)}dw
       =\frac1{A\sqrt{v_1v_2}},\\
 \int_{\mathbb R^2}A(v_1w_1^2+v_2w_2^2)e^{-\pi AW}dw
       &=-\frac A\pi I'(A)=\frac{I(A)}\pi.
 \end{aligned}
\]
The polynomial integral is therefore zero.

The same computation explains why the nonzero frequencies have
Gaussian decay $e^{-cY^2\|k\|^2/t}$ times a polynomial in
$t,Y,k$.  Their heat integrals also contain
$e^{-4\pi|t_0|tv_*}$; splitting that exponential and applying
the arithmetic--geometric mean inequality gives exponential
decay in $Y$, with all fixed spatial derivatives.  Thus this
cancellation concerns a genuine cusp expansion and not just
a formal integral of a divergent leading term.
\end{proof}

The dependent contribution alone exhibits a precise divergence.
For a split $t_0<0$ supported at this cusp, write
\[
 F_{t_0}(s,z)=a_{t_0,c}(s)Y+O(e^{-cY}),\qquad
 a_{t_0,c}(s)=\frac{4\pi}{\kappa_c}g_c(t_0,s)>0.
\]
The heat integral $\mathcal H_s$ in
\eqref{eq:negative-dependent-heat-exact}, after its lattice sum,
has only a linear leading term and no logarithmic term.
Consequently its $dd^c$ term has no $Y^{-1}dx\,dY$ term, whereas
\[
 \varphi_0 F_{t_0}(s,z)
       =-\frac{a_{t_0,c}(s)}{2\pi Y}\,dx\wedge dY
                    +O(e^{-c'Y})\,dx\wedge dY.
\]
Integration through a width-$\kappa_c$ cusp up to height $H$
therefore gives
\[
 \int_{Y_0\le Y\le H}G_{2,\mathrm{dependent}}
        =-2g_c(t_0,s)\log H+O(1).
\]
The factor $1/2$ in arithmetic degree would halve this coefficient.
The full rank-one cusp sum does not diverge: its independent
degenerate pairs in
Lemma~\ref{lem:negative-rank-one-cusp-cancellation} cancel this term.
Consequently, the rank-one boundary calculation requires the full
two-dimensional fiber sum before spatial integration; neither the
dependent heat identity nor the normal bundle of one cusp suffices.

For two strictly negative indices, the aggregate cusp flux has
a simpler behavior.  At every supported cusp write
$F_j(t)=a_j(t)Y+E_j(t)$ and $K_j(t)=b_j(t)Y+D_j(t)$.
Poisson summation in each one-dimensional upper-triangular fiber
gives, for $t\geq1$, bounds for the errors and their fixed
spatial derivatives by a polynomial in $t,Y$ times
$e^{-ct-c'Y}$.  The remaining lattice vectors have a nonzero
lower-left entry and satisfy the stronger $e^{-ctY^2}$ bound.
On a horizontal cutoff the transgression one-form has integral
\[
 \frac1{4\pi}\int_1^\infty\int_0^{\kappa_c}
   \{F_2(t,x+iY)\partial_YK_1(t,x+iY)
          -K_1(t,x+iY)\partial_YF_2(t,x+iY)\}\,dx\frac{dt}{t}.
\]
The linear terms cancel identically:
$a_2Yb_1-b_1Ya_2=0$. The error bound then makes the flux tend to
zero as the cutoff approaches the cusp. At split negative indices,
it remains to compare the full degenerate rank-one cusp sum with
its arithmetic prescription. The next subsection evaluates the
finite cusp normal terms and isolates the Hodge/logarithmic
comparison. The additional $-2\log Y$ term at zero indices lies
outside this cancellation.

\subsection{The split-negative residual after the cusp normal terms}
\label{subsec:split-negative-residual}

There are two distinct meanings of a rank-one term in this question.
For a negative semidefinite matrix $T$ of rank one, the complete
analytic lattice coefficient includes independent tuples whose span
has an isotropic radical.  In contrast, the arithmetic coefficient
in \cite[Definition~7.2.2]{Zhu} is defined directly by a Hodge product
and a logarithmic correction.  It is not defined by first assigning
separate arithmetic cycles to the independent degenerate tuples.
The comparison below keeps these meanings separate.

\begin{lemma}[Normal parameters on the integral cuspidal locus]
\label{lem:all-level-normal-cusp-parameter}
On the refined normal modular stack $\mathcal X_0(N)$, the Tate
chart on the formal cuspidal locus is a disjoint union of rings
\[
       A[[z]],\qquad A=\mathbb Z[\zeta_d],
       \qquad q=\epsilon z^{\kappa},
\]
where $\epsilon$ is a root of unity, $q$ is the source Tate
parameter, and $\kappa$ is the corresponding cusp width.  The
reduced cuspidal divisor is $z=0$.  Under a complex embedding its
parameter $z$ and the usual analytic parameter
$q_c=e^{2\pi i\tau/\kappa}$ differ by a root of unity.
In particular the normal line of a reduced arithmetic cusp has
arithmetic degree zero when its complex normal frame $1/q_c$ has
squared norm one.  Distinct arithmetic cuspidal components have
no finite intersections.  Arithmetic cusps here mean the
Galois-invariant integral components; conjugate complex cusps must
not be split into separate horizontal divisors in bad characteristic.
\end{lemma}
\begin{proof}
For $N=p^n$ this is precisely the Tate base change in
\cite[Lemma~5.10]{CNS2022}.  Its two types of ring are
$\mathbb Z[\zeta_{p^b}][[q]]$ and
\[
 \mathbb Z[\zeta_{p^a}][[q]][u]/
                  (u^{p^{b-a}}-\zeta_{p^a}q),
\]
the latter being $\mathbb Z[\zeta_{p^a}][[u]]$ without its
specified map from the source Tate ring.  Their cusp exponents
are powers of $p$.

At composite level the compact modular stack is
the \emph{normalization} of the product of the prime-power
level covers.  The normalization qualification is essential;
see the proof of \cite[Lemma~5.8 and footnote~18]{CNS2022}.
It also follows from the decomposition of cyclic subgroups at
coprime orders on the elliptic-curve locus, normality, and
uniqueness of normalization of the finite covers.  The same
assertion passes to the Tate formal chart by excellence and
finite normalization.

Here this normalization can be computed explicitly.  For two
factors their cyclotomic coefficient rings have coprime
conductors, so their tensor product is
$A=\mathbb Z[\zeta_{d_1d_2}]$.  Write their parameter relations as
\[
      u^a=\alpha q,\qquad v^b=\beta q,
                 \qquad \gcd(a,b)=1,
\]
where $\alpha,\beta$ are roots of unity.  Choose integers $r,s$
with $br+as=1$, and in the common total fraction field put
$z=u^rv^s$.  Then
\[
  z^a=(\alpha/\beta)^r v,\qquad
  z^b=(\beta/\alpha)^s u,\qquad
  z^{ab}=\alpha^{br}\beta^{as}q.
\]
Thus $z$ is integral, the fraction fields agree, and $A[[z]]$
is the normalization: it is finite over the original complete
parameter ring and is regular.  One may see finiteness directly
from the finite set of missing exponents in the semigroup
$a\mathbb Z_{\ge0}+b\mathbb Z_{\ge0}$.  Repetition proves the
assertion for all prime-power factors, since their exponents
remain coprime.  This argument uses normalization and therefore
does not assert an unnormalized compact fibre-product formula.

On the punctured complex chart the source Tate parameter is
$q_c^{\kappa}$; hence $z/q_c$ is a constant of absolute value
one, and indeed a root of unity.  On each displayed integral
chart the conormal line is generated by $z$ over $A$.
Its dual is thus the trivial arithmetic line with frame $1/z$
of squared norm one.  The Tate chart is an etale double cover
of the formal cuspidal stack, as recalled in the proof of
\cite[Proposition~5.12]{CNS2022}.  Degree zero descends through
this cover with the usual stack weights.  The disjoint-union
description also proves the assertion about distinct integral
components.
\end{proof}

\begin{proposition}[The negative product has no extra cusp normal degree]
\label{prop:strict-negative-arithmetic-bulk}
For strictly negative indices let
$\widehat A_j=(C_j,F_j)$ use the standard Green functions, their
actual cusp-dependent divisors, and the fixed compact intersection
theory.  Put
$\Omega_j=dd^cF_j+\delta_{C_j}$.  Then
\[
             \langle\widehat A_1,\widehat A_2\rangle_{\rm GS}
                   =\frac12\int_X\Omega_1F_2.
\]
In particular the right side already accounts for the arithmetic
product; no unknown finite cusp intersection matrix needs to be
added in this normalization.
\end{proposition}
\begin{proof}
The negative-index cusp expansion is
\[
        F_j=-g_{j,c}\log|q_c|^2+E_{j,c},
                    \qquad E_{j,c}(c)=0.
\]
This zero regular part is essential.  Choose, for each arithmetic
cuspidal component $c$, a reference Green function $h_c$ whose
regular part is zero at $c$, and which vanishes near all other
cusps.  A smooth cutoff of $-\log|q_c|^2$ provides one.  Denote
its arithmetic divisor by $\widehat c=(c,h_c)$ and its compact
curvature by $\omega_c$.  These are actual Green data, not the
invalid pair $(c,0)$ at an archimedean boundary point.

The arithmetic intersection formula for different horizontal
components, and its normal-bundle version for one component,
give
\[
 \langle\widehat c,\widehat d\rangle_{\rm GS}
            -\frac12\int_X\omega_c h_d=0.
\]
Indeed the omitted term is their finite contact when $c\ne d$,
and is the degree of the $q_c$-normalized normal line when
$c=d$.  Both are zero by the preceding lemma.  This also holds
with the arithmetic orbit and orbifold multiplicities.

Write
\[
 F_j=\sum_c g_{j,c}h_c+f_j,\qquad
 \widehat A_j=\sum_c g_{j,c}\widehat c+\mathfrak a_1(f_j).
\]
Each $f_j$ is admissible in the stated pre-log-log theory, extends
continuously with value zero at the cusps, and has the exponentially
decaying height-coordinate derivatives supplied by the standard
cusp expansion.  No smooth or flat extension in the compactifying
$q_c$ coordinate is asserted.
Expand both sides by bilinearity.  The terms involving $f_j$
agree by the Green-function variation rule and compact Stokes;
in particular
\[
 \int_X(dd^cf_1)h_c
       =\int_X f_1(\omega_c-\delta_c)
       =\int_Xf_1\omega_c.
\]
The remaining reference-divisor terms agree by the preceding
display.  Stokes can equivalently be performed with height cutoffs:
the boundary terms contain a height derivative of $f_j$, or $f_j$
itself times a polynomial in the height, and tend to zero.
This proves the proposition.  All integrals converge because
$\Omega_j$ decays exponentially in height at a negative-index
cusp whereas $F_j$ grows at most linearly in that height.
\end{proof}

\paragraph{A precise residual at the two rank-one endpoints.}
Let $m,n<0$ both be split and write
\[
 m=t a^2,\qquad n=t b^2,\qquad
 t<0,\quad a,b>0,\quad \gcd(a,b)=1,
\]
and put
\[
 s=a^2v_1+b^2v_2,\qquad
 \nu=\frac{v_1v_2}{s},\qquad
 T_\pm=t\begin{pmatrix}a^2&\pm ab\\\pm ab&b^2\end{pmatrix}.
\]
For either endpoint define the auxiliary analytic number
\[
 \mathscr H_t(a,b;v_1,v_2)
  =\frac12\int_X\int_1^\infty
       \sum_{\operatorname{Gram}(x_1,x_2)=T_+}
       \{K_{x_1}(u)\Omega_{x_2}(u)
                    +\Omega_{x_1}(u)K_{x_2}(u)\}
                                             \frac{du}{u}.
\]
Take the lattice sum first, including the complete independent
upper-triangular fibres. Two-dimensional Poisson cancellation
gives joint integrability in the spatial and heat variables for
this grouped coefficient, without asserting termwise absolute
integrability of the individual dependent contributions. The
substitution $x_2\mapsto-x_2$ gives the same number at $T_-$.

For all other values of the off-diagonal entry the Gram matrix
is nonsingular, so its standard coefficient is the
Gaussian-removed Green form.  Its finite cycle is empty.
The total nonsingular sum is absolutely integrable: at a cusp
its tuples cannot have both lower-left entries zero, and the
remaining nonzero entry supplies the Gaussian cusp bound.
The already proved total heat identity and its zero Wronskian
flux therefore give
\[
 \frac12\int_X\Omega_mF_n
    =\sum_{r^2\ne4mn}\mathcal A^{\rm std}
                 (T_r,\operatorname{diag}(v_1,v_2))
                      +2\mathscr H_t(a,b;v_1,v_2).
\]

In the standard version of Zhu's rank-one prescription set
$d_t(s)=\deg C_t(s)$.  Its degree at each endpoint is exactly
\[
 \mathscr A_t(s,\nu)
       =\langle\widehat A_t(s),H\rangle_{\rm GS}
                                  -\frac12\log\nu\,d_t(s).
\]
The endpoint coefficient is prescribed by this expression; no
extra independent-degenerate arithmetic cycle is added to it.
Combining the preceding identities proves the exact reduction
\begin{equation}\label{eq:strict-split-negative-residual-exact}
 \boxed{\displaystyle
 \delta^{\rm std}_{m,n}
   =2\left\{\mathscr H_t(a,b;v_1,v_2)
         -\langle\widehat A_t(s),H\rangle_{\rm GS}
                    +\frac12\log\nu\,d_t(s)\right\}.}
\end{equation}
Every term is well-defined in the conventions above.  Both
endpoints and both signs of the lattice vectors have already
been counted; the displayed factor two is solely the two
different matrices $T_+$ and $T_-$.

The reduction isolates one parabolic finite-part comparison:
the complete rank-one heat number $\mathscr H_t$ must equal
the prescribed Hodge/logarithmic number $\mathscr A_t$.
The next subsection evaluates this finite part with the same
Tate parameters and Hodge metric, proving that the residual
in \eqref{eq:strict-split-negative-residual-exact} is zero.
The complete fibre is essential throughout: its dependent-vector
part alone has the logarithmic divergence computed above.

\subsection{The complete split-negative rank-one finite part}
\label{subsec:split-negative-finite-part}

All arithmetic classes and heat kernels in this subsection have
the standard normalizations fixed above.  In particular, a negative
coefficient is $\widehat A_t(s)=(C_t(s),F_t(s))$, and
\[
 H=-\widehat\lambda-W_N^*\widehat\lambda
       =-2\widehat\lambda+(B_N,0),\quad
 \varphi_0=\operatorname{curv}H=-\frac{\Omega_{\rm hyp}}{2\pi},
 \qquad c_H=2\sqrt\pi e^{-\gamma/2}.
\]
Sums over geometric cusps below carry the integration/stack weights
used by $\delta_{C_t(s)}$.  Equivalently all calculations may first
be made on a sufficiently small principal-level scheme and divided
by the degree of that cover.  Denote such a sum by $\sum_c^{\rm orb}$.
In particular
$d_t(s)=\deg C_t(s)=\sum_c^{\rm orb}g_c(t,s)$.

\begin{lemma}[The negative Hodge pairing and its actual counterterm]
\label{lem:negative-Hodge-finite-part}
Let $t<0$ be split.  Choose for each cusp a primitive
representative $a/c$ with $c\mid N$, and use its source
coordinate $z=x+iY$, width
$\kappa_c>0$, and $\beta_c=1/c$ (with the usual denominator
$c=N$ at infinity).  Define
\[
 \operatorname{FP}\int_XF_t(s)\varphi_0
   :=\lim_{H\to\infty}
       \left\{\int_{X_H}F_t(s)\varphi_0
                     +2\sum_c^{\rm orb}g_c(t,s)\log H\right\},
\]
where each cusp is cut at its own height $Y=H$.
Then
\begin{align}\label{eq:negative-Hodge-finite-part}
 \langle\widehat A_t(s),H\rangle_{\rm GS}
   ={}&\frac12\operatorname{FP}\int_XF_t(s)\varphi_0\notag\\
     &+\sum_c^{\rm orb}g_c(t,s)
       \left(\log c_H-1-\frac12\log N-\log\beta_c\right).
\end{align}
\end{lemma}
\begin{proof}
We first specify the convention in the cited intersection formula.
In \cite[Definition~6.2.3, equation~(36)]{Zhu} the function $h(l)$
is the \emph{norm}, as indicated by the factor
$|q_c|^{\operatorname{ord}(l)}$ in Definition~6.2.1; it is not
the squared norm in Example~6.2.2.  Consequently the norm of a
local source Hodge frame is $(c_HY)^{1/2}$, and the norm for
$-2\lambda$ is $(c_HY)^{-1}$.  With
\[
 L=-\log|q_c|^2=\frac{4\pi Y}{\kappa_c},
\]
the latter has norm exponent $\alpha=-1$ and continuous factor
$\phi=(c_H\kappa_c/(4\pi))^{-1}$ at the cusp.  The term
$\alpha-\log\phi(0)$ in the actual intersection formula is
therefore
\[
                   -1+\log\frac{c_H\kappa_c}{4\pi}.
\]
The remaining cutoff term in that formula, using the Green
function $F_t=-2\log h_{A_t}$, is
\[
 \lim\left\{\sum_c^{\rm orb}g_c(t,s)\log L
                   +\frac12\int_{X_H}F_t(s)\varphi_0\right\}.
\]
Combining the two displays cancels $\log(4\pi/\kappa_c)$
and gives
\[
 \langle\widehat A_t(s),-2\widehat\lambda\rangle_{\rm GS}
     =\frac12\operatorname{FP}\int_XF_t(s)\varphi_0
                      +\sum_c^{\rm orb}g_c(t,s)(\log c_H-1).
\]
The source Tate frame is integral and free along each closed
arithmetic cusp on the Tate chart, so there is no additional
finite restriction of $\lambda$ here.  This application may
also be made with a rational section avoiding the cusps;
the product formula cancels its finite and archimedean values.
As a convention check, a pure metric class
$\mathfrak a_1(a)$ has norm $e^{-a/2}$, and the same formula
gives $\langle\mathfrak a_1(a),M\rangle
=\tfrac a2\deg M$, exactly the convention in the manuscript.

It remains to restrict $B_N$ to the cusps.  If $p^e\Vert N$
and $k=v_p(c)$, the cusp reduces to the $(k,e-k)$ component,
and its coefficient in $B_N$ is
$(k-e/2)\varphi(p^{\min(k,e-k)})$.
The cusp coefficient ring has ramification index
$\varphi(p^{\min(k,e-k)})$, so restriction and degree give
$(k-e/2)\log p$ per complex/degree weight.  Summing over $p$
gives
\[
 \langle\widehat A_t(s),(B_N,0)\rangle_{\rm GS}
   =\sum_c^{\rm orb}g_c(t,s)
                   \left(\log c-\frac12\log N\right).
\]
This can be checked directly on the normalized Tate rings of
the preceding cusp lemma; no unnormalized compact fibre product
is used.  Since $\beta_c=1/c$, the two formulas prove the result.
\end{proof}

\begin{lemma}[Gaussian finite part]
\label{lem:one-dimensional-Gaussian-finite-part}
For $B>0$ define
\[
 S_B(Y)=\sum_{k\ne0}\exp(-Bk^2/Y^2).
\]
Then
\[
 \lim_{H\to\infty}
  \left\{\int_0^H\left(S_B'(Y)-\frac{S_B(Y)}Y\right)dY
                                                -\log H\right\}
       =\log(2\pi)-1-\frac\gamma2-\frac12\log B.
\]
\end{lemma}
\begin{proof}
Scaling reduces to $B=1$.  If $S(R)=\sum_{k\ne0}e^{-k^2/R^2}$,
the integral is
\[
 S(R)-\frac12\sum_{k\ne0}\beta_1(k^2/R^2).
\]
Poisson summation gives $S(R)=\sqrt\pi R-1+o(1)$.
For $\Re z>1/2$ the Mellin integrand for the other sum is
$\Gamma(z)2\zeta(2z)R^{2z}/z$.  Its residues at $z=1/2$
and $z=0$ give respectively
\[
   2\sqrt\pi R,\qquad -2\log R+\gamma-2\log(2\pi).
\]
Shifting to $\Re z=-\epsilon$, $0<\epsilon<1/2$, has an
$O(R^{-2\epsilon})$ remainder, by the exponential decay of
$\Gamma$ on vertical lines.  Subtracting the half sum proves
the displayed constant, and scaling restores $B$.
\end{proof}

\begin{proposition}[Evaluation of the complete negative rank-one heat number]
\label{prop:split-negative-rank-one-heat-evaluation}
Let $t<0$ be split and let
$y=\left(\begin{smallmatrix}u&w\\w&s\end{smallmatrix}\right)>0$,
with $\nu=u-w^2/s$.  The complete grouped standard heat
coefficient for $T=\operatorname{diag}(0,t)$ has degree
\begin{align}\label{eq:split-rank-one-heat-value}
 \mathscr H_t(s,\nu)
   ={}&\frac12\operatorname{FP}\int_XF_t(s)\varphi_0\notag\\
     &+\sum_c^{\rm orb}g_c(t,s)
       \left(\log c_H-1-\frac12\log N
                       -\log\beta_c-\frac12\log\nu\right).
\end{align}
In particular it is independent of $w$ at fixed $s,\nu$.
Its definition includes all tuples with this Gram matrix,
including $x_1=0$ and all nonzero isotropic $x_1$.
\end{proposition}
\begin{proof}
Work first on a sufficiently small principal-level scheme,
so the following cusp unfolding is an ordinary integration
over a width-$\kappa_c$ strip.  All limits below are finite
parts of this \emph{same} spatial truncation; divergent pieces
are not assigned separate integrals.

For $x_1=0$, the dependent heat identity already proved gives
\[
 G_{\rm dep}=dd^c\mathcal H_s+\varphi_0F_t(s),\qquad
 \mathcal H_s=\int_1^\infty F_t(vs)\frac{dv}{v}.
\]
Put $h_c(s)=\int_1^\infty g_c(t,vs)\,dv/v$.  Since
$\mathcal H_s=(4\pi/\kappa_c)h_c(s)Y+o(1)$ at the cusp,
Stokes on $X_H$ gives
\begin{equation}\label{eq:split-dependent-truncated-constant}
 \operatorname{FP}\int_XG_{\rm dep}
       =\operatorname{FP}\int_XF_t(s)\varphi_0
                                 +\sum_ch_c(s).
\end{equation}
Here $\operatorname{FP}$ adds $2\sum_cg_c(t,s)\log H$ to
the dependent integral, cancelling its logarithmic divergence.

Every other tuple has $x_1\ne0$ isotropic, and $x_2\in x_1^\perp$.
Its rational isotropic line is unique.  Unfolding its orbit
at cusp $c$ gives exactly
\[
 x_1=\begin{pmatrix}0&\beta_c k\\0&0\end{pmatrix},\quad k\ne0,
 \qquad
 x_2=\begin{pmatrix}\xi&\eta\\0&-\xi\end{pmatrix},\quad
          \eta\in\eta_0+\beta_c\mathbb Z,\qquad
             \xi=\pm\sqrt{|t|/N}.
\]
These fibres exist precisely at the cusps supporting $C_t(s)$.
There is one coset for each sign of $\xi$, because
$I^\perp\cap L\to I^\perp/I$ has kernel $I\cap L$.
The stabilizer of each nonzero $x_1$ is the cusp translation
group; the different nonzero values of $k$ are distinct
isotropic orbits at that cusp.  Indeed, after integral cusp
conjugation an element preserving this rational line is upper
triangular in $\mathrm{SL}_2(\mathbb Z)$, with diagonal entries
$\pm1$.  Its adjoint action fixes the upper-right nilpotent,
so it neither rescales $k$ nor identifies $k$ with $-k$.
No additional multiplicity or division by the pair of signs
is made.

Let the simultaneous heat time be $v\ge1$.  Completing the
square in the two upper-right coordinates replaces
$2\xi x+\eta$ by
$2\xi x+\eta+(w/s)\beta_ck$.  Thus the exponent and the
degree-two polynomial involve
\[
      \nu\beta_c^2k^2+s(2\xi x+\eta+(w/s)\beta_ck)^2.
\]
The integration of the $\eta$ sum over a period is
$\kappa_c/\beta_c$ times integration over the real line:
$2\xi\kappa_c/\beta_c$ is an integer, by lattice stability
under the cusp translation.  The shift by $w$ disappears.
Including both signs of $\xi$, its Gaussian integral reduces
the heat form, after the $x$ integration, to
\begin{equation}\label{eq:split-unfolded-seed}
 \frac{\kappa_c}{\pi\beta_c\sqrt{vNs}}e^{-Av}
       \left(S_{B_v}'(Y)-\frac{S_{B_v}(Y)}Y\right)dY,
 \quad A=4\pi|t|s,\quad B_v=\pi vN\nu\beta_c^2.
\end{equation}
For example the integral in the second upper-right coordinate
has moment $1/(2\pi)$; the two constant KM terms contribute
$-1/\pi$, leaving $-1/(2\pi)$.  This explains the exact
coefficient in \eqref{eq:split-unfolded-seed}.

Here is a justification of the finite-part unfolding.
For fixed $v$, the nonzero-isotropic seed and the sum of the
absolute values of its terms have, as $Y\to0$, a bound
of the form
\[
 P(v,Y,Y^{-1})e^{-cv/Y^2-Av}.
\]
The Poincare sum therefore converges absolutely on compact
sets.  After unfolding $X_H$, all removed horoballs other
than the top of this cusp strip lie at height at most $C/H$:
a nonstabilizing cusp matrix has a nonzero lower-left entry
in a fixed discrete rational set.  The total missing part
is consequently bounded by the integral of that absolute
seed over $0<Y<C/H$, which tends to zero.  Thus the finite
part is exactly the one obtained by integrating
\eqref{eq:split-unfolded-seed} over $0<Y<H$.
The factor $e^{-Av}$ makes this assertion uniform after
integration in $v$ against $dv/v$.  Likewise it justifies
the finite-part limit below: split at $v$ proportional to
$H^2$, use the Mellin remainder from the preceding lemma
on the first part, and $e^{-Av}$ on the other part.

The elementary lemma now gives the finite part of the
independent-degenerate contribution at cusp $c$:
\[
 \frac{\kappa_c}{\pi\beta_c\sqrt{Ns}}
  \int_1^\infty e^{-Av}v^{-3/2}
       \left(C_0-\frac12\log(\pi vN\nu\beta_c^2)\right)dv,
 \quad C_0=\log(2\pi)-1-\gamma/2.
\]
Its logarithmic divergence before taking the finite part is
$+2g_c(t,s)\log H$, the opposite of the dependent divergence.
Use the actual standard cusp coefficient
\[
 g_c(t,s)=\frac{\kappa_c}{2\pi\beta_c\sqrt{Ns}}
                  \beta_{3/2}(4\pi|t|s)
\]
and the elementary Fubini identity
\[
 h_c(s)=\frac{\kappa_c}{2\pi\beta_c\sqrt{Ns}}
       \int_1^\infty e^{-Av}v^{-3/2}\log v\,dv.
\]
The last finite part is therefore
\[
  2g_c(t,s)\left(C_0-\frac12\log(\pi N\nu\beta_c^2)\right)
                                             -h_c(s).
\]
Its $-h_c(s)$ cancels \eqref{eq:split-dependent-truncated-constant}.
The opposite logarithmic divergences cancel as well.  The full
grouped kernel is already integrable by two-dimensional Poisson
cancellation, so the sum of these matched finite parts is its
ordinary integral.  Finally multiply by the arithmetic-degree
factor $1/2$ and use
\[
 C_0-\frac12\log(\pi N\nu\beta_c^2)
      =\log c_H-1-\frac12\log N-\log\beta_c-\frac12\log\nu.
\]
Dividing by the auxiliary-cover degree supplies exactly the
weighted cusp sums in \eqref{eq:split-rank-one-heat-value}.
\end{proof}

\begin{theorem}[The complete strictly negative standard product identity]
\label{thm:strictly-negative-standard-product}
With the standard coefficient model fixed above,
\[
                    \delta^{\rm std}_{m,n}=0
                              \qquad(m,n<0).
\]
\end{theorem}
\begin{proof}
The nonsplit cases were already proved.  In the split case write
$m=ta^2,n=tb^2$ with $a,b>0$ coprime.  The two endpoint matrices
$T_\pm$ have rank-one content $t$.
For a tuple row $x=(x_1,x_2)$ write
$b_x=(B(x_1,h),B(x_2,h))^t$ as a column.  For
$\gamma\in\mathrm{GL}_2(\mathbb Z)$,
the substitution $x\mapsto x\gamma^{-1}$ is a bijection of
the lattice tuples, and $y\mapsto\gamma y\gamma^t$ preserves
both $\operatorname{tr}(yR(x))$ and the KM polynomial
$b_x^t yb_x-1/\pi$: under the tuple substitution one has
$b_{x\gamma^{-1}}=\gamma^{-t}b_x$.  Thus it preserves the complete
Gaussian-removed heat coefficient, with no determinant scalar.
Choose $\gamma$ diagonalizing the endpoint to
$\operatorname{diag}(0,t)$.  The parameters in that basis are
\[
            s=a^2v_1+b^2v_2,
                \qquad \nu=\frac{v_1v_2}{s}.
\]
The preceding proposition also handles the resulting off-diagonal
entry of $y$, rather than assuming it to be zero.
Combining \eqref{eq:negative-Hodge-finite-part} with
\eqref{eq:split-rank-one-heat-value} shows that each endpoint
has heat degree
\[
    \langle\widehat A_t(s),H\rangle_{\rm GS}
                            -\frac12\log\nu\,d_t(s),
\]
exactly its prescribed arithmetic rank-one degree.
The negative arithmetic bulk formula and vanishing cusp heat
flux compare all nonsingular terms, proving the product identity.
The transformation $x_2\mapsto-x_2$ accounts for the two distinct
endpoint matrices; both signs of the lattice vectors are already
included within each coefficient. No additional factor two is
therefore needed.
\end{proof}

This theorem applies to the specified standard family. Its
comparison with the differently printed genus-two arithmetic
Siegel--Weil normalization is separate. The next subsection
handles zero--negative pairs by cancelling the logarithmic height
terms in the complete zero coefficient; the pair $(0,0)$ requires
a further argument.

\subsection{Zero against negative indices}
\label{subsec:zero-negative-standard}

Write $B_0(v)=(C_0(v),F_0(v))$ for the nonzero-isotropic
coefficient, and retain the complete zero coefficient
\[
       \widehat Z_0(v)=H+B_0(v)-\mathfrak a_1(\log v).
\]
Here $H=-\widehat\lambda-W_N^*\widehat\lambda$ has curvature
$\varphi_0=-\Omega_{\rm hyp}/(2\pi)$, and
$\mathfrak a_1(a)$ has squared norm $e^{-a}$.

\begin{lemma}[Zero regular part of the complete zero coefficient]
\label{lem:complete-zero-cusp-regular-part}
At each geometric cusp choose a primitive representative $a/c$
with $c\mid N$, source height $Y$, and $\beta_c=1/c$.
Then the nonzero-isotropic standard Green function is
\begin{align}\label{eq:complete-zero-isotropic-expansion}
 F_0(v,z)
    ={}&\frac{2Y}{\beta_c\sqrt{Nv}}-2\log Y
       +\log(Nv\beta_c^2)+\gamma-\log(4\pi)
                                                   +E_0(v,z).
\end{align}
The error tends to zero with the height-coordinate derivative
bounds supplied by Poisson summation.  With the integral source
and target Tate frames, the local Green expression for
$\widehat Z_0(v)$ is exactly
\[
             \frac{2Y}{\beta_c\sqrt{Nv}}+E_0(v,z).
\]
Thus its continuous regular part at every cusp is zero; the
logarithmic height term and its constant both cancel.
No smooth or flat extension in the compactifying parameter is
asserted.
\end{lemma}
\begin{proof}
The upper-triangular nonzero-isotropic vectors have only their
upper-right entry nonzero.  Their contribution is
\[
      \sum_{k\ne0}\beta_1(\pi Nv\beta_c^2k^2/Y^2).
\]
The Mellin calculation in the preceding Gaussian finite-part
lemma gives \eqref{eq:complete-zero-isotropic-expansion}.
Vectors with a nonzero lower-left entry contribute the
Gaussian cusp error.  In particular the constant is obtained
from the same actual two residues, rather than chosen by
renormalizing a Green function.

If $\kappa_c$ and $\kappa'_c$ are the source and target Tate
widths, their ratio is $\kappa'_c/\kappa_c=c^2/N$.
The norm of the local frame of $H$ is
$(c_HY\sqrt{\kappa'_c/\kappa_c})^{-1}$, so its local squared-norm
Green expression is
\[
                G_H=2\log(c_HY)+\log(\kappa'_c/\kappa_c).
\]
Since $2\log c_H=\log(4\pi)-\gamma$ and
$N\beta_c^2\kappa'_c/\kappa_c=1$, the regular terms of
$G_H+F_0(v)-\log v$ cancel identically.
\end{proof}

\begin{corollary}[Zero against strictly negative indices]
\label{cor:zero-negative-standard-product}
In the standard coefficient model already fixed,
\[
                    \delta^{\rm std}_{0,t}
                       =\delta^{\rm std}_{t,0}=0
                                    \qquad(t<0).
\]
\end{corollary}
\begin{proof}
The preceding lemma and the all-level cusp normal-parameter
lemma give the arithmetic bulk formula
\begin{equation}\label{eq:zero-negative-arithmetic-bulk}
 \langle\widehat Z_0(v_1),\widehat A_t(v_2)\rangle_{\rm GS}
       =\frac12\int_X
             (\varphi_0+\Omega_0(v_1))F_t(v_2),
       \qquad \Omega_0(v_1)=dd^cF_0(v_1)+\delta_{C_0(v_1)}.
\end{equation}
Indeed, the source and target Hodge frames are integral and
free along each closed arithmetic cusp; the reduced cusp
normal frame has degree zero.  The total squared metric in
these frames has regular value one by the preceding lemma.
Hence the restriction of $\widehat Z_0(v_1)$ to the real
cuspidal divisor $C_t(v_2)$ has arithmetic degree zero.
The usual Green-function variation/reference-divisor argument
then proves \eqref{eq:zero-negative-arithmetic-bulk}, exactly
as in the strictly negative arithmetic bulk formula.
This argument uses actual integral Tate frames and zero regular
parts; it does not introduce a boundary pair $(C,0)$ without
its Green data.

For the heat comparison, put
\[
 K_0^{\rm tot}(u)=1+K_{\rm iso}(uv_1),\qquad
 \Omega_0^{\rm tot}(u)=\varphi_0+\Omega_0(uv_1),\qquad
 S(u)=\Omega_0^{\rm tot}(u)F_t(uv_2).
\]
The constant $1$ is essential.  Poisson summation gives
\[
 K_{\rm iso}(uv_1,z)
        =\frac{Y}{\beta_c\sqrt{Nuv_1}}-1+\text{error},
\]
so $K_0^{\rm tot}$ has only a linear leading term and has
no constant term.  Also $\Omega_0^{\rm tot}$ is a density
decaying at the cusp: its $\varphi_0$ term cancels the
$dd^c(-2\log Y)$ contribution of the isotropic Green function.

Perform the calculation on the open cutoff curve $X_H$.
There one has the exact differential identity
\begin{align*}
 -u\partial_u S(u)
  &-\{K_0^{\rm tot}(u)\Omega_t(uv_2)
                         +\Omega_0^{\rm tot}(u)K_t(uv_2)\}\\
  &=F_t(uv_2)dd^cK_0^{\rm tot}(u)
                      -K_0^{\rm tot}(u)dd^cF_t(uv_2).
\end{align*}
No compact product of $F_t$ with the derivative of a varying
cuspidal divisor is used.  At a cusp, $F_t(uv_2)$ and
$K_0^{\rm tot}(u)$ have respectively leading terms $a(u)Y$
and $b(u)Y$, both with zero constant part.  The integrated
Wronskian is therefore zero at leading order.  Its remaining
boundary terms tend to zero after heat integration.

For completeness, the zero-index Poisson errors have Gaussian
factor $e^{-cY^2/u}$ times polynomials in $u,Y,Y^{-1}$,
and the negative factor has an additional $e^{-Au}$ with
$A>0$.  Splitting that exponential gives joint integrability
and decay in $Y$ by the arithmetic--geometric mean inequality.
Terms with a nonzero lower-left entry have the stronger
$e^{-cuY^2}$ cusp factor.  These estimates justify the joint
cutoff/heat limits; they also show
$\int_XS(u)\to0$ as $u\to\infty$, since the region
$Y$ comparable to $\sqrt u$ contributes at most a polynomial
in $u$ times $e^{-Au}$.  The interior compact region has the
same negative heat decay.  Consequently
\begin{equation}\label{eq:zero-negative-total-heat}
 \int_X(\varphi_0+\Omega_0(v_1))F_t(v_2)
       =\int_X\int_1^\infty
          \{K_0^{\rm tot}\Omega_t
                      +\Omega_0^{\rm tot}K_t\}\frac{du}{u}.
\end{equation}

The right side includes all norm-zero vectors in the first
variable, including the zero vector.  For $r\ne0$ its
Gram matrices are nonsingular and have empty finite cycles;
their absolutely convergent sum is the standard genus-two
Green sum.  The only singular matrix is
$\operatorname{diag}(0,t)$.  When $t$ is split its complete
heat degree is the rank-one value proved above, with parameters
$s=v_2$, $\nu=v_1$; hence it is precisely
\[
       \langle\widehat A_t(v_2),H\rangle_{\rm GS}
                         -\frac12\log v_1\,\deg C_t(v_2).
\]
When $t$ is nonsplit there is no nonzero isotropic vector
perpendicular to a norm-$t$ vector, and the dependent heat
identity gives the same statement with $C_t=0$.
Thus both the nonsingular and singular terms in
\eqref{eq:zero-negative-total-heat} have the prescribed
arithmetic degree.  Combining this with
\eqref{eq:zero-negative-arithmetic-bulk} proves the corollary.
\end{proof}

The absence of a strictly negative heat factor changes the
answer at $(0,0)$. The next subsection calculates its dilation
defect and proves that a boundary correction is necessary.

\subsection{The double-zero boundary term and a dilation obstruction}
\label{subsec:double-zero-boundary}

For the double-zero pair, the variation under simultaneous scaling
differs from that of zero--negative pairs and shows that the
unmodified standard product identity is false.
Retain
\[
 \widehat Z_0(v)=H+B_0(v)-\mathfrak a_1(\log v),\qquad
 \operatorname{curv}H=\varphi_0,\qquad
 \widehat{\deg}\mathfrak a_2(T)=\tfrac12\int_XT.
\]
The prescribed rank-zero degree is
\[
 R_0(v_1,v_2)=\widehat{\deg}(H^2)
                    -\tfrac12\log(v_1v_2)\deg H.
\]
Let $S_{\rm ns}(v_1,v_2)$ be the sum of the degrees of the
Gaussian-removed nonsingular coefficients with matrices
$\left(\begin{smallmatrix}0&r/2\\r/2&0\end{smallmatrix}\right)$,
$r\ne0$, and set
\[
 D(v_1,v_2)=\delta^{\rm std}_{0,0}(v_1,v_2)
   =\langle\widehat Z_0(v_1),\widehat Z_0(v_2)\rangle_{\rm GS}
                          -R_0(v_1,v_2)-S_{\rm ns}(v_1,v_2).
\]
Thus $D$ is the complete diagonal coefficient defect, including
all its nonsingular matrices. Write
\[
 d_0=\deg C_0(1)=\sum_c^{\rm orb}
                  \frac{\kappa_c}{2\pi\beta_c\sqrt N}>0.
\]
For $Q(k,l)=v_1k^2+v_2l^2$, let
\[
 \mathcal Z_Q(s)=\sum_{(k,l)\in\mathbb Z^2\setminus\{(0,0)\}}
                                Q(k,l)^{-s}\qquad(\Re s>1)
\]
denote the Epstein zeta function and its meromorphic continuation.
Its value at $s=1/2$ below is not the divergent defining series.

\begin{theorem}[A nonzero double-zero dilation defect]
\label{thm:double-zero-dilation-obstruction}
For the specified standard coefficients, and for $u,v_1,v_2>0$,
\begin{equation}\label{eq:double-zero-dilation-anomaly}
 -u\frac d{du}D(uv_1,uv_2)
             =-\frac{d_0}{4}u^{-1/2}\mathcal Z_Q(1/2).
\end{equation}
Equivalently, for every $\lambda>0$,
\begin{equation}\label{eq:double-zero-two-scale-difference}
 D(v_1,v_2)-D(\lambda v_1,\lambda v_2)
       =-\frac{d_0}{2}(1-\lambda^{-1/2})\mathcal Z_Q(1/2).
\end{equation}
In particular $D(1,1)-D(4,4)>0$. Hence the unmodified standard
arithmetic product identity cannot hold at all imaginary parameters.
\end{theorem}
\begin{proof}
\smallskip\noindent\emph{Arithmetic variation.}
Since \(C_0(v)=v^{-1/2}C_0(1)\), differentiation in the real
arithmetic divisor group gives
\begin{equation}\label{eq:double-zero-arithmetic-derivative}
 -v\frac d{dv}\widehat Z_0(v)
     =\widehat K_0(v)
     :=\bigl(C_0(v)/2,\ 1+K_{\rm iso}(v)\bigr).
\end{equation}
The variation includes the derivative of the nonzero cusp divisor,
so it is not a pure metric variation. The \(1\) comes from
\(-\mathfrak a_1(\log v)\), while the class \(H\) is fixed.

By the already established cusp expansions,
\[
 1+K_{\rm iso}(v,z)
       =\frac{Y}{\beta_c\sqrt{Nv}}+\text{decaying error}.
\]
Its regular part is zero, just as for the complete \(\widehat Z_0\).
The integral Tate frames and the zero degree of the cusp normal
line therefore give, by the same reference-Green argument as
in the zero--negative proof,
\[
 \langle\widehat K_0(v_1),\widehat Z_0(v_2)\rangle_{\rm GS}
       =\frac12\int_X(1+K_{\rm iso}(v_1))
                         (\varphi_0+\Omega_0(v_2)).
\]
Here is an explicit reference-divisor justification. Choose
$\widehat c=(c,h_c)$ with $h_c=-\log|q_c|^2$ near $c$,
zero regular part, and the integral cusp normal frame of norm one.
The actual source and target Tate frames, and the disjointness
of distinct closed arithmetic cusps, give
\[
 \langle\widehat Z_0(v_2),\widehat c\rangle_{\rm GS}
           =\tfrac12\int_Xh_c\,\operatorname{curv}\widehat Z_0(v_2).
\]
Decompose $\widehat K_0(v_1)$ into
$\sum_c(g_c(0,v_1)/2)\widehat c$ plus the metric class of
$1+K_{\rm iso}(v_1)-\sum_c(g_c(0,v_1)/2)h_c$.
The latter function has zero regular part and the required
height-coordinate estimates. The metric variation formula now
gives the asserted bulk integral. In particular the finite
restriction to $C_0(v_1)/2$ has been included. The complete
$\widehat Z_0$ metric has no logarithmic-height exponent, so
no additional logarithmic-metric counterterm occurs here.

Put \(P(u)=\langle \widehat Z_0(uv_1),\widehat Z_0(uv_2)\rangle_{\rm GS}\).
Writing \(K_i=1+K_{\rm iso}(uv_i)\) and
\(\Omega_i=\varphi_0+\Omega_0(uv_i)\), bilinearity gives
\[
             -uP'(u)=\frac12\int_X(K_1\Omega_2+\Omega_1K_2).
\]
The differentiation is legitimate on compact \(u\)-intervals:
the standard Poisson cusp estimates are uniform on such
intervals, the differentiated coefficients have zero regular
parts, and the reference-Green formulas give dominated
integrals.  No singular Green function is multiplied by a
cuspidal Dirac current.

\smallskip\noindent\emph{Separation by Gram matrix.}
The simultaneous heat integrand is the complete lattice sum over
all norm-zero tuples, including the zero vector in either slot.
Separate it by Gram matrix into \(J_{\rm ns}(u)\) and
\(J_{\rm rk0}(u)\).
The nonsingular heat formula gives
\[
 -u\frac d{du}S_{\rm ns}(uv_1,uv_2)
                          =\frac12\int_XJ_{\rm ns}(u).
\]
Here all interchanges are absolute.  Two norm-zero vectors with
nonzero pairing cannot both lie in the same cusp isotropic line.
At every cusp at least one lower-left entry is nonzero, giving
the Gaussian cusp bound.  Also the nonzero integral pairing
gives \(v_1R(x_1,z)+v_2R(x_2,z)\ge c|r|\sqrt{v_1v_2}\)
for a fixed positive normalization constant \(c\).
It supplies exponential large-heat decay and convergence of
the sum over the discrete \(r\ne0\).

The prescribed rank-zero term satisfies
\[
 -u\frac d{du}R_0(uv_1,uv_2)=\deg H=\int_X\varphi_0.
\]
Consequently the exact dilation test is
\begin{equation}\label{eq:double-zero-dilation-before-unfolding}
 -u\frac d{du}D(uv_1,uv_2)
        =\frac12\int_X\bigl(J_{\rm rk0}(u)-2\varphi_0\bigr).
\end{equation}
The right side must not simply be set to zero.
The zero tuple contributes \(2\varphi_0\), but Gram rank zero
also includes every nonzero pair on a rational isotropic line,
including pairs with one entry zero.

\smallskip\noindent\emph{Unfolding the rank-zero orbits.}
Two nonzero isotropic vectors which are perpendicular in this
ternary space are proportional.  Hence every rank-zero tuple
other than the zero tuple has a unique rational isotropic line.
At cusp \(c\), its unfolded entries are exactly
\[
 x_1=\begin{pmatrix}0&\beta_c k\\0&0\end{pmatrix},\qquad
 x_2=\begin{pmatrix}0&\beta_c l\\0&0\end{pmatrix},\qquad
                (k,l)\in\mathbb Z^2\setminus\{(0,0)\}.
\]
After integral cusp conjugation the stabilizer acts by upper
triangular integral matrices with diagonal entries \(\pm1\).
Their adjoint action fixes the nilpotent generator.  Thus
different ordered pairs \((k,l)\) are distinct orbits within
this cusp, and \((k,l)\) is not identified with \((-k,-l)\).
There is no extra factor two or a division by two here.
As throughout the paper, calculation on a small principal-level
scheme followed by division by its degree supplies precisely
the orbifold cusp weights.

Put \(Q(k,l)=v_1k^2+v_2l^2\).  The unfolded integral at cusp \(c\)
is
\begin{equation}\label{eq:double-zero-cusp-seed}
 I_c(u)=\kappa_c\int_0^\infty
  \sum_{(k,l)\ne(0,0)}
    \left(\frac{uN\beta_c^2Q(k,l)}{Y^2}-\frac1\pi\right)
      e^{-\pi uN\beta_c^2Q(k,l)/Y^2}\,\frac{dY}{Y^2}.
\end{equation}
The sum is taken before integration.
At \(Y\to0\) it and its absolute lattice sum are Gaussian small.
At \(Y\to\infty\), two-dimensional Poisson summation cancels
the full Gaussian zero mode: the same sum with \((0,0)\)
included decays exponentially, while the omitted zero term
is \(-1/\pi\).  Thus the displayed grouped integrand is
\(1/(\pi Y^2)\) plus an integrable exponentially small error.
The grouped seed is absolutely integrable on the strip.
Its Poincare unfolding is therefore justified by absolute
integration of the grouped seed.  Termwise absolute integration
of the individual \((k,l)\) summands is neither needed nor true.
It follows that
\[
 \int_X(J_{\rm rk0}(u)-2\varphi_0)
                              =\sum_c^{\rm orb} I_c(u).
\]
In particular \(I_c(u)=u^{-1/2}I_c(1)\), by \(Y\mapsto\sqrt u\,Y\).

\smallskip\noindent\emph{Evaluation of the cusp integral.}
Poisson summation in the theta Mellin integral gives the
meromorphic continuation of $\mathcal Z_Q(s)$; it is regular
at $s=1/2$.
The convergent integral \eqref{eq:double-zero-cusp-seed} is
\begin{equation}\label{eq:double-zero-Epstein-value}
              I_c(u)
                =-\frac{\kappa_c}{4\pi\beta_c\sqrt{uN}}\,
                                      \mathcal Z_Q(1/2).
\end{equation}
Indeed, put \(p=\beta_c\sqrt{\pi uN}/Y\) and
\(\Theta(p)=\sum_{k,l}e^{-p^2Q(k,l)}\).
Then the integral is
\[
 \frac{\kappa_c}{\pi\beta_c\sqrt{\pi uN}}
    \int_0^\infty
        \left\{-\frac p2\Theta'(p)-(\Theta(p)-1)\right\}dp.
\]
Since \(\Theta(p)=\pi/(p^2\sqrt{\det Q})+o(1)\) at zero,
integration by parts on \([\epsilon,\infty)\) cancels that
single divergent term and gives
\[
 \int_0^\infty\left\{-\frac p2\Theta'(p)-(\Theta(p)-1)\right\}dp
       =-\frac12\operatorname{FP}\int_0^\infty(\Theta(p)-1)\,dp
       =-\frac{\sqrt\pi}{4}\mathcal Z_Q(1/2).
\]
This proves \eqref{eq:double-zero-Epstein-value} by evaluating
the preceding convergent grouped integral; no arbitrary
regularization of the original coefficient is used.

Combining \eqref{eq:double-zero-dilation-before-unfolding} and
\eqref{eq:double-zero-Epstein-value}, with the definition of $d_0$,
proves \eqref{eq:double-zero-dilation-anomaly}. Integration over
a finite interval in $u$ proves
\eqref{eq:double-zero-two-scale-difference}; no limit of $D$ at
infinity is used.

Finally the sum-of-two-squares identity gives
\[
 \mathcal Z_{k^2+l^2}(s)=4\zeta(s)\beta_{\rm Dir}(s),
 \qquad
 \beta_{\rm Dir}(s)=\sum_{j\ge0}(-1)^j(2j+1)^{-s}.
\]
For $\Re s>1$ this follows from
$r_2(n)=4\sum_{d\mid n}\chi_{-4}(d)$, and continuation gives
it at $s=1/2$. The alternating series for $\beta_{\rm Dir}(s)$
is positive for real $s>0$. Likewise
$\zeta(s)=\eta(s)/(1-2^{1-s})<0$ for $0<s<1$, since the
alternating eta series is positive. Thus
$\mathcal Z_{k^2+l^2}(1/2)<0$, and setting $\lambda=4$
in \eqref{eq:double-zero-two-scale-difference} proves the strict
inequality.
\end{proof}

\begin{corollary}[The necessary scale-dependent boundary term]
\label{cor:double-zero-required-boundary}
The function
\begin{equation}\label{eq:double-zero-scale-invariant-remainder}
 D(v_1,v_2)+\frac{d_0}{2}\mathcal Z_{\operatorname{diag}(v_1,v_2)}(1/2)
\end{equation}
is invariant under simultaneous positive scaling of $v_1,v_2$.
Consequently it is a function of $v_1/v_2$ alone.
\end{corollary}
\begin{proof}
Use \eqref{eq:double-zero-two-scale-difference} and
$\mathcal Z_{\lambda Q}(1/2)=\lambda^{-1/2}\mathcal Z_Q(1/2)$.
\end{proof}

This determines the scale-dependent part of the correction.
If one keeps the other standard coefficients fixed, the rank-zero
degree would have to acquire
$-\tfrac12d_0\mathcal Z_{\operatorname{diag}(v_1,v_2)}(1/2)$,
up to a scale-invariant remainder. The next subsection proves
that this remainder is zero, completing the double-zero
calculation. The automorphic interpretation of the corrected
coefficients is a separate issue.

\subsection{The exact double-zero boundary correction}
\label{subsec:double-zero-absolute}

To evaluate the scale-invariant remainder in
Corollary~\ref{cor:double-zero-required-boundary}, we control the
escaping region where \(Y\) is comparable to \(\sqrt u\).
Convergence of degrees or of Green functions on compact sets
alone would not justify continuity of heights.

Fix \(v_1,v_2>0\).  Write
\[
 B_i(u)=B_0(uv_i),\quad F_i(u)=F_0(uv_i),\quad
 d_i=\deg C_0(v_i),\quad
 Z_i(u)=H+B_i(u)-\mathfrak a_1(\log(uv_i)).
\]
Thus \(\deg B_i(u)=u^{-1/2}d_i\).  All cusp sums retain the
same arithmetic-orbit/orbifold weights as before.

\begin{lemma}[A uniform profile for the isotropic boundary layer]
\label{lem:zero-boundary-profile-large-scale}
At a cusp \(c\), with source coordinate \(z=x+iY\), width
\(\kappa_c>0\), and isotropic step \(\beta_c\), put
\[
 \Phi_{i,c}(r)=\sum_{k\ne0}
        \beta_1\left(\frac{\pi Nv_i\beta_c^2k^2}{r^2}\right),
 \qquad
 a_{i,c}=\frac{2}{\beta_c\sqrt{Nv_i}},\qquad
 g_{i,c}=\frac{\kappa_ca_{i,c}}{4\pi}.
\]
For \(u\ge1\), \(Y\ge Y_0\),
\[
       F_i(u,x+iY)=\Phi_{i,c}(Y/\sqrt u)+E_{i,c}(u,x,Y),
\]
where the error and any fixed number of height-coordinate
derivatives are bounded by a polynomial in \(u,Y,Y^{-1}\)
times \(e^{-cuY^2}\).  On the fixed compact complement of the
cusps, \(F_i(u)\) and its fixed derivatives are
\(O(u^A e^{-cu})\).  Moreover,
\[
 \begin{aligned}
 \Phi_{i,c}(r)&=O(r^{-A}e^{-b/r^2})&& (r\to0),\\
 \Phi_{i,c}(r)&=a_{i,c}r-2\log r+b_{i,c}
                    +O(r^Ae^{-br^2})&& (r\to\infty),
 \end{aligned}
\]
with corresponding bounds after any specified fixed number
of derivatives.  In particular
\[
       \Phi_{i,c}''(r)-2/r^2=O(r^Ae^{-br^2})
                                                  \quad(r\to\infty).
\]
\end{lemma}
\begin{proof}
In the transformed lattice, an upper-triangular norm-zero vector
has zero diagonal entries.  Its remaining entry is
\(\beta_ck\), \(k\ne0\), giving exactly the displayed profile.
Every other norm-zero vector has a nonzero lower-left entry
in a fixed discrete rational set.  Its majorant contains
\(Nc^2Y^2/2\); the usual Gaussian summation, with a fixed
fraction of the exponent retained, proves the error bounds.
The compact bounds follow from uniform positive definiteness
of the majorant and the exclusion of the zero vector.

The expansion of \(\Phi\) is the already proved one-dimensional
theta/Mellin expansion.  For a direct check on its differentiated
remainder,
\[
 \Phi_{i,c}'(r)=\frac2r\sum_{k\ne0}
       e^{-\pi Nv_i\beta_c^2k^2/r^2}.
\]
Poisson summation gives \(a_{i,c}-2/r\) plus a Gaussian
remainder at infinity; differentiating gives the asserted
second-derivative cancellation.  At zero, the original sum
and all its fixed derivatives have the stated Gaussian
bound.  These are bounds in \(r,Y\), not a claim of a smooth
flat extension in the compactifying coordinate.
\end{proof}

\begin{lemma}[Large-scale Hodge and isotropic cross-pairings]
\label{lem:zero-isotropic-pairing-large-scale}
There are real constants \(A_i\) such that
\begin{equation}\label{eq:zero-Hodge-large-scale}
 \langle H,B_i(u)\rangle_{\rm GS}
       =u^{-1/2}\left(\frac{d_i}{2}\log u+A_i\right)
                                      +O(u^A e^{-cu}).
\end{equation}
Also
\begin{equation}\label{eq:zero-complete-isotropic-bulk}
 \langle Z_1(u),B_2(u)\rangle_{\rm GS}
       =\frac12\int_X F_2(u)\operatorname{curv}Z_1(u)
       =O(u^{-1/2}),
\end{equation}
and consequently
\begin{equation}\label{eq:zero-isotropic-cross-large-scale}
                       \langle B_1(u),B_2(u)\rangle_{\rm GS}
                                                   =O(u^{-1/2}).
\end{equation}
\end{lemma}
\begin{proof}
We first keep the actual Hodge finite part.  The norm convention
in K\"uhn's pairing, as recorded in Definition~6.2.3 of \cite{Zhu},
gives
\begin{equation}\label{eq:zero-isotropic-Hodge-finite-part}
 \langle B_0(v),H\rangle_{\rm GS}
     =\frac12\operatorname{FP}\int_XF_0(v)\varphi_0
        +\sum_c^{\rm orb}g_c(0,v)
             \left(\log c_H-1-\frac12\log N-\log\beta_c\right).
\end{equation}
Indeed one uses the cuspidal section for the first factor and
the source/target Tate frames for the second factor \(H\).
The exponent in the norm of \(H\) is \(-1\), so the same
\(-1\) counterterm as in the negative Hodge pairing remains.
The additional \(-2\log Y\) in \(F_0\) is integrable against
\(\varphi_0\) at infinity and adds no cutoff divergence.
Its own norm exponent is that of the first factor and is
not an extra cusp term in this application of the pairing.
The integral Tate frames account for the finite restrictions
exactly as before.  Thus the formula is not being inferred
from degree convergence.

Here the finite part adds \(2\sum_cg_c(0,v)\log H\) to the
height-\(H\) integral.  The profile lemma shows that
\[
 L_{i,c}:=\lim_{R\to\infty}
       \left(\int_0^R\frac{\Phi_{i,c}(r)}{r^2}\,dr
                                               -a_{i,c}\log R\right)
\]
is finite.  On the cusp, \(\varphi_0=-dx\,dY/(2\pi Y^2)\).
Changing \(Y=\sqrt u\,r\) and using the fixed compact
complement therefore gives
\[
 \operatorname{FP}\int_XF_i(u)\varphi_0
       =u^{-1/2}\left(
                d_i\log u-\sum_c^{\rm orb}
                            \frac{\kappa_c}{2\pi}L_{i,c}\right)
                                      +O(u^A e^{-cu}).
\]
The missing interval \(0<r<Y_0/\sqrt u\) is Gaussian small,
as are all non-upper-triangular and compact errors.
Substitution into
\eqref{eq:zero-isotropic-Hodge-finite-part} proves
\eqref{eq:zero-Hodge-large-scale}.

For \eqref{eq:zero-complete-isotropic-bulk}, the complete
coefficient \(Z_1(u)\) has zero logarithmic norm exponent
and regular norm one in its integral cusp frames, by
Lemma~\ref{lem:complete-zero-cusp-regular-part}.  Its
restriction to \(C_0(uv_2)\) has degree zero.  Apply the
pairing with \(B_2\) as the first factor and \(Z_1\) as the
second.  Its cusp exponent and continuous regular-value
terms are both zero; hence the pairing is precisely the
displayed convergent bulk integral.  The logarithmic norm
of \(B_2\) does not create an additional counterterm.
This uses the complete \(Z_1\), not \(H\) or \(B_1\) alone.

The bound for that integral controls the entire escaping
layer explicitly.  On a cusp set \(r=Y/\sqrt u\).  Up to the
Gaussian errors of the profile lemma, one has
\[
 \operatorname{curv}Z_1(u)
    =\frac1{4\pi u}\left(\Phi_{1,c}''(r)-\frac2{r^2}\right)
                                                 dx\wedge dY.
\]
Consequently the cusp contribution to its half integral
against \(F_2(u)\) is
\[
 \frac{\kappa_c}{8\pi\sqrt u}
       \int_0^\infty
          \Phi_{2,c}(r)
            \left(\Phi_{1,c}''(r)-\frac2{r^2}\right)dr
                                      +O(u^A e^{-cu}).
\]
This ordinary integral is absolute.  At zero \(\Phi_2\)
is Gaussian small and absorbs \(r^{-2}\); at infinity
\(\Phi_1''-2/r^2\) is Gaussian small and absorbs the
linear/logarithmic growth of \(\Phi_2\).
The compact contribution and all differentiated error
terms are \(O(u^A e^{-cu})\).  This proves
\eqref{eq:zero-complete-isotropic-bulk}.

Finally the exact class identity for \(Z_1\) gives
\[
 \langle B_1(u),B_2(u)\rangle_{\rm GS}
   =\langle Z_1(u),B_2(u)\rangle_{\rm GS}
      -\langle H,B_2(u)\rangle_{\rm GS}
                 +\frac12\log(uv_1)\,u^{-1/2}d_2.
\]
The \(\tfrac12u^{-1/2}d_2\log u\) terms cancel by
\eqref{eq:zero-Hodge-large-scale}.  All the remaining
terms are \(O(u^{-1/2})\), proving the last assertion.
\end{proof}

\begin{lemma}[The nonsingular off-diagonal sum is exponentially small]
\label{lem:double-zero-nonsingular-large-scale}
For fixed \(v_1,v_2>0\), the standard nonsingular heat sum
with zero diagonal satisfies
\[
                  S_{\rm ns}(uv_1,uv_2)=O(u^Ae^{-cu})
                                                    \qquad(u\to\infty).
\]
\end{lemma}
\begin{proof}
For a nonempty tuple in this sum its pairing is a nonzero
integer \(r\).  The positive majorant satisfies
\[
 v_1R(x_1,z)+v_2R(x_2,z)
                  \ge c_1|r|\sqrt{v_1v_2}\ge c_2>0
\]
uniformly in \(z\).  At any cusp the two vectors cannot both
be in its isotropic line, so at least one lower-left entry
is nonzero in a fixed discrete rational set.  Its majorant
also supplies a fixed multiple of \(Y^2\).
Split the heat Gaussian into fixed positive fractions:
one extracts \(e^{-cu}\), one extracts the cusp factor
\(e^{-cuY^2}\), and one retains a positive Gaussian for
the absolute lattice summation.  The KM polynomial and
fixed derivatives are absorbed by a further fixed fraction
of that Gaussian.  The residual lattice sum has at most
polynomial growth in \(u,Y\).  On the fixed compact
complement its growth is also polynomial.
These bounds are integrable in space and in the defining
heat variable, and prove the stated estimate for the
absolute sum of the nonsingular coefficients.
\end{proof}

\begin{theorem}[Exact double-zero defect]
\label{thm:double-zero-absolute-defect}
For the specified standard model,
\begin{equation}\label{eq:double-zero-absolute-defect}
 \boxed{\displaystyle
 D(v_1,v_2)
       =-\frac{d_0}{2}
                    \mathcal Z_{\operatorname{diag}(v_1,v_2)}(1/2).}
\end{equation}
Thus the scale-invariant remainder in
\eqref{eq:double-zero-scale-invariant-remainder} is zero.
\end{theorem}
\begin{proof}
Put \(P(u)=\langle Z_1(u),Z_2(u)\rangle_{\rm GS}\).
Expanding the classes and using the constant-metric rule
gives the exact identity
\[
 P(u)-R_0(uv_1,uv_2)
   =\langle H,B_1(u)\rangle_{\rm GS}
       +\langle Z_1(u),B_2(u)\rangle_{\rm GS}
                -\frac12\log(uv_2)\,u^{-1/2}d_1.
\]
The first and last logarithmic terms cancel by
\eqref{eq:zero-Hodge-large-scale}; the middle term is
\(O(u^{-1/2})\) by the full boundary-layer calculation.
Therefore \(P(u)-R_0(uv_1,uv_2)=O(u^{-1/2})\).
The nonsingular estimate then proves
\[
                        D(uv_1,uv_2)\longrightarrow0.
\]
The limit thus includes both the moving cusp layer and the
cross-pairing, without interchanging a compact limit with an
unbounded-height integral.

Apply the exact two-scale formula
\eqref{eq:double-zero-two-scale-difference} and let its
finite scaling parameter tend to infinity.  The proved
limit gives precisely
\eqref{eq:double-zero-absolute-defect}.
\end{proof}

\begin{corollary}[A coefficientwise boundary correction]
\label{cor:corrected-standard-support}
For every positive definite real matrix $y$, set
\[
 \mathcal A^{\rm corr}(T,y)=
 \begin{cases}
 \mathcal A^{\rm std}(0_2,y)-\dfrac{d_0}{2}\mathcal Z_y(1/2),&T=0_2,\\[4pt]
 \mathcal A^{\rm std}(T,y),&T\ne0_2.
 \end{cases}
\]
The correction at rank zero is invariant under
$y\mapsto\gamma y\gamma^t$ for $\gamma\in\mathrm{GL}_2(\mathbb Z)$.
If $\delta^{\rm corr}$ is defined by
\eqref{eq:standard-defect-definition} with these degree coefficients,
then
\begin{equation}\label{eq:corrected-standard-support}
 \delta^{\rm corr}_{m,n}(v_1,v_2)=0
               \qquad(m,n\in\mathbb Z,\ v_1,v_2>0).
\end{equation}
\end{corollary}
\begin{proof}
The bijection $\mathbf n\mapsto\gamma^t\mathbf n$ of
$\mathbb Z^2\setminus\{0\}$ proves
$\mathcal Z_{\gamma y\gamma^t}(s)=\mathcal Z_y(s)$ in the
half-plane of absolute convergence, and then at $s=1/2$ by
continuation. The rank-zero degree changes only the double-zero
diagonal coefficient. The change is exactly
$D(v_1,v_2)$ by \eqref{eq:double-zero-absolute-defect}, so the
new double-zero defect is zero. Apply
Proposition~\ref{prop:standard-residual-support} away from the
positive square-product sector, and
Theorem~\ref{thm:standard-positive-square-product} in that sector.
\end{proof}

This is a coefficientwise degree statement. Rank-zero invariance
under $\mathrm{GL}_2(\mathbb Z)$ does not establish Siegel
modularity of the corrected series; comparison with an arithmetic
Siegel--Weil kernel remains separate.

\section{Common CM components and the global product comparison}
\label{sec:CM-products}

This section develops the regularized analytic self-intersection,
integral adjunction and canonical vertical correction. An intrinsic
trace formula combines the different, conductor and contacts on
a reduced finite-flat Cartier sheet. The ordinary charts, a
two-morphism conormal calculation and the cyclic difference-divisor
theorem give the finite comparison at every prime. After accounting
for contacts introduced by integral scaling, the excess normal
line and complete Green-current expansion assemble these local
identities into the global common-component product formula.

\subsection{The local heat finite part at a common CM point}
\label{subsec:local-CM-self-finite-part}

The calculation concerns one simple common CM point, with the
tangent coordinate or tangent metric included in the definition.
It assigns no value to the product of a Green function and its
own point delta current.

\paragraph{Local model and its normalization.}
Let $p$ be a point on a complex curve, let $w$ be a holomorphic
coordinate with $w(p)=0$, and let $\chi$ be the norm square of a
holomorphic section with a simple zero at $p$, times a fixed positive
constant.  Equivalently,
\[
 \chi(w)=|w|^2 h(w),\qquad h>0\text{ smooth},\qquad \kappa=h(0)>0.
\]
Use $dd^c=(i/2\pi)\partial\bar\partial$ and define the smooth form
\[
 \varphi_0=-dd^c\log h,
 \qquad dd^c\log\chi=\delta_p-\varphi_0.
\]
For $c,t>0$, set
\begin{align*}
 f_c(t)&=\beta_1(tc\chi),& k_c(t)&=e^{-tc\chi},\\
 \Omega_c(t)&=dd^cf_c(t)+\delta_p,&
 \beta_1(x)&=\int_1^\infty e^{-xs}\,\frac{ds}{s}.
\end{align*}
When the heat parameter is omitted, it equals $1$.  The form
\[
 J_\chi=\frac{i}{2\pi}\frac{\partial\chi\wedge\bar\partial\chi}{\chi}
\]
extends smoothly across $p$, with
$J_\chi(p)=(\kappa/\pi)\,d\Re w\wedge d\Im w$.  Direct
differentiation away from $p$, followed by the Green equation, gives
\begin{equation}\label{eq:self-heat-curvature}
       \Omega_c(t)=e^{-tc\chi}(tcJ_\chi+\varphi_0).
\end{equation}
Thus $\Omega_c(t)$ is smooth.  Also
\begin{equation}\label{eq:self-regular-value}
 f_c(t,w)=-\log|w|^2+h_c(t)+o(1),\qquad
 h_c(t)=-\gamma-\log(tc\kappa),
\end{equation}
where $\gamma$ is Euler's constant.

The standard orthogonal CM Green function is this model with
$\chi=2\pi R(u,z)$, where $Q(u)=m>0$ and
$R(u,z)=\tfrac12B(u,h_z)^2-2m$.  In this case
$\varphi_0=-\omega_{\rm hyp}/(2\pi)$.

\begin{definition}[A fixed-coordinate regularized local star]
\label{def:self-regularized-star}
For $c,d>0$ define the current
\begin{equation}\label{eq:self-star-defined}
 \mathcal S^w_{c,d}(t)
      =h_c(t)\delta_p+\Omega_c(t)f_d(t).
\end{equation}
The second summand is locally integrable.  The first summand uses the
regular value in \eqref{eq:self-regular-value}, not the undefined
expression $f_c(t)\delta_p$.
\end{definition}

There is a direct coalescence interpretation of this definition.
Let $p_\varepsilon\ne p$ tend to $p$, and let
$f_{d,\varepsilon}$ be a smooth family of Green functions with its
simple logarithmic pole at $p_\varepsilon$ and with
$f_{d,\varepsilon}\to f_d$ locally in $L^1$.  The standard functions
obtained by smoothly moving a positive vector of fixed norm have
these properties.  The proper local star currents satisfy
\begin{equation}\label{eq:self-coalescence-definition}
 \begin{split}
 \lim_{\varepsilon\to0}\bigl\{
    f_c(p_\varepsilon)\delta_{p_\varepsilon}
      +\Omega_c f_{d,\varepsilon}
      +\log|w(p_\varepsilon)|^2\,\delta_p\bigr\}
       =\mathcal S^w_{c,d}(1).
 \end{split}
\end{equation}
Indeed the second summand converges in $L^1$, and in the atomic
summand the difference between evaluation of a smooth test function
at $p_\varepsilon$ and at $p$ is $O(|w(p_\varepsilon)|)$; its product
with $\log|w(p_\varepsilon)|^2$ tends to zero.  This also proves
independence of the direction of approach.

\begin{lemma}[The concentrated heat constant]
\label{lem:self-concentrated-constant}
For $c,d>0$ put $s=c+d$.  As currents on a sufficiently small
neighborhood of $p$,
\begin{align}
 \Omega_c(t)f_d(t)&\longrightarrow
                             \log(s/d)\,\delta_p,
                                     \label{eq:self-heat-mass}\\
 \mathcal S^w_{c,d}(t)+\log t\,\delta_p
       &\longrightarrow C^w_{c,d}\,\delta_p,
                                     \label{eq:self-star-asymptotic}\\
 C^w_{c,d}&=-\gamma-\log\kappa-\log\frac{cd}{c+d}.
                                     \label{eq:self-constant}
\end{align}
In particular the limiting constant is symmetric in $c,d$.
\end{lemma}

\begin{proof}
Outside $p$ all the relevant heat functions decrease exponentially
in $t$. On a sufficiently small fixed disk $U$ about $p$,
$\chi\asymp |w|^2$ and $|d\chi|\leq C|w|$. Write
$J_\chi=j_\chi(w)\,d^2w$. With $\zeta=\sqrt t\,w$,
\begin{align*}
 \chi(\zeta/\sqrt t)
  &=t^{-1}\{\kappa|\zeta|^2+O(t^{-1/2}|\zeta|^3)\},\\
 j_\chi(\zeta/\sqrt t)
  &=\kappa/\pi+O(t^{-1/2}|\zeta|),
 \qquad d^2w=t^{-1}d^2\zeta.
\end{align*}
The disk can be chosen so that the quadratic lower bound is
uniform. The inequalities
\[
 \beta_1(x)\leq C(1+|\log x|)\quad(0<x\leq1),\qquad
 \beta_1(x)\leq e^{-x}/x\quad(x\geq1)
\]
provide an integrable logarithmic Gaussian majorant.
The exponential mean-value inequality absorbs the Taylor errors
into $t^{-1/2}$ times a polynomial Gaussian. For the $\beta_1$
error, $\beta_1'(x)=-e^{-x}/x$ shows that the cubic argument
error divided by its quadratic argument contributes at most
$Ct^{-1/2}|\zeta|e^{-c|\zeta|^2}$. For a test function $\phi$,
the error in replacing $\phi(\zeta/\sqrt t)$ by $\phi(0)$ is
at most $t^{-1/2}|\zeta|\|\phi\|_{C^1}$. Each resulting
majorant is integrable, including the factor
$1+|\log|\zeta||$ near $0$. The $\varphi_0$ summand in
\eqref{eq:self-heat-curvature} has an additional factor $t^{-1}$.
The contribution outside a fixed smaller disk is exponentially small.

The two limiting masses are
\[
 s\int_0^\infty e^{-sr}\,dr=1,\qquad
 c\int_0^\infty e^{-cr}\beta_1(dr)\,dr
    =\int_1^\infty\frac{c}{u(c+du)}\,du
    =\log\frac{c+d}{d},
\]
where Tonelli's theorem justifies the second evaluation.
Thus, for every smooth test function supported in $U$,
\begin{align*}
 |\langle\Omega_s(t)-\delta_p,\phi\rangle|
       &\leq Ct^{-1/2}\|\phi\|_{C^1},\\
 |\langle\Omega_c(t)f_d(t)-\log(s/d)\delta_p,\phi\rangle|
       &\leq Ct^{-1/2}\|\phi\|_{C^1}.
\end{align*}
These are weak current estimates; they do not assert convergence
in total variation to an atom. Inserting
$h_c(t)=-\log t-\gamma-\log(c\kappa)$ gives
\eqref{eq:self-star-asymptotic} with the same $O(t^{-1/2})$
bound in the dual $C^1$ norm, and yields
\eqref{eq:self-constant}.
\end{proof}

As a direct check on the coalescence scale, in the flat model
$\chi_1=\kappa|w|^2$, $\chi_2=\kappa|w-\varepsilon|^2$,
$\varphi_0=0$, the two-point heat form integrates over the plane to
\[
 \int_{\mathbb C}
       \{k_c(t)\Omega_{d,\varepsilon}(t)
                  +\Omega_c(t)k_{d,\varepsilon}(t)\}
       =\exp\left(-t\kappa\frac{cd}{c+d}|\varepsilon|^2\right).
\]
This follows by completing the square
$c|w|^2+d|w-\varepsilon|^2
 =(c+d)|w-d\varepsilon/(c+d)|^2+cd|\varepsilon|^2/(c+d)$.
Consequently the integrated proper genus-two Green form is
$\beta_1(\kappa cd|\varepsilon|^2/(c+d))$, and its finite part
after adding $\log|\varepsilon|^2$ is exactly
$-\gamma-\log(\kappa cd/(c+d))$.  The curved calculation in
Lemma~\ref{lem:self-concentrated-constant} establishes the same
constant without replacing the metric by a flat one.

\paragraph{The rank-one heat kernel.}
By \cite[Proposition~2.4.5(a)]{GarciaSankaran}, the normalized
form $\nu^o$ gives the local two-vector heat form
\[
 \mathcal J_{c,d}(t)
       =k_c(t)\Omega_d(t)+\Omega_c(t)k_d(t).
\]
Since $k_c k_d=k_{c+d}$, equation
\eqref{eq:self-heat-curvature} gives the exact identity
\begin{equation}\label{eq:self-rankone-heat-decomposition}
 \mathcal J_{c,d}(t)
        =\Omega_s(t)+\varphi_0 k_s(t),\qquad s=c+d.
\end{equation}
The two-vector kernel therefore depends only on $s$; the
complementary parameter will occur in its comparison with the
regularized star, not in this kernel itself.

The following distributional finite part is unambiguous:
\begin{equation}\label{eq:self-rankone-finite-part}
 \mathcal G_s^{\rm fp}
   =\lim_{T\to\infty}
       \left\{\int_1^T\mathcal J_{c,d}(t)\,\frac{dt}{t}
                                      -\log T\,\delta_p\right\}.
\end{equation}
In fact, put
\[
 H_s=\int_1^\infty f_s(t)\,\frac{dt}{t}.
\]
This function is locally integrable: $f_s(t)$ has local $L^1$
norm $O(t^{-1})$ near $p$, and exponential decay on compact
subsets excluding $p$.  Integrating
$\Omega_s(t)=dd^cf_s(t)+\delta_p$ gives directly
\begin{equation}\label{eq:self-rankone-G-explicit}
          \mathcal G_s^{\rm fp}=dd^cH_s+\varphi_0 f_s.
\end{equation}
Here $dd^cH_s$ is the distributional derivative of an $L^1$
function.  It must not be replaced by an unregularized ordinary
Laplacian on the punctured disk.  Near $p$, $H_s=O(1+\log^2|w|)$.

This is also the Mellin constant term
\[
 \operatorname{CT}_{\rho=0}
     \int_1^\infty\mathcal J_{c,d}(t)t^{-\rho}\,\frac{dt}{t}.
\]
Indeed the quantitative estimate in the proof of
Lemma~\ref{lem:self-concentrated-constant} gives
$\Omega_s(t)-\delta_p=O(t^{-1/2})$ in the dual $C^1$ test norm,
locally in this coordinate.  Also $k_s(t)$ has local
$L^1$ norm $O(t^{-1})$.  Thus the Mellin integral has the form
$\delta_p/\rho$ plus a distribution holomorphic for
$\Re\rho>-1/2$, and its constant term is
\eqref{eq:self-rankone-G-explicit}.  This is precisely the
Gaussian-removed rank-one regularization of
\cite[Definition~2.6.7]{GarciaSankaran}, with subtraction
$\delta_p/\rho$ and no extra finite counterterm.

\begin{theorem}[A regularized local self-star identity]
\label{thm:self-local-identity}
Fix the coordinate $w$ and the finite part
\eqref{eq:self-rankone-finite-part}.  The locally $L^1$ one-forms
\begin{align}
 P_{c,d}
   &=\frac{i}{2\pi}\int_1^\infty
                     f_d(t)\bar\partial k_c(t)\,\frac{dt}{t},
                                      \label{eq:self-P}\\
 Q_{c,d}
   &=\frac{i}{2\pi}\int_1^\infty
                     k_c(t)\partial f_d(t)\,\frac{dt}{t}
                                      \label{eq:self-Q}
\end{align}
satisfy the exact current identity
\begin{equation}\label{eq:self-current-identity}
 \mathcal S^w_{c,d}(1)-\mathcal G_s^{\rm fp}
      =C^w_{c,d}\delta_p+\partial P_{c,d}+\bar\partial Q_{c,d}.
\end{equation}
Off $p$ the primitives have the closed formulas
\begin{align}
 P_{c,d}
   &=-\frac{i}{2\pi}\bar\partial\log\chi\,
                                  (k_c f_d-f_s),
                                      \label{eq:self-P-closed}\\
 Q_{c,d}
   &=-\frac{i}{2\pi} f_s\,\partial\log\chi.
                                      \label{eq:self-Q-closed}
\end{align}
Near $p$ one has $P_{c,d}=O(|w|^{-1})$ and
$Q_{c,d}=O((1+|\log|w||)/|w|)$; these bounds are locally integrable.
Together with \eqref{eq:self-rankone-G-explicit}, the identity
can equivalently be written
\begin{equation}\label{eq:self-local-adjunction-form}
 \mathcal S^w_{c,d}(1)
      =\varphi_0 f_s+C^w_{c,d}\delta_p
         +dd^cH_s+\partial P_{c,d}+\bar\partial Q_{c,d}.
\end{equation}
\end{theorem}

\begin{proof}
The primitives converge strongly in local $L^1$. On a fixed small
disk $U$, the change of variables $\zeta=\sqrt t\,w$ and the
logarithmic Gaussian bounds in the preceding proof give
\[
 \|f_s(t)\|_{L^1(U)}\leq C/t,\qquad
 \|k_c(t)\partial f_d(t)\|_{L^1(U)}
   +\|f_d(t)\bar\partial k_c(t)\|_{L^1(U)}\leq C/\sqrt t.
\]
Here $\partial f_d(t)=-e^{-td\chi}\partial\chi/\chi$,
and $|d\chi|/\chi=O(|w|^{-1})$; the other integrand is
controlled by the same rescaling and the local logarithmic bound
for $\beta_1$. If a superscript $[T]$ denotes truncation of a
defining heat integral at $T$, integration with $dt/t$ yields
\[
 \|H_s-H_s^{[T]}\|_{L^1(U)}\leq C/T,\qquad
 \|P_{c,d}-P_{c,d}^{[T]}\|_{L^1(U)}
   +\|Q_{c,d}-Q_{c,d}^{[T]}\|_{L^1(U)}\leq C/\sqrt T.
\]
Thus distributional differentiation commutes with all three limits.
In particular, for finite $T$, the defining curvature identity gives
\[
 \int_1^T\mathcal J_{c,d}(t)\,\frac{dt}{t}-\log T\,\delta_p
   =dd^cH_s^{[T]}+\varphi_0\{f_s(1)-f_s(T)\}.
\]
The same strong bounds prove its limit
\eqref{eq:self-rankone-G-explicit} directly.

The formula for $Q$ follows immediately by integrating
$-e^{-t(c+d)\chi}\partial\log\chi$.  For $P$, Tonelli's
theorem for its scalar integral and the identity
\[
 \frac{c}{u(c+du)}=\frac1u-\frac{d}{c+du}
\]
give
\[
 c\chi\int_1^\infty e^{-tc\chi}\beta_1(td\chi)\,dt
       =e^{-c\chi}\beta_1(d\chi)-\beta_1((c+d)\chi).
\]
This proves the other closed formula.  Since
$k_cf_d-f_s\to\log(s/d)$ at $p$, it also proves the displayed
singularity bounds.

Heat differentiation now involves only well-defined currents.
Since $-t\partial_t h_c(t)=1$,
\[
 -t\partial_t\mathcal S^w_{c,d}(t)
      =\delta_p+(dd^ck_c(t))f_d(t)+\Omega_c(t)k_d(t).
\]
Subtract \eqref{eq:self-rankone-heat-decomposition}, using
$\Omega_d=dd^cf_d+\delta_p$ and $k_c(p)=1$.  The difference is
\[
 f_d(t)dd^ck_c(t)-k_c(t)dd^cf_d(t)
   =\frac{i}{2\pi}
       \{\partial(f_d(t)\bar\partial k_c(t))
                   +\bar\partial(k_c(t)\partial f_d(t))\}.
\]
This is the distributional Leibniz rule: the mixed first derivatives
cancel.  Its terms are defined because $k_c$ is smooth and $f_d$
and its first weak derivatives are locally integrable.

Integrate from $1$ to $T$.  By
Lemma~\ref{lem:self-concentrated-constant},
$\mathcal S^w_{c,d}(T)=-\log T\,\delta_p+
C^w_{c,d}\delta_p+o(1)$, whereas the integral of the two-vector
heat form is $\log T\,\delta_p+\mathcal G_s^{\rm fp}+o(1)$.
The two logarithmic divergences cancel.  The $L^1$ convergence of
the primitives permits distributional differentiation after
passing to the limit, which proves
\eqref{eq:self-current-identity}.  The final form follows from
\eqref{eq:self-rankone-G-explicit}.
\end{proof}

\paragraph{Coordinate dependence and a tangent metric.}
If $w'=\eta w+O(w^2)$ with $\eta\ne0$, then
\[
 \kappa'=\kappa/|\eta|^2,\qquad
 C^{w'}_{c,d}=C^w_{c,d}+\log|\eta|^2,
 \qquad
 \mathcal S^{w'}_{c,d}=\mathcal S^w_{c,d}
                                  +\log|\eta|^2\delta_p.
\]
Hence the unqualified fixed-coordinate finite part is not an
intrinsic self-intersection.  If the tangent line has a Hermitian
norm with $\|\partial/\partial w\|^2=\ell_p$, replace
$\log|w(p_\varepsilon)|^2$ in
\eqref{eq:self-coalescence-definition} by
$\log(\ell_p|w(p_\varepsilon)|^2)$.  The resulting intrinsic
finite part has
\begin{equation}\label{eq:self-metric-constant}
        C^{\rm tan}_{c,d}
             =-\gamma-\log(\kappa/\ell_p)-\log\frac{cd}{c+d}.
\end{equation}
The combination is invariant because
$\ell'_p=\ell_p/|\eta|^2$.  If
$\chi=2\pi\|s_u\|^2$ in a Hermitian line bundle, then
$\kappa/\ell_p=2\pi\|ds_u:T_pX\to\mathcal E_p\|^2$.
This explicitly identifies the missing normal-metric input in an
arithmetic adjunction comparison.

For the standard CM vector $u$ with $Q(u)=m$ and the hyperbolic
tangent norm, defined explicitly by
$\|\partial/\partial z\|^2=(\Im z)^{-2}$, one has
\[
 R(u,z)=2m\sinh^2d_{\rm hyp}(p,z),\qquad
 \kappa/\ell_p=4\pi m,
\]
and therefore
\begin{equation}\label{eq:self-hyperbolic-constant}
          C^{\rm hyp}_{c,d}
             =-\gamma-\log(4\pi m)-\log\frac{cd}{c+d}.
\end{equation}

\paragraph{Comparison with the rank-one parameters.}
Suppose the two positive vectors are $a u$ and $b u$, where
$a,b\in\mathbb Q^\times$, with imaginary parameters $v_1,v_2>0$.
Then
\[
 c=a^2v_1,\qquad d=b^2v_2,\qquad
 s=v^*=a^2v_1+b^2v_2,\qquad
 v_\perp=\frac{v_1v_2}{v^*}.
\]
The local normalized rank-one two-vector kernel depends exactly on
$v^*$.  Its regularized star comparison contains
\begin{equation}\label{eq:self-rankone-parameter-comparison}
 C^{\rm tan}_{a^2v_1,b^2v_2}
    =-\log v_\perp-\gamma-\log(\kappa/\ell_p)
                                      -\log(a^2b^2).
\end{equation}
In particular when $a=b=1$ the parameter correction is precisely
$-\log v_\perp$, with remaining universal/metric constant
$-\gamma-\log(\kappa/\ell_p)$.  In the hyperbolic tangent metric
the latter is $-\gamma-\log(4\pi m)$.

The rational factors in the general formula are required in a
fixed-coordinate calculation: they result from replacing the two
local sections by $a s_u$ and $b s_u$. Adding the constants
$\log a^2$ and $\log b^2$ to the respective Green functions
increases the concentrated comparison constant by
$\log(a^2b^2)$, giving
$-\log v_\perp-\gamma-\log(\kappa/\ell_p)$. This changes the
arithmetic divisor data and therefore requires the corresponding
rational-section/finite-place comparison.
For the primitive integral parameters of a rank-one matrix,
Corollary~\ref{cor:proportionality-scaling-contact} below proves
the matching finite contact, with the original Green functions
unchanged. Rescaling the base vector $u$ and the coefficients
$a,b$ consistently leaves
the complete constant in
\eqref{eq:self-rankone-parameter-comparison} invariant.

\paragraph{What this proves and what it does not.}
The calculation proves the fixed-coordinate coalescence limit,
the exact local current identity, the rank-one Mellin subtraction
and all logarithmic constants above. It determines the analytic
$-\log v_\perp$ dependence and the additional tangent/section
normalization. An arithmetic self-intersection also requires
applying the integral adjunction and metric formula of
Lemma~\ref{lem:integral-CM-adjunction} to the actual branches and
including the other finite-place intersections with the same
rational-section normalization. Local current exactness alone
does not justify discarding primitives on the boundary of a
coordinate disk. Throughout, the heat kernel is the normalized
$\nu^o$; integrating the ordinary Gaussian $\nu$ instead gives
a different finite part.

\subsection{Integral adjunction for a horizontal Cartier branch}
\label{subsec:integral-branch-adjunction}

\begin{lemma}[Normal bundle, conductor and metric correction]
\label{lem:integral-CM-adjunction}
Let $X$ be a regular flat arithmetic surface with relative dualizing
line bundle $\omega_X=\omega_{X/\mathbb Z}$.  Let
$i:D=\operatorname{Spec}A\hookrightarrow X$ be an integral effective
Cartier divisor finite and flat over $\mathbb Z$, and let
$\nu:\operatorname{Spec}\mathcal O_F\to D$ be its normalization.
Write $j=i\nu$, let $\mathfrak d_F$ be the different of $F/\mathbb Q$,
and put
\[
 \mathfrak f=\{x\in\mathcal O_F:x\mathcal O_F\subset A\}.
\]
There are canonical integral isomorphisms
\begin{align}
 \mathcal O_X(D)|_D
    &\simeq\omega_{D/\mathbb Z}\otimes(i^*\omega_X)^{-1},
                         \label{eq:integral-Cartier-adjunction}\\
 j^*\mathcal O_X(D)
    &\simeq\mathfrak f^{-1}\mathfrak d_F^{-1}
                                   \otimes(j^*\omega_X)^{-1}.
                         \label{eq:integral-normalized-adjunction}
\end{align}

Choose a Hermitian metric on $\omega_X$ and a logarithmic Green
function $g_D$ determining the metric
$\|s_D\|^2=e^{-g_D}$ on $\mathcal O_X(D)$.
For every embedding $\sigma:F\hookrightarrow\mathbb C$, choose
a local coordinate $w_\sigma$ at $p_\sigma=j(\sigma)$ and write
\[
 g_D=-\log|w_\sigma|^2+h_\sigma+o(1),\qquad
 \ell_\sigma=\|\partial/\partial w_\sigma\|^2,
\]
where the tangent norm is dual to the chosen metric on $\omega_X$.
Give the rational trace element of
$\mathfrak f^{-1}\mathfrak d_F^{-1}$ squared norm $1$ at every
embedding.  With the arithmetic-degree convention
\[
 \widehat{\deg}\,\overline I
    =-\log\operatorname{Norm}(I)
                -\frac12\sum_{\sigma:F\hookrightarrow\mathbb C}
                                      \log\|1\|_\sigma^2
\]
for a fractional ideal, the normal-bundle height is
\begin{equation}\label{eq:integral-adjunction-metric-degree}
 \begin{split}
 \widehat{\deg}\,j^*\overline{\mathcal O_X(D)}
   ={}&\log|\operatorname{Disc}(F)|
          +\log\operatorname{Norm}(\mathfrak f)
          -\widehat{\deg}\,j^*\overline\omega_X\\
      &+\frac12\sum_\sigma(h_\sigma+\log\ell_\sigma).
 \end{split}
\end{equation}
All embeddings are counted, including both members of a complex
conjugate pair.  For normal $D$ the conductor is $\mathcal O_F$;
for a quadratic order of conductor $f$ its term is $2\log f$.
\end{lemma}

\begin{proof}
A Cartier divisor in a regular scheme is Gorenstein.  Applying
$\mathcal Hom(-,\omega_X)$ to
$0\to\mathcal O_X(-D)\to\mathcal O_X\to\mathcal O_D\to0$
gives
$\mathcal Ext^1_X(\mathcal O_D,\omega_X)=\omega_X(D)|_D$.
Relative duality identifies this with $\omega_{D/\mathbb Z}$,
proving \eqref{eq:integral-Cartier-adjunction}.

Finite duality identifies
$\omega_A=\operatorname{Hom}_{\mathbb Z}(A,\mathbb Z)$ with
the trace codifferent in $F$ and gives
\[
 \mathfrak d_F^{-1}
   =\operatorname{Hom}_A(\mathcal O_F,\omega_A)
   =\operatorname{Hom}_A(\mathcal O_F,A)\otimes_A\omega_A
   =\mathfrak f\,\nu^*\omega_A.
\]
Here $\omega_A$ is invertible because $A$ is Gorenstein.
An $A$-linear map $\mathcal O_F\to A$ is generic multiplication
by $x\in F$, and its image lies in $A$ precisely when
$x\mathcal O_F\subset A$; hence its Hom module is $\mathfrak f$.
The conductor is an invertible $\mathcal O_F$-ideal, so the last
identity can be solved for $\nu^*\omega_A$.  This proves
\eqref{eq:integral-normalized-adjunction}.

At $p_\sigma$, the residue form of adjunction sends
$(1/w_\sigma)\otimes dw_\sigma$ to the trace element $1$.
The Green metric gives
$\|1/w_\sigma\|^2=e^{-h_\sigma}$, whereas the metric on the
right of \eqref{eq:integral-normalized-adjunction} gives this
frame squared norm $\ell_\sigma$.  Their ratio is therefore
$e^{-h_\sigma}/\ell_\sigma$.  Changing a line metric by this
ratio changes its arithmetic degree by
$\tfrac12\sum_\sigma(h_\sigma+\log\ell_\sigma)$.
Finally
$\operatorname{Norm}(\mathfrak d_F)=|\operatorname{Disc}(F)|$
and the displayed fractional-ideal convention give the two
positive finite terms in
\eqref{eq:integral-adjunction-metric-degree}.
\end{proof}

The correction $h_\sigma+\log\ell_\sigma$ is coordinate
independent: replacing $w$ by $\alpha w+O(w^2)$ adds
$\log|\alpha|^2$ to $h$ and subtracts it from $\log\ell$.
For the individual local Green function
$f_c=\beta_1(c\chi)$ of
Section~\ref{subsec:local-CM-self-finite-part},
$\chi=\kappa|w|^2+O(|w|^3)$ gives
\[
 h_c+\log\ell=-\gamma-\log(c\kappa/\ell),\qquad
 C^{\rm tan}_{c,d}
      =(h_c+\log\ell)+\log\frac{c+d}{d}.
\]
Thus the normal-bundle metric term and the two-vector heat finite
part are related but are not the same constant.  In the stated
hyperbolic tangent norm the first is
$-\gamma-\log(4\pi Q(u)c)$.  If the complete Green function has
an additional smooth term at the CM point, its value must also
be included in $h_\sigma$.

\paragraph{Application to special cycles.}
The next subsection supplies actual Cartier branches on auxiliary
level covers and compares the relative dualizing metric with the
modified Hodge class at every level.  Applying the normal-bundle
formula then requires the finite branch orders and their conductor
ideals.  The arithmetic self-product also retains the
Green-curvature integral, intersections between distinct branches,
and the explicit vertical contact terms.  Rational factors in the
special vectors must be tracked at finite as well as infinite places.

\subsection{Actual Cartier branches and the all-level canonical correction}
\label{subsec:all-level-canonical-correction}

Write $\mathcal X=\mathcal X_0(N)$ for the Deligne--Rapoport
modular stack with the refined compactification of
\cite[Theorem~5.13(a)]{Cesnavicius2017}, and let $C$ be its
reduced cuspidal divisor. Regularity here refers to the stack;
the coarse moduli scheme need not be regular at every elliptic
point in bad characteristic.

\begin{lemma}[Cartier branches of a positive special divisor]
\label{lem:actual-CM-Cartier-branches}
Fix a prime $p$, put $R=\mathbb Z_{(p)}$, and let $m>0$.
For the positive special-moduli divisor in
\cite[Proposition~4.3.7]{Zhu23}, its map to $\mathcal X_R$ is
finite unramified and generalized Cartier.  After an auxiliary
prime-to-$p$ full-level cover of the elliptic-curve locus, every
reduced irreducible horizontal branch of every Cartier sheet is
an effective Cartier divisor on a regular scheme.  It is finite
and flat over $R$, misses the cusps, and its finite dualizing
module is invertible.  Consequently the integral adjunction
isomorphisms of Lemma~\ref{lem:integral-CM-adjunction} apply to
these actual branches, locally over $R$.
\end{lemma}

\begin{proof}
The stated proposition supplies both generalized Cartierness and
finiteness, and its proof supplies flatness over $R$ and avoidance
of the cusps.  Choose auxiliary full level $M\geq3$ prime to $p$
on both the source and target elliptic curves, with $M^2>4m$,
and extend this finite etale cover of the elliptic-curve locus
by normalization over the proper stack $\mathcal X_R$.
This extended cover is finite over $\mathcal X_R$ and proper
over $R$.  Its elliptic-curve locus is a regular scheme $U$,
since the cover there is etale and $\mathcal X$ is regular.
The full levels make $j|_{E[M]}$ a matrix over $\mathbb Z/M$,
so its finitely many matrix values decompose the special-moduli
scheme into open-and-closed subschemes.  Each such piece has at
most one geometric point over each point of $U$.  Indeed, if
$j_1,j_2$ of degree $m$ have the same matrix, then
$j_1-j_2=Mh$ for a homomorphism $h$.  The degree parallelogram
identity gives
\[
 M^2\deg h=\deg(j_1-j_2)
       \leq2\deg j_1+2\deg j_2=4m,
\]
so $h=0$.  A finite unramified map with these geometric fibres
is a closed immersion: its fibre algebra is either zero or
the residue field, and Nakayama's lemma makes the finite
algebra a quotient of the base algebra.  Generalized
Cartierness now says that the ideal of each matrix piece is
invertible.  Thus these are actual globally defined Cartier
sheets on $U$, not a presumed gluing of local sheets with
monodromy.  If special morphisms have a fixed denominator,
first clear it, replace $m$ by the resulting integral degree,
and choose $M$ prime to that denominator as well.

In a regular local ring each reduced height-one prime is
principal.  Therefore each irreducible horizontal branch of
one of these Cartier sheets is again Cartier on $U$.

For completeness, flatness can also be checked without assuming
that the special fibre is reduced.  A local equation $f$ of a
horizontal sheet has no irreducible factor in common with $p$:
otherwise the sheet would contain a vertical component.  This
is excluded in the proof of Zhu's proposition by the ordinary
generic point of each special-fibre component.  In the regular
local unique-factorization domain, $f\mid pg$ implies $f\mid g$
whenever $f$ and $p$ have no common prime factor.  Thus $p$ is
a nonzerodivisor modulo $f$, proving flatness over $R$.  The same
argument applies to each reduced horizontal factor.  A branch
closure on the extended auxiliary-level cover is proper over
$R$.  It misses the cusps and has no vertical component, so
its fibres over $R$ are zero-dimensional.  It is consequently
proper and quasi-finite, hence finite over $R$.  Such a Cartier
quotient of a regular surface is Gorenstein, so its finite
dualizing module is invertible, as required for adjunction.
\end{proof}

The argument establishes Cartierness but does not identify a
branch order or its conductor. Distinct moduli sheets must retain
their multiplicities before they are combined. At each $p$, choose
the auxiliary level prime to that $p$ and descend degrees by the
degree of the finite etale cover; equivalently, retain the usual
stack stabilizer weights.

Put $\lambda=\omega_N$ for the source Hodge line and
$K=\omega_{\mathcal X/\mathbb Z}$ for the relative dualizing
line.  If $p^e\Vert N$, index the reduced irreducible components
of the $p$-fibre by $E_{p;a,b}$, where $a+b=e$ and the generic
cyclic $p^e$-subgroup has connected part of order $p^a$ and
etale part of order $p^b$.  Set
\[
 r_{a,b}=\varphi(p^{\min(a,b)}),\qquad
 \operatorname{div}(p)=\sum_{a+b=e}r_{a,b}E_{p;a,b},
 \qquad \varphi(1)=1.
\]
The component corresponding to Zhu's index $j$ has $j=a-b$.
Define
\begin{align}
 D_N&=\sum_{p\mid N}\sum_{a+b=e}d_p(a,b)E_{p;a,b},
                       \label{eq:explicit-canonical-different}\\
 d_p(a,b)&=
 \begin{cases}
 b,&a=0,\\
 p^{\min(a,b)-1}\bigl((p-1)b-1\bigr),&a,b\geq1,\\
 0,&b=0,
 \end{cases} \notag\\
 U_N&=\sum_{p^e\Vert N}\ \sum_{\substack{a+b=e\\a,b\geq1}}
                  p^{\min(a,b)-1}E_{p;a,b}.
                       \label{eq:internal-canonical-correction}
\end{align}
Empty sums in $U_N$ are zero.  All formulas include $p=2$.

\begin{proposition}[All-level dualizing--Hodge comparison]
\label{prop:all-level-dualizing-Hodge}
In $\operatorname{Pic}(\mathcal X)\otimes\mathbb Q$ one has
\begin{equation}\label{eq:all-level-canonical-algebraic}
 K\simeq\lambda^{\otimes2}(-C+D_N),\qquad
 -[K]=H+[C]+[U_N],
\end{equation}
where $H=-[\lambda]-[W_N^*\lambda]$ is the algebraic modified
Hodge class used in the rank-one coefficient.

More precisely, put $Y=\mathcal X\setminus C$, equip $K|_Y$
with the hyperbolic cotangent metric
$\|d\tau\|^2=(\operatorname{Im}\tau)^2$, and equip $\lambda$
with Zhu's metric.  Use the Tate-normalized Kodaira--Spencer map
$\mathrm{KS}((ds/s)^{\otimes2})=dq/q$ and the standard
$q$-expansion identification of modular forms with Hodge
sections.  There is an isometry of metrized rational line
bundles on $Y$,
\begin{equation}\label{eq:all-level-canonical-arithmetic}
 \overline K_{\rm hyp}^{-1}|_Y
  \simeq
  \left(\overline{\mathcal H}\otimes
        \overline{\mathcal O}(U_N)\otimes
          \mathbf1(\gamma+\log\pi-\log N)\right)|_Y.
\end{equation}
Here $\overline{\mathcal H}$ represents
$\widehat H=-\widehat\omega_N-W_N^*\widehat\omega_N$,
$\overline{\mathcal O}(U_N)$ has canonical section of norm one,
and $\mathbf1(c)$ is the trivial line with generator of squared
norm $e^{-c}$.  Its constant-metric class is
$\mathfrak a_1(c)=(0,c)$.  We do not regard $(C,0)$ as a global
arithmetic divisor in the fixed extended Chow theory.
\end{proposition}

\begin{proof}
The first isomorphism is the integral Kodaira--Spencer formula
of \cite[Proposition~5.6]{CNS2022}; its different multiplicities
are exactly \cite[Proposition~5.12]{CNS2022}.  Its regular-locus
qualification causes no loss here, because the modular stack
$\mathcal X$ is regular.  The different is computed in strict
henselian discrete valuation rings at generic fibre components;
it is not the different of a number-field CM branch.

Write
\[
 B_N=\sum_{p^e\Vert N}\sum_{a+b=e}
             \frac{a-b}{2}r_{a,b}E_{p;a,b}.
\]
Zhu's source--target comparison gives
$\widehat\omega_N-W_N^*\widehat\omega_N=(B_N,0)$
\cite[Corollary~6.4.4]{Zhu}.  In particular,
$-\widehat H=2\widehat\omega_N-(B_N,0)$.
An elementary calculation of each component gives
\begin{equation}\label{eq:different-plus-B}
 D_N+B_N=\frac12\operatorname{div}(N)-U_N.
\end{equation}
Indeed, for $\min(a,b)=0$ the coefficient on the left is $e/2$.
For $k=\min(a,b)\geq1$ it is
\[
 p^{k-1}\left((p-1)b-1+\frac{(a-b)(p-1)}2\right)
   =\frac e2\varphi(p^k)-p^{k-1}.
\]
This proves \eqref{eq:different-plus-B}, including its dyadic
case, and proves the second formula in
\eqref{eq:all-level-canonical-algebraic} after discarding the
principal divisor $\operatorname{div}(N)/2$.

We spell out the metric correction.  Let
$c_H=2\sqrt\pi e^{-\gamma/2}$, so that a Hodge section with
scalar modular function $f$ has squared norm
$c_H\operatorname{Im}(\tau)|f(\tau)|^2$
\cite[Example~6.2.2]{Zhu}.  The $q$-expansion frame is
$e_T=ds/s$; analytically $s=e^{2\pi iz}$ on the elliptic curve,
and $q=e^{2\pi i\tau}$ on its parameter space.  Thus
\[
 \mathrm{KS}(e_T^{\otimes2})=2\pi i\,d\tau,
 \qquad
 \frac{\|\mathrm{KS}(e_T^{\otimes2})\|_{\rm hyp}^2}
      {\|e_T^{\otimes2}\|_{\lambda^{\otimes2}}^2}
       =\frac{4\pi^2}{c_H^2}.
\]
The Tate normalization of KS may for example be checked in
\cite[Section~1.1.3--1.1.4]{NekovarPlater2000}.  With
$\mathfrak a_1(c)=(0,c)$ a squared-metric ratio $r$ contributes
$\mathfrak a_1(-\log r)$.  Therefore, on $Y$, where the
canonical cusp section trivializes $\mathcal O(C)$,
\[
 \widehat K_{\rm hyp}|_Y
   =\left(2\widehat\omega_N+(D_N,0)
                    +\mathfrak a_1(-\gamma-\log\pi)\right)|_Y.
\]
Finally
$(\operatorname{div}(N),-2\log N)=0$ implies
$(\operatorname{div}(N)/2,0)=\mathfrak a_1(\log N)$.
Substitution of \eqref{eq:different-plus-B} gives
\eqref{eq:all-level-canonical-arithmetic}.
\end{proof}

For comparison with the local heat calculation, use the hyperbolic
tangent metric and write the two special vectors as $au,bu$, with
$Q(u)>0$. Combining \eqref{eq:self-rankone-parameter-comparison}
with the metric constant in
\eqref{eq:all-level-canonical-arithmetic} gives
\begin{equation}\label{eq:canonical-heat-constant-cancellation}
 C^{\rm hyp}_{a^2v_1,b^2v_2}+\gamma+\log\pi-\log N
   =-\log v_\perp-\log\bigl(4NQ(u)a^2b^2\bigr).
\end{equation}
Here $C^{\rm hyp}$ denotes the previously defined tangent-metric
constant with the hyperbolic norm. Thus the Euler constant and
$\pi$ cancel in the combined local comparison. This numerical
identity leaves the field different, conductor and finite
intersection contributions in the following formula to be evaluated.

The terms in this proposition can be used directly in integral
adjunction.  On a representable number-field branch as in
Lemma~\ref{lem:integral-CM-adjunction}, let
$j:\operatorname{Spec}\mathcal O_F\to\mathcal X$ be its
normalization map and $\mathfrak f$ its conductor.  Since a
positive branch misses $C$, the open-curve isometry
\eqref{eq:all-level-canonical-arithmetic} pulls back directly
to this proper branch.  Equivalently, $j^*\mathcal O(C)$ is
canonically trivialized by the cusp section; no arithmetic
Chow class $(C,0)$ is needed.  Hence the normal-bundle formula becomes
\begin{equation}\label{eq:branch-adjunction-with-explicit-level}
 \begin{split}
 \widehat{\deg}\,j^*\overline{\mathcal O(D)}
  ={}&\log|\operatorname{Disc}F|
          +\log\operatorname{Norm}\mathfrak f
          +\widehat{\deg}\,j^*\widehat H
          +\widehat{\deg}\,j^*(U_N,0)\\
    &+\frac12\sum_{\sigma:F\hookrightarrow\mathbb C}
       \left(h_\sigma+\log\ell_\sigma
                            +\gamma+\log\pi-\log N\right).
 \end{split}
\end{equation}
As before, both complex-conjugate embeddings are included.
The finite correction is concrete: if $t_{p;a,b}$ is a local
equation of $E_{p;a,b}$ in a regular chart, then
\begin{equation}\label{eq:CM-internal-component-valuation}
 \widehat{\deg}\,j^*(U_N,0)
 =\sum_{p^e\Vert N}\sum_{\substack{a+b=e\\a,b\geq1}}
    p^{\min(a,b)-1}
    \sum_{\mathfrak q\mid p}
       \operatorname{ord}_{\mathfrak q}
                       (j^*t_{p;a,b})\log\operatorname{Norm}\mathfrak q.
\end{equation}
Unit changes of $t_{p;a,b}$ do not affect these valuations.
The formula is local on auxiliary etale covers and descends
with the stack degree convention.  If the branch avoids all
interior components, the term is zero; in particular $U_N=0$
for squarefree $N$.  If its reduction lies on an interior
component, each corresponding summand is nonnegative and may
not be omitted.  At level $N=p^2$ the correction is the middle
reduced component with coefficient one.

\begin{lemma}[Vanishing on a split ordinary CM branch]
\label{lem:split-ordinary-CM-no-interior}
At a place above $p^e\Vert N$, suppose a CM branch has ordinary
good reduction and carries the full action
$\mathcal O_K\otimes\mathbb Z_p\simeq
\mathbb Z_p\times\mathbb Z_p$ on its elliptic curve.  Suppose
also that the cyclic $p^e$-isogeny is linear for this action.
Then the branch reduces to one of the two extreme components
$E_{p;e,0}$ or $E_{p;0,e}$, away from all interior components.
Its contribution to \eqref{eq:CM-internal-component-valuation}
is zero.
\end{lemma}
\begin{proof}
The two idempotents split the ordinary $p$-divisible group into
its connected height-one and etale height-one factors.  On the
characteristic-zero fibre the Tate module is free of rank one
over $\mathbb Z_p\times\mathbb Z_p$.  Therefore a finite
subgroup stable under this action has abstract geometric group
$\mathbb Z/p^r\mathbb Z\oplus\mathbb Z/p^s\mathbb Z$ for
some $r+s=e$.  Cyclicity implies $\min(r,s)=0$.  The finite
flat closure of the subgroup is consequently wholly in one
idempotent factor: the projection to the other factor is zero
generically, hence zero by flatness.  Its connected--etale
type in the ordinary special fibre is thus $(e,0)$ or $(0,e)$.
The ordinary loci of the distinct components are disjoint,
so no interior component meets this specialization.
\end{proof}

The vertical term therefore vanishes under these restricted
hypotheses. The full split maximal-order action must be integral
and compatible with the isogeny; its existence cannot be inferred
from a complex CM point alone. Cycles whose conductor is divisible
by $p$, or which lack that integral action, need not satisfy the
hypotheses.

\paragraph{What this proves and what is still required.}
The Cartier hypothesis and the relative-dualizing/Hodge
comparison at arbitrary level are supplied by the preceding
lemmas and formulas.  The bad-level correction is now the
explicit component valuation \eqref{eq:CM-internal-component-valuation},
not an unspecified canonical-line comparison.  To obtain the
full common-component product formula one still has to
identify the finite CM branch orders and their conductors,
evaluate these component valuations and intersections between
distinct sheets, and compare their weighted sum with the
singular special-cycle multiplicity.  The smooth contribution
to each complete Green function belongs in $h_\sigma$;
rational proportionality factors of the two special vectors
also retain their finite-place contributions.  None of these
data follows just from generalized Cartierness.  Thus
\eqref{eq:branch-adjunction-with-explicit-level} is a new
all-level adjunction input, and is not a claim that the full
arithmetic doubling identity has been established.

\subsection{A finite-place contact calculation at an ordinary CM branch}
\label{subsec:ordinary-split-contact-five}

First we identify precisely which discriminants enter the
remaining comparison.  Let $D=\operatorname{Spec}A$ be a
representable horizontal Cartier branch, with normalization
$j_D:\operatorname{Spec}\mathcal O_F\to D$ and conductor
$\mathfrak f_A$.  At its generic point let
$\pi:E\to E'$ and $j:E\to E'$ be the two special isogenies
of degrees $N$ and $m$, respectively.  Put
\[
 \beta=\pi^\vee j\in\operatorname{End}_F(E),\qquad
 K=\mathbb Q(\beta),\qquad
 f_\beta=[\mathcal O_K:\mathbb Z[\beta]].
\]
Orthogonality gives $\beta+\bar\beta=0$ and
$\beta\bar\beta=Nm$, whence
\begin{equation}\label{eq:beta-order-discriminant}
 \beta^2=-Nm,\qquad
 \operatorname{disc}\mathbb Z[\beta]
      =-4Nm=\operatorname{Disc}(K)f_\beta^2.
\end{equation}
The action on the one-dimensional $F$-space $\operatorname{Lie}E$
embeds $K$ in $F$: the derivative of a nonzero endomorphism in
characteristic zero cannot vanish.  Write $r=[F:K]$, so that
$[F:\mathbb Q]=2r$.

\begin{lemma}[The finite package left by a single branch]
\label{lem:single-branch-discriminant-package}
Suppose the pullback of $U_N$ to $D$ is zero.  The sum of the
conductor/different contribution from integral adjunction and
the universal tangent/KS constant for the primitive pair
$a=b=1$, after removing $-\log v_\perp$, is
\begin{align}
 \mathcal R(D,\beta)
  &:=\log|\operatorname{Disc}(F)|
       +\log\operatorname{Norm}_{F/\mathbb Q}\mathfrak f_A
                                  -r\log(4Nm)\notag\\
  &=\log\operatorname{Norm}_{F/\mathbb Q}\mathfrak D_{F/K}
       +\log\operatorname{Norm}_{F/\mathbb Q}\mathfrak f_A
                                  -2r\log f_\beta.
                         \label{eq:single-branch-residual-package}
\end{align}
Here $\mathfrak D_{F/K}$ is the relative different ideal in
$\mathcal O_F$.  Its norm equals the absolute norm of the
relative discriminant ideal in $\mathcal O_K$.
This expression is only the indicated matching package; it is
not the complete arithmetic product defect.
\end{lemma}

\begin{proof}
The universal constant proved above, combined with the
dualizing--Hodge metric correction, is
$-\log v_\perp-\log(4Nm)$ at each complex point.
The arithmetic degree counts all embeddings with factor
$1/2$, so the parameter-independent contribution is
$-r\log(4Nm)$.  Add the two finite terms from
Lemma~\ref{lem:integral-CM-adjunction}.  Finally use the
discriminant tower formula
\[
 |\operatorname{Disc}(F)|
  =|\operatorname{Disc}(K)|^r
             \operatorname{Norm}_{F/\mathbb Q}\mathfrak D_{F/K}
\]
and \eqref{eq:beta-order-discriminant}.
\end{proof}

The three objects must be distinguished: the number-field
different of $F/\mathbb Q$, the conductor ideal of the branch
order $A\subset\mathcal O_F$, and the integer conductor of
$\mathbb Z[\beta]\subset\mathcal O_K$. Maximal CM on $E$ neither
implies $f_\beta=1$ nor identifies the first two objects with
the quadratic-order conductor.

The following example computes a contact term omitted from the
single-branch package. Take
\begin{equation}\label{eq:five-CM-data}
 N=m=5,\qquad E:\ y^2=x^3-x,\qquad
 \pi=[2+i],\quad j=[-1+2i]=[i]\pi.
\end{equation}
These maps are defined after adjoining $i$.  They both have
degree $5$, the map $\pi$ is cyclic, and
\[
 \pi^\vee j=[5i],\qquad j^\vee\pi=[-5i].
\]
Thus this is a point of the positive special-moduli divisor,
with full maximal CM order $\mathbb Z[i]$.  Its reduction at
$5$ is ordinary, since the displayed curve has trace $-2$
over $\mathbb F_5$ (equivalently, it has eight points there).
In particular $5$ is split and unramified in $K=\mathbb Q(i)$.

\begin{proposition}[An actual canonical--quasi-canonical contact]
\label{prop:ordinary-five-contact}
Fix the reduction of \eqref{eq:five-CM-data} at either prime
above $5$, and add prime-to-$5$ auxiliary full level on both
elliptic curves to remove automorphisms and select the sheet
of the fixed special morphism $j$.  Work after strict
unramified base change, with $W=W(\overline{\mathbb F}_5)$.
Over $W$ the formal local surface
$\mathcal X_0(5)$ is $\operatorname{Spf}W[[T]]$, and the
specified special-moduli sheet has exact equation, up to a unit,
\begin{equation}\label{eq:five-special-local-equation}
 (1+T)^5-1
      =T\Phi_5(1+T)=0.
\end{equation}
The maximal-order canonical branch is $D_0:(T=0)$, and the
other horizontal branch is $D_1:(\Phi_5(1+T)=0)$.
They are distinct and reduced; the second is a ramified
degree-four branch whose characteristic-zero points are
quasi-canonical lifts of CM conductor $5$.  Their local
intersection multiplicity is
\begin{equation}\label{eq:five-contact-length}
 i(D_0,D_1)=
 \operatorname{length}_W
       W[[T]]/(T,\Phi_5(1+T))=1.
\end{equation}
Consequently the restriction of the full sheet to $D_0$
differs from the restriction of $D_0$ itself by this nonzero
contact term.  No interior vertical component is involved.
The length is unchanged by this faithfully flat unramified
base change, so it is the geometric intersection multiplicity
on the original auxiliary-level chart as well.
\end{proposition}

\begin{proof}
Use Serre--Tate parameters $q_1,q_2\in1+\mathfrak m$ for the
source and target ordinary elliptic curves, with bases chosen
compatibly with their polarizations.  For a homomorphism $f$
write $f_{\rm et}$ for its scalar on the physical Tate module
and $f_{\rm et}^\vee$ for the scalar of its dual on the dual
physical Tate module.  The exact homomorphism-lifting criterion
of \cite[Theorem~2.1(4)]{KatzST81} is
\begin{equation}\label{eq:ST-hom-lift-criterion}
             q_1^{f_{\rm et}^\vee}=q_2^{f_{\rm et}}.
\end{equation}
These are equations of functors on Artinian local $W$-algebras,
not merely equations on characteristic-zero points.  The
lift of a homomorphism, when it exists, is unique.

Let $\iota\in\mathbb Z_5$ be the scalar by which $i$ acts on
the physical Tate module, so $\iota^2=-1$.  Put
$r_0=2+\iota$, $s_0=2-\iota$.  Then
$r_0s_0=5$, and one of $r_0,s_0$ is a unit and the other
has valuation one.  The pair $(\pi,j)$ lifts exactly when
\begin{align}
 q_1^{s_0}&=q_2^{r_0},\notag\\
 q_1^{-\iota s_0}&=q_2^{\iota r_0}.
                       \label{eq:five-two-isogeny-equations}
\end{align}
The first equation gives the complete ordinary deformation
space of $\pi$.  Because one exponent is a unit, it eliminates
one of the two variables, leaving a smooth formal curve over
$W$.  Inside that curve write the common value in the first
equation as $z_0$.  The second equation becomes
$z_0^{2\iota}=1$.  The endomorphism $z\mapsto z^{2\iota}$
of the formal multiplicative group is an automorphism, since
$2\iota\in\mathbb Z_5^\times$; hence it is equivalent to
$z_0=1$.

If $r_0$ is a unit, eliminate $q_2$ using
$q_2=q_1^{s_0/r_0}$ and put
$q=q_1^{s_0/(5r_0)}$.  This is a valid formal coordinate
change by a $5$-adic unit, and $q_2=q^5$.
The remaining equation $z_0=1$ is exactly $q^5=1$.
If $s_0$ is a unit, interchange source and target and use
$q=q_2^{r_0/(5s_0)}$.  The same equation results.
Setting $T=q-1$ proves \eqref{eq:five-special-local-equation}.
Degree and the orthogonality relation impose no further
equations: they hold for the fixed special-fibre maps, and
rigidity of homomorphisms gives their identities on every
infinitesimal lift.  The degree is locally constant.

Both CM parameters equal one on $T=0$, the canonical lift.
At a nontrivial fifth root of unity one of $q_1,q_2$ is
primitive of order $5$ and the other equals one.  The source
or target, respectively, has conductor-$5$ CM.  To check the
order directly, the special-fibre endomorphism ring is
$\mathbb Z[i]$, which is already maximal.  The endomorphism
$a+bi$ lifts at a parameter $q$ precisely when
$q^{2b\iota}=1$, by \eqref{eq:ST-hom-lift-criterion}.
For a primitive fifth root this is equivalent to $5\mid b$;
rigidity embeds every lifted endomorphism in the special-fibre
ring.  The lifted order is therefore exactly
$\mathbb Z+5\mathbb Z[i]$.
The polynomial $\Phi_5(1+T)$ is Eisenstein over $W$, so its
branch is integral and reduced, and its normalization is the
cyclotomic DVR $W[\zeta_5]$.  Finally $\Phi_5(1)=5$ proves
\eqref{eq:five-contact-length}.
\end{proof}

This proves a finite compensation, with its multiplicity fixed.
For the data \eqref{eq:five-CM-data},
$\mathbb Z[\beta]=\mathbb Z[5i]$ has conductor $f_\beta=5$,
whereas the branch $D_0$ is unramified and normal at $5$.
Thus the $5$-part of
\eqref{eq:single-branch-residual-package} is $-2r\log5$.
After retaining finite auxiliary level, let $F$ be the residue
field of this canonical branch.  Its completions above $5$
are unramified, by the equation $T=0$ and the prime-to-$5$
etaleness of the auxiliary level.  Summing the contact once
at each such place gives
\begin{equation}\label{eq:five-contact-compensation}
 \sum_{\mathfrak q\mid5}
    i_{\mathfrak q}(D_0,D_1)\log\operatorname{Norm}\mathfrak q
       =[F:\mathbb Q]\log5=2r\log5.
\end{equation}
This cancels the indicated $5$-part exactly.  The coefficient
``once'' is intrinsic to the restriction identity
\[
 \mathcal O(D_0+D_1)|_{D_0}
   =\mathcal O(D_0)|_{D_0}\otimes\mathcal O(D_1)|_{D_0}.
\]
In a square of the complete divisor the two ordered cross
terms must of course both be counted.  The sheets belonging
to $j$ and $-j$ are separate moduli sheets even when their
images coincide, and retain their own multiplicities and
stack weights; they are not an additional factor inserted
into \eqref{eq:five-contact-length}.

\paragraph{Precise scope.}
The proposition proves the local $5$-adic conductor/contact
compensation for this ordinary example. Other horizontal branches
meet the canonical branch, so maximal CM and absence of internal
vertical contact do not suffice for a single-branch product
formula. The calculation does not evaluate the CM ramification
prime $2$, prove a global common-component product formula, or
compare every sheet of an arbitrary-level special cycle. That
global comparison must also include the complete Green regular
values and the finite contributions of rational proportionality
factors. It does not replace $\operatorname{Disc}F$ by a quadratic
CM discriminant.

\begin{corollary}[All split ordinary canonical-branch contacts]
\label{cor:ordinary-split-general-contact}
Let $p$ be any prime, including $2$, and consider an actual
ordinary CM point of the special-moduli problem
$\pi:E\to E'$, $j:E\to E'$ with
\[
 \deg\pi=N,\qquad \deg j=m>0,\qquad
                 \pi^\vee j+j^\vee\pi=0.
\]
Assume the common quadratic CM field $K$ is split at $p$,
the canonical branch carries the full maximal split
$p$-adic CM action, and $\pi$ is linear for that action.
Assume its $p$-primary kernel is cyclic.  Put
\[
 \beta=\pi^\vee j,\qquad
 f_\beta=[\mathcal O_K:\mathbb Z[\beta]],\qquad
                       k=\operatorname{ord}_p(f_\beta).
\]
After strict unramified base change to
$W=W(\overline{\mathbb F}_p)$ and prime-to-$p$ auxiliary level
which removes automorphisms and selects the special-morphism
sheet, its formal equation on the ordinary modular curve is
\begin{equation}\label{eq:general-split-sheet-equation}
 \operatorname{Spf}W[[T]],\qquad
              (1+T)^{p^k}-1=0,
\end{equation}
up to multiplication of the equation by a unit.  The canonical
branch is $D_0:(T=0)$.  Its contact with all other branches of
this sheet, counted with their scheme multiplicities, is
\begin{equation}\label{eq:general-split-contact-length}
 i(D_0,D_{\rm rest})=k,\qquad
 D_{\rm rest}:\quad
       \Psi_k(T)=\frac{(1+T)^{p^k}-1}{T}=0.
\end{equation}
For $k=0$ take $\Psi_0=1$ and the contact to be zero.
For $k>0$ the other branches are precisely the reduced
cyclotomic branches $\Phi_{p^h}(1+T)=0$, $1\leq h\leq k$,
each contributing contact one with $D_0$.
\end{corollary}

\begin{proof}
Take compatible bases of the source and target physical Tate
modules and their polarized duals.  Write $r,s$ for the
scalars of $\pi,\pi^\vee$ and $u,v$ for those of
$j,j^\vee$, respectively.  Thus
\[
 rs=N,\qquad su+rv=0.
\]
All four scalars are integral in $\mathbb Z_p$.  The split
maximal CM action decomposes each $p$-divisible group into
height-one factors.  Cyclicity of the CM-linear $p$-kernel
implies at least one of $r,s$ is a unit: otherwise its
characteristic-zero geometric group would contain two
nonzero cyclic $p$-primary direct summands.  Prime-to-$p$
level imposes only an etale cover of the deformation problem.

By the exact lifting criterion of
\cite[Theorem~2.1(4)]{KatzST81}, with source and target
Serre--Tate parameters $q_1,q_2$, the two lifting equations are
\[
                     q_1^s=q_2^r,\qquad q_1^v=q_2^u.
\]
Suppose first $r$ is a unit.  The first equation eliminates
$q_2=q_1^{s/r}$, identifying the ambient ordinary deformation
space with $\operatorname{Spf}W[[q_1-1]]$.  The second equation
then becomes
\[
 q_1^{v-su/r}=1,\qquad
                v-\frac{su}{r}=\frac{rv-su}{r}.
\]
The two split components of $\beta$ are $su$ and $rv$.
Under $\mathcal O_K\otimes\mathbb Z_p\simeq
\mathbb Z_p\oplus\mathbb Z_p$, the suborder
$\mathbb Z_p[\beta]$ has the $\mathbb Z_p$-basis
$(1,1),(su,rv)$, because $\beta$ satisfies a monic quadratic
equation.  Its index has valuation
\[
 k=\operatorname{ord}_p(rv-su).
\]
This is the determinant index formula, not an identification
of a number-field discriminant with a quadratic discriminant.
Since $r$ is a unit, the remaining exponent has exactly
valuation $k$.  A unit exponent is an automorphism of the
formal multiplicative group, so its equation is equivalent,
as an Artinian deformation functor, to $q_1^{p^k}=1$.
Putting $T=q_1-1$ proves
\eqref{eq:general-split-sheet-equation}.  If $s$ is a unit,
eliminate $q_1$ instead; the remaining exponent is
$(su-rv)/s$ and has the same valuation.  This proves the
assertion in both orientations.

The canonical CM lift is $q_1=q_2=1$.  All equations above
are exact formal equations, since lifts of homomorphisms are
unique and their degree and orthogonality identities persist
by rigidity.  Factoring the remaining equation gives
\[
 (1+T)^{p^k}-1
    =T\prod_{h=1}^k\Phi_{p^h}(1+T),\qquad
                      \Psi_k(0)=p^k.
\]
Each factor $\Phi_{p^h}(1+T)$ is Eisenstein over $W$.
Consequently the branches are reduced, and
\[
 W[[T]]/(T,\Psi_k(T))=W/(p^k)
\]
has length $k$.  This geometric intersection length is
unchanged by the strict unramified base change.

No step requires $p$ odd.  In particular, at $p=2$ it would
be wrong to remove the factor $2$ from $rv-su=-2su$:
the determinant above already includes it in $k$.  The
degree-one branch $\Phi_2(1+T)=T+2$ has contact one, just as
the higher cyclotomic branches do.
\end{proof}

For a proper canonical branch with residue field $F\supset K$
and $r_F=[F:K]$, the local equation $T=0$ is unramified and
normal at $p$, so its branch conductor and relative
different have no $p$-part.  If the hypotheses hold at every
place of this branch above $p$, the package
\eqref{eq:single-branch-residual-package} has $p$-part
$-2r_F k\log p$.  The contact term in the restriction of
the complete sheet to its canonical branch is
\[
 \sum_{\mathfrak q\mid p}
       k\log\operatorname{Norm}\mathfrak q
           =[F:\mathbb Q]k\log p=2r_F k\log p,
\]
and cancels it exactly, giving the full local conductor/contact
compensation for these split ordinary canonical branches. The
divisor-restriction identity counts each other horizontal branch
once; distinct moduli sheets and both ordered cross terms in a
complete square retain their multiplicities. This calculation
covers the canonical outer ordinary branch. The following
subsections prove the all-branch and all-prime finite comparisons.
The global arithmetic doubling comparison is proved in
Theorem~\ref{thm:full-standard-arithmetic-doubling}.

\subsection{All branches in an ordinary split chart}
\label{subsec:ordinary-split-all-branches}

\begin{proposition}[Compensation on every cyclotomic branch]
\label{prop:ordinary-split-all-branch-compensation}
Assume the hypotheses and exact formal chart of
Corollary~\ref{cor:ordinary-split-general-contact}.  Thus, after
strict unramified base change and auxiliary level, the chosen
special-morphism sheet in the regular ordinary surface is
\[
 R=W[[T]],\qquad q=1+T,\qquad f(q)=q^{p^k}-1,
       \qquad W=W(\overline{\mathbb F}_p).
\]
Retain all the reduced horizontal branches
\[
 D_0:(q-1=0),\qquad
 D_h:(\Phi_{p^h}(q)=0),\quad 1\leq h\leq k.
\]
The special point is in the extreme ordinary locus of the
modular special fibre.  Its formal neighborhood misses the
interior components, so $U_N$ restricts to zero on every
$D_h$, not just on $D_0$.
Put $e_0=1,d_0=0$, and for $h\geq1$ put
\[
 e_h=\varphi(p^h),\qquad d_h=h\varphi(p^h)-p^{h-1}.
\]
Then the branch $D_h$ is normal.  Its arithmetic branch-order
conductor is trivial at $p$, its relative different has valuation
$d_h$, and its intersection with the complete union of the other
branches is
\begin{equation}\label{eq:ordinary-all-branch-contact}
 i(D_h,D_{\ne h})=k e_h-d_h.
\end{equation}
Consequently the finite different/conductor--tangent package of
Lemma~\ref{lem:single-branch-discriminant-package}, localized at
this branch and with $U_N=0$, is exactly the negative of its
contact with all the other branches.  This holds for every
$0\leq h\leq k$ and for every prime $p$, including $2$.
\end{proposition}

\begin{proof}
For $h\geq1$, the polynomial
$P_h(T)=\Phi_{p^h}(1+T)$ is Eisenstein over $W$: its reduction
is $T^{e_h}$ and its constant coefficient is $p$.
Thus
\[
 A_h=R/(P_h)=W[\zeta_{p^h}]
\]
is a complete discrete valuation ring, with uniformizer
$\zeta_{p^h}-1$, ramification index $e_h$, and no branch-order
conductor.  For $h=0$ the branch ring is $A_0=W$.
Here normality of $A_h$ must not be confused with maximality of
the elliptic curve's quadratic CM order.  The latter can be
nonmaximal on $D_h$.

The different of a monogenic separable extension $W[\theta]/W$
is generated by $P'(\theta)$, where $P$ is the monic minimal
polynomial.  One elementary justification is Lagrange
interpolation: the functional
$a\mapsto\operatorname{Tr}(a/P'(\theta))$ extracts the
coefficient of $\theta^{e-1}$ in the remainder of $a$ modulo
$P$.  Its pairing on the power basis is integral and has
unit determinant, so the trace-dual lattice is
$P'(\theta)^{-1}W[\theta]$.

At $\zeta=\zeta_{p^h}$, differentiation of
\[
 \Phi_{p^h}(q)=\frac{q^{p^h}-1}{q^{p^{h-1}}-1}
\]
gives
\[
 \Phi_{p^h}'(\zeta)
       =\frac{p^h\zeta^{p^h-1}}
                    {\zeta^{p^{h-1}}-1}.
\]
The denominator is $\zeta_p-1$ and has $A_h$-valuation
$p^{h-1}$: $p$ is a unit times $(\zeta_p-1)^{p-1}$.
The numerator has valuation $h e_h$.  Thus the different
valuation is precisely $d_h$.

Write $f=P_hG_h$ for $h\geq1$, with the harmless coordinate
shift understood.  Since $P_h$ has a simple root in
characteristic zero,
\[
 G_h(\zeta)=\frac{f'(\zeta)}{\Phi_{p^h}'(\zeta)},
          \qquad v_{A_h}(f'(\zeta))=k e_h.
\]
The intersection length with the other branches is
$\operatorname{length}_{A_h}A_h/(G_h(\zeta))$, hence is
$k e_h-d_h$.  This is also the geometric $W$-length since
the residue fields agree.  For $h=0$ the same assertion is
$G_0(1)=p^k$, giving length $k$.

For a direct check of every individual contact, for $h,l>0$
distinct one has
\[
 i(D_h,D_l)=\varphi(p^{\min(h,l)}),
                    \qquad i(D_0,D_h)=1.
\]
Indeed, for $l<h$ the quotient of
$\zeta_{p^h}^{p^l}-1$ by
$\zeta_{p^h}^{p^{l-1}}-1$ has valuation
$p^l-p^{l-1}$; for $l>h$ the cyclotomic sum evaluates to
$p$, of valuation $e_h$.  The lower-branch contacts sum to
$1+\sum_{l=1}^{h-1}\varphi(p^l)=p^{h-1}$, and the upper
ones sum to $(k-h)e_h$, agreeing with
\eqref{eq:ordinary-all-branch-contact}.

To relate these local integers to the number-field package,
let $\mathfrak q$ be a place of a normalized arithmetic
branch with residue degree $f_{\mathfrak q}$ over $p$ whose
strict unramified completion is the displayed $A_h$.
The CM field $K$ is split at $p$, so its local different is
trivial.  The local number-field different contributes
$d_h f_{\mathfrak q}\log p$, and the branch-order conductor
contributes zero.  Since $v_p(f_\beta)=k$, the local share
of $-2[F:K]\log f_\beta$ is
$-k e_h f_{\mathfrak q}\log p$.  Thus the package is
\[
 (d_h-k e_h)\log\operatorname{Norm}\mathfrak q,
\]
whereas restriction of the other branches contributes
\[
 (k e_h-d_h)\log\operatorname{Norm}\mathfrak q.
\]
They cancel at this place.  Strict unramified descent
preserves both lengths and different exponents.  In
particular there is no replacement of the number-field
different by a quadratic-order discriminant.
The contact here occurs once in restriction of the complete
selected sheet to $D_h$.  Other special-morphism sheets,
coincident images after forgetting level, and ordered cross
terms in the full product retain their own multiplicities.

All arguments apply at $2$.  For example $h=1$ has
$\Phi_2(1+T)=T+2$, $e_1=1$, and $d_1=0$; its contact
with all other branches is $k$, exactly as the formula says.
\end{proof}

\paragraph{Scope and an example.}
The earlier canonical-branch calculation is the case $h=0$.
The proposition also treats all the noncanonical cyclotomic,
or quasi-canonical, branches in the same exact ordinary
chart.  These include lifts with nonmaximal quadratic CM
order.  For instance at $p=5,k=1$, the branch
$\Phi_5(q)=0$ has ramification index $4$, different exponent
$3$, and contact one with the canonical branch: its local
package is $(3-4)\log5$, canceled by that contact.
The branch arithmetic ring itself is a DVR; the conductor
of its quadratic endomorphism order is a different object.

\subsection{A finite comparison for reduced Cartier sheets}
\label{subsec:general-Cartier-sheet-different}

The preceding compensation has an intrinsic extension that does
not require a monogenic branch order or a smooth curve over the
base discrete valuation ring. The input is the trace element in
the finite dualizing module of the complete selected sheet.

\begin{theorem}[Different, conductor and contact on a Cartier sheet]
\label{thm:Cartier-sheet-different-contact}
Let $W$ be a complete discrete valuation ring of characteristic
zero with perfect residue field $k$. Let $R$ be a complete
regular local ring of dimension two, flat over $W$, with residue
field $k$. Suppose
\[
 f=\prod_{i=1}^{s}f_i\in R,\qquad A=R/(f)
\]
is finite flat over $W$ and generically etale, where the $f_i$
are distinct, pairwise nonassociate prime elements. Put
$A_i=R/(f_i)$, let $O_i$ be its finite normalization, and set
\[
 g_i=\prod_{j\ne i}f_j,\qquad
 \mathfrak c_i=\{x\in O_i:xO_i\subset A_i\},\qquad
 \mathfrak d_i=\mathfrak D_{O_i/W}.
\]
Each $O_i$ is a complete discrete valuation ring. Write $v_i$
for its normalized valuation, $c_i=v_i(\mathfrak c_i)$,
$d_i=v_i(\mathfrak d_i)$, and
$b_i=[O_i/\mathfrak m_i:k]$.

The finite dualizing module $\omega_{A/W}=\operatorname{Hom}_W(A,W)$
is free of rank one over $A$. Choose a generator $\eta$ and
write the trace functional as
$\operatorname{Tr}_{A/W}=J\eta$. The image $J_i$ of $J$ in
$O_i$ is nonzero and satisfies
\begin{equation}\label{eq:Cartier-sheet-different-ideal}
             J_i O_i=g_i\mathfrak c_i\mathfrak d_i.
\end{equation}
Consequently
\begin{align}
 v_i(J_i)&=d_i+c_i+v_i(g_i),
                         \label{eq:Cartier-sheet-different-contact}\\
 \operatorname{length}_W R/(f_i,g_i)
      &=b_i v_i(g_i)
        =\sum_{j\ne i}\operatorname{length}_W R/(f_i,f_j).
                         \label{eq:Cartier-sheet-contact-length}
\end{align}
For an empty set of other branches, $g_i=1$ and the contact
length is zero. Changing $\eta$ changes $J$ by a unit and
leaves these identities unchanged.
\end{theorem}

\begin{proof}
The quotients $A_i$ are finite torsion-free $W$-modules, hence
finite flat. Their fraction fields are finite separable
extensions of $\operatorname{Frac}(W)$. Completeness and
normalization give the rings $O_i$ stated above. The rings
$A$ and $A_i$ are hypersurfaces in the regular ring $R$, so
are Gorenstein. Their finite dualizing modules are therefore
invertible, and are free because these rings are local.
This also follows from the finite Gorenstein criterion in
\cite[Section~49.16]{StacksDifferent}.

Put $B=\prod_i O_i$, the normalization of $A$. Inside the total
quotient ring, let $\mathfrak C=\operatorname{Hom}_A(B,A)$;
this is the conductor of the complete sheet. Since
\[
 \operatorname{Ann}_A(f_i)=g_iA,\qquad
                 A_i\xrightarrow{\ \cdot g_i\ }g_iA
                 \text{ is an isomorphism},
\]
its $i$th component is
\begin{equation}\label{eq:whole-sheet-conductor}
                    \mathfrak C_i=g_i\mathfrak c_i.
\end{equation}
Indeed an $A$-linear map $O_i\to A$ takes values in
$\operatorname{Ann}_A(f_i)$; under the displayed isomorphism
it is an $A_i$-linear map $O_i\to A_i$, hence multiplication
by an element of $\mathfrak c_i$.

The trace pairing on $A\otimes_W\operatorname{Frac}(W)$
is nondegenerate. It identifies $\eta$ with multiplication
by $J^{-1}$ followed by trace; in particular $J$ is a
nonzerodivisor. Finite duality now gives, as fractional
$B$-modules,
\[
 \operatorname{Hom}_W(B,W)
   =\operatorname{Hom}_A(B,\omega_{A/W})
   =\mathfrak C\eta=J^{-1}\mathfrak C.
\]
The left side has $i$th component $\mathfrak d_i^{-1}$.
Combining this with \eqref{eq:whole-sheet-conductor} gives
\eqref{eq:Cartier-sheet-different-ideal}, including the
conductor with its indicated positive sign.

To justify the contact length also when $A_i$ is nonnormal,
set $Q_i=O_i/A_i$, a finite-length $W$-module. Multiplication
by the nonzero element $g_i$ on $0\to A_i\to O_i\to Q_i\to0$
gives an exact sequence
\[
 0\longrightarrow Q_i[g_i]\longrightarrow A_i/g_iA_i
 \longrightarrow O_i/g_iO_i\longrightarrow Q_i/g_iQ_i
 \longrightarrow0.
\]
The kernel and cokernel of an endomorphism of a finite-length
module have the same length. Thus
$\operatorname{length}_W A_i/g_iA_i
 =\operatorname{length}_W O_i/g_iO_i=b_i v_i(g_i)$.
Apply this also to each $f_j$ and use additivity of the
valuation of their product. This proves
\eqref{eq:Cartier-sheet-contact-length}.
\end{proof}

\begin{corollary}[Computing the trace element]
\label{cor:Cartier-trace-Jacobian}
In the theorem, if $R=W[[T]]$ and $f$ is a distinguished
polynomial, one may take $J=f'(T)$. More generally, if the
finite algebra has a complete-intersection presentation
\[
 A=\bigl(W[x_1,\ldots,x_n,1/h]/(F_1,\ldots,F_n)\bigr)_{\mathfrak m}
\]
with the stated finite-flat hypotheses, its different is
generated by
\[
                    J=\det\left(\frac{\partial F_a}
                                           {\partial x_b}\right).
\]
The same assertion holds for the completed local factor of
such a finite algebra.
\end{corollary}
\begin{proof}
For a monic separable polynomial the trace-dual lattice of
$W[T]/(f)$ is $f'(T)^{-1}W[T]/(f)$; this follows by the
power-basis interpolation argument used in
Proposition~\ref{prop:ordinary-split-all-branch-compensation}.
The distinguished condition identifies this algebra with
$W[[T]]/(f)$. For the displayed complete intersection,
shrink to a quasi-finite neighborhood of $\mathfrak m$;
this does not change the local algebra. Then
\cite[Lemma~49.12.4]{StacksDifferent} supplies a generator
of the dualizing module for which the trace is the stated
Jacobian. A finite algebra over complete $W$ splits into
its completed local factors. Both the trace and its finite
dualizing module restrict to these factors, proving the
last assertion. No unstated monogenic presentation is needed.
\end{proof}

\paragraph{The finite CM comparison as one valuation.}
Apply the theorem to a reduced selected CM sheet on an
actual regular chart. Resolve its normalization into the
branches $O_i$, using this same decomposition for the
conductor and contact terms. At a finite place $\mathfrak q$
above $p$, take $W$ unramified over $\mathbb Z_p$ and let
$O_i$ be the corresponding normalized branch ring.
Let $u_i$ be the valuation on this branch of the pullback
of $U_N$. For the primitive pair $a=b=1$, retain
$\beta=\pi^\vee j$, so that $\beta^2=-Nm$.
The local finite package from adjunction, the contacts with
the other branches of this selected sheet, and the universal
CM constant is exactly
\begin{equation}\label{eq:CM-trace-different-residual}
 \begin{split}
 &(d_i+c_i+v_i(g_i)+u_i-v_i(2\beta))
                                  \log\operatorname{Norm}\mathfrak q\\
 &\hspace{15mm}=(v_i(J_i)+u_i-v_i(2\beta))
                                  \log\operatorname{Norm}\mathfrak q.
 \end{split}
\end{equation}
Indeed, the conductor, different and contact terms are those
of the theorem, and the vertical contribution is
\eqref{eq:CM-internal-component-valuation}. If the global
branch field is $F\supset K=\mathbb Q(\beta)$, then
\[
 |\operatorname{Norm}_{F/\mathbb Q}(2\beta)|
     =(4Nm)^{[F:K]}.
\]
Thus the parameter-independent term $-[F:K]\log(4Nm)$
in Lemma~\ref{lem:single-branch-discriminant-package}
is precisely the negative sum of the local $v_i(2\beta)$
terms with their residue-field weights. This argument uses
the absolute number-field different, so it applies equally
to split, inert and ramified quadratic CM fields.

Formula \eqref{eq:CM-trace-different-residual} reduces the
vanishing of this finite package to the specific equality
\begin{equation}\label{eq:CM-trace-different-target}
                         v_i(J_i)+u_i=v_i(2\beta).
\end{equation}
For the ordinary charts already computed,
$J=p^kq^{p^k-1}$, $u_i=0$, and split maximal-order
coordinates give $v_i(2\beta)=e_i k$; hence this recovers
the proven compensation on every cyclotomic branch.
The intrinsic trace theorem alone does not establish
\eqref{eq:CM-trace-different-target}. The next subsections
compute the actual cyclic charts and then prove this comparison
at every prime by a paired conormal calculation and the
cyclic difference-divisor identity.

All conductors here refer to the chosen formal branch
orders. More explicitly, if one global branch completes to
the union of formal branches indexed by $I$, its conductor
on $O_i$, for $i\in I$, has valuation
\[
 c_i+\sum_{\substack{j\in I\\j\ne i}}v_i(f_j).
\]
Only contacts with $j\notin I$ may then be added; the total
is $c_i+v_i(g_i)$, as in \eqref{eq:whole-sheet-conductor}.
This proves compatibility with the global branch decomposition
and prevents counting an internal contact twice. Distinct
special-morphism sheets and rational proportionality factors
also retain their own contributions.

For example, at an odd prime, the single branch
$A_1=W[T]/(T^2-p^3)$ has normalization $O_1=W[\pi]$ with
$\pi^2=p$ and $T=p\pi$. Its conductor is $pO_1$, of
valuation $2$, and its different has valuation $1$.
The derivative $2p\pi$ has valuation $3$, exactly their
sum; omitting the branch conductor would already fail in
this simple nonnormal example. More generally regularity
and Cartierness do not imply a monogenic presentation:
$R=W[[x,y]]/(xy-p)$ is regular, but
$R/(x^2+y^2)$ has special fibre
$k[[x,y]]/(xy,x^2+y^2)$ of embedding dimension two and
cannot be a local monogenic $W$-algebra. The intrinsic
trace formulation covers this case as well.

\subsection{The primitive cyclotomic chart at an inner ordinary point}
\label{subsec:inner-ordinary-cyclic-chart}

\begin{theorem}[The actual inner ordinary cyclic chart]
\label{thm:inner-ordinary-cyclic-chart}
Let $k$ be an algebraically closed field of characteristic $p>0$,
let $W=W(k)$, and let
$\pi_0:E_0\longrightarrow E'_0$ be an ordinary cyclic isogeny
of degree $N$, in the group-scheme sense defining $\mathcal Y_0(N)$.
Fix the prime-to-$p$ data, and add prime-to-$p$ auxiliary level
if needed to remove automorphisms. Choose compatible bases of
the physical Tate modules and their polarized duals.
Let $r,s\in\mathbb Z_p$ be the respective scalars of
$\pi_0,\pi_0^\vee$, and suppose
\[
 r=p^aR_0,\qquad s=p^bS_0,\qquad
 1\leq a\leq b,\qquad R_0,S_0\in\mathbb Z_p^\times.
\]
In particular $a+b=v_p(N)$.
Write $q_1,q_2$ for the source and target Serre--Tate parameters.
Set
\[
 x=q_1^{S_0},\qquad y=q_2^{R_0},\qquad
 z=x^{p^{b-a}}/y,\qquad X=x-1,\quad Z=z-1.
\]
Then the complete local ring of the actual cyclic deformation
problem is
\begin{equation}\label{eq:inner-cyclic-local-ring}
 R_{\rm cyc}
  =W[[X,Z]]/\bigl(\Phi_{p^a}(1+Z)\bigr)
  \simeq V_a[[X]],\qquad V_a=W[\zeta_{p^a}],
\end{equation}
where the isomorphism sends $z$ to $\zeta_{p^a}$.
Here $V_a$ is a complete discrete valuation ring, totally
ramified over $W$ of degree $\varphi(p^a)$.
In this chart
\begin{equation}\label{eq:inner-cyclic-source-target}
 q_1=(1+X)^{S_0^{-1}},\qquad
 q_2=\left(\zeta_{p^a}^{-1}
                 (1+X)^{p^{b-a}}\right)^{R_0^{-1}}.
\end{equation}
Thus the chart is regular and its special fibre has one
reduced component with multiplicity $\varphi(p^a)$.
The assertion includes $p=2$.
\end{theorem}

\begin{proof}
We separate the deformation equation, generic cyclicity, and
identification of the integral cyclic locus.

By \cite[Theorem~2.1(1), (2), (4)]{KatzST81}, the marked
deformations of the two ordinary elliptic curves are classified
by $q_1,q_2$, and $\pi_0$ lifts if and only if
\[
                         q_1^s=q_2^r .
\]
Such a lift is unique. It is an isogeny of the same degree:
the identities with its dual lift by rigidity, or equivalently
its finite locally free kernel is the deformation of that
of $\pi_0$. Therefore the complete local ring of the
selected degree-$N$ homomorphism problem before imposing
cyclicity is
\[
 S=W[[q_1-1,q_2-1]]/(q_1^s-q_2^r).
\]
Unit powers are automorphisms of the formal multiplicative
group. The changes to $x,y$, and then to $x,z$, are
invertible changes of formal coordinates. Since $y$ is a unit,
\begin{equation}\label{eq:inner-all-isogeny-ring}
 S\simeq W[[X,Z]]/\bigl((1+Z)^{p^a}-1\bigr).
\end{equation}

Next we determine which characteristic-zero components of
\eqref{eq:inner-all-isogeny-ring} have cyclic $p$-primary kernel. Over $k$, the ordinary
$p$-divisible groups split into their connected and etale
factors. The physical scalar of $\pi_0$ is divisible by $p$,
as is the scalar of its dual; hence $\pi_0$ vanishes on both
factors of $E_0[p]$. Consequently it kills $E_0[p]$ and
factors uniquely as
\[
                \pi_0=\rho_0\circ[p]_{E_0}
\]
for an actual homomorphism $\rho_0:E_0\to E'_0$.
This uses the quotient $E_0/E_0[p]\simeq E_0$, and does
not merely assert divisibility in a $p$-adic Hom space.
The physical scalars of $\rho_0,\rho_0^\vee$ are $r/p,s/p$.
Katz's criterion therefore says that $\rho_0$ lifts exactly
on the closed formal locus
\begin{equation}\label{eq:inner-division-locus}
 q_1^{s/p}=q_2^{r/p}
 \quad\Longleftrightarrow\quad
                    z^{p^{a-1}}=1 .
\end{equation}
If it lifts, uniqueness implies that its composition with
$[p]$ is the selected lift of $\pi_0$.
Conversely any factorization of that lift through $[p]$
reduces to $\rho_0$ and hence satisfies
\eqref{eq:inner-division-locus}.

This also detects generic divisibility; it
does not assume that a generic homomorphism extends without
proof. The reduced factors of \eqref{eq:inner-all-isogeny-ring} are
\[
 S_h=W[[X,Z]]/\bigl(\Phi_{p^h}(1+Z)\bigr)
       \simeq V_h[[X]],\qquad 0\leq h\leq a,
\]
where $\Phi_1(1+Z)=Z$, $V_0=W$, and $z=\zeta_{p^h}$.
These rings are integral and $W$-flat. Suppose that the
generic lift on $S_h$ kills $E[p]$. The scheme $E[p]$ is
finite flat over $S_h$, so its characteristic-zero part
is schematically dense. The equalizer of the two morphisms
$\pi|_{E[p]}$ and $0$ is closed because $E'$ is separated.
Their generic equality therefore implies equality over
$S_h$. Equivalently, in the finite locally free coordinate
modules no nonzero section can be killed by a power of $p$.
The quotient by $E[p]$ then gives a factorization over
$S_h$, so its reductions to all Artinian quotients satisfy
\eqref{eq:inner-division-locus}. Thus
$\zeta_{p^h}^{p^{a-1}}=1$.
The converse follows directly by lifting $\rho_0$ using
the same equation. Passing from the Artinian homomorphisms
to the complete local ring is legitimate by uniqueness and
effectivity of homomorphisms of proper elliptic curves
\cite[Theorem~3.4.2(b)]{Cesnavicius2017};
one can equally make the equalizer argument on the actual
completed isogeny chart.

Over a characteristic-zero field, a finite subgroup of
an elliptic curve of order $p^{a+b}$ is noncyclic if and
only if it contains the full subgroup $E[p]$: this is the
elementary classification of subgroups of
$(\mathbb Q_p/\mathbb Z_p)^2$. Its isogeny then factors
through $[p]$, and conversely such a factorization makes
the kernel noncyclic. Hence the generic factor $S_h$ is
cyclic precisely when $h=a$.
The type of a finite etale group scheme is locally constant
on its base, so this is the entire characteristic-zero
cyclic locus, rather than just a statement at one generic
point. The fixed prime-to-$p$ kernel is already cyclic
and is unaffected by this ordinary deformation.

Finally, cyclicity is a universally closed condition for
finite subgroup schemes of elliptic curves
\cite[Proposition~4.1.4, citing Katz--Mazur
Theorem~6.4.1]{Zhu23}. Therefore the cyclic deformation
functor is represented by a quotient $S/I_{\rm cyc}$.
It is the completed local ring of the actual
$\mathcal Y_0(N)$ with the auxiliary level, not a newly
defined candidate model. The source morphism from
$\mathcal Y_0(N)$ to the elliptic-curve stack is finite
flat \cite[Lemma~4.2.7]{Zhu23}. After the indicated etale
auxiliary level and completion, it follows that
$S/I_{\rm cyc}$ is $W$-flat. Thus $I_{\rm cyc}$ is
$p$-saturated:
\[
           I_{\rm cyc}=S\cap I_{\rm cyc}[1/p].
\]
The generic calculation gives
$I_{\rm cyc}[1/p]=(\Phi_{p^a}(z))[1/p]$.
The ideal $(\Phi_{p^a}(z))\subset S$ is itself
$p$-saturated, since its quotient is $V_a[[X]]$,
which is $W$-flat. Taking contractions proves
$I_{\rm cyc}=(\Phi_{p^a}(z))$, proving
\eqref{eq:inner-cyclic-local-ring}.
In particular no additional embedded or vertical
thickening has been inferred from generic points alone.

The polynomial $\Phi_{p^a}(1+Z)$ is Eisenstein over $W$.
Its root $\zeta_{p^a}-1$ is a uniformizer of $V_a$,
and
$p\in V_a^\times(\zeta_{p^a}-1)^{\varphi(p^a)}$.
This proves regularity and the asserted special-fibre
multiplicity. Solving for $q_1,q_2$ gives
\eqref{eq:inner-cyclic-source-target}.
\end{proof}

\begin{remark}[What the lower cyclotomic factors mean]
The lower factors $0\leq h<a$ are genuine deformations of
the same special isogeny in the degree-$N$ homomorphism
problem, but their generic isogenies have a common
$[p]$ factor. They are excluded from the actual cyclic
model. It would be incorrect to use the full equation
$z^{p^a}=1$ as the completed local ring of
$\mathcal Y_0(N)$ at this inner ordinary point.
In characteristic $p$, the group-scheme definition of
cyclicity may hold even though the isogeny factors
through $[p]$; only the characteristic-zero equivalence
was used in the proof.
\end{remark}

\begin{remark}[Orientation of source and target]
Here $a=v_p(r)$ is the valuation on the physical Tate
module of $\pi$, and $b=v_p(s)$ is that for $\pi^\vee$.
For $a\leq b$, the source parameter is a unit power of
$1+X$, whereas the target parameter has an additional
$p^{b-a}$ power. Thus the source forgetful map on this
chart has cyclotomic different exponent
$a\varphi(p^a)-p^{a-1}$; the target has the corresponding
additional power contribution. Component labels in an
external source must be identified with these scalars
before comparing ordered different formulas. The fibre
multiplicity $\varphi(p^{\min(a,b)})$ is symmetric and
has no such ambiguity.
In the manuscript's convention that $E_{p;A,B}$ has
connected kernel of order $p^A$ and etale kernel of
order $p^B$, the chart here is $E_{p;b,a}$.
Thus its source different is
$p^{a-1}((p-1)a-1)
 =a\varphi(p^a)-p^{a-1}$, in agreement with that
ordered component formula.
\end{remark}

\subsection{The finite comparison on an inner ordinary sheet}
\label{subsec:inner-ordinary-finite-compensation}

\begin{theorem}[Compensation on the actual inner ordinary chart]
\label{thm:inner-ordinary-finite-compensation}
Use the actual cyclic chart of
Theorem~\ref{thm:inner-ordinary-cyclic-chart}, with
\[
 r=p^aR_0,\quad s=p^bS_0,\quad
 1\leq a\leq b,\quad R_0,S_0\in\mathbb Z_p^\times,
 \qquad V_a=W[\zeta],\quad \zeta=\zeta_{p^a}.
\]
Let $j_0:E_0\to E'_0$ be a positive special isogeny,
with physical scalars $u,v$ for $j_0,j_0^\vee$.
Thus $uv=m>0$ and $su+rv=0$.
Put
\[
 \rho=u/R_0=p^l\epsilon,\quad
 l=v_p(u)\geq0,\quad \epsilon\in\mathbb Z_p^\times,
 \qquad \delta=v_p(2),\quad n=b-a+l+\delta.
\]
There is a unit-power change of the source Serre--Tate
coordinate for which the complete selected special-morphism
sheet is
\begin{equation}\label{eq:inner-special-sheet-polynomial}
 A=V_a[X]/(X^{p^n}-\zeta^{p^l}),\qquad X\equiv1
                                                  \pmod{\mathfrak m}.
\end{equation}
It is finite flat and generically etale over $W$.
Its trace element can be chosen as
\begin{equation}\label{eq:inner-special-trace-element}
             J=\Phi_{p^a}'(\zeta)\,p^nX^{p^n-1}.
\end{equation}
For every normalized branch $O_i$ of this complete sheet,
put $e_i=v_i(p)$ and let $u_i$ be the valuation of the
pullback of $U_N$. Then
\begin{equation}\label{eq:inner-special-trace-cancellation}
 \begin{split}
 v_i(J_i)&=e_i(b+l+\delta)-\frac{e_i}{p-1},\\
 u_i&=\frac{e_i}{p-1},\\
 v_i(2\beta)&=e_i(b+l+\delta),\qquad \beta=\pi^\vee j.
 \end{split}
\end{equation}
Consequently the finite different/conductor/contact package
\eqref{eq:CM-trace-different-residual} is zero on every
branch. The conclusion includes $p=2$ and uses the whole
selected sheet, with all its internal contacts.
\end{theorem}

\begin{proof}
Write $x=q_1^{S_0}$ and $y=q_2^{R_0}$, as in the
actual cyclic chart. Its equation is
$y=\zeta^{-1}x^{p^{b-a}}$.
The equation for lifting $j_0$ is $q_1^v=q_2^u$.
Orthogonality gives $v=-p^{b-a}S_0u/R_0$, hence this
equation is exactly
\[
                   x^{2p^{b-a}\rho}=\zeta^\rho.
\]
Set $\rho=p^l\epsilon$. Raising to the unit
$\epsilon^{-1}$ is an automorphism of the formal
multiplicative group, so the equation is equivalent
as an Artinian deformation functor to
\[
                   x^{2p^{b-a+l}}=\zeta^{p^l}.
\]
The unit-power automorphism acts on the entire ratio of the two
sides, so no unit exponent remains on the right-hand side.
Put $X=x^{2/p^\delta}$. Since $2/p^\delta$ is a
$p$-adic unit, this is a valid formal coordinate, and
the equation becomes \eqref{eq:inner-special-sheet-polynomial}.
The exponent $n$ is nonnegative. In $T=X-1$ its
polynomial is distinguished, including the linear case
$n=0$. Thus its formal quotient is the displayed finite
algebra, free of rank $p^n$ over $V_a$. Its generic
derivative $p^nX^{p^n-1}$ is nonzero at every root,
so it is generically etale and reduced.

Over $W$, a finite polynomial complete-intersection
presentation uses the equations
\[
 \Phi_{p^a}(Z)=0,\qquad X^{p^n}-Z^{p^l}=0.
\]
The triangular Jacobian is exactly
\eqref{eq:inner-special-trace-element}; the
trace--Jacobian formula of
Corollary~\ref{cor:Cartier-trace-Jacobian} applies.
The different exponent of $V_a/W$ is
$a\varphi(p^a)-p^{a-1}$. Therefore
\[
 v_i\!\left(\Phi_{p^a}'(\zeta)\right)
   =\frac{e_i}{\varphi(p^a)}
                      \bigl(a\varphi(p^a)-p^{a-1}\bigr)
   =e_i a-\frac{e_i}{p-1}.
\]
Since $X$ is a unit, adding $v_i(p^n)=e_i n$
proves the first formula in
\eqref{eq:inner-special-trace-cancellation}.

The physical valuation is $a$ and the connected
valuation is $b$, so this chart is $E_{p;b,a}$ in
the manuscript's component convention. The reduced
component has local equation $\zeta-1$ and coefficient
$p^{a-1}$ in $U_N$. Consequently
\[
 u_i=p^{a-1}v_i(\zeta-1)
      =p^{a-1}\frac{e_i}{\varphi(p^a)}
      =\frac{e_i}{p-1}.
\]
This interior contribution must therefore be retained under
ordinary reduction.
Finally the two split components of $\beta$ on the
special-fibre Tate factors are $su$ and $rv=-su$.
Their valuations are both $b+l$. Under either
compatible characteristic-zero tangent embedding
the valuation of $2\beta$ on $O_i$ is therefore
$e_i(b+l+\delta)$. The same valuation follows
from the fixed equation $\beta^2=-Nm$.
Adding the first two formulas proves the target
$v_i(J_i)+u_i=v_i(2\beta)$ on every branch.
Multiplication by each residue-degree weight in
$\log\operatorname{Norm}\mathfrak q$ preserves the
equality. The whole-sheet different theorem accounts
for all branch conductors and contacts exactly once.
\end{proof}

If $a>b$, interchange source and target. The involution sends
$(\pi,j)$ to $(\pi^\vee,j^\vee)$, preserves the norms and
orthogonality, and exchanges $a,b$. The divisor $U_N$ is
symmetric in the two component indices, so the same valuation
identity follows. Together with
Proposition~\ref{prop:ordinary-split-all-branch-compensation},
this proves the finite comparison at every ordinary point.

\paragraph{A middle-component check.}
Take $p=5$, $N=25$, $m=1$, and the ordinary curve
$E_0:y^2=x^3-x$ over $\overline{\mathbb F}_5$.
With $i^2=-1$, take $\pi_0=[5]$ and $j_0=[i]$.
Its kernel $E_0[5]\simeq\mu_5\times\mathbb Z/5$
is cyclic in the special-fibre Katz--Mazur sense:
an etale generator traverses the five components,
each five times as a divisor. Thus the actual point
lies on the inner ordinary chart with $a=b=1$.
Here $l=0$, $\delta=0$, $n=0$, and the selected
sheet is $A=V_1=W[\zeta_5]$. It is normal with
different exponent $3$, no other horizontal branch
of this selected sheet, and $U_N$ valuation $1$.
The target is $v(2\beta)=v(10i)=4$, giving
$3+1=4$. This example exhibits the necessity of
the interior term. 

\subsection{The paired conormal calculation at every finite prime}
\label{subsec:all-prime-paired-conormal}

The following argument uses actual morphism-lifting functors.
In particular the two equations before imposing cyclicity are
the equations for lifting the fixed morphisms, with all their
deformation multiplicities.

\begin{lemma}[The two-morphism trace element]
\label{lem:paired-hom-trace-character}
Let $k$ be algebraically closed of characteristic $p$,
$W=W(k)$, and let
\[
                 R=W[[t_1,t_2]]
\]
be the marked deformation chart of a pair of elliptic curves
$(E_0,E'_0)$. Fix isogenies
$\pi_0,j_0:E_0\longrightarrow E'_0$ of positive degrees
$N,m$, satisfying
\[
 \pi_0^\vee j_0+j_0^\vee\pi_0=0 .
\]
Suppose the actual individual lifting functors are Cartier
divisors with equations $F_\pi,F_j\in R$. Put
\[
 I=(F_\pi,F_j),\qquad A_{\rm all}=R/I .
\]
Then $A_{\rm all}$ is finite flat and generically etale over
$W$. Let $\beta=\pi^\vee j$ on its universal source curve,
and let $\alpha$ be its tangent scalar. In a generator of
$\operatorname{Hom}_W(A_{\rm all},W)$, the trace element
$J_{\rm all}$ satisfies
\begin{equation}\label{eq:paired-all-hom-trace}
 (J_{\rm all})=
 \left(\det\frac{\partial(F_\pi,F_j)}
                         {\partial(t_1,t_2)}\right)
                         =(2\alpha)
                 \quad\text{as ideals of }A_{\rm all}.
\end{equation}
The assertion includes $p=2$.
\end{lemma}

\begin{proof}
The relations satisfied by the two special-fibre morphisms
continue on their simultaneous lifting locus by rigidity.
Thus $\beta^\vee=-\beta$ and $\beta^2=-Nm$ there.
The special fibre of this locus is zero-dimensional.
Indeed a positive-dimensional component of the actual
two-curve chart would have a generic source with
transcendental $j$-invariant; its endomorphism ring is
$\mathbb Z$. The nonzero trace-zero endomorphism $\beta$
cannot occur on such a curve. A fixed source admits only
finitely many degree-$N$ targets and morphisms, so a
positive-dimensional component cannot have constant source
$j$ either. Hence $A_{\rm all}/p$ is Artinian.
Completeness makes $A_{\rm all}$ finite over $W$.
The ideal $(p,F_\pi,F_j)$ is therefore a parameter ideal
in the three-dimensional regular ring $R$. Its generators,
in any order, form a regular sequence. Consequently
$F_\pi,F_j$ is a regular sequence, and $p$ is a
nonzerodivisor on $A_{\rm all}$, proving flatness.

Set $S=R/I^2$. Since $I/I^2\simeq A_{\rm all}^{\,2}$,
both $S$ and $A_{\rm all}$ are finite free over $W$.
For $n\geq1$ put
\[
 S_n=S/p^n,\quad A_n=A_{\rm all}/p^n,\quad
 I_n=\ker(S_n\longrightarrow A_n)\simeq A_n^{\,2}.
\]
Equip the square-zero ideal $I_n$ with its trivial
nilpotent divided powers. The local nilpotence required in
Grothendieck--Messing theory concerns this kernel, whose
divided-power filtration is zero from degree two onward.
This applies in residue characteristic two as well. Flatness
over $W$ makes the intersection of the square-zero ideal with
the ideal generated by $p$ equal to its $p$-multiple. On that
intersection, all divided powers of degree at least two vanish,
so the trivial structure is compatible with the canonical
divided powers on the base. The elliptic curves over $R$
already provide lifts to $S_n$; thus the local
liftability condition in
\cite[Chapter~V, Theorem~1.6]{Messing1972} is satisfied.
The obstruction to lifting each of $\pi,j$ to these
chosen curves is its failure to preserve the Hodge
filtration. After choosing generators of the two
relevant Hom lines over $A_{\rm all}$, write the
compatible obstructions as
\[
 \begin{pmatrix}h_{\pi,n}\\h_{j,n}\end{pmatrix}
       =C_n
          \begin{pmatrix}F_\pi\\F_j\end{pmatrix},
                    \qquad C_n\in M_2(A_n).
\]

The matrix $C_n$ is invertible. Put
$J_n=(h_{\pi,n},h_{j,n})\subset I_n$,
$\overline S_n=S_n/J_n$, and $\overline I_n=I_n/J_n$.
The quotient need not be flat over $W/p^nW$, but $A_n$ is.
Tensoring the exact sequence
\[
 0\longrightarrow\overline I_n\longrightarrow\overline S_n
   \longrightarrow A_n\longrightarrow0
\]
with $W/pW$ over $W/p^nW$ is therefore left exact. Hence
\[
 \overline I_n\cap p\overline S_n=p\overline I_n.
\]
The remaining kernel is square-zero and carries trivial divided
powers. On the displayed intersection these agree with the
canonical base divided powers: for $px$ with
$x\in\overline I_n$, all divided powers of degree at least two
vanish because $x^2=0$. Thus the quotient used to test the
obstructions satisfies the same PD compatibility, also for $p=2$.
Both obstructions now
vanish, so Grothendieck--Messing and the general
Serre--Tate equivalence
\cite[Theorem~1.2.1]{KatzST81} lift both actual
elliptic-curve morphisms. The representing property of
their simultaneous lifting functor forces the resulting
map from $R$ to kill $I$. Hence
$(h_{\pi,n},h_{j,n})=I_n$.
Since $I_n$ is free on the displayed two equations,
$C_n$ is surjective and therefore invertible.
Choose the generators of the two Hom lines compatibly over
$A_{\rm all}$. Crystalline base change makes the obstruction
elements compatible under reduction in $n$. Since the classes
of $F_\pi,F_j$ form the specified basis of $I_n$ over $A_n$,
the coefficient matrices $C_n$ are compatible as well.
Their inverse limit is a matrix $C$ over $A_{\rm all}$. Its
determinant is a unit, since its reduction is the determinant
of the invertible matrix $C_1$. Hence
$C\in\operatorname{GL}_2(A_{\rm all})$.

Map $S_n$ to
$A_n[\epsilon_1,\epsilon_2]/(\epsilon_1,\epsilon_2)^2$
by $t_i\mapsto t_i+\epsilon_i$.
This is a map of the chosen PD thickenings.
Relative to the constant lifts of the two curves,
their Hodge lines vary by
$\kappa_1\epsilon_1$ and $\kappa_2\epsilon_2$,
where the two Kodaira--Spencer coefficients
$\kappa_i$ are units. Their reductions at the closed
point are the isomorphisms from the universal elliptic
deformation tangent lines to the corresponding Hodge
deformation lines; Nakayama's lemma gives the assertion.

Let $r,s$ be the tangent scalars of $\pi,\pi^\vee$,
and $u,v$ those of $j,j^\vee$, in compatible frames.
For the contravariant de Rham map of $\pi$, the
scalar on the Hodge line is $r$ and on its quotient
is $s$; similarly the two scalars for $j$ are $u,v$.
The first-order obstruction rows are therefore
\[
              (-r\kappa_1,\ s\kappa_2),
       \qquad (-u\kappa_1,\ v\kappa_2).
\]
The possible upper-right matrix entries contribute
only quadratic infinitesimals. Moreover
$su+rv=0$ and $\alpha=su$. Taking determinants
after the Taylor pullback consequently gives
\[
 \det(C)\det\frac{\partial(F_\pi,F_j)}
                         {\partial(t_1,t_2)}
        =\kappa_1\kappa_2(su-rv)
        =2\alpha\kappa_1\kappa_2
                       \quad\text{in }A_{\rm all}.
\]
As an ideal identity, this is unaffected by changes of frames or
by the sign introduced by ordering the variables. The factor $2$
is retained, including dyadically. Since $\alpha^2=-Nm$, the
Jacobian becomes a unit after inverting $p$. The Jacobian criterion
therefore gives an etale finite generic fibre without assuming
that the CM lifting equations are reduced in advance.

The trace--Jacobian formula can be applied before completion,
as the following algebraic presentation shows.
The two individual Hom divisors are actual
generalized Cartier divisors on the product of
the elliptic moduli stacks
\cite[Proposition~4.3.3 and Definition~4.3.2]{Zhu23}.
After auxiliary level and an etale scheme cover,
choose their fixed closed-immersion sheets.
A smooth affine chart is etale over a polynomial
chart after adding auxiliary variables and equations
whose Jacobian block is a unit.
The two-morphism intersection is quasi-finite over
$W$ near the chosen point, by the fibre argument
above. Applying
\cite[Lemma~49.12.4]{StacksDifferent} to this
polynomial complete intersection gives its trace
element as its full Jacobian.
Eliminating the etale variables multiplies this
Jacobian by a unit, leaving the displayed
two-variable determinant. One may equivalently
use a finite neighborhood supplied by Zariski's
Main Theorem and the local factor over the
henselian ring $W$.
Its completion is the finite flat algebra
$A_{\rm all}$; trace and the finite dualizing
module restrict to this factor and commute with
completion. This proves the first equality in
\eqref{eq:paired-all-hom-trace}, and the obstruction
calculation proves the second.
\end{proof}

\begin{lemma}[The divisor removed by cyclicity]
\label{lem:cyclic-removed-divisor-U}
On the elliptic locus $\mathcal Y_0(N)$ of the actual
regular cyclic model, let $G_p$ denote
the local equation for the locus on which the
universal $\pi$ is divisible by $[p]$, with $G_p=1$
where this locus is empty. Then it is an effective
Cartier divisor, supported in the special fibre,
and
\begin{equation}\label{eq:cyclic-removed-divisor-U}
 \operatorname{div}(G_p)=
   \sum_{\substack{A+B=v_p(N)\\A,B\geq1}}
                        p^{\min(A,B)-1}E_{p;A,B}
                         =U_{N,p}.
\end{equation}
Locally in the product deformation chart one has
\[
               F_\pi=F_{\rm cyc}G_p
                       \quad\text{up to a unit}.
\]
\end{lemma}

\begin{proof}
Divisibility means that there is a morphism
$\rho$ with $\pi=[p]\rho$; it is unique, since
Hom groups have no $p$-torsion. Thus this is an
intrinsic closed lifting condition and local
equations differ only by units.
At a supersingular point, Proposition~5.2.3
and Theorems~6.2.2--6.2.3 of \cite{Zhu23} give
precisely
\[
 F_\pi=f_\pi,\qquad
 F_{\rm cyc}=d_\pi=f_\pi/f_{\pi/p},\qquad
 G_p=f_{\pi/p}.
\]
Here $f_{\pi/p}=1$ if $\pi/p$ is not integral.
These results include every prime $p$; the
difference divisor is the actual cyclic
deformation space and is regular, and its
equation is coprime to all the removed factors.
At ordinary points the same statements follow
from the actual Serre--Tate cyclic chart proved
in Theorem~\ref{thm:inner-ordinary-cyclic-chart};
the cases with one slope zero have $G_p$ a unit.
Thus restriction of $G_p$ to the regular
cyclic surface is everywhere a nonzerodivisor.
There is no horizontal divisor: a characteristic-zero
cyclic kernel cannot contain the whole $E[p]$.

The vertical multiplicities can be computed at the generic
ordinary points of the components. Write the physical and dual
slopes as $a,b$ with $1\leq a\leq b$. In the actual chart
$V_a[[X]]$, $V_a=W[\zeta_{p^a}]$, the equation
for lifting $\pi/p$ is, up to a unit,
\[
                    \zeta_{p^a}^{\,p^{a-1}}-1 .
\]
Its valuation at the reduced component
$\zeta_{p^a}-1=0$ is $p^{a-1}$.
For example this follows by multiplying
$1-\zeta_{p^a}^{\,j}$ over the primitive
$p$-th roots, or from
$v_p(1-\zeta_{p^h})=1/\varphi(p^h)$.
The chart is $E_{p;b,a}$ because $a$ is the
etale-kernel exponent and $b$ the connected-kernel
exponent. Interchanging source and target handles
$a>b$, yielding the coefficient in
\eqref{eq:cyclic-removed-divisor-U}.
The components with one exponent zero have
coefficient zero. Equality of all height-one
valuations on a regular surface proves the
Cartier-divisor identity, including its
supersingular specializations. No additional
divisor can be supported at isolated closed points.
\end{proof}

\begin{theorem}[All finite selected-sheet packages]
\label{thm:all-prime-finite-CM-compensation}
Fix an actual positive selected special-isogeny
sheet on the regular cyclic model, with
$\beta=\pi^\vee j$ and $\beta^2=-Nm$.
At every prime and every normalized formal
branch of the complete selected sheet, let
$J_i$ be its trace element and $u_i$ the
valuation of the pullback of $U_N$.
Then
\begin{equation}\label{eq:all-prime-trace-U-target}
                         v_i(J_i)+u_i=v_i(2\beta).
\end{equation}
Consequently the finite selected-sheet package
\eqref{eq:CM-trace-different-residual}, for the
primitive pair $a=b=1$, vanishes at every finite
place. The complete-sheet trace includes its
different, conductor, and internal contact terms.
\end{theorem}

\begin{proof}
At a supersingular point use the actual product
chart and the equations of the preceding lemmas.
The cyclic selected sheet is
$A_{\rm cyc}=R/(F_{\rm cyc},F_j)$.
The positive-CM special-fibre argument again
makes it finite flat. Its trace
element is its complete-intersection Jacobian,
with the algebraic justification given in
Lemma~\ref{lem:paired-hom-trace-character};
the cyclic equation also algebraizes, since
it represents the actual closed cyclic
morphism functor.
Modulo $(F_{\rm cyc},F_j)$, differentiating
$F_\pi=F_{\rm cyc}G_p$ gives
\[
 \left.\det\frac{\partial(F_\pi,F_j)}
                    {\partial(t_1,t_2)}\right|_{A_{\rm cyc}}
      =G_p
         \det\frac{\partial(F_{\rm cyc},F_j)}
                    {\partial(t_1,t_2)} .
\]
The unit identity in
\eqref{eq:paired-all-hom-trace} restricts along
$A_{\rm all}\to A_{\rm cyc}$, so
\[
                         (G_pJ_{\rm cyc})=(2\alpha)
                                             \quad\text{in }A_{\rm cyc}.
\]
After inverting $p$, $G_p$ and $2\alpha$
are units, so this Jacobian identity also
proves that the cyclic generic fibre is etale.
Lemma~\ref{lem:cyclic-removed-divisor-U} identifies
$v_i(G_p)$ with $u_i$. The action on the tangent
line embeds the characteristic-zero quadratic
field generated by $\beta$ by one of its two
embeddings; hence $v_i(\alpha)=v_i(\beta)$.
Taking valuations proves
\eqref{eq:all-prime-trace-U-target}.
At ordinary points the same identity follows
from the explicit complete-sheet calculations,
including Theorem~\ref{thm:inner-ordinary-finite-compensation}
and the already proved outer ordinary case.
The whole-sheet different/contact theorem then
identifies this equality with the vanishing
of the stated finite package. It introduces
no second copy of the internal contacts.
\end{proof}

\begin{corollary}[Integral scaling of a common CM branch]
\label{cor:proportionality-scaling-contact}
Let $j_0$ be a fixed positive special morphism on an actual
cyclic chart. Let $a,b$ be nonzero coprime integers.
Write $f_j,f_{aj},f_{bj}$ for the actual closed
lifting-sheet equations of the fixed morphisms
$j_0,aj_0,bj_0$, in the same regular local cyclic
surface ring. Then
\begin{equation}\label{eq:Bezout-morphism-sheet}
 (f_{aj},f_{bj})=(f_j),\qquad
 f_{aj}=f_jg_a,\quad f_{bj}=f_jg_b,\quad
                         (g_a,g_b)=(1).
\end{equation}
On every normalized formal branch of the complete
$j$-sheet,
\begin{equation}\label{eq:scaling-contact-valuation}
                    v_i(g_a)=v_i(a),\qquad
                    v_i(g_b)=v_i(b).
\end{equation}
Thus the total finite scaling-contact contribution
on a global branch with field $F$ is
\begin{equation}\label{eq:scaling-contact-global}
 \sum_{\mathfrak q}
   \bigl(v_{\mathfrak q}(g_a)+v_{\mathfrak q}(g_b)\bigr)
                         \log\operatorname{Norm}\mathfrak q
       =[F:\mathbb Q]\log|ab|.
\end{equation}
If $K=\mathbb Q(\pi^\vee j)\subset F$, this equals
$[F:K]\log(a^2b^2)$.
\end{corollary}

\begin{proof}
A lift of $j_0$ gives lifts of its integral multiples.
Conversely choose integers $c,d$ with $ca+db=1$.
When $aj_0$ and $bj_0$ lift, the morphism
$c(aj)+d(bj)$ lifts $j_0$; rigidity ensures that
its multiples are the two specified morphisms.
This proves equality of the actual lifting
functors, hence of their defining ideals.
The individual inclusions give the displayed
factorizations. Since $f_j$ is a nonzerodivisor,
cancellation proves $(g_a,g_b)=(1)$.

Use the same ambient dualizing generator for
all three hypersurfaces. Equivalently, in the
product elliptic chart, use the fixed first
equation $F_{\rm cyc}$ for the cyclic surface
and compute their two-equation trace Jacobians.
Modulo $f_j$, differentiation of
$f_{aj}=f_jg_a$ gives
\[
                      J_{aj}|_{A_j}=g_aJ_j.
\]
This formula is unaffected by changing the
chosen equations by units; it is an identity
of their generated fractional ideals.
Theorem~\ref{thm:all-prime-finite-CM-compensation},
applied to $j$ and $aj$, gives on the fixed
normalized $j$-branch
\[
 v_i(J_j)+u_i=v_i(2\beta),\qquad
 v_i(J_{aj})+u_i=v_i(2a\beta),
                        \qquad\beta=\pi^\vee j.
\]
Subtracting proves $v_i(g_a)=v_i(a)$.
The same argument proves the assertion for $b$.
The trace identities include the complete
selected sheets, so no branch conductor or
contact has been omitted in taking their ratio.

For a nonzero rational integer $c$, the product
formula for its finite divisor is
\[
 \sum_{\mathfrak q}v_{\mathfrak q}(c)
                         \log\operatorname{Norm}\mathfrak q
                     =[F:\mathbb Q]\log|c|.
\]
Apply this to $a$ and $b$ and use $[F:K]=[F:\mathbb Q]/2$.
This proves \eqref{eq:scaling-contact-global}.
\end{proof}

\begin{lemma}[The integral rank-one pair functor]
\label{lem:rank-one-pair-moduli}
Let $t>0$ and let $a,b$ be coprime integers. For
$T=t(a,b)^t(a,b)$, the map
\[
             j\longmapsto(aj,bj)
\]
identifies the actual special-morphism functor $\mathcal Z(t)$
with the ordered-pair functor $\mathcal Z(T)$, over arbitrary
bases, with the level data fixed. In particular it preserves
the scheme multiplicities and stabilizers.
\end{lemma}
\begin{proof}
For a pair $(j_1,j_2)$ with Gram matrix $T$, the degree
quadratic form gives
\[
 \deg(bj_1-aj_2)=b^2a^2t+a^2b^2t-2ab(abt)=0.
\]
A homomorphism of elliptic curves of degree zero is zero;
fibrewise this is immediate, and rigidity gives the assertion
also over nonreduced bases. Thus $bj_1=aj_2$.
Choose integers $c,d$ with $ca+db=1$, and put
$j=cj_1+dj_2$. The preceding relation gives $aj=j_1$
and $bj=j_2$, while the degree quadratic form gives
$\deg j=t$. Orthogonality to $\pi$ is preserved by the
integral linear combination. These constructions are inverse
and commute with base change. An automorphism preserves
$j$ if and only if it preserves the ordered pair, since
$j$ is their indicated integral linear combination.
\end{proof}

For rank-one coefficients with both diagonal entries positive
(so that $ab\ne0$), the primitive integral
parameters of Lemma~\ref{lem:standard-rank-one-intrinsic-parameters}
therefore place the scaling comparison in exactly the setting
of Corollary~\ref{cor:proportionality-scaling-contact}.
On a branch of degree $[F:K]$, its finite contribution
$[F:K]\log(a^2b^2)$ cancels the parameter-independent scaling
term in \eqref{eq:self-rankone-parameter-comparison}, with the
factor $1/2$ over all complex embeddings. This cancellation
uses the original Green functions and their actual finite
sheets.

\subsection{The common-component comparison for positive coefficients}
\label{subsec:global-positive-rankone-assembly}

\begin{lemma}[Removing the proportional ordered pairs]
\label{lem:positive-independent-subsum}
For fixed positive indices $m,n$ and parameters $v_1,v_2$,
remove from all two-vector heat sums the ordered pairs
$(x,y)$ with $x,y$ rationally proportional.
The remaining proper-star currents, two-vector heat
integrals, and transgression primitives converge
absolutely in the same compact-current sense as in
Theorem~\ref{thm:compact-proper-green-star}.
Their compact transgression has zero integral.
Writing $S_{\rm ind}(U)$ for the upper star endpoint at heat time
$U\geq1$, one has
\[
 \|S_{\rm ind}(U)\|_{\rm TV}\leq C(1+U)^Ae^{-cU}.
\]
This holds also when the full divisors have common CM points.
\end{lemma}

\begin{proof}
On a compact geometric chart apply
Lemma~\ref{lem:positive-heat-finite-window} in each variable. Products
involving a remainder have exponential decay in heat time. Their
products with a finite singular term are integrable by the displayed
logarithmic and \(r^{-1}\) bounds, or by the stronger local \(L^1\)
estimates in that lemma.

It remains to consider two vectors in the finite windows.
Independent positive vectors cannot have the same zero, because a
point corresponds to one positive line and a norm fixes the rational
proportionality factor. Their finitely many zero sets can therefore
be separated pairwise. On a disk about the zero of one factor, the
other factor and all fixed derivatives decay exponentially; the
first has at most polynomial heat growth. This gives absolute
integrability of the heat form and both primitives. It also gives
the stated total-variation bound for
\(F_x(U)\delta_{D_y}+\Omega_x(U)F_y(U)\): the atomic value is
exponentially small and the continuous term is bounded using the
local logarithmic \(L^1\) estimate. The bounds are uniform because
only finitely many window pairs remain.

There is no issue from different lifts of a quotient CM point.
The estimates concern actual lattice vectors in a fixed chart.
Alternatively, discreteness of \(B(x,y)\) and
\(\cosh d(p_x,p_y)=|B(x,y)|/(2\sqrt{mn})\) give a positive uniform
separation for every independent pair before quotienting.

At a cusp the positive-index estimate is termwise absolute.
The bound for the whole double sum is the product of the two
one-variable bounds; a subsum is bounded by this absolute sum.
The same statement holds for the first derivatives in the
primitives. Integrating \(t\) from \(U\) to infinity preserves a
bound \(C(1+U)^A(1+Y)^Be^{-c'UY^2}\). This is integrable in the
joint spatial and heat variables. Conversion of a fixed number
of spatial derivatives to the local cusp parameter costs only
\(e^{CY}\), which is absorbed by this Gaussian. The primitives
extend flatly across the cusp, so no boundary distribution appears.

The finite-\(T\) Leibniz identity is legitimate for each independent
pair. The absolute bounds justify summing it, passing to the heat
limit, and taking its distributional derivatives. Orbifold charts
or a complex torsion-free cover preserve these estimates with
their finite weights. Absolute convergence also permits grouping by $r=B(x,y)$;
all real nonempty independent terms satisfy $r^2>4mn$.
On the compact curve the exact current pairs to zero against $1$. The interior \(r^{-1}\) terms are
derivatives of globally defined \(L^1\) currents, not discarded
punctured-disk boundary fluxes.
\end{proof}

\begin{lemma}[Globalization of the dependent current identity]
\label{lem:positive-dependent-global-current}
Fix a positive rank-one half-integral matrix
$T=t(a,b)^t(a,b)$, where $t>0$, $(a,b)$ is
primitive, and $ab\ne0$. Choose $a>0$.
Put
\[
 c=a^2v_1,\quad d=b^2v_2,\quad
 s=c+d,\quad \nu=v_1v_2/s,\qquad
 \varphi_0=\operatorname{curv}H.
\]
Let $\mathcal S_T^{\rm hyp}$ denote the sum of the
regularized stars of the ordered pairs $(au,bu)$
over all norm-$t$ vectors $u$, with the original
moduli and orbifold weights and the hyperbolic
tangent metric. The orbit sums $H_T,P_T,Q_T$ of the local
potentials and primitives are global $L^1$ currents on the compact
curve. At heat cutoff $U\geq1$, their tails have $L^1$ bounds
$C/U$, $C/\sqrt U$, and $C/\sqrt U$, respectively. Then
\begin{equation}\label{eq:positive-dependent-global-integral}
 \int_{\overline X}\mathcal S_T^{\rm hyp}
   =\int_{\overline X}\varphi_0 F_{t,s}
      +C^{\rm hyp}_{c,d}\deg_{\mathbb Q} Z(t),
 \qquad
 C^{\rm hyp}_{c,d}
   =-\log\nu-\gamma-\log(4\pi t)-\log(a^2b^2).
\end{equation}
All smooth regular values of the complete Green
functions that come from nonproportional pairs
belong to Lemma~\ref{lem:positive-independent-subsum};
they have not been omitted from the arithmetic product.
\end{lemma}

\begin{proof}
For each positive vector $u$ use $\chi_u=2\pi R(u,z)$.
This function and the closed formulas for its primitives in
Theorem~\ref{thm:self-local-identity} are equivariant under the
arithmetic group, so their orbit sums descend. Their extensions
across the CM support are canonical $L^1$ extensions. On a compact
chart Lemma~\ref{lem:positive-heat-finite-window} leaves finitely
many singular vectors. The local bounds in the proof of
Theorem~\ref{thm:self-local-identity} give the stated tail rates for
each; the remaining vectors have a normally convergent exponential
remainder. This proves the tail bounds on the compact core.

At a cusp, $\chi_u\geq cY^2$ on the positive norm shell.
Derivatives of $\log\chi_u$ have polynomial bounds in the majorant
and height, since the exponent is bounded away from zero there.
The closed formulas
\[
 P_u=-\frac{i}{2\pi}\bar\partial\log\chi_u
       (k_cf_d-f_s),\qquad
 Q_u=-\frac{i}{2\pi}f_s\,\partial\log\chi_u
\]
and $\beta_1(x)\leq e^{-x}/x$ give termwise Gaussian bounds.
Lemma~\ref{lem:positive-heat-cusp-decay} therefore applies to the
orbit sums and all fixed spatial derivatives. Integrating their
heat expressions absolutely from $U$ to infinity yields a bound
$C(1+U)^A(1+Y)^Be^{-c'UY^2}$. The cusp tails are exponentially
small, including after conversion to the local cusp parameter.
Consequently the compact extensions have no cusp atom or flux,
and the global $L^1$ tail bounds follow.

The tangent-regularized atom is intrinsic. For an explicit check,
if $p=\alpha+i\beta$ is the CM point of $u$, the hyperbolic distance
formula gives
\[
 \chi_u(z)=\frac{\pi t}{\beta^2Y^2}
                      |z-p|^2|z-\bar p|^2.
\]
Hence $-dd^c\log(\chi_u/|z-p|^2)
=-\omega_{\rm hyp}/(2\pi)$ and
$\kappa/\ell_p=4\pi t$ for $\ell_p=\beta^{-2}$.
The local atom therefore has the same coefficient on every orbit,
\[
 C^{\rm hyp}_{c,d}
   =-\gamma-\log(4\pi t)-\log\frac{cd}{c+d}
   =-\log\nu-\gamma-\log(4\pi t)-\log(a^2b^2).
\]

We can now sum the finite-cutoff identities absolutely and take
their distributional limit. The strong primitive tails and the
weak concentrated-atom estimates in the proof of
Lemma~\ref{lem:self-concentrated-constant} give
\[
 \mathcal S_T^{\rm hyp}
   =\varphi_0F_{t,s}
      +C^{\rm hyp}_{c,d}\delta_{Z(t)}
      +dd^cH_T+\partial P_T+\bar\partial Q_T .
\]
No cutoff of a local primitive is inserted. Pairing with $1$ kills
all three exact currents on the compact curve and proves
\eqref{eq:positive-dependent-global-integral}. The atomic mass is
the actual orbifold degree of $Z(t)$. Choosing $a>0$ makes $(a,b)$
unique and retains both $u$ and $-u$ in the norm shell; this
introduces no additional sign quotient or factor $2$.
\end{proof}

\begin{lemma}[The finite package descends to the actual CM stack]
\label{lem:CM-stack-principal-finite-package}
On the normalization of the actual rank-one
CM stack $\mathcal Z(t)$, the total finite
normal/contact package for the pair $(aj,bj)$
is the principal divisor of the tangent scalar
$2ab\alpha$, where $\alpha$ is the action of
$\beta=\pi^\vee j$ on the source tangent line.
In particular its degree is
\begin{equation}\label{eq:CM-stack-finite-package-degree}
 \frac{\deg_{\mathbb Q}\mathcal Z(t)}2
                              \log(4Nta^2b^2).
\end{equation}
At each finite prime, the statement uses an etale chart prime
to that prime; it requires no single auxiliary full-level cover
that is etale everywhere.
\end{lemma}

\begin{proof}
The tangent endomorphism is a scalar intrinsically,
so $\alpha$ is a regular function on the CM stack
and $\alpha^2=-Nt$. At a finite prime choose the
actual etale closed-immersion sheets used in
the local calculations. The different,
conductor, and ordered internal contacts are
the vanishing valuations of the trace section
in the finite dualizing line, by
Theorem~\ref{thm:Cartier-sheet-different-contact}.
This description is unchanged by etale
transition maps: relative duality identifies
the dualizing lines and their trace sections.
Equivalently, after strict completion the
transition maps identify these actual
deformation rings, and a change of generator
only multiplies the trace element by a unit.
Thus this finite divisor contribution
descends to the normalization of the
actual CM stack.

Adding the pullback of $U_N$ gives exactly
$\operatorname{div}(2\alpha)$, by
Theorem~\ref{thm:all-prime-finite-CM-compensation}.
The two additional scaling contacts give
$\operatorname{div}(ab)$ by
Corollary~\ref{cor:proportionality-scaling-contact}.
These identities are local at every prime
and therefore identify their full finite
Weil divisors on the actual CM stack.

Apply the principal-divisor product formula
for stack arithmetic degrees. One can verify
it without any everywhere-etale auxiliary
level: on each generic component choose a
number-field object and, after a finite
extension, extend it to all finite places by
properness of the positive CM stack.
Only finitely many places require extension,
since the generic object spreads out.
Normalization gives a finite generically
finite scheme cover
$\operatorname{Spec}\mathcal O_F\to
\widetilde{\mathcal Z(t)}$ of that component.
The map is proper, quasi-finite, and
representable; denote its degree by $d$.
Pullback of a principal divisor and the
projection formula give its degree as
$d^{-1}$ times the usual number-field
principal-divisor degree. Ramification of
this cover is included by that projection
formula and has not been discarded.
In particular
$[F:\mathbb Q]/d$ is the generic stack
degree of the component.
At every embedding,
$|2ab\alpha|^2=4Nta^2b^2$.
The number-field product formula, summed
over these components, proves
\eqref{eq:CM-stack-finite-package-degree}.
\end{proof}

\begin{theorem}[Positive common-component coefficients]
\label{thm:standard-positive-square-product}
For the standard data and degree prescriptions of
Definition~\ref{def:standard-coefficients}, one has
\[
                 \delta^{\rm std}_{m,n}(v_1,v_2)=0
                            \qquad(m,n>0).
\]
The identity includes the positive square-product sector and
is coefficientwise for the standard data. It neither asserts a
modular transformation law for the corrected genus-two family
nor identifies these data with a different arithmetic
Siegel--Weil normalization.
\end{theorem}

\begin{proof}
The nonsquare case was proved already.
For the remaining case expand the actual product
on auxiliary level covers where the two positive
special-morphism moduli spaces split into
closed-immersion Cartier sheets. One can use
different prime-to-$p$ auxiliary levels at the
finite primes: the resulting degree identities
descend by their finite etale degrees. Retain
the ordered morphisms and their stabilizers.

The Gram entries are locally constant on the
pair moduli. Each nonsingular Gram component
has a proper finite intersection, even when
the two full divisors have other common components.
Indeed a common horizontal component would
give two independent special morphisms on a
characteristic-zero elliptic curve; its
special Hom space is one-dimensional, so their
Gram matrix would instead have rank one.
The positive sheets have no vertical component.
Thus the same regular-sequence and scheme-length
argument as in
Lemma~\ref{lem:standard-proper-finite-coefficients}
applies separately on every nonsingular Gram
component. Its finite term is exactly the
prescribed finite genus-two coefficient.
Its archimedean term is the independent
subsum of Lemma~\ref{lem:positive-independent-subsum}.
The sum is absolutely convergent, so these
terms can be removed from both sides.

Consider one of the two remaining rank-one
matrices $T=t(a,b)^t(a,b)$, with $a>0$.
Lemma~\ref{lem:rank-one-pair-moduli} identifies
its actual ordered-pair functor with
$j\mapsto(aj,bj)$ on $\mathcal Z(t)$,
over arbitrary bases and with unchanged
scheme multiplicities and stabilizers.
For its common Cartier sheet the actual local
equations satisfy
$f_{aj}=f_jg_a$, $f_{bj}=f_jg_b$,
$(g_a,g_b)=(1)$.
The arithmetic line-bundle product along that
sheet is therefore its normal restriction
with the contacts from $g_a$ and $g_b$,
together with its Green-curvature integral.
This is the usual Cartier restriction rule:
$\mathcal O(D+A)|_D=
\mathcal O(D)|_D\otimes\mathcal O(A)|_D$.
There is no additional intersection of the
two remainder divisors at this local
rank-one locus, since $(g_a,g_b)=(1)$.
More intrinsically, write $N_t,N_{aj},N_{bj}$
for the Cartier normal lines and $P=\mathcal Z(t)$.
There is an exact sequence on $P$
\[
 0\longrightarrow E^\vee
 \longrightarrow N_{aj}^{\vee}|_P\oplus N_{bj}^{\vee}|_P
 \xrightarrow{\ (g_a,g_b)\ }N_t^\vee
 \longrightarrow0.
\]
Surjectivity is precisely the displayed
comaximality. Taking determinants gives
\[
 E=N_{aj}|_P\otimes N_{bj}|_P\otimes N_t^{-1}
       =N_t\otimes\mathcal O(G_a+G_b)|_P.
\]
These formulas are intrinsic and descend
from the closed-immersion sheets. They
show directly which normal line and
which two contact divisors occur in
the common contribution.
One can check the excess assertion directly:
the Koszul complex of $(fg_a,fg_b)$ has
$H_0=\mathcal O_P$ and $H_1=E^\vee$.
Indeed every relation is a multiple of
$(-g_b,g_a)$, and its Koszul boundaries
are precisely its multiples by $f$.
Thus its excess line is $E$ with the
indicated sign. This also verifies the
Cartier intersection grouping without
combining distinct moduli sheets.

To be explicit about the metric, an individual
singular Green summand has local expansion
$-\log|w|^2+h_c$. The normal frame $1/w$
has squared norm $e^{-h_c}$.
Integral adjunction in the hyperbolic tangent
metric contributes precisely
$\frac12(h_c+\log\ell)$ for this embedding.
The other archimedean contribution is
$\frac12\int\Omega_cf_d$.
These are exactly the two terms in
$\frac12\int\mathcal S_T^{\rm hyp}$.
Every additional smooth summand in a complete
Green function is assigned to the corresponding
nonproportional pair. Hence the local regular
value has entered once, not twice.
Here is a precise way to perform this separation.
At every quotient CM point choose a cutoff
reference Green function $\rho_p$ equal to
$-\log|w|^2$ near $p$, so that its local
regular part is zero. Subtracting the
appropriate sum of these reference functions
from each full Green coefficient gives a
smooth function on the compact curve.
Expand the arithmetic product using
bilinearity for these smooth metric changes.
The regular value of each full coefficient
is the sum of $h_c$ from its singular
proportional summand and the finite values
of all the other vector lifts. The latter
are exactly the atomic terms
$f_x(p_y)\delta_{p_y}$ of the independent
star subsum. The curvature-times-Green
terms have the same ordered-tuple expansion.
Absolute convergence from the two preceding
lemmas permits its limit and regrouping.
Thus the normal metric takes $h_c$ once,
while all other smooth regular values
are retained in the independent contribution.

Integral adjunction and the
dualizing--Hodge isometry give the Hodge
restriction, the branch different/conductor,
the internal contacts, and the $U_N$ term,
with the remaining metric constant
$\frac12(\gamma+\log\pi-\log N)
\deg_{\mathbb Q}\mathcal Z(t)$.
Lemma~\ref{lem:CM-stack-principal-finite-package}
identifies the total finite contribution,
including scaling contacts, with
\[
 \frac{\deg_{\mathbb Q}\mathcal Z(t)}2
                                      \log(4Nta^2b^2).
\]
The trace element counts every internal contact of the complete
selected $j$-sheet once. If a global branch is normalized first,
its internal formal contacts are already included in the enlarged
conductor and must not be added again.

By \eqref{eq:positive-dependent-global-integral},
the archimedean part not contained in the
Hodge restriction is its curvature integral
and the constant
\[
 \frac{\deg_{\mathbb Q}\mathcal Z(t)}2
 \{C^{\rm hyp}_{c,d}+\gamma+\log\pi-\log N\}
 =-\frac{\deg_{\mathbb Q}\mathcal Z(t)}2
                  \{\log\nu+\log(4Nta^2b^2)\}.
\]
All finite parameter-independent terms thus
cancel exactly. The two remaining Hodge
terms combine, by the arithmetic
line-bundle pairing rule and
$\operatorname{curv}H=\varphi_0$, into
$\langle\widehat Z_{\rm std}(t,s),H\rangle_{\rm GS}$.
The remaining parameter term is
$-\frac12\log\nu\,\deg_{\mathbb Q} Z(t)$.
Summing over the global branches, with the
unchanged moduli weights, therefore gives
exactly $\mathcal A^{\rm std}(T,y)$.
The same argument applies to the other
endpoint, including its sign of $b$.
Together with the nonsingular components
this is the full original product and
proves the asserted defect vanishing.
\end{proof}

\section{Comparison with arithmetic Siegel--Weil formulas}
\label{sec:remaining}
\label{sec:height-reduction}

We fix a genus-two Eisenstein series and compare its derivative
with the standard arithmetic coefficients. The nonsingular terms
are followed by the rank-one terms and the rank-zero middle orbit.
An explicit modular completion and Green-energy bounds justify
the full Fourier expansion and double Petersson pairing. We also
specify the data needed to express the resulting height through
an adelic doubling formula. That analytic identity is proved in
Appendix~\ref{sec:analytic-doubling}.

\subsection{The nonsingular arithmetic Siegel--Weil coefficients}
\label{subsec:all-level-nonsingular-ASW}

Fix the standard additive character used in the preceding
Weil representation. Let $\mathcal E_N(z,s)$ be the classical
incoherent genus-two Eisenstein series whose finite section
is obtained by the Rallis map from $1_{\widehat L_N^{\,2}}$,
and whose archimedean standard section has scalar
metaplectic $K$-type $\det^{3/2}$ and value $1$ at the
identity. Use the classical factor $\det(y)^{-3/4}$.
Write
\[
 \mathcal E_N(z,s)=
  \sum_T C_{N,T}(y,s)q^T,\qquad
 q^T=e^{2\pi i\operatorname{tr}(Tz)}.
\]
This is the normalization in
\cite[Definition~2.13 and equations~(90)--(94)]{SSY2022}.
It is also Zhu's classical standard-section normalization
after the isometry
\[
 \begin{pmatrix}b&a/N\\c&-b\end{pmatrix}
       \longmapsto
 \begin{pmatrix}Nb&a\\c&-b\end{pmatrix}
                 \in\Delta(N).
\]
Both quadratic forms of this map are $-Nb^2-ac$;
in particular no rescaling of the additive character
or of the index $T$ is involved in this identification.

\begin{proposition}[All nonsingular standard coefficients]
\label{prop:all-level-nonsingular-ASW}
For every $N\geq1$, every nonsingular half-integral
matrix $T$, and every $y>0$, the standard coefficients
of Definition~\ref{def:standard-coefficients} satisfy
\begin{equation}\label{eq:all-level-nonsingular-ASW}
 \mathcal A^{\rm std}(T,y)
    =\frac{\psi(N)}{24}\,
                  \left.\frac{\partial}{\partial s}
                      C_{N,T}(y,s)\right|_{s=0}.
\end{equation}
For the non-positive-definite matrices the Green form
here is the normalized $g_2^o$, constructed from the
intrinsic Garcia--Sankaran Schwartz form $\nu^o$.
It is not the literal heat integral of $\nu$ before
removing its quadratic Gaussian.
\end{proposition}

\begin{proof}
Suppose first that $T>0$. There is no real tuple
with this Gram matrix in signature $(1,2)$.
Consequently the standard coefficient is exactly
the finite degree of the actual ordered
special-morphism stack. The actual moduli
problem, scheme lengths, residue degrees, and
stabilizer factors are the ones in
\cite[Definition~4.3.5 and Theorem~4.4.1]{Zhu23}.
The stated characteristic function is
$T$-admissible. That all-level theorem gives
\eqref{eq:all-level-nonsingular-ASW} directly.
This part does not use any Green-function
or boundary normalization.

For the remaining nonsingular $T$ we use
the actual Schwartz-form calculation in
\cite[Section~4.2]{SSY2022}. We describe why
its archimedean argument applies to arbitrary
$N$ and fixes the normalization required here.
The normalized form is $\nu^o$, with
\[
 \nu(x)=e^{-2\pi(Q(x_1)+Q(x_2))}\nu^o(x),
 \qquad
 g_2^o(x)=\int_1^\infty
                       \nu^o(\sqrt{\lambda}\,x)\,
                                      \frac{d\lambda}{\lambda},
\]
as in equations~(59)--(63) of that paper.
The identification with our heat expression
was proved above from the intrinsic
Garcia--Sankaran forms.

More explicitly the real Schwartz form used
in its Siegel--Weil argument is fixed by
\cite[equations~(178)--(179)]{SSY2022}.
For
$e=\operatorname{diag}(1,-1)$ and
$f_1=\left(\begin{smallmatrix}0&0\\1&0\end{smallmatrix}\right)$,
these equations give
\[
 R(ae+xf_1,i)=\frac N2(x^2+4a^2).
\]
This is exactly the full majorant
$R=\frac12B(x,h_z)^2-2Q(x)$ used here.
In particular we fix $R$ by this
intrinsic Schwartz construction and this
explicit calculation, rather than
transplanting a potentially ambiguous
projection notation from another displayed formula.

Let $\phi_\infty$ be the scalar coefficient
of that real Schwartz form at the base point.
The only nontrivial convergence input in
Lemma~4.7 of \cite{SSY2022} is the vanishing
of its mixed-model restriction at $(A,B)=(0,0)$.
For the present lattice it is
\[
 e^{-2\pi N(a_1^2+a_2^2)}
 \int_{\mathbb R^2}
  \left\{-\frac1\pi+N(x_1^2+x_2^2)\right\}
             e^{-\pi N(x_1^2+x_2^2)}\,dx_1dx_2=0.
\]
Indeed the Gaussian integral is $1/N$,
and its $N(x_1^2+x_2^2)$ moment is
$1/(\pi N)$, canceling the constant term.
This calculation holds for every $N>0$.
It uses no squarefreeness or parity.
The finite Schwartz function in that
argument is arbitrary; it may therefore
be $1_{\widehat L_N^{\,2}}$ at every level.
The Kudla--Rallis mixed-model criterion
and the convergent Siegel--Weil identity
of Lemma~4.7 and Corollary~4.8 consequently
apply with these finite data.
Here arbitrariness of the finite Schwartz
function is used for convergence; in the
coefficient identity we choose precisely
$1_{\widehat L_N^{\,2}}$.
There is no additional genus sum at higher
level. The Clifford-group identification
of \cite[Lemma~4.2.16]{Zhu23} gives
$\operatorname{GSpin}(V)\simeq\operatorname{GL}_2$
and the finite level $K_N=\Gamma_0(N)(\widehat{\mathbb Z})$,
which stabilizes that characteristic function.
Since $\det K_N=\widehat{\mathbb Z}^{\,\times}$,
strong approximation for $\operatorname{SL}_2$
and the ideal class number one of $\mathbb Q$
give one finite adelic component.
The normalized theta average is therefore
the normalized integral on the stated
$\Gamma_0(N)$ modular curve. Passing from
this level to a larger orthogonal stabilizer,
if desired, changes its integral and volume
by the same finite index.

The transgression calculation of
\cite[Theorem~4.10, equation~(182)]{SSY2022},
using the intrinsic Garcia--Sankaran forms,
then reads in this same normalization
\[
 \int_{\mathcal X_0(N)(\mathbb C)}g_2^o(T,y)
    =\operatorname{vol}(\mathcal X_0(N)(\mathbb C),d\Omega)
        \left.\partial_s C_{N,T}(y,s)\right|_{s=0},
 \qquad d\Omega=\frac{dx\,dY}{2\pi Y^2}.
\]
No squarefree local-density computation
enters this archimedean identity. For
the actual modular stack the volume is
\[
 \operatorname{vol}(\mathcal X_0(N)(\mathbb C),d\Omega)
   =\frac12[\operatorname{SL}_2(\mathbb Z):\Gamma_0(N)]
       \int_{\operatorname{SL}_2(\mathbb Z)\backslash\mathbb H}
                          \frac{dx\,dY}{2\pi Y^2}
   =\frac{\psi(N)}{12}.
\]
The first factor $1/2$ accounts for the
ineffective central automorphisms in
this stack integration convention.
The arithmetic degree of the purely
archimedean class has the further factor
$1/2$. It is therefore $\psi(N)/24$
times the displayed Fourier coefficient
derivative, proving
\eqref{eq:all-level-nonsingular-ASW}.
This is also the unrestricted-level
extension expressly invoked at the
end of \cite[Section~8.3]{Zhu23}.
\end{proof}

This proposition covers all nonsingular $T$.
Section~\ref{subsec:positive-rankone-comparison} treats the
positive rank-one terms using the Du--Yang Green construction.
We first express the singular comparison in terms of three
functions. These are evaluated below, together with the complete
modular family required by the pairing criterion in
Section~\ref{subsec:RS-positive-support}.

\subsection{An explicit singular-coefficient comparison}
\label{subsec:singular-ASW-comparison}

Put $\mathfrak c_N=\psi(N)/24$ and
$\mathcal A^{\rm SW}(T,y)=\mathfrak c_N\partial_z C_{N,T}(y,z)|_{z=0}$,
using the fixed Eisenstein series of
Section~\ref{subsec:all-level-nonsingular-ASW}.
For the genus-one Eisenstein series with the same finite
characteristic function and the standard weight-$3/2$ section, set
\[
 e_{N,t}(u,z)=e^{2\pi tu}E^{(1)}_t(iu,z;\Delta(N)),\qquad t\in\mathbb Z.
\]
Here $E^{(1)}_t$ denotes the full Fourier term, so the exponential
removes its $q^t$ factor. The spectral parameter is $z$ and the
positive imaginary parameter is $u$. Let
\[
 \kappa_{N,t}=2+\log N+4\frac{\Lambda'(2)}{\Lambda(2)}
                  -\sum_{p\mid N}\beta'_{p,t}(0),
 \qquad \Lambda(z)=\pi^{-z/2}\Gamma(z/2)\zeta(z),
\]
where $\beta_{p,t}$ is the local factor of
\cite[equation~(28)]{Zhu}; the dependence on the nonzero index
$t$ is included in this notation. Define two comparison functions
\begin{align}
 h_t(u)&=\langle\widehat Z_{\rm std}(t,u),H\rangle_{\rm GS}
       -\mathfrak c_N\{2\partial_z e_{N,t}(u,z)|_{z=1/2}
                              +\kappa_{N,t}e_{N,t}(u,1/2)\},
                  \label{eq:rankone-Hodge-comparison-function}\\
 d_t(u)&=\deg Z^*_{\rm std}(t,u)
                              +2\mathfrak c_N e_{N,t}(u,1/2).
                  \label{eq:rankone-degree-comparison-function}
\end{align}
These are explicit genus-one comparisons; their vanishing is
not assumed.
At rank zero retain the comparison function
\begin{align}\label{eq:rankzero-base-comparison-function}
 b_N(y)&=\mathcal A^{\rm std}(0_2,y)-\mathcal A^{\rm SW}(0_2,y)\notag\\
 &=\widehat{\deg}(H^2)-\frac12\log\det(y)\deg H
       -\mathfrak c_N\left.\partial_z C_{N,0_2}(y,z)\right|_{z=0}.
\end{align}
This definition makes no rank-zero Siegel--Weil identification.

\begin{proposition}[The singular comparison functions]
\label{prop:singular-ASW-comparison}
For a nonzero rank-one $T$, let
$t$, $s=t^{-1}\operatorname{tr}(Ty)$ and $\nu=\det(y)/s$
be its intrinsic parameters. Then
\begin{equation}\label{eq:rankone-ASW-comparison}
 \mathcal A^{\rm corr}(T,y)-\mathcal A^{\rm SW}(T,y)
                =h_t(s)-\frac12\log\nu\,d_t(s).
\end{equation}
At rank zero the comparison is
\begin{equation}\label{eq:rankzero-ASW-comparison}
 \mathcal A^{\rm corr}(0_2,y)-\mathcal A^{\rm SW}(0_2,y)
                    =b_N(y)-\frac{d_0}{2}\mathcal Z_y(1/2).
\end{equation}
Together with Proposition~\ref{prop:all-level-nonsingular-ASW},
these formulas express the comparison with the fixed Eisenstein
derivative in terms of $h_t,d_t$, and the retained function $b_N$.
\end{proposition}
\begin{proof}
The analytic singular-coefficient formula
\cite[Corollary~4.2.5]{Zhu}, with its $q^T$ factor removed, is
\[
 \mathcal A^{\rm SW}(T,y)
 =\mathfrak c_N\{2\partial_z e_{N,t}(s,z)|_{z=1/2}
                   +(\log\nu+\kappa_{N,t})e_{N,t}(s,1/2)\}.
\]
One first applies that formula at
$T=\operatorname{diag}(0,t)$ and
$y=\operatorname{diag}(\nu,s)$; integral change of basis and
upper triangular real shear give the displayed intrinsic
parameters for general rank-one $T$.
Subtract this expression from the rank-one degree prescription
in Definition~\ref{def:standard-coefficients}.
The coefficient of $\log\nu$ is
$-\tfrac12\deg Z^*_{\rm std}(t,s)-\mathfrak c_Ne_{N,t}(s,1/2)$,
which proves \eqref{eq:rankone-ASW-comparison}, including its sign.

At rank zero the definition of $b_N$, together with
$\mathcal A^{\rm corr}(0_2,y)=\mathcal A^{\rm std}(0_2,y)
-d_0\mathcal Z_y(1/2)/2$, proves
\eqref{eq:rankzero-ASW-comparison}. In particular, this step
does not set $b_N$ equal to zero.
\end{proof}

\begin{corollary}[The double-zero comparison]
\label{cor:doublezero-SW-comparison}
For the standard genus-one arithmetic coefficients, the
double-zero pullback comparison with the fixed Eisenstein
derivative is
\begin{align}\label{eq:SW-doublezero-obstruction}
 &\langle\widehat Z_{\rm std}(0,v_1),
               \widehat Z_{\rm std}(0,v_2)\rangle_{\rm GS}
 -\sum_{r\in\mathbb Z}\mathcal A^{\rm SW}
      \left(\begin{pmatrix}0&r/2\\r/2&0\end{pmatrix},
                                \operatorname{diag}(v_1,v_2)\right)
 \notag\\
 &\hspace{18mm}=b_N(\operatorname{diag}(v_1,v_2))-\frac{d_0}{2}
                  \mathcal Z_{\operatorname{diag}(v_1,v_2)}(1/2).
\end{align}
\end{corollary}
\begin{proof}
The term $r=0$ contributes the difference $b_N(y)$ between
the standard and Siegel--Weil degrees. Every $r\ne0$ gives
a nonsingular matrix, where
Proposition~\ref{prop:all-level-nonsingular-ASW} applies.
Replacing those nonzero terms in the convergent coefficient
sum and using Theorem~\ref{thm:double-zero-absolute-defect}
gives the result.
\end{proof}

\begin{remark}[Scope of the rank-zero comparison]
The standard-product defect and the Eisenstein discrepancy
are different quantities: the former is the computed Epstein
term, whereas the latter also contains $b_N(y)$. Thus its
nonvanishing cannot be inferred from that Epstein term alone.
To evaluate $b_N$, the rank-zero Fourier expansion must include
the middle Weyl orbit, with its actual stabilizer, local sections,
and central factors. This orbit appears explicitly in
\cite[equations~(239), (253), (255)--(257)]{SSY2022}.
In particular, a vanishing argument using a factor at a prime
dividing $N$ cannot be extended to $N=1$ without a separate
calculation. The rank-zero formula in \cite[Section~9.3]{Zhu}
is therefore not used to identify the fixed data here.
Section~\ref{subsec:rankzero-middle-orbit} computes this orbit
and evaluates $b_N$ directly. The earlier paired argument uses
the entire modular residual and its positive coefficients, so
it is independent of that rank-zero evaluation.
\end{remark}

\subsection{The positive rank-one comparison}
\label{subsec:positive-rankone-comparison}

\begin{proposition}[Positive-index Kronecker-limit comparison at every level]
\label{prop:positive-DY-KLF}
Use the actual positive special-morphism divisor and its
Du--Yang Kudla Green function with
\[
 R(x,z)=\tfrac12 B(x,h_z)^2-2Q(x),\qquad
 F_{t,u}(z)=\sum_{Q(x)=t}\beta_1(2\pi uR(x,z)),\qquad t>0.
\]
Retain both vectors \(x\) and \(-x\), with the original moduli
weights. Integrate on the actual modular stack, so integration
of an invariant form is one half its effective
\(\Gamma_0(N)\backslash\mathbb H\) integral. Let
\[
 E_N^{(1)}(\tau,z)=\sum_t e_{N,t}(\operatorname{Im}\tau,z)q^t
\]
be the classical genus-one Eisenstein series obtained from
\(1_{\widehat L_N}\), with standard archimedean section of
weight \(3/2\) and value \(1\) at the identity, and with the
full rational parabolic in its adelic summation.
Write
\[
 \mathfrak c_N=\frac{\psi(N)}{24},\qquad
 H=-\widehat\omega_N-W_N^*\widehat\omega_N.
\]
Then, for every positive integer \(t\) and \(u>0\),
\begin{align}
 \deg_{\mathbb Q}\mathcal Z(t)
   &=-2\mathfrak c_N e_{N,t}(u,\tfrac12),
       \label{eq:positive-DY-degree}\\
 \langle\widehat Z_{\rm std}(t,u),H\rangle_{\rm GS}
   &=\mathfrak c_N\left\{
      2\left.\partial_z e_{N,t}(u,z)\right|_{z=1/2}
                  +\kappa_{N,t}e_{N,t}(u,\tfrac12)\right\},
       \label{eq:positive-DY-Hodge}
\end{align}
where
\[
 \kappa_{N,t}
   =2+\log N+4\frac{\Lambda'(2)}{\Lambda(2)}
                    -\sum_{p\mid N}\beta'_{p,t}(0).
\]
Here \(\beta_{p,t}\) has exactly the local-density
normalization of Zhu's equation (28).
Consequently the comparison functions \(h_t\) and \(d_t\)
in \eqref{eq:rankone-Hodge-comparison-function} and
\eqref{eq:rankone-degree-comparison-function} vanish for \(t>0\).
\end{proposition}

\begin{proof}
All the analytic input below is in Part I of
\cite{DY1}, and is valid for arbitrary \(N\): the authors
expressly restrict squarefreeness to the arithmetic part
at the end of their introduction.

We first record both normalization changes rather than
combining them into an unexplained factor. Let \(E_L^{\rm b}(\tau,s)\)
be Du--Yang's unnormalized vector-valued Eisenstein series,
whose seed is \(u^{(s-1)/2}e_0\), and put
\[
 C_N(s)=-\frac{s}{4\pi}\Lambda(2s)
              \prod_{p\mid N}(1-p^{-2s})N^{(3s+1)/2}.
\]
Their normalized series is
\(\mathcal E_L^{\rm DY}(\tau,s)=C_N(s)E_L^{\rm b}(\tau,s)\).
For the scalar zero component one has
\begin{equation}\label{eq:DY-central-spectral-bridge}
 (E_L^{\rm b}(\tau,s))_0
             =2E_N^{(1)}(\tau,s-\tfrac12).
\end{equation}
The parameter shift follows by comparing the classical
seeds \(u^{(s-1)/2}\) and \(u^{z/2-1/4}\).
For the factor \(2\), the primitive-row sum in the
Poisson-summation proof in Du--Yang Section 2 retains both
\((c,d)\) and \((-c,-d)\). The adelic parabolic contains
the full center and identifies these rows. The finite
section and its value at the identity are otherwise
unchanged. Thus this factor is distinct from the
ineffective-central-automorphism factor in stack integration.
In particular
\begin{equation}\label{eq:DY-normalized-scalar-bridge}
 (\mathcal E_L^{\rm DY}(\tau,s))_0
           =2C_N(s)E_N^{(1)}(\tau,s-\tfrac12).
\end{equation}

Let \(\Theta_{\rm KM}\) be the standard full-majorant
Kudla--Millson form and let \(I_{\rm eff}\) denote its
integral on the effective modular curve. Du--Yang
Theorem 1.6 says
\[
 I_{\rm eff}(\tau,1)=\frac2{\varphi(N)}
                    \mathcal E_L^{\rm DY}(\tau,1),\qquad
 I_{\rm eff}(\tau,\log\|\Delta_N\|)
                    =-12(\mathcal E_L^{\rm DY})'(\tau,1).
\]
The positive coefficient of \(\Theta_{\rm KM}\) is the
compact Green curvature of \(F_{t,u}\). Its stack
integral is the degree of the actual divisor, by the
compact Green equation. Taking one half of the first
effective integral and using
\eqref{eq:DY-normalized-scalar-bridge} gives
\[
 \deg_{\mathbb Q}\mathcal Z(t)
       =\frac{2C_N(1)}{\varphi(N)}
                         e_{N,t}(u,\tfrac12).
\]
Since \(C_N(1)=-\varphi(N)\psi(N)/24\), this proves
\eqref{eq:positive-DY-degree}.

For the height calculation use the rational section
\(\Delta_N\) of weight \(12\varphi(N)\) with the metric
\[
 \|f(z)\|^2=(c_H\operatorname{Im}z)^k|f(z)|^2,\qquad
 c_H=2\sqrt\pi e^{-\gamma/2}.
\]
This is both Du--Yang's equation (1.9) and the
Petersson normalization used here. The function
\(\Delta_N\) and its vertical divisor are the ones in
Zhu Section 6.4. Thus
\[
 \widehat\Delta_N
  =(\psi(N)\varphi(N)P_\infty,-\log\|\Delta_N\|^2)
  =12\varphi(N)\widehat\omega_N
                         -\sum_{p\mid N}(f_p(N),0).
\]
The positive special divisor misses every arithmetic
cusp. Its Green function and the required height
derivatives decay at every geometric cusp; meanwhile
\(\log\|\Delta_N\|=O(Y+\log Y)\). The cutoff
integration-by-parts boundary therefore tends to zero.
The proper arithmetic star product gives
\[
 \langle\widehat Z_{\rm std}(t,u),\widehat\Delta_N\rangle_{\rm GS}
    =-\int_{\mathcal X_0(N)(\mathbb C)}
                       \log\|\Delta_N\|\,\Omega_{t,u}
    =-\frac12[I_{\rm eff}(\tau,\log\|\Delta_N\|)]_t/q^t .
\]
In the first equality the factor \(1/2\) in arithmetic
degree cancels the \(2\) in \(-\log\|\Delta_N\|^2\).
The second \(1/2\) is stack integration. Hence
\begin{equation}\label{eq:positive-Delta-height-DY}
 \langle\widehat Z_{\rm std}(t,u),\widehat\Delta_N\rangle_{\rm GS}
     =12\left.\frac{d}{ds}
             \{C_N(s)e_{N,t}(u,s-\tfrac12)\}\right|_{s=1}.
\end{equation}
There is no extra factor in the derivative of the
spectral parameter.

The finite calculation uses only
positive actual moduli and does not transfer any
nonpositive Green prescription. In Zhu Corollary 8.2.3
and Proposition 8.2.4 the metrics on the vertical classes
are zero. The positive divisor is the same actual
special-morphism divisor; thus those finite
intersection identities apply unchanged. In the
notation \(n_p=v_p(N)\), let \(g_{p,t}\) be Zhu's
equation (25) and set
\[
 \rho_{p,t}=\frac{p^{-1}g_{p,t}(0)}
                         {1-p^{-1}g_{p,t}(0)},\qquad
 V_{N,t}=\sum_{p\mid N}
       \left(\frac{n_p-1-n_pp}{p-1}-\rho_{p,t}\right)\log p .
\]
These denominators are nonzero. Indeed,
\(D_p=\operatorname{Den}(H_4^+,\langle t,N\rangle)\)
is finite and positive: embed the two rank-one
forms in the two separate hyperbolic planes.
The unit coordinates give a smooth integral
representation. Zhu's local difference identity
in the proof of (27), evaluated at \(k=0\),
identifies
\(D_p(1-p^{-1}g_{p,t}(0))\), up to a nonzero
normalizing factor, with the density for
\(\Delta(N)\) to represent \(t\).
That density is strictly positive: in our
coordinates \(Q=-ac-Nb^2\), the representation
\((a,b,c)=(1,0,-t)\) has
\(\partial Q/\partial c=-1\).
This works also at \(p=2\). Thus neither the
definition of \(\rho_{p,t}\) nor the displayed
derivative of \(\beta_{p,t}\) requires assigning
a value to an indeterminate quotient here.
Equivalently, taking the positive coefficient in his
finite computation preceding equation (59) gives
\[
 \frac1{12}\sum_{p\mid N}
       \langle\widehat Z_{\rm std}(t,u),(f_p(N),0)\rangle_{\rm GS}
          =C_N(1)V_{N,t}e_{N,t}(u,\tfrac12).
\]
The appearances of lower-level degree series in that
calculation are justified by
\eqref{eq:positive-DY-degree} at those levels.
Combining this with \eqref{eq:positive-Delta-height-DY}
proves
\[
 \varphi(N)\langle\widehat Z_{\rm std}(t,u),
                                     \widehat\omega_N\rangle_{\rm GS}
 =C_N(1)\left[
    \left.\partial_z e_{N,t}(u,z)\right|_{z=1/2}
       +\left\{\frac{C_N'(1)}{C_N(1)}+V_{N,t}\right\}
                              e_{N,t}(u,\tfrac12)\right].
\]
The positive actual moduli and their full-majorant
Green function are \(W_N\)-invariant, so pairing with
\(H\) is minus twice pairing with \(\widehat\omega_N\).
Finally Zhu's elementary equations (60)--(61), or direct
differentiation, give
\[
 \frac{C_N'(1)}{C_N(1)}
  =1+2\frac{\Lambda'(2)}{\Lambda(2)}
       +2\sum_{p\mid N}\frac{\log p}{p^2-1}
       +\frac32\log N,\qquad
 \beta'_{p,t}(0)
       =\left(\frac2{p+1}+2\rho_{p,t}\right)\log p .
\]
Using
\[
 \frac{n-1-np}{p-1}+\frac2{p^2-1}
                         =-n-\frac1{p+1}
\]
shows that \(2(C_N'(1)/C_N(1)+V_{N,t})=\kappa_{N,t}\).
The value of \(C_N(1)/\varphi(N)\) now proves
\eqref{eq:positive-DY-Hodge}.
\end{proof}

\begin{corollary}[Positive rank-one Siegel--Weil coefficients]
\label{cor:positive-rankone-ASW}
For every \(N\), every nonzero positive semidefinite
half-integral rank-one matrix \(T\), and every \(y>0\),
\[
 \mathcal A^{\rm std}(T,y)
       =\mathfrak c_N
                     \left.\partial_z C_{N,T}(y,z)\right|_{z=0}.
\]
Consequently the standard arithmetic product kernel
and the fixed Eisenstein-derivative pullback have the
same Fourier coefficient at every pair of strictly
positive diagonal indices.
\end{corollary}
\begin{proof}
The intrinsic rank-one parameter \(t\) is positive.
Insert \(h_t=d_t=0\) in
\eqref{eq:rankone-ASW-comparison}; the rank-zero
correction does not affect \(T\ne0\).
For the last assertion, use the proved standard
positive-product identity and the nonsingular
Siegel--Weil comparison. A symmetric matrix with both
diagonal entries positive is either nonsingular or
nonzero positive semidefinite of rank one, so these
cases exhaust its off-diagonal sum.
\end{proof}

\paragraph{Strict scope of the result.}
This argument fixes the positive coefficients for the
canonical Du--Yang full-majorant Green family, independently
of any printed negative-index boundary correction. It
does not identify the negative or zero coefficients of
Zhu's separately modular arithmetic family with these
standard coefficients. Nor does it silently resolve the
ambiguous displayed \(R\), index, and \(u\)-dependence in
Zhu Section 7.1.2. To identify the positive coefficients
of that modular family, one must use its actual
Du--Yang Green construction: equality follows if its
positive coefficient has the standard compact
Kudla--Millson curvature and cusp limit zero. Indeed,
subtracting the two Green functions then gives a
globally harmonic function on the compact curve;
the removable cusp limit is zero, so the function
is identically zero. This criterion needs the actual
standard curvature, not merely a formally similar
Green equation with an unverified normalization.

\subsection{Positive Fourier support and the double cuspidal pairing}
\label{subsec:RS-positive-support}

\begin{theorem}[Double Rankin--Selberg support criterion]
\label{thm:RS-positive-positive-cusp-criterion}
For $i=1,2$ let $\Gamma_i$ be a finite-index modular group with a
unitary scalar multiplier of weight $k_i>0$, whose values on cusp
translations have finite order. Let $f\in S_{k_1}(\Gamma_1)$
and $h\in S_{k_2}(\Gamma_2)$. Suppose that $K(z,w)$ is a smooth kernel
with weight $k_1$ in $z$ and conjugate weight $k_2$ in $w$, with the
same respective multipliers as $f$ and $\overline h$.

Assume joint moderate growth: in every pair of cusp coordinates,
and also when one variable ranges over a compact core, its invariant
norm is bounded by a fixed polynomial in the two cusp heights.
At a chosen cusp in each variable use width-one scaling coordinates.
Suppose that the Fourier coefficients of $K$ with frequencies
$(m,-n)$ vanish identically for all $m,n>0$. Then
\[
 \int_{(\Gamma_1\backslash\mathbb H)
                  \times(\Gamma_2\backslash\mathbb H)}
 K(z,w)\overline{f(z)}h(w)
       (\operatorname{Im}z)^{k_1}(\operatorname{Im}w)^{k_2}
                                        d\mu(z)d\mu(w)=0.
\]
The integral is absolutely convergent. The result includes
$k_1=k_2=3/2$ and requires no holomorphy of $K$.
\end{theorem}

\begin{proof}
First suppose the chosen cusp coordinates have width one without a
change of variable. The general statement follows by applying the
two scaling matrices and the corresponding weight slash operators.
Pass to a finite-index subgroup if necessary to make the multipliers
trivial on the cusp translations; this multiplies the relevant
integrals by finite indices and does not affect the assertion.

Put $y=\operatorname{Im}z$, $v=\operatorname{Im}w$ and
\[
 A(z,w)=K(z,w)\overline{f(z)}h(w)y^{k_1}v^{k_2}.
\]
The weight factors cancel, so $A$ is invariant in both variables.
The invariant norms of $f$ and $h$ decrease exponentially at every
cusp. Consequently the assumed polynomial bound on $K$ makes $A$
rapidly decreasing jointly in the two quotient heights, including
the regimes where only one variable enters a cusp. In particular it
is absolutely integrable and bounded on the product of quotients.

Let $E_i(z,s_i)$ be the weight-zero Eisenstein series of the selected
cusp, normalized by a width-one scaling coordinate. For
$\operatorname{Re}s_1,\operatorname{Re}s_2>1$ consider
\[
 J(s_1,s_2)=\int A(z,w)E_1(z,s_1)E_2(w,s_2)
                                             d\mu(z)d\mu(w).
\]
On this region $|E_i(z,s_i)|\leq E_i(z,\operatorname{Re}s_i)$.
Polynomial Eisenstein growth and the rapid decrease of $A$ make the
integral of the absolute value finite. Tonelli's theorem therefore
justifies unfolding both Poincare sums. The result is
\[
 J(s_1,s_2)=\int_0^\infty\int_0^\infty
    y^{s_1-2}v^{s_2-2}
       \int_0^1\int_0^1 A(x+iy,u+iv)\,dx\,du\,dy\,dv.
\]
This is an absolutely convergent integral; no assertion of absolute
integrability of every individual Fourier term was used.

For fixed positive $y,v$, the Fourier expansions of the two cusp
forms converge absolutely and uniformly in $x,u$. Their frequencies
are strictly positive. The horizontal integral consequently pairs
only the $(m,-n)$ Fourier coefficients of $K$ with $m,n>0$:
\[
 \int_0^1\int_0^1
 K(x+iy,u+iv)\overline{f(x+iy)}h(u+iv)\,dx\,du=0.
\]
One can justify this equality using uniform convergence of the two
cusp expansions and the boundedness of the smooth kernel on this
fixed compact torus. The hypothesis makes each selected coefficient
zero. It follows that $J(s_1,s_2)=0$ throughout its initial domain.

Each $E_i$ continues meromorphically and has a simple pole at $s_i=1$
with a nonzero constant residue $r_i$. Near $1$, the functions
$(s_i-1)E_i(z,s_i)$ satisfy uniform polynomial cusp bounds. Multiplying
by the rapidly decreasing $A$ gives a common integrable majorant.
Hence $J$ continues near $(1,1)$, with at most a simple pole in each
variable, and dominated convergence gives
\[
 \lim_{s_1,s_2\to1}(s_1-1)(s_2-1)J(s_1,s_2)
                =r_1r_2\int A(z,w)d\mu(z)d\mu(w).
\]
Analytic continuation of the zero identity makes the left side zero.
As $r_1r_2\ne0$, the desired pairing vanishes.

For effective quotient measure and width-one normalized Eisenstein
series, $r_i$ is the reciprocal of the hyperbolic volume. An orbifold
or metaplectic central-volume convention changes the unfolding
measure compatibly; only nonvanishing of the two residues is used.
\end{proof}

\paragraph{Spectral derivatives and negative frequencies.}
The two Eisenstein parameters in this proof are auxiliary
Rankin--Selberg parameters. A derivative already present in the
fixed kernel $K$ introduces logarithmic cusp factors, which are still
bounded by a polynomial. It does not introduce a higher pole in
either auxiliary parameter. No holomorphic-projection formula of
weight $3/2$ is required. Negative or zero Fourier coefficients may
have substantial nonholomorphic imaginary-part dependence; the
proof never estimates them separately after unfolding.

\begin{corollary}[From positive coefficient comparison to paired doubling]
\label{cor:RS-positive-comparison-doubling}
Suppose two complete kernels with the same separate automorphy laws
are smooth and jointly of moderate growth. If their Fourier
coefficients agree for every pair of strictly positive indices, then
their double Petersson pairings with holomorphic cusp forms agree.
\end{corollary}
\begin{proof}
Apply the theorem to their difference.
\end{proof}

\paragraph{Application and the necessary distinction.}
The theorem applies to a fully modular residual. It does not apply
to the isolated Epstein zero mode, an isolated boundary Eichler
series, or a coefficientwise Fourier cut that has not itself been
proved modular. Thus the earlier warning against simply integrating
nonpositive modes over a fundamental-domain compact core remains
valid, but it does not rule out the double Rankin--Selberg argument.

For the arithmetic application, positive rank-one matching is
Corollary~\ref{cor:positive-rankone-ASW}. The next subsection
establishes the needed Green-energy bound, and
Theorem~\ref{thm:explicit-modular-completion} identifies the positive
arithmetic classes of the complete modular family. These results
supply the hypotheses of the support criterion for that family.
The corollary removes the complete nonpositive discrepancy after
double cuspidal pairing; it makes no assertion about the pairing
of an individual nonpositive coefficient in isolation.

\begin{lemma}[Strong parameter regularity of the full curvature]
\label{lem:KM-strong-parameter-regularity}
Let $X$ be the compact geometric modular orbifold and let
$K\Subset\mathbb H$ be a compact set in the arithmetic variable
$\tau=u+iv$. For every component of the full Kudla--Millson
curvature, the map
\[
 K\longrightarrow C^\infty(X,\Omega_X^2),\qquad
 \tau\longmapsto\Theta_\mu(\tau)
\]
is smooth, with all derivatives bounded in every smooth seminorm
on $K$. At every geometric cusp the extensions of these forms
and all their arithmetic derivatives are flat in the local
compactifying parameter.

Fix a smooth positive probability form $\nu$ on $X$. Write
\[
 H^1_\nu(X)=\{u\in H^1(X):\textstyle\int_Xu\nu=0\},\qquad
 H^{-1}_0(X)=\{\rho\in H^{-1}(X):\rho(1)=0\}.
\]
The inverse of $dd^c$ is a continuous linear map
\[
 G_\nu:H^{-1}_0(X)\longrightarrow H^1_\nu(X).
\]
Consequently
\[
 u_\mu(\tau)=G_\nu\!\left(\Theta_\mu(\tau)
                     -\left(\int_X\Theta_\mu(\tau)\right)\nu\right)
\]
is a smooth $H^1_\nu(X)$-valued function, and all arithmetic
parameter derivatives commute with $G_\nu$.
\end{lemma}
\begin{proof}
On the compact geometric core, the majorant
$\mathcal R(\lambda,z)+Q(\lambda)$ is uniformly positive
definite. Any fixed collection of derivatives in $u,v,z$
multiplies the summand of \eqref{eq:KM-proof-kernel} by a
polynomial in the lattice coordinates, with coefficients
bounded when $\tau\in K$. A slightly smaller Gaussian
therefore dominates every differentiated series uniformly.

At a geometric cusp write $z=x+iY$, with $x$ in a fixed
cylinder period and $Y\geq Y_0>0$. It suffices to treat one
of the finitely many rectangular rational cosets used in
\eqref{eq:KM-Poisson-majorant}. Before taking absolute values,
its Poisson terms are a fixed constant times
\[
 \frac{Y^3}{v^{3/2}}(m-c\overline\tau)^2
 e^{-\pi NY^2|m-c\tau|^2/v}
 e^{-2\pi Nv(b-cx)^2}
\]
multiplied by the unit-modulus phases from that calculation.
Those phases have polynomial exponents in the indicated
variables. In particular all their derivatives produce
polynomial factors. The pair $(m,c)=(0,0)$ is identically
zero before differentiation and remains zero after any
arithmetic or geometric derivative.

There is $a_K>0$ such that, for every real pair $(m,c)$,
\[
 \frac{|m-c\tau|^2}{v}
       =\frac{(m-cu)^2}{v}+c^2v
       \geq a_K(m^2+c^2)\qquad(\tau\in K).
\]
All frequencies belong to a fixed rational grid. Each mixed
derivative of the Poisson terms, including the factor $Y^{-2}$
from the area form, is bounded by a finite sum of expressions
\[
 C_K(1+Y)^A(1+|m|+|c|+|b-cx|)^B
 e^{-a_KY^2(m^2+c^2)}e^{-a_K(b-cx)^2},
 \qquad (m,c)\ne(0,0),
\]
where $A,B$ depend on the orders of differentiation. Powers
of $v^{-1}$ and $Y^{-1}$ are harmless on the specified ranges.
The shifted Gaussian moments in $b$ are uniformly bounded.
Absorbing the remaining polynomial factors into smaller
Gaussians and using the nonzero minimum length of the frequency
grid gives, for the cylinder density $k_\mu$,
\begin{equation}\label{eq:KM-mixed-compact-Gaussian}
 \left|\partial_u^a\partial_v^b
          \partial_x^c\partial_Y^d k_\mu(\tau,x,Y)\right|
     \leq C_{K,a,b,c,d}(1+Y)^M e^{-\epsilon_KY^2}.
\end{equation}
The same argument provides a summable majorant before taking
any of these derivatives, so termwise differentiation is
justified, not merely estimated after the fact.

For a cusp parameter $q=e^{2\pi iz/\kappa}$, rewriting the
two-form introduces a factor proportional to $|q|^{-2}$.
Each fixed number of derivatives in $q,\overline q$ introduces
only further factors bounded by a fixed exponential in $Y$
times a polynomial. These are absorbed by the Gaussian in
\eqref{eq:KM-mixed-compact-Gaussian}. Thus the density extends
by zero with every derivative zero at $q=0$, uniformly on $K$.
Finite orbifold charts give the same conclusion at elliptic
points. This proves the stated strong smoothness; in particular
it implies smoothness with values in $H^{-1}(X)$.

On the compact connected orbifold, the Dirichlet form is coercive
on $H^1_\nu(X)$ by the Poincare inequality. The weak formulation
of $dd^c$ and Lax--Milgram therefore identify this space
continuously with $H^{-1}_0(X)$. The kernel consists precisely
of constants, removed by the mean-zero condition. Applying
this fixed continuous inverse to the smooth mass-zero family
proves the final assertions.
\end{proof}

\subsection{A finite-energy growth estimate for the smooth coordinate}
\label{subsec:KM-Green-energy-growth}

Let $X$ be the compact complex modular orbifold. Fix a smooth positive
probability form $\nu$ and its mean-zero Green operator $G_\nu$.
Smoothness and the Green construction here are in orbifold charts;
finite stabilizer weights are retained in all integrals. No
coarse-coordinate smoothness assertion at elliptic points is needed.
Let $\Theta_\mu(\tau)$ be any component of the full standard
Kudla--Millson theta curvature, viewed as a measure in the geometric
variable, and put
\[
 d_\mu(\tau)=\int_X\Theta_\mu(\tau),\qquad
 \rho_{\mu,\tau}=\Theta_\mu(\tau)-d_\mu(\tau)\nu,
 \qquad dd^c u_{\mu,\tau}=\rho_{\mu,\tau},\quad
                         \int_Xu_{\mu,\tau}\nu=0.
\]
Thus $u_{\mu,\tau}$ is the smooth coordinate supplied by $G_\nu$.
For each fixed $\tau$ these are smooth compact currents: the full
Kudla--Millson theta curvature has Gaussian decay at geometric cusps.

\begin{lemma}[Weighted total variation of the full KM curvature]
\label{lem:KM-weighted-TV-growth}
For $\tau=a+iv$ with $v\geq v_0>0$, the following bounds hold
uniformly in $a$ and in the finitely many discriminant components:
\[
 \|\Theta_\mu(\tau)\|_{\rm TV}\ll1,\qquad
 \int_X H(z)\,|\Theta_\mu(\tau)|(z)\ll1+\log(1+v).
\]
Here $H(z)$ is the geometric cusp height $Y$ in each cusp strip and
is bounded on the compact core. In each cusp cylinder, the density
with respect to $dx\,dY$ is uniformly bounded. The same assertions
hold for $\rho_{\mu,\tau}$.
\end{lemma}

\begin{proof}
On a fixed compact core the ordinary majorant Gaussian gives a
uniform bound for the full theta series and its curvature density:
the polynomial factor $vB(x,h_z)^2$ is absorbed into a slightly
smaller Gaussian, and the zero vector contributes a fixed form.

At a geometric cusp use a rational cusp matrix, with cylinder
period $\kappa$ in $x$. The Poisson identity \eqref{eq:KM-Poisson-majorant}
for the full KM kernel bounds its density by
\[
 C\frac{Y}{v^{3/2}}
       \sum_{m\in D^{-1}\mathbb Z\setminus\{0\}}
           m^2e^{-cY^2m^2/v}+Ce^{-cvY^2}.
\]
The fixed denominator $D$ may depend on the cusp and discriminant
coset. To obtain this bound, the zero isotropic-frequency pair is
annihilated exactly by the squared-frequency polynomial in the
Poisson formula. The anisotropic one-dimensional Gaussian sum is
uniformly bounded for $v\geq v_0$. Terms with nonzero second
isotropic frequency have $|m-c\tau|^2\geq c^2v^2$; summing the first
frequency and absorbing their polynomial factors gives the last
Gaussian. These statements are uniform under the finitely many
rational coset shifts. Absolute values are taken after the original
Poisson grouping; no absolute estimate on the uncancelled original
isotropic sum is asserted.

The elementary bounds for the remaining one-dimensional Gaussian
moment imply the piecewise majorant
\[
 C\begin{cases}
 Y^{-2},&Y\leq\sqrt v,\\
 Yv^{-3/2}e^{-cY^2/v},&Y\geq\sqrt v,
 \end{cases}
       +Ce^{-cvY^2}.
\]
If $\sqrt v$ is below the bottom of the cusp strip the second bound
alone suffices. Its integral in $Y$ is uniformly bounded. Multiplying
by $Y$ gives $\int_{Y_0}^{\sqrt v}dY/Y=O(1+\log(1+v))$, and the
integral of the remaining Gaussian tail is bounded. This also gives
the uniform cylinder-density bound. The degree is bounded by total
variation, and a fixed smooth compact form has finite height moment
and bounded cylinder density. Subtracting $d_\mu(\tau)\nu$ therefore
preserves all the bounds.
\end{proof}

\begin{proposition}[Green energy and joint continuity]
\label{prop:KM-smooth-energy-growth}
For the above mean-zero smooth coordinate,
\[
 \mathcal E(u_{\mu,\tau})
       :=-\int_X\overline{u_{\mu,\tau}}\,dd^cu_{\mu,\tau}
        \ll1+\log(1+v).
\]
Consequently the Hermitian smooth-coordinate arithmetic pairing
satisfies
\[
 |\langle\mathfrak a_1(u_{\mu,\tau_1}),
                 \mathfrak a_1(u_{\eta,\tau_2})\rangle|
 \ll\sqrt{(1+\log(1+v_1))(1+\log(1+v_2))}.
\]
These estimates hold at every arithmetic cusp after the weight
slash operation, since such a slash is a finite linear combination
of the discriminant components of the KM theta kernel.
\end{proposition}

\begin{proof}
Use the norm on $H^1_\nu(X)$ equivalent to the square root of
\[
 \mathcal E(u)=-\int_X\overline u\,dd^cu
              =c_{dd^c}\int_X|\nabla u|^2
\]
in conformal charts, where $c_{dd^c}>0$ is the fixed constant
in our convention for $dd^c$. Orbifold weights are included in the
integral. The coefficient of the $c=0$ part of the Poisson
majorant, as a density with respect to $dx\,dY$, is at most
\[
 \frac{CY}{v^{3/2}}
 \sum_{m\in D^{-1}\mathbb Z\setminus\{0\}}
                    m^2e^{-cY^2m^2/v}.
\]
The terms with $c\ne0$ contribute at most $Ce^{-cvY^2}$,
after the Gaussian moments and their polynomial factors are
absorbed. The elementary moment estimate gives
\[
 |\Theta_\mu(u+iv)|\leq
 C\begin{cases}
 Y^{-2},&Y\leq\sqrt v,\\
 Yv^{-3/2}e^{-cY^2/v},&Y\geq\sqrt v,
 \end{cases}
 \quad+Ce^{-cvY^2}
\]
on the cusp cylinder, uniformly for $v\geq v_0$ and all real $u$.
This is an estimate on the Poisson-grouped density. The compact
core has a uniform bound from its positive-definite majorant.
It follows that the density is uniformly bounded on the
cylinders, its total variation is bounded, and
\[
 \int_X H(z)\,|\Theta_\mu(u+iv)|(z)
      \leq C(1+\log(1+v)).
\]
Indeed the first range contributes
$\int_{Y_0}^{\sqrt v}dY/Y$, when nonempty, while the second
weighted integral becomes a bounded Gaussian integral under
$Y=\sqrt v\,t$. Subtracting $d_\mu(\tau)\nu$ preserves these
bounds because $|d_\mu(\tau)|$ is bounded by total variation
and the fixed smooth form $\nu$ has finite height moment.

Write $\rho=\Theta_\mu-d_\mu\nu$. The compact Green kernel
satisfies, in one cusp cylinder,
\[
 |G_\nu(z,w)|\leq
 C+C\min(Y,Y')+C\log^+\!\frac1{d_{\rm cyl}(z,w)}.
\]
To obtain this bound, factor $q-q'$ into the larger of
$|q|,|q'|$ and the remaining difference of two unit-scale
exponentials. Their logarithms give respectively the
$\min(Y,Y')$ term and the ordinary diagonal logarithm.
Cross terms between separated cusp neighborhoods are bounded;
on the compact core the usual diagonal logarithmic estimate
applies. These estimates hold as well in finite orbifold charts.

The bounded cylinder density controls the diagonal logarithm:
its integral over a unit cylinder ball is uniformly finite.
Bounded total variation then controls the second integration.
The height term is bounded using $\min(Y,Y')\leq Y$ and the
weighted estimate above. Hence
\[
 \iint_{X\times X}|G_\nu(z,w)|\,
                     d|\rho|(z)d|\rho|(w)
       \leq C(1+\log(1+v)).
\]
For the smooth fixed-parameter density, Green representation
and integration by parts identify its energy with the
corresponding signed double integral. Taking absolute values
therefore proves the asserted energy bound. Poincare also gives
\begin{equation}\label{eq:KM-H1-growth-explicit}
 \|u_\mu(u+iv)\|_{H^1}
                \leq C\sqrt{1+\log(1+v)}.
\end{equation}
The smooth--smooth arithmetic pairing is a fixed multiple
of the Dirichlet pairing; Cauchy--Schwarz gives the joint bound.
At an arithmetic cusp, the weight slash replaces the curvature
by a fixed finite linear combination of its Weil components.
The same estimates therefore hold after every required slash.
\end{proof}

\begin{lemma}[Pairing with a fixed logarithmic Hodge metric]
\label{lem:Hodge-smooth-energy-continuity}
Choose a smooth reference metric on the same Hodge line, with
curvature $\deg(\lambda)\nu$. The difference between the fixed
Petersson metric and this reference metric is a function $b$ with
$b=c\log Y+b_0$ at each cusp, where $b_0$ extends smoothly in the
local cusp parameter. In particular $b$ has finite Dirichlet energy.
Therefore
\[
 |\langle\mathfrak a_1(u_{\mu,\tau}),\widehat\lambda\rangle|
                       \ll\sqrt{1+\log(1+v)}.
\]
The same statement holds for the fixed modified line $H$.
\end{lemma}

\begin{proof}
In an integral Tate frame the Petersson metric is a constant
multiple of the geometric height. Relative to a smooth reference
metric on the same line, its logarithmic metric difference is
$b=c\log Y+b_0(q)$, with $b_0$ smooth. In a cusp cylinder,
\[
 \int_{Y_0}^\infty\int_0^\kappa
               |\nabla(c\log Y)|^2\,dx\,dY
       =\frac{\kappa c^2}{Y_0}<\infty.
\]
The $L^2$ norm for a smooth compact volume form is finite as
well, since that volume decays exponentially in $Y$. Thus
$b\in H^1(X)$. Smooth truncations which agree with $b$ for
$Y\leq R$ and replace $\log Y$ by a smoothed constant for
$Y\geq2R$ converge in $H^1$: their tail Dirichlet error is
$O(R^{-1})$, and their compact-volume $L^2$ error tends to zero.
The normal flux of $\log Y$ at height $R$ is $O(R^{-1})$,
so its distributional $dd^c$ has no additional cusp atom.

The reference curvature contributes zero on mean-zero $u$.
For smooth $u$, the remaining pairing is the fixed multiple
of the Dirichlet pairing with $b$ dictated by the arithmetic
metric-variation formula. Its bound is
\[
 C\,\mathcal E(u)^{1/2}\mathcal E(b)^{1/2}.
\]
The truncations prove agreement with the extended arithmetic
pairing; continuity gives the same functional on all of
$H^1_\nu(X)$. Equivalently the curvature of the fixed Hodge
metric belongs to $H^{-1}(X)$. Substituting
\eqref{eq:KM-H1-growth-explicit} proves the desired estimate.
The source and target Hodge factors in $H$ satisfy the same
argument; their cusp-width ratios only add constants to $b$.
\end{proof}

\paragraph{Use in a direct-sum modular completion.}
In the arithmetic Chow decomposition, assume the remaining
Mordell--Weil, degree, vertical, and fixed-reference-line coordinates
have been constructed as their known finite-dimensional modular
forms and Eisenstein series or derivatives. Their products have
joint polynomial growth. The proposition handles the full
infinite-dimensional smooth--smooth pairing, and the last lemma
handles the fixed Hodge functional used when solving the constant
coordinate. Thus the resulting complete scalar arithmetic product
kernel has joint moderate growth. This conclusion uses the explicit
direct-sum construction and the finite-energy topology, rather than
inferring a bilinear growth bound from separate weak modularity
against individual linear functionals.

\subsection{An explicit modular completion and paired arithmetic doubling}
\label{subsec:explicit-modular-completion}

A complete modular arithmetic family is constructed from fixed
coordinates. Its strictly positive coefficients are identified
with the standard arithmetic special divisors in this subsection.
The nonpositive coefficients are initially defined by the
construction and are compared with the standard prescription in
Section~\ref{subsec:nonpositive-completion-comparison}.

We work on the regular modular stack and its compact complex
orbifold. Smoothness throughout this subsection means orbifold
smoothness. At an elliptic point, an invariant smooth two-form
need not be smooth in the branched coordinate of the coarse curve.
Fix a positive orbifold form \(\nu\) of
stack mass one. Normalize the generic cuspidal divisor
\[
 B_\infty=\frac{P_\infty}{\deg_{\mathbb Q}P_\infty}
\]
to have degree one. Choose its arithmetic representative
\(\widehat D_\infty\) with curvature \(\nu\), and make the
vertical corrections so that its pairing with the chosen vertical
quotient is zero. Use this cuspidal degree-one class in the
admissible decomposition. Throughout this construction, the
admissible Mordell--Weil representatives have zero compact
curvature, pair trivially with the vertical quotient, and are
normalized to be orthogonal to \(\widehat D_\infty\).
The last condition fixes their remaining constant because
\(\langle\one,\widehat D_\infty\rangle_{\rm GS}=1/2\).
The series \(Y\) and its vector-valued lift use this same realization.

For each bad prime choose a basis of the vertical quotient modulo
the full fibre. Write its arithmetic representatives as
\[
 Y_i=\widehat V_i+c_i\one,\qquad
 \widehat V_i=(V_i,0),\quad \one=(0,1),\qquad
 \langle\widehat D_\infty,Y_i\rangle_{\rm GS}=0.
\]
Let \(M=(\langle Y_i,Y_j\rangle_{\rm GS})\), with its primewise
blocks. It is invertible. The following uniqueness statement
will be used in both coefficient comparisons.

\begin{lemma}[Determination by admissible coordinates]
\label{lem:arithmetic-coordinate-uniqueness}
Let \(\mathscr H\) be the chosen admissible direct sum, with its
mean-zero metric summand completed in \(H^1(X)\). Fix a continuous
linear functional \(\ell\) on this space with \(\ell(\one)\ne0\).
Two classes in \(\mathscr H\) are equal if they have the same
Mordell--Weil projection, degree, pairings with all \(Y_i\),
distributional curvature, and value under \(\ell\).
\end{lemma}
\begin{proof}
By \eqref{eq:arith-decomposition}, the difference of two such
classes is \(V+c\one+\mathfrak a_1(u)\), where \(V\) belongs
to the chosen vertical quotient and \(u\in H^1(X)\) has mean zero.
Equality of curvature gives \(dd^c u=0\). A weakly harmonic
function on the compact connected orbifold is constant, so the
mean-zero condition gives \(u=0\). The remaining pairings with
the \(Y_i\) are obtained by multiplying the coefficient vector
of \(V\) by \(M\), since constant metric classes pair trivially
with vertical classes.
Invertibility of \(M\) therefore gives \(V=0\). Finally
\(0=\ell(c\one)=c\ell(\one)\) gives \(c=0\).
\end{proof}

\paragraph{The modular input functions.}
All the following functions have the same scalar weight \(3/2\)
and multiplier on \(H_N\).
\begin{enumerate}
\item Let \(Y(\tau)=R\widehat\Phi_{\rm MW}(\tau)\) be the
admissible lift of the cuspidal-base Heegner series in the Jacobian,
whose modularity and actual positive generic-cycle comparison
were fixed above. Its positive coefficient is
\[
 \operatorname{pr}_{\rm MW}\widehat Z_{\rm std}(t,v)
        =[\mathcal Z(t)-(\deg\mathcal Z(t))B_\infty]^{\rm adm}.
\]
These trace Heegner classes lie in the Mordell--Weil group over
a fixed number field; thus the relevant real or complex
Mordell--Weil space is finite dimensional.
\item Let \(\Theta_0(\tau)\) be the full zero-component
Kudla--Millson curvature on the geometric orbifold, including
its zero vector, and put
\[
 d_N(\tau)=\int_X\Theta_0(\tau)
           =-\frac{\psi(N)}{12}E_N^{(1)}(\tau,\tfrac12).
\]
The equality is the full analytic Du--Yang identity with the
two central normalizations recorded above. It is not a claim
about a nonpositive arithmetic divisor prescription.
\item Define full scalar functions \(F_{p,a}^N(\tau)\) by
the finite component formulas, replacing each degree series
there by the just-defined \(d_M(\tau)\). Explicitly, at an
outer component, \(n=v_p(N)\geq1\),
\[
 F_{p,\pm n}^N
   =\frac{\log p}{2}
       \left\{d_N-(p-1)
                 \sum_{j=1}^{\lfloor n/2\rfloor}d_{Np^{-2j}}\right\}.
\]
For \(-n<a<n\), set
\(F_{p,a}^N=F_{p,a}^{Np^{-2}}\) recursively; when \(p\)
no longer divides the level the unique component is the
full fibre and its function is \((\log p)d_M\).
These are exactly the positive-index finite identities of
Corollary 8.2.3 and Proposition 8.2.4 of \cite{Zhu}.
They are now definitions at all other indices. They respect
the full-fibre relation
\(\widehat F_p=2\log p\,\one\). Since \(N/M\) is a square,
the lower-level characters coincide on \(H_N\).
For a linear combination \(V_i\) of components, use the same
linear combination \(F_i\) of these functions and set
\[
 v_i(\tau)=F_i(\tau)+\frac{c_i}{2}d_N(\tau),\qquad
 V(\tau)=\sum_{i,j}(M^{-1})_{ij}v_j(\tau)Y_i.
\]
Thus \(\langle V(\tau),Y_i\rangle=v_i(\tau)\).
\item Fix the continuous mean-zero Green inverse \(G_\nu\):
\[
 dd^cG_\nu(\eta)=\eta,\qquad
 \int_XG_\nu(\eta)\nu=0,\qquad \int_X\eta=0.
\]
Set
\[
 u(\tau)=G_\nu\{\Theta_0(\tau)-d_N(\tau)\nu\}.
\]
Compact orbifold elliptic theory supplies this operator.
For each fixed \(\tau\) the full theta curvature is
smooth, including at the geometric cusps by its Poisson
Gaussian bound. On compact sets of arithmetic \(\tau\)
that bound is uniform with derivatives, so taking its
Fourier coefficients commutes with this fixed inverse.
\item Finally put
\[
 J_N(\tau)=12\left.\frac{d}{ds}
       \{C_N(s)E_N^{(1)}(\tau,s-\tfrac12)\}\right|_{s=1},
\qquad
 C_N(s)=-\frac{s}{4\pi}\Lambda(2s)
          \prod_{p\mid N}(1-p^{-2s})N^{(3s+1)/2}.
\]
This is a fixed Eisenstein derivative and value. It is
modular and has polynomial growth with logarithmic factors.
The positive-index Kronecker-limit comparison proves
\[
 (J_N)_t(v)
      =\langle\widehat Z_{\rm std}(t,v),
                                      \widehat\Delta_N\rangle_{\rm GS},
                      \qquad t>0,
\]
after removing the common \(q^t\) factor.
\end{enumerate}

Let
\[
 \ell(z)=\langle z,\widehat\Delta_N\rangle_{\rm GS},\qquad
 \ell(\one)=\tfrac12\deg\Delta_N\ne0.
\]
Although the Petersson metric has logarithmic height
singularities, this is a well-defined fixed functional in
the present extended arithmetic pairing. On a pure smooth
metric it is
\(\ell(\mathfrak a_1(f))=\tfrac12\int_Xf\,\operatorname{curv}
\widehat\Delta_N\). That curvature is a finite-mass
hyperbolic density, so the functional is continuous for
uniform convergence. Its finite-energy continuity and
growth on \(u(\tau)\) also follow from
Lemma~\ref{lem:Hodge-smooth-energy-continuity}.

Define, explicitly,
\begin{align}
 Q(\tau)&=Y(\tau)+d_N(\tau)\widehat D_\infty
                                      +V(\tau)+\mathfrak a_1(u(\tau)),
                         \label{eq:completion-preconstant}\\
 c(\tau)&=\frac{J_N(\tau)-\ell(Q(\tau))}{\ell(\one)},\\
 \widehat\Phi_{\rm cmp}(\tau)&=Q(\tau)+c(\tau)\one.
                         \label{eq:explicit-modular-completion}
\end{align}
These equations define every coefficient, rather than
postulating a matching nonpositive arithmetic family.

\begin{theorem}[Modularity and positive coefficient compatibility]
\label{thm:explicit-modular-completion}
The family \(\widehat\Phi_{\rm cmp}\) is a weight-\(3/2\)
modular arithmetic family on \(H_N\) in the direct-sum
topology just specified. It has jointly moderate scalar
arithmetic product kernel. Its positive coefficients
are strong arithmetic-class equalities:
\begin{equation}\label{eq:completion-positive-class-equality}
 (\widehat\Phi_{\rm cmp})_t(v)
                         =\widehat Z_{\rm std}(t,v),\qquad t>0.
\end{equation}
Its Mordell--Weil projection is \(Y(\tau)\).
\end{theorem}
\begin{proof}
Each finite-dimensional input is modular: Heegner
modularity gives \(Y\); the analytic genus-one
Siegel--Weil formula gives \(d_N\); the finite
component formulas give \(v_i,V\) as lower-level
Eisenstein combinations; and the Kronecker-limit
normalization gives \(J_N\).
The fixed continuous Green inverse is linear and
commutes with the arithmetic modular action, which
acts trivially on the geometric variable. Thus \(u\)
is modular. Applying a fixed linear functional and
dividing by the nonzero constant \(\ell(\one)\)
shows that \(c\) and \(\widehat\Phi_{\rm cmp}\)
are modular.

Fix \(t>0\). The positive generic-cycle comparison
gives the same Mordell--Weil coordinate on both
sides of \eqref{eq:completion-positive-class-equality}.
The positive degree formula gives the same degree.
Because the metrics on the raw vertical components
are zero, the actual positive finite intersection
formulas give the same \(Y_i\)-pairings, including
the \(c_i d_N/2\) correction. The compact Green
equation gives the same curvature
\(\Omega_{t,v}\); the mean-zero normalization
therefore gives the same smooth coordinate.
The positive Kronecker-limit comparison and the definition
of \(c\) also give equal \(\ell\)-values. Since
\(\ell(\one)\ne0\), Lemma~\ref{lem:arithmetic-coordinate-uniqueness}
proves \eqref{eq:completion-positive-class-equality}.
This comparison uses only the positive coefficients.

The only infinite-dimensional contribution to
the scalar product is the smooth coordinate.
Proposition~\ref{prop:KM-smooth-energy-growth}
gives its joint logarithmic energy bound at every
arithmetic cusp, and
Lemma~\ref{lem:Hodge-smooth-energy-continuity}
controls the fixed functional \(\ell\).
All other coordinates are finite-dimensional
modular forms or Eisenstein values and derivatives,
with polynomial cusp growth. The mixed terms
are fixed continuous pairings. These estimates
prove joint moderate growth of the complete
scalar product, including one-variable cusp
regimes; no conclusion is drawn merely from
separate weak modularity. The direct-sum formula
also gives the asserted Mordell--Weil projection.
\end{proof}

\begin{corollary}[Existence of the full lift and its Mordell--Weil value]
\label{cor:completion-full-lift-MW}
For \(f\in S_{3/2}(H_N,\chi_N)\), the integral
\[
 \widehat\theta_{\rm cmp}(f)
    =\int_{H_N\backslash\mathbb H}
       \widehat\Phi_{\rm cmp}(\tau)\overline{f(\tau)}
                                      v^{3/2}\,d\mu(\tau)
\]
exists as a Bochner integral in the finite-dimensional
admissible summands plus the mean-zero finite-energy
metric space. In fact it is in the original
Mordell--Weil summand and equals
\[
 \widehat\theta_{\rm cmp}(f)
       =\int_{H_N\backslash\mathbb H}
                       Y(\tau)\overline{f(\tau)}
                                      v^{3/2}\,d\mu(\tau).
\]
The continuous arithmetic product commutes with the
two integrals. Thus \eqref{eq:completion-paired-doubling}
is also the arithmetic product identity for these full
completed lifts. The previously established Petersson
adjunction identifies this lift with the vector
Mordell--Weil lift at \(Tf\), with the same
unnormalized measures.
\end{corollary}
\begin{proof}
Give the finite-dimensional summands any fixed
Euclidean norms and the mean-zero metric summand
its \(H^1\), or equivalently Dirichlet, norm.
The mean-zero Poincare inequality makes this a
Hilbert space. The arithmetic product is a continuous
bilinear form there: its metric--metric part is
one half the Dirichlet pairing; the mixed parts
are the fixed continuous functionals in the
admissible decomposition. In particular the
functional \(\ell\) is continuous by the
finite-energy Hodge-metric estimate.

The finite-dimensional coordinate functions have polynomial
arithmetic cusp growth, with logarithmic factors. Equation
\eqref{eq:KM-H1-growth-explicit} controls the metric coordinate,
and the fixed Hodge functional controls its contribution to
the constant coordinate. At each arithmetic cusp there are
constants $C,M,A$ such that, in its slash coordinate,
\[
 \|\widehat\Phi_{\rm cmp}|_{3/2}\sigma(u+iv)\|_{\mathscr H}
       \leq C(1+v)^M(1+\log(1+v))^A.
\]
A cusp form satisfies
$|(f|_{3/2}\sigma)(u+iv)|\leq C_f e^{-\epsilon_fv}$.
Thus the norm of the Petersson integrand has tail bounded by
\[
 C\int_{v_0}^\infty
       (1+v)^M(1+\log(1+v))^A e^{-\epsilon_fv}v^{-1/2}\,dv
       <\infty.
\]
There are only finitely many cusps. On the remaining compact
region, continuity gives integrability. The integrand is
strongly measurable because it is continuous with values
in the separable Hilbert space $\mathscr H$. This proves
Bochner existence.

Write $\Phi=\widehat\Phi_{\rm cmp}$. Write $\Phi=\widehat\Phi_{\rm cmp}$. If $B$ is the continuous
arithmetic bilinear pairing on this
direct sum, then $|B(x,y)|\leq C_B\|x\|\|y\|$. The product
of the two preceding integrable norm majorants dominates
the double Petersson integrand. Consequently scalar Fubini
and the defining continuity property of the Bochner integral
give
\[
 B\!\left(\int\Phi\overline f\,d\mu_{3/2},
          \int\Phi\overline h\,d\mu_{3/2}\right)
   =\iint B(\Phi(\tau_1),\Phi(\tau_2))
                 \overline{f(\tau_1)h(\tau_2)}
                    \,d\mu_{3/2,1}\,d\mu_{3/2,2}.
\]
Here $d\mu_{3/2}=v^{3/2}d\mu$.
The continuous map $dd^c:H^1\to H^{-1}$ likewise commutes
with the integral. These convergence statements do not by
themselves prove the subsequent vanishing of any coordinate;
that uses the orthogonality results cited in the full-lift proof.

The degree and vertical coordinate integrals
vanish by ordinary Eisenstein cusp orthogonality.
For the smooth coordinate, distributional
\(dd^c\) is continuous from \(H^1\) to \(H^{-1}\).
Its integral is
\[
 \int_{H_N\backslash\mathbb H}
       \{\Theta_0(\tau)-d_N(\tau)\nu\}
                    \overline{f(\tau)}v^{3/2}d\mu(\tau)=0.
\]
The first term vanishes by the proved
Kudla--Millson raising/Stokes lemma and compact
testwise Fubini, including the spatial cusps;
the second is Eisenstein orthogonality.
The mean-zero integral of \(u\) is consequently
a harmonic \(H^1\) function, hence zero.

The remaining integral is \(x+c_f\one\),
where \(x\) is the displayed Mordell--Weil lift.
Its \(\ell\)-value is zero because \(J_N\)
is an Eisenstein derivative. Thus even before
using the next orthogonality argument one has
the exact formula
\[
 \widehat\theta_{\rm cmp}(f)
          =x-\frac{\ell(x)}{\ell(\one)}\one.
\]
Such a constant correction does not change
any product of these degree-zero lifts:
\(\langle x,\one\rangle=\deg(x)/2=0\)
and \(\langle\one,\one\rangle=0\).
In the stated admissible normalization, moreover,
\(\ell(x)=0\). To justify this last statement
with the cusp metric, the generic divisor of
\(\widehat\Delta_N\) is cuspidal, so its
degree-zero Jacobian class is zero by
Manin--Drinfeld. Relative to the same
cuspidal degree-one class it therefore has
only degree, vertical, constant, and metric
coordinates. The metric difference from a
smooth reference metric has a logarithm of
the cusp height and finite Dirichlet energy.
The admissible Mordell--Weil class is
orthogonal to each of these coordinates;
the metric assertion follows by finite-energy
approximation from its zero compact curvature.
Thus \(0=\ell(x+c_f\one)=c_f\ell(\one)\)
and \(c_f=0\). In particular the final
Bochner value lies in the original finite-dimensional
arithmetic space despite having been integrated
in its energy completion.
\end{proof}

\begin{remark}[Other choices of degree-one class]
If the fixed admissible degree-one class is
\(\widehat D_0\) with generic class \(D_0\),
replace \(Y(\tau)\) above by
\[
 Y(\tau)+d_N(\tau)[B_\infty-D_0]^{\rm adm}
\]
and replace \(\widehat D_\infty\) by
\(\widehat D_0\). The sign follows from
\([Z(t)-d_tD_0]=[Z(t)-d_tB_\infty]+d_t[B_\infty-D_0]\).
The other coordinates are adjusted with the
same fixed representatives. This preserves
the construction and the positive coefficients.
We use the cuspidal choice to keep the
Mordell--Weil coordinate equal to \(Y\)
without an Eisenstein multiple of a fixed
Jacobian class.
\end{remark}

\begin{theorem}[Paired arithmetic doubling for the explicit completion]
\label{thm:explicit-completion-paired-doubling}
Let \(f,h\in S_{3/2}(H_N,\chi_N)\). Use the
unnormalized invariant Petersson measure and
the arithmetic degree convention of the preceding
coefficient comparisons. For the complex-bilinear
extension of the arithmetic product, one has
\begin{align}
 &\int_{H_N\backslash\mathbb H}
       \int_{H_N\backslash\mathbb H}
  \langle\widehat\Phi_{\rm cmp}(\tau_1),
                 \widehat\Phi_{\rm cmp}(\tau_2)\rangle_{\rm GS}
       \overline{f(\tau_1)h(\tau_2)}
           v_1^{3/2}v_2^{3/2}\,d\mu_1d\mu_2
                       \label{eq:completion-paired-doubling}\\
 &\qquad=\frac{\psi(N)}{24}
  \int_{H_N\backslash\mathbb H}
       \int_{H_N\backslash\mathbb H}
    \left.\partial_s\mathcal E_N
              (\operatorname{diag}(\tau_1,\tau_2),s)\right|_{s=0}
       \overline{f(\tau_1)h(\tau_2)}
              v_1^{3/2}v_2^{3/2}\,d\mu_1d\mu_2 .
 \notag
\end{align}
Both integrals converge absolutely. The Hermitian
version is obtained by conjugating the second
arithmetic argument and using the reflected
second analytic variable with the same convention.
\end{theorem}
\begin{proof}
The two complete scalar kernels are separately
modular with the same weight and multipliers.
The preceding theorem proves joint moderate
growth of the arithmetic kernel. Standard
Eisenstein growth and its parameter derivative
give the corresponding bound for the analytic
kernel, including simultaneous cusps. Their
difference therefore satisfies the full-kernel
hypotheses of the double Rankin--Selberg
support criterion.

For each \(m,n>0\), strong positive class
compatibility and the proved standard
positive-product identity give its arithmetic
coefficient as
\[
 \sum_{r\in\mathbb Z}
   \mathcal A^{\rm std}
       \left(\begin{pmatrix}m&r/2\\r/2&n\end{pmatrix},
                       \operatorname{diag}(v_1,v_2)\right).
\]
Indeed the two fixed-height Fourier integrals
are Bochner integrals on compact tori, and
continuity of the arithmetic bilinear form
allows it to pass through those two integrals.
Thus the arithmetic coefficient is the product
of the two class coefficients; no unproved
termwise product of divisor representatives
is required.
Every matrix in this sum is nonsingular or
nonzero positive semidefinite of rank one.
The all-level nonsingular comparison and the
positive rank-one Kronecker-limit comparison
identify each term with the matching Fourier
coefficient of \(\psi(N)\partial_s\mathcal E_N/24\).
The proved proper and dependent estimates
justify the original off-diagonal sum. Thus
all strictly positive diagonal Fourier
coefficients of the complete residual vanish.

Apply the bilinear version of the double
Rankin--Selberg criterion: unfold two auxiliary
weight-zero Eisenstein series, and take their
two constant residues at \(1\). The horizontal
integral selects precisely those positive
diagonal coefficients, so it is zero.
This proves \eqref{eq:completion-paired-doubling}.
Joint polynomial growth and the exponential
decay of \(f,h\) give absolute convergence
throughout, as required in that criterion.
\end{proof}

\paragraph{Scope of the paired argument.}
This theorem proves the paired identity for the explicitly
constructed modular family. The comparison with the standard
nonpositive coefficients is proved next and yields the pointwise
identity in Section~\ref{subsec:full-standard-arithmetic-doubling}.
The Rankin--Selberg argument itself uses the entire modular
residual, including all nonpositive terms. It does not require
the isolated Epstein correction to be modular or separately
orthogonal to cusp forms.

\subsection{The nonpositive coefficients of the explicit completion}
\label{subsec:nonpositive-completion-comparison}

We use the actual stack degree and
$\langle\mathfrak a_1(1),L\rangle=\deg L/2$ throughout.
Write $A_t(v)=\widehat Z_{\rm std}(t,v)$, so that
\[
 A_t(v)=(C_t(v),F_t(v))\quad(t<0),\qquad
 A_0(v)=H+(C_0(v),F_0(v))-\mathfrak a_1(\log v).
\]
The standard compact Green equations identify the curvature of
the complete family with the full KM theta curvature, including
its zero vector. All cusp and finite restrictions below are on
the normalized integral model, using its integral Tate frames.
The precise previously proved inputs are
Lemma~\ref{lem:all-level-normal-cusp-parameter} and
Lemma~\ref{lem:complete-zero-cusp-regular-part}.

\begin{lemma}[Restriction of a complete standard coefficient to a cusp]
\label{lem:standard-all-index-cusp-restriction}
For every $t\leq0$ and every closed arithmetic cusp $P$, the
restriction of $A_t(v)$ to $P$, with the extended cusp metric,
has arithmetic degree zero. Its logarithmic-height exponent is
zero. This assertion concerns the complete $A_0(v)$ when $t=0$;
its three displayed summands must not be restricted separately
with their logarithmic-height metrics omitted.
\end{lemma}
\begin{proof}
For $t<0$, the proved cusp expansion is
\[
 F_t(v,z)=-g_P(t,v)\log|q_P|^2+E_P(v,z),\qquad E_P(v,z)\longrightarrow0.
\]
The same expansion with $g_P=0$ holds at a cusp not supporting
$C_t(v)$. Thus the normal frame of $\mathcal O(C_t(v))$ has
limiting norm one. There is no logarithm of the cusp height.
Lemma~\ref{lem:all-level-normal-cusp-parameter} shows that its finite
restriction has degree zero: after an integral Tate chart
extension, the cusp ideal is generated by its normalized
parameter, and changes of the normal frame are integral units
with the corresponding product-formula cancellation. Distinct
closed arithmetic cusps are treated before summing embeddings;
complex embeddings are not counted as distinct arithmetic
closures. Descent divides the resulting degree by the chart
degree. This is exactly the arithmetic normal restriction,
not the invalid global cycle $(P,0)$.

For $t=0$, Lemma~\ref{lem:complete-zero-cusp-regular-part} says that
the metric of $H+(C_0(v),F_0(v))-\mathfrak a_1(\log v)$, in the
product of the source and target integral Tate frames and the
normal frame, has logarithmic-height exponent zero and continuous
regular norm one. The constants in the Mellin expansion of $F_0$
cancel the two Hodge metric constants and $\log v$ exactly.
The source and target Tate frames restrict to free integral
lines; the normal frame has the same zero finite degree as above.
The complete restricted metric is one, so its arithmetic degree
is zero. This proves both assertions, including at bad primes.
For real cusp coefficients, apply the assertion first to the
actual rational line generators and smooth metric changes, then
extend by real linearity. No single-valued rational section of
a hypothetical irrational tensor power is used.
\end{proof}

\begin{lemma}[A cusp-section height identity]
\label{lem:cusp-section-height-identity}
Let $L$ be an extended arithmetic class whose metric extends
with logarithmic-height exponent zero at the cusps. Suppose its
restricted arithmetic degree at $P_\infty$ is zero. Then
\[
 \langle L,\widehat\Delta_N\rangle
       =-\int_{\mathcal X_0(N)(\mathbb C)}
                    \log\|\Delta_N\|\,\operatorname{curv}(L).
\]
The integral is assumed convergent; for $L=A_t(v)$ it is
absolutely convergent by the full KM Gaussian cusp estimate.
\end{lemma}
\begin{proof}
The arithmetic line $\widehat\Delta_N$ is represented by the
cusp section with divisor $\psi(N)\varphi(N)P_\infty$ and
Green function $-2\log\|\Delta_N\|$. Choose rational sections
of the actual rational line generators of $L$ whose generic
divisors avoid $P_\infty$, treating the remaining smooth metric
functions by the usual Green variation formula. Then extend
the resulting identity by real or complex linearity.
K\"uhn's intersection formula is then applied to disjoint
generic divisors, with $\widehat\Delta_N$ as the first factor.
Its finite restriction and the finite continuous cusp value
combine to
$\psi(N)\varphi(N)\widehat\deg(L|_{P_\infty})=0$.
There is no logarithmic-height counterterm because the second
factor has exponent zero. The remaining term is
$-\int\log\|\Delta_N\|\operatorname{curv}(L)$, as claimed.
Independence of the auxiliary rational section is the product
formula. In particular this proof does not apply the disjoint-
section formula directly to two canonical sections having a
common cusp, and does not drop a finite self-restriction term.
\end{proof}

\begin{proposition}[All nonpositive discriminant heights]
\label{prop:nonpositive-discriminant-heights}
For every $t\leq0$,
\[
 \langle A_t(v),\widehat\Delta_N\rangle
   =12\left.\partial_s
      \{C_N(s)e_{N,t}(v,s-1/2)\}\right|_{s=1}=(J_N)_t(v).
\]
Here $J_N$ is the fixed full modular function used to construct
$\widehat\Phi_{\rm cmp}$; the common factor $q^t$ is removed.
\end{proposition}
\begin{proof}
Apply the two preceding lemmas to the complete standard class.
Its curvature is the $t$-coefficient of $\Theta_{L,0}$, and
$\log\|\Delta_N\|$ grows only linearly and logarithmically at
geometric cusps. The full KM estimate therefore permits taking
this Fourier coefficient under the geometric integral. Thus
\[
 \sum_t\langle A_t(v),\widehat\Delta_N\rangle q^t
       =-\frac12 I_{\rm eff,0}(\tau,\log\|\Delta_N\|)
\]
at the level of the relevant coefficient; no modularity of the
standard arithmetic family is assumed in this calculation.
Du--Yang's full analytic Kronecker-limit identity is
$I_{\rm eff}(\log\|\Delta_N\|)=-12\mathcal E_{L}^{\prime}(1)$.
The full-parabolic scalar component relation
$\mathcal E_{L,0}(s)=2C_N(s)E_N^{(1)}(s-1/2)$ and the actual
stack half give the stated factor $12$. This argument includes
$t=0$ and uses no zero-index height identity from Zhu24.
\end{proof}

\paragraph{The other coordinates.}
The degree and curvature comparisons at all nonpositive indices
follow by integrating, or retaining, the complete compact Green
equation and the full analytic genus-one Siegel--Weil identity.
The generic divisors for $t<0$ are cuspidal; the generic Hodge
class at $t=0$ is also rationally cuspidal by $\Delta_N$.
Manin--Drinfeld consequently makes their degree-zero Jacobian
coordinates zero, as in the cuspidal-base Heegner series.

The finite vertical comparison can be checked
directly from normalized cusps, without importing a zero-index
arithmetic formula. Write $N=p^nM$, $p\nmid M$, and index the
reduced components by $a=2k-n$, $0\leq k\leq n$. Put
$h=\min(k,n-k)$. A cusp with denominator $p^kd$, $d\mid M$,
has cyclotomic $p$-ramification index $\varphi(p^h)$ and
$\varphi(p^h)$ complex embeddings over its prime-to-$p$ data.
Its intersection with the reduced component therefore has its
generic degree divided by $\varphi(p^h)$. The standard cusp
coefficient is
\[
 g_{p^kd}^{N}(t,v)=g_{d}^{Mp^{|a|}}(t,v)
\]
with the appropriate outer denominator on the right; this holds
both for $t<0$ (including the support condition) and for $t=0$.
Indeed its common factors are
$\sqrt N/p^h=\sqrt{Mp^{|a|}}$ and
$N/p^{2h}=Mp^{|a|}$.
Thus every interior component reduces exactly to the corresponding
outer component at level $Np^{-2h}$.

Let $d_L^{\rm cusp}(t,v)$ denote the degree of the cusp part.
At an outer component, direct summation of the two extreme
cusp sets gives
\[
 F_{p,\pm n}^{\rm cusp}
   =\frac{\log p}{2}\left\{d_N^{\rm cusp}
            -(p-1)\sum_{j=1}^{\lfloor n/2\rfloor}
                              d_{Np^{-2j}}^{\rm cusp}\right\}.
\]
To verify the equality coefficientwise, group cusp denominators
by $h$. The two extreme sets have weight one; the other pairs
have weights $\varphi(p^h)=p^{h-1}(p-1)$, with a single middle
set when $n$ is even. The displayed subtraction telescopes these
weights. More explicitly, let $A_r$ be the mass of one outer
cusp set at level $Mp^r$, with $A_0=d_M^{\rm cusp}$. Then
\[
 d_{Mp^n}^{\rm cusp}
  =2\sum_{h=0}^{\lfloor(n-1)/2\rfloor}\varphi(p^h)A_{n-2h}
    +\begin{cases}\varphi(p^{n/2})A_0,&n\text{ even},\\0,&n\text{ odd}.
      \end{cases}
\]
In the stated difference, the coefficient of every $A_{n-2h}$
with $h\geq1$ is zero by
$\sum_{r=0}^{h-1}\varphi(p^r)=p^{h-1}$; the coefficient of
$A_n$ is $2$. This verifies the prescribed finite formula at every
nonpositive index, including nonsplit negative indices where
all its cusp terms are zero.

At $t=0$ one must also keep $H$. Its pairing with a reduced
component is
\[
 \langle H,\widehat{\mathcal X}_{p}^{a}(N)\rangle
    =-\frac{\psi(M)}{24}(p^{|a|}+1)\log p.
\]
This follows from the source and target maps of that component,
whose degrees are $p^{\max(-a,0)}$ and $p^{\max(a,0)}$, and
$\deg\lambda_M=\psi(M)/24$. At the outer component it equals
the constant part of the same finite formula, since
\[
 \psi(N)-(p-1)\sum_{j=1}^{\lfloor n/2\rfloor}\psi(Np^{-2j})
                          =\psi(M)(p^n+1).
\]
The interior equality is the already noted reduction in level.
The constant metric $-\mathfrak a_1(\log v)$ has zero pairing
with vertical classes. Therefore the complete $t=0$ finite
coordinate agrees as well. The full-fibre relation is preserved
before passing to the vertical quotient.

\begin{theorem}[The explicit completion has the standard nonpositive coefficients]
\label{thm:nonpositive-standard-completion}
With the cusp data of
Lemmas~\ref{lem:all-level-normal-cusp-parameter} and
\ref{lem:complete-zero-cusp-regular-part}, for every $t\leq0$
and $v>0$ one has the strong arithmetic-class equality
\[
       (\widehat\Phi_{\rm cmp})_t(v)=A_t(v).
\]
\end{theorem}
\begin{proof}
The preceding comparisons identify the degree, Mordell--Weil
coordinate, every vertical pairing, and curvature. Proposition~
\ref{prop:nonpositive-discriminant-heights} identifies the values
of $\ell=\langle\,\cdot\,,\widehat\Delta_N\rangle$ by evaluating
the standard coefficients from their Tate restrictions and the
analytic Kronecker-limit formula. As
$\ell(\mathfrak a_1(1))=\deg\Delta_N/2\ne0$,
Lemma~\ref{lem:arithmetic-coordinate-uniqueness} proves the equality.
\end{proof}

\begin{corollary}[Strong Fourier convergence of the standard family]
\label{cor:standard-family-strong-Fourier}
Let $\mathscr H$ be the finite-dimensional admissible arithmetic
summands plus the mean-zero $H^1$ metric summand used for the
explicit completion. For each fixed $v>0$, the standard Fourier
series
\[
       \sum_{t\in\mathbb Z} A_t(v)e^{2\pi it(x+iv)}
\]
converges in $\mathscr H$ to $\widehat\Phi_{\rm cmp}(x+iv)$,
uniformly in $x$. The same convergence is locally uniform for
$v$ in compact subintervals of $(0,\infty)$.
\end{corollary}
\begin{proof}
Work in the Hilbert direct sum $\mathscr H$ of the finite
admissible coordinates and $H^1_\nu(X)$. The finite modular
inputs are smooth. Lemma~\ref{lem:KM-strong-parameter-regularity}
gives strong smoothness of the metric coordinate. The fixed
Hodge functional is continuous on $H^1$ by
Lemma~\ref{lem:Hodge-smooth-energy-continuity}, so the final
constant coordinate is strongly smooth as well.

Put $F_v(x)=\widehat\Phi_{\rm cmp}(x+iv)$ and
\[
 b_t(v)=\int_0^1 F_v(x)e^{-2\pi itx}\,dx.
\]
These are Bochner integrals. The already proved positive and
nonpositive class comparisons give
$b_t(v)=e^{-2\pi tv}A_t(v)$. For every compact interval
$I\Subset(0,\infty)$ and every integer $r\geq2$, integration
by parts in the Hilbert space gives
\[
 \|b_t(v)\|_{\mathscr H}
    \leq\frac{\sup_{x\in[0,1],\,v\in I}
                    \|\partial_x^rF_v(x)\|_{\mathscr H}}
                 {(2\pi|t|)^r}\qquad(t\ne0).
\]
The right side is finite by strong smoothness. Therefore, for $M\geq1$,
\[
 \sup_{x\in\mathbb R,\,v\in I}
 \left\|\sum_{|t|>M}b_t(v)e^{2\pi itx}\right\|_{\mathscr H}
       \leq C_{I,r}M^{1-r}.
\]
The uniformly convergent Fourier series equals $F_v$:
applying any continuous linear functional reduces this
assertion to uniqueness of scalar Fourier coefficients, and
such functionals separate points of $\mathscr H$. Substituting
the coefficient identification proves the stated standard
series convergence. The same proof applies to a fixed larger
period or a translation twist, with the corresponding
frequency lattice. No assertion about convergence of
uncancelled singular Green representatives is used.
\end{proof}

\begin{corollary}[All nonzero rank-one comparison functions vanish]
\label{cor:all-nonzero-rankone-comparison-zero}
For every $t\ne0$ and $v>0$ one has $d_t(v)=h_t(v)=0$ in the
fixed rank-one comparison notation. No division by a possibly
zero Fourier value $e_{N,t}(v,1/2)$ is made.
\end{corollary}
\begin{proof}
The complete degree comparison gives $d_t=0$ for every nonzero $t$.
For $h_t$, the positive proof's algebra applies to all nonzero
$t$ after the discriminant-height proposition and the direct
finite vertical comparison above. In detail, put
\[
 \rho_{p,t}=\frac{p^{-1}g_{p,t}(0)}{1-p^{-1}g_{p,t}(0)},\qquad
 V_{N,t}=\sum_{p\mid N}
      \left(\frac{n_p-1-n_pp}{p-1}-\rho_{p,t}\right)\log p.
\]
The same local density denominator is nonzero for negative $t$.
Indeed $Q=-ac-Nb^2$ represents any integral $t$ at
$(a,b,c)=(1,0,-t)$ with unit partial derivative in $c$.
The hyperbolic rank-four denominator for the pair $\langle t,N\rangle$
also has a smooth representation: with $Q=ab+cd$, take
$(1,t,0,0)$ and $(0,0,1,N)$. The derivatives of the two norms
and their mutual pairing have a unit triangular minor.
The exact local Whittaker-density relation therefore gives
$1-p^{-1}g_{p,t}(0)\ne0$, including at $p=2$.

Use the full discriminant height and the finite divisor of
$\Delta_N$ to obtain the identity of Fourier values
\[
 \varphi(N)\langle A_t(v),\widehat\lambda_N\rangle
  =C_N(1)\left[e'_{N,t}(v,1/2)
       +\{C_N'(1)/C_N(1)+V_{N,t}\}e_{N,t}(v,1/2)\right].
\]
This is a linear identity; it remains valid if the last Fourier
value is zero. The standard negative cusp divisors and Green
functions are $W_N$-invariant: this follows from the lattice
isometry, with $c\leftrightarrow N/c$ preserving the quotient
integer $\gcd(c,N/c)$ and the cusp coefficients. Hence
\[
 \langle A_t,H\rangle=-2\langle A_t,\widehat\lambda_N\rangle.
\]
Finally $C_N(1)=-\varphi(N)\psi(N)/24$, and the same elementary
logarithmic derivative calculation as in the positive proof gives
\[
 2\{C_N'(1)/C_N(1)+V_{N,t}\}
     =2+\log N+4\Lambda'(2)/\Lambda(2)
                                  -\sum_{p\mid N}\beta'_{p,t}(0).
\]
Substitution is precisely $h_t=0$. This computation only divides
by nonzero normalization and local density factors, never by
$e_{N,t}(v,1/2)$.
\end{proof}


\subsection{The middle orbit in the actual rank-zero Fourier coefficient}
\label{subsec:rankzero-middle-orbit}

Let \(E_N^{(1)}(\tau,\sigma)\) and
\(\mathcal E_N^{(2)}(z,s)\) have the full rational
parabolic summation convention, standard archimedean
weight-\(3/2\) section of value one, and finite sections
obtained from \(1_{\widehat L_N}\) and
\(1_{\widehat L_N^2}\). Set
\[
 \mathcal Z_y(w)=\sum_{\mathbf n\in\mathbb Z^2\setminus\{0\}}
                   (\mathbf n^ty\mathbf n)^{-w}
                 \quad(\Re w>1),\qquad
 d_0=\deg C_0(1).
\]
The Epstein series retains both signs.
All formulas below are meromorphic identities; at
\(w=1/2\) this is its analytic continuation.

Use the quadratic form
\[
 q_N(a,b,c)=-ac-Nb^2
\]
and the global additive character $\psi$ of $\mathbb A/\mathbb Q$
with $\psi_\infty(x)=e^{2\pi ix}$ and
$\psi_p(x)=e^{-2\pi i\{x\}_p}$, where $\{x\}_p$ is the
$p$-adic fractional part. Each finite character has conductor
$\mathbb Z_p$.
The Fourier transform is self-dual and
\(\omega_p(w)\phi(0)=\gamma_p(q_N)^{-1}\int\phi(x)\,dx\),
where \(w=\left(\begin{smallmatrix}0&-1\\1&0\end{smallmatrix}\right)\).
Here \(\gamma_p\) denotes the raw Weil index in this convention.
At infinity use the standard positive-definite weight-\(3/2\)
section of value one. These conventions give the full rational
parabolic Eisenstein series in
Lemma~\ref{lem:rankzero-actual-middle-orbit}.

\begin{lemma}[All-prime zero densities]\label{lem:all-prime-truncated-zero-densities}
For \(n\geq0\), \(\epsilon\in\mathbb Z_p^\times\), set
\[
 \alpha_{n,r}
  =p^{-2r}\#\{(a,b,c)\bmod p^r:
                      ac+\epsilon p^n b^2=0\bmod p^r\},
 \qquad \alpha_{n,0}=1.
\]
The limit \(\alpha_n=\lim_{r\to\infty}\alpha_{n,r}\)
exists and, at every prime including \(2\), is
\[
 \alpha_0=1+p^{-1},\qquad
 \alpha_1=2,\qquad
 \alpha_n=2(1-p^{-1})+\alpha_{n-2}
          =n+1-(n-1)p^{-1}\quad(n\geq2).
\]
\end{lemma}
\begin{proof}
For \(n=0\), the equation over \(\mathbb F_p\) has \(p^2\)
solutions: \(b=0\) gives \(2p-1\), and each of the \(p-1\)
nonzero \(b\)'s gives \(p-1\). Every nonzero zero has
\(a\ne0\) or \(c\ne0\), hence is smooth even at \(p=2\).
The smooth contribution to \(\alpha_{0,r}\), \(r\ge2\), is
\(1-p^{-2}\). The remaining solutions have all three
coordinates divisible by \(p\), and give
\[
 \alpha_{0,r}=1-p^{-2}+p^{-1}\alpha_{0,r-2}.
\]
In particular
\(\alpha_{0,2j}=\alpha_{0,2j+1}
=1+p^{-1}-p^{-j-1}\) for \(j\ge0\).

For \(n=1\), the smooth zeros modulo \(p\) are the
\(2p(p-1)\) zeros with \(a\ne0\) or \(c\ne0\).
If \(a,c\in p\mathbb Z_p\) and \(b\) is a unit, the expression
has valuation exactly one, so it cannot vanish modulo
\(p^r\) for \(r\ge2\). The remaining singular contribution
has all coordinates divisible by \(p\). Thus
\[
 \alpha_{1,r}=2(1-p^{-1})+p^{-1}\alpha_{1,r-2}
                   \qquad(r\ge2).
\]
For \(n\ge2\), the singular contribution is \(a,c\in p\mathbb Z_p\);
after writing \(a=pa'\), \(c=pc'\) and dividing the equation
by \(p^2\), its normalized density is
\(\alpha_{n-2,r-2}\). Hence
\[
 \alpha_{n,r}=2(1-p^{-1})+\alpha_{n-2,r-2}
                   \qquad(r\ge2).
\]
These recurrences prove convergence, exponential convergence
in \(r\) for fixed \(n\), and the claimed limits.
No step divides by \(2\).
\end{proof}

\begin{lemma}[The normalized genus-one zero Whittaker germ]
\label{lem:middle-seed-local-normalization}
Let \(w_N(s)=W^{(1)}_0(1,s+\tfrac12,\Phi_{1,N})\).
There are functions \(a_{p,n}(s)\), holomorphic near zero,
such that
\[
 w_N(s)=-
 \frac{s\Lambda(2s+1)}{(s+1)\Lambda(2s+2)}
       \prod_{p^n\Vert N}a_{p,n}(s),
 \qquad
 a_{p,n}(0)=p^{-n/2}
             \frac{n+1-(n-1)/p}{1+1/p}.
\]
Consequently \(w_N\) is holomorphic at zero and
\[
 w_N(0)=-\frac3\pi
       \prod_{p^n\Vert N}
         p^{-n/2}\frac{n+1-(n-1)/p}{1+1/p}.
\]
\end{lemma}
\begin{proof}
At a finite prime $p$, write $N=p^n\epsilon$ with
$\epsilon\in\mathbb Z_p^\times$. The polar Gram determinant
of $q_N$ is $2N$, so the
self-dual volume of \(L_{N,p}\) is \(|2N|_p^{1/2}\).
Fourier inversion over \(p^{-r}\mathbb Z_p\) therefore gives
\[
 \int_{p^{-r}\mathbb Z_p}
      \Phi_{1,N,p}(wn(t),\tfrac12)\,dt
    =\gamma_p(q_N)^{-1}|2N|_p^{1/2}\alpha_{n,r}.
\]
The integration domain $p^{-r}\mathbb Z_p$ has volume $p^r$;
thus the normalization here is exactly the density
\(p^{-2r}\#\{q_N=0\bmod p^r\}\), rather than twice that density.

On the annulus \(|t|_p=p^r>1\), changing the induction parameter
from \(1/2\) to \(s+1/2\) multiplies the standard section by
\(p^{-rs}\). It is unchanged on \(\mathbb Z_p\).
Thus the local zero Whittaker germ is
\[
 \gamma_p(q_N)^{-1}|2N|_p^{1/2} Z_{p,n}(s),\qquad
 Z_{p,n}(s)=1+
    \sum_{r\ge1}p^{-rs}(\alpha_{n,r}-\alpha_{n,r-1}).
\]
Lemma~\ref{lem:all-prime-truncated-zero-densities} shows that this series is holomorphic
for \(s\) in a neighborhood of zero, and
\(Z_{p,n}(0)=\alpha_n\).
For \(n=0\), the explicit density sequence gives
\[
 Z_{p,0}(s)
   =1+(1-p^{-1})\sum_{j\ge1}p^{-j(1+2s)}
   =\frac{\zeta_p(2s+1)}{\zeta_p(2s+2)}.
\]
This also holds at \(p=2\). Define
\[
 a_{p,n}(s)
   =p^{-n/2}\frac{Z_{p,n}(s)}{Z_{p,0}(s)}.
\]
The denominator is nonzero at zero and
\(a_{p,n}(0)\) is as stated.

For the real section there is a direct beta-integral check.
The branch of \((t+i)^{-3/2}\) has argument in \((0,\pi)\), and
\[
 \begin{split}
 W^{(1)}_{0,\infty}(1,s+\tfrac12)
 &=\int_{\mathbb R}(t+i)^{-3/2}(1+t^2)^{-s/2}\,dt\\
 &=e^{-3\pi i/4}
   \frac{\pi\Gamma(s+\tfrac12)}
        {2^{s-1/2}\Gamma(s/2+\tfrac32)\Gamma(s/2)}\\
 &=e^{-3\pi i/4}\sqrt2\,
    \frac{s\zeta_\infty(2s+1)}
         {(s+1)\zeta_\infty(2s+2)}.
 \end{split}
\]
The substitution \(t=\tan\theta\) proves the middle equality;
duplication of the gamma function proves the last one.
Initially these integrals are evaluated in an absolute
half-plane, then continued meromorphically.

It remains to retain the global phase rather than discard it.
The coherent rational form \(q_N\) has signature \((1,2)\);
its real raw Weil index is \(e^{-\pi i/4}\). The Weil-index
product formula gives
\(\prod_{p<\infty}\gamma_p(q_N)^{-1}=e^{-\pi i/4}\).
The archimedean positive-definite section contributes
\(e^{-3\pi i/4}\), so their product is \(-1\).
Also \(\prod_{p<\infty}|2|_p^{1/2}=2^{-1/2}\),
which cancels the displayed \(\sqrt2\).
Multiplication in the initial Euler-product half-plane now
gives the asserted expression for \(w_N(s)\).

Finally
\[
 \lim_{s\to0}s\Lambda(1+2s)=\tfrac12,
 \qquad \Lambda(2)=\pi/6.
\]
The quotient is therefore \(3/\pi\), with the negative
phase already fixed above. This proves both holomorphy
and the stated value, without a genus-one scalar/vector
central-factor convention.
\end{proof}

\begin{lemma}[Cusp count and the coefficient \(d_0\)]
\label{lem:cusp-sum-middle-seed}
With the arithmetic-stack convention and the standard
isotropic boundary divisor,
\[
 d_0=\frac1{4\pi\sqrt N}
      \prod_{p^n\Vert N}
            p^{n-1}\bigl((n+1)p-(n-1)\bigr).
\]
In particular the preceding seed is
\(w_N(0)=-12d_0/\psi(N)\).
\end{lemma}
\begin{proof}
Choose cusp representatives \(a/c\) with \(c\mid N\).
There are \(\varphi(\gcd(c,N/c))\) effective geometric cusps
with this denominator. In the existing analytic/orbifold cusp
sum defining \(d_0\), the effective cusp integral carries the
generic ineffective-center weight \(1/2\). This is the
normalization obtained by integration on the actual curve,
or on a local scheme chart followed by its stack degree.
It is not a separate assertion that the residual gerbe of
every compactified cusp has automorphism group of order two:
if a ramified stack parameter is used, its parameter index
must be retained in the divisor coefficient. The compact
Green equation fixes their combined contribution. Closed
arithmetic cusps are summed with their field degrees before
summing complex embeddings. The effective cusp width and
isotropic step are
\[
 \kappa_c=\frac N{\gcd(N,c^2)},\qquad \beta_c=\frac1c.
\]
Thus
\[
 d_0=\frac1{4\pi\sqrt N}
      \sum_{c\mid N}\varphi(\gcd(c,N/c))
                         \frac{Nc}{\gcd(N,c^2)}.
\]
This sum is multiplicative. For \(N=p^n\), its two
endpoints contribute \(2p^n\); each of the \(n-1\)
interior denominators contributes \(p^{n-1}(p-1)\).
It is therefore
\(p^{n-1}((n+1)p-(n-1))\).
Comparing with
\(\psi(p^n)=p^{n-1}(p+1)\) proves the assertion.
\end{proof}

\begin{lemma}[The actual rank-one Bruhat orbit and its seed]
\label{lem:rankzero-actual-middle-orbit}
Let
\[
 w_N(s)=W_0^{(1)}(1,s+\tfrac12,\Phi_{1,N})
\]
be the global genus-one zero Whittaker germ with
the same self-dual additive measures as the
Eisenstein series. The classical middle-orbit
contribution to the genus-two zero coefficient is
\begin{equation}\label{eq:actual-middle-orbit}
 \mathcal M_N(y,s)
   =(\det y)^{s/2}w_N(s)
          \frac{\mathcal Z_y(s+\tfrac12)}
                              {2\zeta(2s+1)}.
\end{equation}
The germ \(w_N\) is holomorphic at zero and satisfies
\begin{equation}\label{eq:middle-seed-from-DY}
                 w_N(0)=-\frac{12d_0}{\psi(N)}.
\end{equation}
Consequently
\begin{equation}\label{eq:middle-derivative-Epstein}
 \left.\partial_s\mathcal M_N(y,s)\right|_{s=0}
                  =-\frac{12d_0}{\psi(N)}
                                      \mathcal Z_y(1/2).
\end{equation}
\end{lemma}

\begin{proof}
We first check the discrete orbit, independently of
any primitive-row convention in a later formula.
Let \(U\) be the standard rational Lagrangian in
the four-dimensional symplectic space. Left
cosets of the full Siegel parabolic identify
with rational Lagrangians \(L=g^{-1}U\).
The middle Bruhat cell is characterized by
\(\dim(L\cap U)=1\). Its discrete invariant
under the Siegel unipotent, which fixes \(U\)
pointwise, is the rational line \(L\cap U\).
The Levi acts transitively on these lines,
and the stabilizer is the full line
parabolic. Thus the residual discrete
sum is indexed by \(\mathbb P^1(\mathbb Q)\).
In primitive integral bottom rows its
index is
\[
 \{(c,d)\in\mathbb Z^2:\gcd(c,d)=1\}/\{\pm1\}.
\]
There is no choice of orientation.
Explicitly, representatives \(\gamma\)
and \(-\gamma\) give the same original
Siegel-parabolic coset: \(m(-I_2)=-I_4\)
is rational central, belongs to that
parabolic, and commutes with the rank-one
Weyl representative. The full rational
splitting makes this an equality of
adelic metaplectic cosets as well.
No factor from the kernel of the cover
is being added at this step.

Unfold the affine unipotent variables
for each fixed line. Their compact
\(\mathbb Q\backslash\mathbb A\) quotients
have volume one, leaving the genus-one
zero Whittaker functional. This is the
orbit formula underlying \cite[equations~(238)--(239)]{SSY2022}.
The pullback section is precisely
\[
 \eta^*\Phi_{2,N}(s)=\Phi_{1,N}(s+\tfrac12).
\]
At finite places this follows either
from \cite[Lemma~4.12]{SSY2022} or directly from the
Rallis map: restricting the tuple
characteristic function to \((0,x)\)
gives \(1_{L_N}(x)\). At infinity,
restriction has the same standard
weight-\(3/2\) compact type and value
one at the identity. These data determine
the standard section, so no scalar
remains to be chosen. Finite Levi
matrices in \(\mathrm{SL}_2(\mathbb Z_p)\)
preserve \(1_{L_{N,p}^2}\). Indeed, on the canonical lift of the
determinant-one Levi subgroup, the Weil action substitutes an
invertible integral linear combination for the two lattice vectors.
The inverse substitution is also integral, so it preserves their
characteristic function. The restriction of the standard section
to the compact subgroup is independent of the spectral parameter.
This invariance therefore holds for the whole family of sections,
including at the dyadic prime, without a self-duality assumption on
the lattice. The finite factors are consequently independent of
the primitive row.

For completeness the archimedean
dependence can be obtained without
choosing an absolute Weil-index phase.
Write \(y=aa^t\), \(\det a>0\), and
apply a real QR decomposition to
\(\gamma a\). If \(r^2=(c,d)y(c,d)^t\),
its upper triangular factor has
diagonal \((\det a/r,r)\). The first
diagonal coordinate contributes
\((\det a/r)^{s+3/2}\); the genus-one
zero Whittaker scaling law contributes
\(r^{1/2-s}\). The orthogonal factor
lies in the real \(\mathrm{SO}_2\)
inside the rank-two compact subgroup
with determinant one, and hence acts
trivially on its scalar determinant
type. The upper triangular unipotent
is in the stabilizer of this zero
Whittaker functional. Multiplying
the two diagonal factors gives
\[
 (\det y)^{s/2+3/4}
                 ((c,d)y(c,d)^t)^{-s-1/2}.
\]
After the classical factor
\((\det y)^{-3/4}\), the seed is
therefore \((\det y)^{s/2}w_N(s)\).

Grouping nonzero integral vectors
as a positive integer multiple of
an unoriented primitive vector gives
\[
 \mathcal Z_y(w)
     =2\zeta(2w)
        \sum_{(c,d)\ {\rm primitive}/\{\pm1\}}
                            ((c,d)y(c,d)^t)^{-w}.
\]
This proves \eqref{eq:actual-middle-orbit}
in its initial convergence region,
and then meromorphically.

Lemmas~\ref{lem:middle-seed-local-normalization} and
\ref{lem:cusp-sum-middle-seed} give the seed value
\eqref{eq:middle-seed-from-DY} and its holomorphy, with the
phase and every bad-prime factor included. This also agrees with
the genus-one degree identity: its zero coefficient is
$-\psi(N)/12+d_0/\sqrt v$, whereas the zero coefficient of
$E_N^{(1)}(\tau,1/2)$ is $1+w_N(0)v^{-1/2}$. The analytic
regularity is also consistent with the normalized theta-lift
expression in \cite[Theorem~1.4]{DY1}.

Finally \(\mathcal Z_y(w)\) is regular
at \(w=1/2\), whereas
\[
                 \frac1{2\zeta(1+2s)}=s+O(s^2).
\]
Differentiating \eqref{eq:actual-middle-orbit}
now proves \eqref{eq:middle-derivative-Epstein}.
Neither the derivative of \((\det y)^{s/2}\)
nor that of \(w_N(s)\) contributes.
\end{proof}

\begin{proposition}[The complete rank-zero analytic comparison]
\label{prop:actual-rankzero-ASW-comparison}
For the actual full-parabolic Eisenstein series just
fixed,
\begin{equation}\label{eq:actual-rankzero-ASW}
 \frac{\psi(N)}{24}
       \left.\partial_s C_{N,0_2}(y,s)\right|_{s=0}
       =\mathcal A^{\rm std}(0_2,y)
                          -\frac{d_0}{2}\mathcal Z_y(1/2).
\end{equation}
Thus the specified coefficient correction
\(\mathcal A^{\rm corr}(0_2,y)\) is the actual
rank-zero Eisenstein coefficient in this normalization.
Equivalently, the base-comparison function
\(b_N(y)=\mathcal A^{\rm std}(0_2,y)
                -\mathcal A^{\rm SW}(0_2,y)\)
is \(d_0\mathcal Z_y(1/2)/2\).
\end{proposition}

\begin{proof}
We check the other two orbits explicitly,
without using any assertion that the middle
one vanishes. Put
\[
 A_{p,n}(s)=\frac{p^{-n}}{1+p^{-1-s}}
  \left\{\sum_{j=0}^{n}p^{j(1-s)}
             -p^{-2s}\sum_{j=0}^{n-2}p^{j(1-s)}\right\},
                    \qquad n\geq1,
\]
where the last sum is zero for \(n=1\).
The open-orbit zero Whittaker function,
with the classical factor removed, is
\[
 \mathcal W_N(y,s)
  =(\det y)^{-s/2}
       \frac{(s-1)\Lambda(2s-1)}
                    {(s+1)\Lambda(2s+2)}
           \prod_{p\mid N}A_{p,v_p(N)}(s).
\]
These are the bad-prime local factors in
\cite[Lemma~4.1.1]{Zhu}, obtained directly from
the density limit in its proof. At an
unramified prime the relative factor is
\(1\). In particular a separate
\((1+p^{-1-s})^{-1}\) must not be
inserted at every good prime.

For precision, with \(X=p^{-s}\)
the displayed density limit multiplied
by its normalizing factor is
\[
 p^{-n}(1-p^{-1}X)
 \frac{1-X^2+(p^{n-1}-p^{n+1})X^{n+1}}
                     {(1-pX)(1-pX^2)}.
\]
Dividing by the unramified factor
\((1-p^{-2}X^2)/(1-pX^2)\) gives
\[
 p^{-n}
 \frac{1-X^2+(p^{n-1}-p^{n+1})X^{n+1}}
                       {(1-pX)(1+p^{-1}X)},
\]
which is exactly \(A_{p,n}(s)\)
for \(n\geq1\) and is \(1\) for
\(n=0\). This algebra also fixes
the good-prime normalization.

Each bad factor has \(A_{p,n}(0)=1\) and
\[
 A'_{p,n}(0)=L_{p,n},\qquad
 L_{p,n}=
   \frac{-np^{n+1}+2p^n+np^{n-1}-2}
                       {p^{n-1}(p^2-1)}\log p .
\]
The identity orbit is \((\det y)^{s/2}\).
The functional equation of \(\Lambda\)
gives \(\mathcal W_N(y,0)=-1\),
and direct differentiation therefore gives
\begin{equation}\label{eq:rankzero-two-extreme-orbits}
 \left.\partial_s\{(\det y)^{s/2}
                              +\mathcal W_N(y,s)\}\right|_{s=0}
  =\log\det y+2-4\frac{\Lambda'(-1)}{\Lambda(-1)}
                                -\sum_{p\mid N}L_{p,v_p(N)}.
\end{equation}

The arithmetic comparison of these
two extreme orbits is independent
of the analytic middle orbit. The
actual Hodge self-intersection and
actual vertical intersection matrix
give
\[
 \widehat\lambda_N^2
       =\frac{\psi(N)}{24}
                    \left(\tfrac12-\frac{\Lambda'(-1)}{\Lambda(-1)}\right),
 \qquad
 (W_N^*\widehat\lambda_N-\widehat\lambda_N)^2
       =\frac{\psi(N)}{24}\sum_{p\mid N}L_{p,v_p(N)}.
\]
Here the first identity is the
K\"uhn--Bost Hodge self-intersection
in the fixed Petersson normalization
(also \cite[Example~6.2.5]{Zhu}), and its
all-level form follows by the degree
\(\psi(N)\) source-forgetful projection
formula. The second follows from
the actual component intersection
matrix \cite[Corollary~6.3.2]{Zhu} and the
antisymmetry of each local difference
under \(W_N\). Neither identity
uses Zhu Proposition 4.1.2 or its
rank-zero arithmetic Siegel--Weil
conclusion.
Since \(H=-\widehat\lambda_N-W_N^*\widehat\lambda_N\),
bilinearity and \(W_N\)-invariance give
\[
 H^2=4\widehat\lambda_N^2
            -(W_N^*\widehat\lambda_N-\widehat\lambda_N)^2.
\]
Moreover \(\deg H=-\psi(N)/12\).
Thus multiplication of
\eqref{eq:rankzero-two-extreme-orbits}
by \(\psi(N)/24\) gives exactly
\[
 \widehat{\deg}(H^2)-\tfrac12\log\det y\,\deg H
                         =\mathcal A^{\rm std}(0_2,y).
\]
There are precisely the three Bruhat
ranks \(0,1,2\). Adding the middle
derivative from
\eqref{eq:middle-derivative-Epstein}
proves \eqref{eq:actual-rankzero-ASW}.
\end{proof}

\paragraph{Normalization checks and scope.}
The factor \(1/2\) in the primitive-row
sum comes from an unoriented rational
line. It is independent of stack
integration and of arithmetic degree.
The sign of the seed is fixed by the
full genus-one theta integral, not
by choosing an archimedean phase to
match the desired answer. For \(N=1\)
one obtains \(w_1(0)=-3/\pi\) and
\(\partial_s\mathcal M_1(y,0)
                    =-3\mathcal Z_y(1/2)/\pi\).
These are checks of the general proof.
The conclusion concerns the precisely
specified full-parabolic sections;
it is not an assertion about an
unverified re-normalization of a
printed formula.

\subsection{The full standard arithmetic doubling identity}
\label{subsec:full-standard-arithmetic-doubling}

\begin{proof}[Proof of Theorem~\ref{thm:full-standard-arithmetic-doubling}]
The equality of the two genus-one families, including all
nonpositive coefficients and convergence in the stated topology,
is Theorem~\ref{thm:nonpositive-standard-completion} together with
positive compatibility and
Corollary~\ref{cor:standard-family-strong-Fourier}.

We first establish \eqref{eq:all-coefficients-ASW}. For nonsingular
$T$ this is Proposition~\ref{prop:all-level-nonsingular-ASW}.
For nonzero rank-one $T$, apply
Proposition~\ref{prop:singular-ASW-comparison} and
Corollary~\ref{cor:all-nonzero-rankone-comparison-zero}.
For $T=0_2$, use
Proposition~\ref{prop:actual-rankzero-ASW-comparison}.
These three cases exhaust the indices; the last includes the
middle orbit, so no separate zero-mode residual remains.

Fix $v_1,v_2>0$. By strong Fourier convergence and continuity
of the arithmetic product, the coefficient of $q_1^m q_2^n$
on the left of \eqref{eq:full-standard-arithmetic-doubling} is
the product of the corresponding standard classes.
Corollary~\ref{cor:corrected-standard-support} gives this product
as the convergent sum
\[
 \sum_{r\in\mathbb Z}\mathcal A^{\rm corr}
    \left(\begin{pmatrix}m&r/2\\r/2&n\end{pmatrix},
                         \operatorname{diag}(v_1,v_2)\right).
\]
For the positive common-component case, its proof uses
Theorem~\ref{thm:standard-positive-square-product}, including
the complete finite adjunction and Green-current contributions.
For the double-zero case, its correction is the explicitly
computed Epstein term. Thus the product formula used here
includes all integer pairs $(m,n)$.
Replace each summand by its equal value in
\eqref{eq:all-coefficients-ASW}.

These sums are exactly the diagonal Fourier coefficients of
the analytic pullback. To justify evaluating the off-diagonal
real variable at zero, observe that for fixed positive definite
$y$ the Eisenstein derivative is a smooth periodic function
of all three independent entries of its real symmetric part.
Repeated integration by parts bounds its full Fourier
coefficients faster than any polynomial. Its Fourier series
therefore converges absolutely on that three-dimensional torus,
and restriction and summation over $r$ commute. The same bounds
are locally uniform in $y$.
Both sides of \eqref{eq:full-standard-arithmetic-doubling} are
continuous on the two real-variable torus and have the same
Fourier coefficients for every $(m,n)\in\mathbb Z^2$.
Fourier uniqueness proves the pointwise identity, using the
nonpositive comparisons as well as the positive ones.

Finally the complete arithmetic kernel has joint moderate
growth, and the cusp forms decay exponentially at every cusp.
Corollary~\ref{cor:completion-full-lift-MW} supplies the full
Bochner lifts and permits both integrals to pass through the
continuous arithmetic product. This proves the paired assertion,
with the same measures and arithmetic degree convention.
The Hermitian form follows by complex conjugating the second
arithmetic factor and reflecting its analytic variable as in
\eqref{eq:red-arithmetic-kernel}.
\end{proof}

\begin{corollary}[The height of the full standard lift]
\label{cor:full-standard-lift-height}
The full standard cusp lift equals the admissible
Mordell--Weil lift. Consequently the scalar height theorem
\ref{thm:main-projected-height}, with its specified Jacobi
newform line and factor $\|Tf\|_\Gamma^2$, holds for the
full standard lift as well. No replacement of that trace norm
by an uncorrected scalar norm is made.
\end{corollary}
\begin{proof}
Use the equality of the two complete families in the theorem
and Corollary~\ref{cor:completion-full-lift-MW}, then apply
the stated projected-height theorem.
\end{proof}

\subsection{Height radicals and minimal scalar representatives}
\label{subsec:height-radical}

The next result describes the degeneracy of the full lift on a
specified newform line. It is a consequence of the height formula
and trace adjunction, with attainability kept explicit.

\begin{corollary}[The height radical and the scalar new-line quotient]
\label{cor:newline-height-radical}
Retain the hypotheses on $G$ and $g_G$ of
Theorem~\ref{thm:BY-newline-height}. Put
\[
 \mathcal E=S_{3/2}(H_N,\chi_N),\qquad
 \mathcal E_G=T^{-1}(\mathbb Cg_G),\qquad
 a_G=\frac{2\sqrt N}{\pi}L'(G,1),
\]
and let $\widehat\theta_{\rm std}$ denote the full standard cusp lift.
On $\mathcal E_G$ define the Hermitian form, linear in its first
argument, by
\[
 B_G(f,h)=H_{\rm NT}\bigl(\widehat\theta_{\rm std}(h),
                              \widehat\theta_{\rm std}(f)\bigr).
\]
Then
\begin{equation}\label{eq:newline-cross-height}
 B_G(f,h)=a_G\langle Tf,Th\rangle_\Gamma
         =a_G\langle f,RTh\rangle_{H_N}.
\end{equation}
In particular,
\begin{equation}\label{eq:newline-height-radical}
 \operatorname{rad}B_G
   =\ker(\widehat\theta_{\rm std}|_{\mathcal E_G})
   =\begin{cases}
       \ker T,&L'(G,1)\ne0,\\
       \mathcal E_G,&L'(G,1)=0.
     \end{cases}
\end{equation}
Thus the full lift vanishes on $\mathcal E_G$ when $L'(G,1)=0$.
For every $G$ there is a canonical linear isomorphism
\begin{equation}\label{eq:newline-quotient-image}
 \mathcal E_G/\ker T
    \xrightarrow{\ \overline T\ }
       \operatorname{im}T\cap\mathbb Cg_G.
\end{equation}
Its dimension is either zero or one. It is one precisely when
$g_G\in\operatorname{im}T$, equivalently when $g_G$ takes values
in the cyclic coefficient space of $e_0$.
If this attainability condition holds and $L'(G,1)\ne0$, then
$a_G>0$, and the induced full lift is a conjugate-linear
isomorphism from this one-dimensional quotient onto a
Mordell--Weil line. For the quotient product
\[
 \langle[f],[h]\rangle_T=\langle Tf,Th\rangle_\Gamma,
\]
the lift scaled by $a_G^{-1/2}$ is a conjugate-linear isometry.

If $g_G\in\operatorname{im}T$, let $S_N^\dagger$ be the inverse
of the coefficient frame on its cyclic subspace and zero on
its orthogonal complement. Then
\begin{equation}\label{eq:newline-minimal-representative}
 f_G=R(S_N^\dagger g_G),\qquad
 Tf_G=g_G,\qquad
 \|f_G\|_{H_N}^2=\langle S_N^\dagger g_G,g_G\rangle_\Gamma.
\end{equation}
Every $f\in\mathcal E_G$ is uniquely $f=k+b f_G$, with
$k\in\ker T$ and $b\in\mathbb C$, and
\[
 \widehat\theta_{\rm std}(f)
      =\overline b\,\widehat\theta_{\rm std}(f_G).
\]
The vector $b f_G$ is the unique representative of $[f]$ of
minimal scalar Petersson norm.
\end{corollary}

\begin{proof}
Corollary~\ref{cor:full-standard-lift-height} identifies the full
lift with the Mordell--Weil lift. The latter factors through
$T$ by \eqref{eq:red-MW-adjunction} and is conjugate-linear in
the test form. Theorem~\ref{thm:main-projected-height} gives
$B_G(f,f)=a_G\|Tf\|_\Gamma^2$ for every $f\in\mathcal E_G$.
The constant $a_G$ is real by \eqref{eq:BY-factor-and-height}.
Hermitian polarization therefore gives the first equality in
\eqref{eq:newline-cross-height}; Lemma~\ref{lem:red-adjunction}
gives the second. The reversed order in the definition of
$B_G$ accounts for the conjugate-linearity of the lift.

The trace Heegner classes lie over a fixed number field $K$,
as in the construction preceding
Corollary~\ref{cor:completion-full-lift-MW}. Their admissible
Mordell--Weil span is contained in
$(J_0(N)(K)/J_0(N)(K)_{\rm tors})\otimes_\mathbb Z\mathbb C$.
The Hermitian N\'eron--Tate pairing is positive definite on
this finite-dimensional space. Consequently a lift has zero
self-height if and only if it is zero. Formula
\eqref{eq:newline-cross-height} now proves
\eqref{eq:newline-height-radical}, including the case
$L'(G,1)=0$. The first isomorphism theorem gives
\eqref{eq:newline-quotient-image}. Corollary~\ref{cor:all-level-norm}
identifies $\operatorname{im}T$ with the cyclic-valued vector
cusp forms. If $g_G$ is attainable and $a_G\ne0$, applying
the height formula to a preimage of $g_G$ shows $a_G>0$;
the asserted conjugate-linear isomorphism and isometry
follow from \eqref{eq:newline-cross-height}.

Finally $TR=S_N$ and
$S_NS_N^\dagger g_G=g_G$, so the displayed $f_G$ maps to $g_G$.
It lies in $\operatorname{im}R=(\ker T)^\perp$ by adjunction.
The adjunction identity also gives
\[
 \|R(S_N^\dagger g_G)\|_{H_N}^2
   =\langle S_N^\dagger g_G,TR(S_N^\dagger g_G)\rangle_\Gamma
   =\langle S_N^\dagger g_G,g_G\rangle_\Gamma.
\]
Every preimage of $b g_G$ differs from $b f_G$ by $\ker T$.
The resulting orthogonal decomposition yields
\[
 \|k+b f_G\|_{H_N}^2
      =\|k\|_{H_N}^2+|b|^2\|f_G\|_{H_N}^2,
\]
which proves uniqueness and minimality. Factoring the
conjugate-linear lift through $T$ proves its stated value.
\end{proof}

\begin{remark}
The conclusion that the quotient is one-dimensional requires
attainability by the specified scalar trace.
When $\operatorname{im}T\cap\mathbb Cg_G=0$, one has
$\mathcal E_G=\ker T$ and the restricted lift is zero,
regardless of $L'(G,1)$. A zero restricted lift then gives
no information about this derivative.
\end{remark}
\subsection{The cuspidal Eisenstein kernel as a height Gram form}
\label{subsec:cuspidal-height-Gram}

The full doubling identity identifies the double cuspidal
projection of the reflected Eisenstein derivative with a height
Gram kernel. We extend the N\'eron--Tate pairing Hermitianly,
linearly in its first argument, and retain the positive constant
\(\eta\) in \eqref{eq:red-FH-normalization}. This fixes the sign
of the Gram form. The image of the full lift determines its rank.

Put
\[
 S=S_{3/2}(H_N,\chi_N),\qquad
 \Theta(f)=\widehat\theta_{\rm std}(f),\qquad
 c_{\rm ht}=\frac{24\eta}{\psi(N)}>0.
\]
The lift \(\Theta\) is conjugate-linear. Write
\[
 \mathscr K_N(z,w)
   =\left.\partial_s\mathcal E_N
              (\operatorname{diag}(z,-\overline w),s)\right|_{s=0}
\]
for the reflected kernel with exactly the
section, real structure, and measures of
\eqref{eq:red-arithmetic-kernel}. Define
\begin{equation}\label{eq:cuspidal-Eis-Gram-form}
 \mathscr D_N(f,h)
  =\int\!\!\int \mathscr K_N(z,w)
       \overline{f(z)}h(w)
          y_z^{3/2}y_w^{3/2}\,d\mu(z)d\mu(w).
\end{equation}
The integrals are over \(H_N\backslash\mathbb H\).
This convention is conjugate-linear in \(f\)
and linear in \(h\). To use the manuscript's
linear-first convention on the scalar space,
one instead uses \(\mathscr D_N(h,f)\).

\begin{corollary}[Sign, radical, and arithmetic rank]
\label{cor:cuspidal-Eis-height-Gram}
For all \(f,h\in S\),
\begin{equation}\label{eq:cuspidal-Eis-NT-Gram}
 \mathscr D_N(f,h)
       =-\frac{24\eta}{\psi(N)}
                       H_{\rm NT}(\Theta(f),\Theta(h)).
\end{equation}
Both sides retain the conjugate-linearity of the lift.

Suppose the lift values belong to the admissible realization of
\[
 V_K=(J_0(N)(K)/J_0(N)(K)_{\rm tors})
                                      \otimes_{\mathbb Z}\mathbb C
\]
for one number field \(K\), with the usual positive
N\'eron--Tate pairing for the fixed polarization, in the
normalization of \eqref{eq:red-FH-normalization}.
The construction preceding
Corollary~\ref{cor:completion-full-lift-MW}
uses such a fixed field of definition.
Then \(\mathscr D_N\) is Hermitian negative semidefinite,
and
\begin{align}
 \mathscr D_N(f,f)=0
      &\ \Longleftrightarrow\ \Theta(f)=0,\\
 \operatorname{rad}\mathscr D_N
      &=\ker\Theta,\label{eq:cuspidal-Eis-radical}\\
 \operatorname{rank}\mathscr D_N
      &=\dim_{\mathbb C}\operatorname{im}\Theta
       \ \leq\ \operatorname{rank}J_0(N)(K).
                    \label{eq:cuspidal-Eis-rank}
\end{align}
In particular, if \(J_0(N)(K)\) has rank zero,
the entire double cuspidal pairing vanishes.
For arbitrary \(f,h\), one also has
\[
 |\mathscr D_N(f,h)|^2
                 \leq\mathscr D_N(f,f)\mathscr D_N(h,h),
\]
with equality precisely when \(\Theta(f)\) and
\(\Theta(h)\) are complex linearly dependent,
including the case of a zero lift.
\end{corollary}

\begin{proof}
All standard arithmetic coefficients are real
arithmetic classes, so
\(\widehat\Phi_{\rm std}(-\overline w)
                   =\overline{\widehat\Phi_{\rm std}(w)}\).
Reflecting the second variable in
\eqref{eq:full-standard-arithmetic-doubling} gives
\[
 \frac{\psi(N)}{24}\mathscr K_N(z,w)
   =H_{\rm GS}(\widehat\Phi_{\rm std}(z),
                                  \widehat\Phi_{\rm std}(w)).
\]
The full theorem and
Corollary~\ref{cor:completion-full-lift-MW}
give absolute convergence and the Bochner/Fubini
interchange. Double pairing therefore gives
\[
 \frac{\psi(N)}{24}\mathscr D_N(f,h)
              =H_{\rm GS}(\Theta(f),\Theta(h)).
\]
The full lifts are in the admissible
Mordell--Weil subspace. Equation
\eqref{eq:red-FH-normalization} now gives
\eqref{eq:cuspidal-Eis-NT-Gram}, with the stated sign and height normalization.

On \(V_K\) the Hermitian N\'eron--Tate
pairing is positive definite after tensoring
away torsion. The first equivalence and
the radical formula follow. The kernel
of the conjugate-linear map \(\Theta\)
is a complex subspace; equivalently,
view \(\Theta\) as a linear map from
the conjugate complex space \(\overline S\).
The first isomorphism theorem gives the
rank identity and the dimension bound.
Finally Hermitian Cauchy--Schwarz on
\(V_K\) gives the displayed inequality
and its equality criterion.
\end{proof}

\begin{corollary}[Finite-rank double cuspidal projection]
\label{cor:cuspidal-Eis-finite-rank-kernel}
Retain the fixed-field and positive-height
hypotheses of Corollary~\ref{cor:cuspidal-Eis-height-Gram}.
Let \(u_1,\ldots,u_d\) be a Petersson
orthonormal basis of \(S\). Define the
double cuspidal projection by its finite
coefficient expansion
\begin{equation}\label{eq:finite-cusp-projection-definition}
 \mathscr K_N^{\rm cusp}(z,w)
    =\sum_{a,b=1}^d
          \mathscr D_N(u_a,u_b)u_a(z)\overline{u_b(w)}.
\end{equation}
This definition uses the convergent cusp
pairings and does not assume that the
unprojected kernel is globally \(L^2\).
It is independent of the chosen
orthonormal basis.

Let \(W=\operatorname{im}\Theta\), of
dimension \(r_\Theta\), and choose a
N\'eron--Tate orthonormal basis
\(e_1,\ldots,e_{r_\Theta}\) of \(W\).
There are unique cusp forms
\(G_1,\ldots,G_{r_\Theta}\in S\) such that
\[
 H_{\rm NT}(\Theta(f),e_j)
                               =\langle G_j,f\rangle_{H_N}
                  \qquad(f\in S).
\]
They are linearly independent, and
\begin{equation}\label{eq:finite-cusp-height-feature-kernel}
 \mathscr K_N^{\rm cusp}(z,w)
    =-\frac{24\eta}{\psi(N)}
         \sum_{j=1}^{r_\Theta}G_j(z)\overline{G_j(w)}.
\end{equation}
Consequently its Petersson integral operator on \(S\) is
self-adjoint, negative semidefinite, and has
rank exactly \(r_\Theta\) and kernel \(\ker\Theta\).

Alternatively, if the whole finite Mordell--Weil
coefficient span is written in an orthonormal
basis as
\[
 Y(\tau)=R\widehat\Phi_{\rm MW}(\tau)
                           =\sum_{j=1}^{r}F_j(\tau)e_j,
\]
then one may use
\[
 \mathscr K_N^{\rm cusp}(z,w)
   =-\frac{24\eta}{\psi(N)}
      \sum_{j=1}^{r}(P_SF_j)(z)\overline{(P_SF_j)(w)},
 \quad
 P_SF_j=\sum_{a=1}^d\langle F_j,u_a\rangle_{H_N}u_a.
\]
These pairings converge by the growth statements
already proved in the manuscript. No assertion
that \(F_j\) is cuspidal at every cusp is
needed; when that assertion is independently
available, \(P_SF_j=F_j\).
\end{corollary}

\begin{proof}
The coordinate functional
\(f\mapsto H_{\rm NT}(\Theta(f),e_j)\)
is conjugate-linear and continuous on
the finite-dimensional Petersson space
\(S\). Riesz representation in the
linear-first convention gives \(G_j\)
with exactly the displayed order.
Thus
\[
 \Theta(f)=\sum_j\langle G_j,f\rangle_{H_N}e_j.
\]
Substitution in
\eqref{eq:cuspidal-Eis-NT-Gram}
and then in
\eqref{eq:finite-cusp-projection-definition}
proves \eqref{eq:finite-cusp-height-feature-kernel}.
The coefficient identity also proves
basis independence.

Since \(\Theta(S)=W\), its coordinate
functionals are linearly independent.
Their Riesz representatives \(G_j\)
are therefore independent. More
explicitly, if \(\sum_jc_jG_j=0\),
then \(\sum_jc_jH_{\rm NT}(x,e_j)=0\)
for every \(x\in W\), which forces
every \(c_j=0\).
The integral operator with this
kernel acts by
\[
 \mathcal A_Nf
      =-\frac{24\eta}{\psi(N)}
                    \sum_j\langle f,G_j\rangle_{H_N}G_j.
\]
It has precisely the stated sign,
range, rank, and kernel. In particular
\(\mathscr D_N(f,h)=\langle\mathcal A_Nh,f\rangle_{H_N}\);
the reversed order is required by
the scalar convention in
\eqref{eq:cuspidal-Eis-Gram-form}.
For the final formulation, the full
lift equals the lift of \(Y\), so
its coordinates are
\(\langle F_j,f\rangle=\langle P_SF_j,f\rangle\).
The same calculation applies, allowing
zero or dependent projected functions
when the chosen coefficient span is
larger than \(W\).
\end{proof}

\begin{corollary}[Determinant minors and rank certificates]
\label{cor:cuspidal-Eis-Gram-minors}
Retain the same fixed-field and positive-height hypotheses.
For \(f_1,\ldots,f_m\in S\), put
\[
 D_{ij}=\mathscr D_N(f_i,f_j),\qquad
 M_{ij}=H_{\rm NT}(\Theta(f_i),\Theta(f_j)).
\]
Then \(D=-c_{\rm ht}M\), and
\[
 \operatorname{rank}D
       =\dim_{\mathbb C}\operatorname{span}
                                \{\Theta(f_1),\ldots,\Theta(f_m)\}.
\]
Its nullspace is exactly
\[
 \left\{\mathbf a\in\mathbb C^m:
                 \Theta\!\left(\sum_i a_if_i\right)=0\right\}.
\]
For any subset \(I\) of size \(\ell\),
\begin{equation}\label{eq:cuspidal-Eis-principal-minor}
 (-1)^\ell\det D_{I,I}
       =c_{\rm ht}^{\,\ell}\det M_{I,I}\geq0,
\end{equation}
with strict inequality exactly when
the corresponding \(\ell\) lift classes
are complex linearly independent.
All minors of size greater than
\(\operatorname{rank}J_0(N)(K)\) vanish.
In particular a strictly positive
sign-adjusted principal minor of size
\(\ell\) certifies
\(\operatorname{rank}J_0(N)(K)\geq\ell\).

For two index sets \(I,J\) of equal
size \(\ell\), one also has
\[
 \det D_{I,J}=(-c_{\rm ht})^\ell\det M_{I,J},
 \qquad
 |\det D_{I,J}|^2
    \leq\det(-D_{I,I})\det(-D_{J,J}).
\]
When the two principal determinants
are nonzero, equality in the latter
inequality holds precisely when
the two \(\ell\)-tuples of lifts
span the same complex subspace.
\end{corollary}

\begin{proof}
The entrywise identity is
\eqref{eq:cuspidal-Eis-NT-Gram}.
Positivity of the N\'eron--Tate Gram
matrix proves the rank and determinant
claims. To check the coefficients in
the nullspace assertion, use the
conjugate-linearity of \(\Theta\):
\[
 \mathbf a^*M\mathbf a
     =H_{\rm NT}\!\left(
           \Theta(\sum_i a_if_i),
           \Theta(\sum_i a_if_i)\right).
\]
The nonnegative Gram matrix has
zero quadratic value exactly on its
kernel. Principal Gram determinants
are positive exactly for independent
vectors. Finally the mixed-minor
inequality is Cauchy--Schwarz on
the \(\ell\)-th exterior power of
the finite-dimensional height space.
For nonzero exterior vectors,
equality is equivalently equality
of their \(\ell\)-dimensional spans.
\end{proof}

\begin{remark}
The finite-rank statement applies to the double cuspidal
projection. It does not imply finite rank or a pointwise sign
for the unprojected Eisenstein kernel. No eigenpacket condition
on \(f,h\) is needed. The rank bounds use the fixed-field
realization and positive-height hypotheses stated above.
\end{remark}

\subsection{The precise normalization problem}

For a nonzero scalar cusp form \(g\) in a specified Shimura packet,
put \(r(g)=\|Tg\|_\Gamma^2/\|g\|_H^2\). Under the arithmetic,
analytic, and normalization hypotheses of
Assumption~\ref{ass:red-inputs}, the height reduction gives the
following coefficient in front of \(\|g\|_H^2L'(G,1)\) when
\(L(G,1)=0\).
\begin{equation}\label{eq:final-coefficient}
 C(g,N)=-\frac{c_{\rm SW}A_g(0)}{\eta}\,r(g)
\end{equation}
For the projected new line, Theorem~\ref{thm:main-projected-height}
identifies the scalar coefficient as $(2\sqrt N/\pi)r(g)$.
Corollary~\ref{cor:BY-scalar-normalizations} gives the precise
obstruction to omitting $r(g)$, and the canonical isometric
normalization for which that factor becomes one. A change to $T$
or to the scalar measure must be propagated through the arithmetic
lift as well.

The nonsingular and positive rank-one comparisons fix the
normalization relevant to the completed paired identity as
\begin{equation}\label{eq:zhu-normalization}
 c_{\rm SW}=\frac{\psi(N)}{24},
 \qquad \psi(N)=N\prod_{p\mid N}(1+p^{-1}),
\end{equation}
for the fixed section of
Section~\ref{subsec:all-level-nonsingular-ASW}.
Consequently $c_{\rm SW}$ in
\eqref{eq:red-ASW-input} cannot simply be set to $1$ while retaining
this Eisenstein normalization. For the explicit completion,
Theorem~\ref{thm:explicit-completion-paired-doubling} proves the
required equality after double cuspidal pairing. This application
uses neither a raw coefficientwise rank-zero identification nor
an assertion that each boundary mode separately pairs to zero.

\subsection{The arithmetic kernel and the sign of the height}

Let $\mathcal A_\mathbb R$ be the real arithmetic Chow space in which the
arithmetic theta series takes values, and let
$[\ ,\ ]_{\rm GS}$ be its real symmetric arithmetic intersection pairing.
On its complexification use the Hermitian extension
\[
 H_{\rm GS}(x,y)=[x,\overline y]_{\rm GS}.
\]
Suppose the scalar arithmetic theta series $\widehat\phi(z)$ has weight
$k$ and character $\chi$ for $H$.  Define
\begin{equation}\label{eq:red-arithmetic-lift}
 \widehat\theta(f)=\int_{H\backslash\mathbb H}
           \widehat\phi(z)\overline{f(z)}y_z^k\,d\mu(z).
\end{equation}
With our convention this map is conjugate-linear in $f$.  Let
$\widetilde{\rm MW}_\mathbb R\subset\mathcal A_\mathbb R$ be the
admissible degree-zero Mordell--Weil subspace.  Fix the height convention
by specifying a constant $\eta>0$ such that
\begin{equation}\label{eq:red-FH-normalization}
 H_{\rm GS}(x,y)=-\eta H_{\rm NT}(x,y)
       \quad(x,y\in\widetilde{\rm MW}_\mathbb C).
\end{equation}
This is where the Faltings--Hriljac comparison and the precise
normalization of the N\'eron--Tate pairing enter.  For example, defining
$H_{\rm NT}=-H_{\rm GS}$ on this subspace amounts to choosing $\eta=1$;
one must not simultaneously import a different height convention.

Set
\begin{equation}\label{eq:red-arithmetic-kernel}
 K_{\rm ar}(z,w)=
 H_{\rm GS}(\widehat\phi(z),\widehat\phi(w))
 =[\widehat\phi(z),\overline{\widehat\phi(w)}]_{\rm GS}.
\end{equation}
The conjugation in the second variable is essential.  For a theta series
with real arithmetic Fourier coefficients it is equivalent to evaluating
its second factor at $-\overline w$.

\begin{assumption}[Arithmetic and analytic inputs]
\label{ass:red-inputs}
For the chosen form $f$, assume the following statements.
\begin{enumerate}
 \item The lift \eqref{eq:red-arithmetic-lift} exists and lies in
       $\widetilde{\rm MW}_\mathbb C$.
 \item There are a specified constant $c_{\rm SW}$ and a kernel $R$ such
       that
       \begin{equation}\label{eq:red-ASW-input}
        K_{\rm ar}(z,w)=c_{\rm SW}
          \left.\frac{\partial}{\partial s}K_{00}(z,w;s)\right|_{s=0}
          +R(z,w),
       \end{equation}
       where the double pairing of $R$ with
       $\overline{f(z)}f(w)y_z^ky_w^k$ is zero.
 \item The kernels and their derivatives are sufficiently integrable to
       interchange the arithmetic pairing, both integrals, and the
       derivative at $s=0$.  For instance, locally uniform polynomial
       growth at every cusp, together with locally integrable
       singularities and a common integrable majorant, suffices because
       $f$ is cuspidal.
 \item The classical kernel and the adelic data are matched as in
       Proposition~\ref{prop:paired-identity-sufficient}, so that the
       paired doubling identity \eqref{eq:red-doubling-consequence} holds,
       and its multiplier $A_f$ is holomorphic at $0$.
       The form $G$ is the normalized weight-two newform of the precise
       Waldspurger correspondent for this inducing character and
       additive character.
\end{enumerate}
\end{assumption}

For $\widehat\phi=\widehat\Phi_{\rm cmp}$, the first three
arithmetic and interchange requirements follow from
Corollary~\ref{cor:completion-full-lift-MW} and
Theorem~\ref{thm:explicit-completion-paired-doubling}, with
$c_{\rm SW}=\psi(N)/24$ and the reflected second-variable
normalization of that theorem. The standard uniform polynomial
bounds for Eisenstein series in a regular parameter neighborhood
justify differentiation under the cuspidal integral. The fourth
requirement specifies the additional match to the chosen adelic
data and its local multiplier; it is retained when those data
are varied.

\begin{proposition}[Formal height reduction]\label{prop:red-height}
Under Assumption~\ref{ass:red-inputs} and the normalization
\eqref{eq:red-FH-normalization},
\begin{align}
 H_{\rm NT}(\widehat\theta(f),\widehat\theta(f))
  &=-\frac{c_{\rm SW}}{\eta}\,\|Tf\|_\Gamma^2
       \bigl(A_f(0)L'(G,1)+A_f'(0)L(G,1)\bigr).
       \label{eq:red-general-height}
\end{align}
If $L(G,1)=0$, this simplifies to
\begin{equation}\label{eq:red-central-height}
 H_{\rm NT}(\widehat\theta(f),\widehat\theta(f))
   =-\frac{c_{\rm SW}A_f(0)}{\eta}\,
        \|Tf\|_\Gamma^2 L'(G,1).
\end{equation}
If in addition \eqref{eq:red-projection-hypothesis} has been proved,
then
\begin{equation}\label{eq:red-projected-height}
 H_{\rm NT}(\widehat\theta(f),\widehat\theta(f))
   =-\frac{c_{\rm SW}A_f(0)c_f}{\eta}\,
        \|f\|_H^2 L'(G,1).
\end{equation}
\end{proposition}

\begin{proof}
Expanding \eqref{eq:red-arithmetic-lift} in both arguments of the
Hermitian intersection pairing gives
\[
 H_{\rm GS}(\widehat\theta(f),\widehat\theta(f))
 =\int\!\!\int K_{\rm ar}(z,w)\overline{f(z)}f(w)
                    y_z^ky_w^k\,d\mu(z)d\mu(w).
\]
The arithmetic kernel identity and the stated vanishing of $R$ turn the
right-hand side into $c_{\rm SW}D_f'(0)$.  Apply
\eqref{eq:red-FH-normalization}, then
\eqref{eq:red-doubling-consequence}, and differentiate the product
$A_f(s)L(G,s+1)$.  This proves
\eqref{eq:red-general-height}.  The remaining assertions follow from
$L(G,1)=0$ and \eqref{eq:red-cf}.
\end{proof}

\begin{corollary}[Odd functional equation]\label{cor:red-sign}
Suppose $G$ has level $N$ and trivial character, and use the finite
Hecke $L$-function and completion
\begin{equation}\label{eq:red-completed-L}
 \Lambda(G,s)=N^{s/2}(2\pi)^{-s}\Gamma(s)L(G,s).
\end{equation}
If $\Lambda(G,s)=-\Lambda(G,2-s)$, then $L(G,1)=0$ and
\[
 \Lambda'(G,1)=\frac{\sqrt N}{2\pi}L'(G,1).
\]
In this case the trace-corrected height formula
\begin{equation}\label{eq:red-target}
 H_{\rm NT}(\widehat\theta(f),\widehat\theta(f))
  =\frac{2\sqrt N}{\pi}\,\|Tf\|_\Gamma^2L'(G,1)
  =4\|Tf\|_\Gamma^2\Lambda'(G,1)
\end{equation}
follows from Proposition~\ref{prop:red-height} if the independent
normalization calculation proves
\begin{equation}\label{eq:red-target-constant}
 -\frac{c_{\rm SW}A_f(0)}{\eta}=\frac{2\sqrt N}{\pi}.
\end{equation}
In scalar norms the coefficient in \eqref{eq:red-target} is
$(2\sqrt N/\pi)r(f)$.  If $c_f=r(f)=\nu/2$ is separately
established, it is $\nu\sqrt N/\pi$.
The uncorrected scalar coefficient and the isometric alternative
are distinguished in Corollary~\ref{cor:BY-scalar-normalizations}.
\end{corollary}

\begin{proof}
The odd functional equation implies $\Lambda(G,1)=0$.  The multiplier
$N^{s/2}(2\pi)^{-s}\Gamma(s)$ is nonzero at $s=1$, proving
$L(G,1)=0$.  Differentiating \eqref{eq:red-completed-L} at $1$ gives
the derivative identity because the term proportional to $L(G,1)$
vanishes.  Substitute \eqref{eq:red-target-constant} into
\eqref{eq:red-central-height}.
\end{proof}

\subsection{Boundary shadows and cuspidal pairing}
\label{subsec:boundary-defect-note}

Write $N=N_0h^2$, where $N_0$ is positive and square-free, and
put $\kappa_N=\sqrt N/(2\pi)$. Zhu's equation~(48) gives the
coefficient of the cuspidal divisor in $Z^*(t,v)$. Its generating
function is
\begin{equation}\label{eq:boundary-eichler-series}
 b_N(\tau)=\kappa_Nv^{-1/2}
    \left(1+\sum_{j\geq1}
          \beta_{3/2}(4\pi N_0j^2v)q^{-N_0j^2}\right),
 \qquad
 \beta_s(x)=\int_1^\infty e^{-xu}u^{-s}\,du.
\end{equation}
The equality $-Nt>0$ square is equivalent to $t=-N_0j^2$.
This identifies the boundary coefficient explicitly without
claiming that it is an independently modular summand of the full
arithmetic theta series.

\begin{lemma}[The boundary coefficient has unary shadow]
\label{lem:boundary-shadow}
The series \eqref{eq:boundary-eichler-series} and its derivatives
converge locally uniformly on $\mathbb H$. For
$\xi_kF=2iv^k\overline{\partial F/\partial\overline\tau}$ one has
\begin{equation}\label{eq:boundary-shadow}
 \xi_{3/2}b_N=-\frac{\kappa_N}{2}\,
          \vartheta_{N_0}(\tau),
 \qquad
 \vartheta_{N_0}(\tau)=\sum_{j\in\mathbb Z}q^{N_0j^2}.
\end{equation}
In particular, $b_N$ is annihilated by the weight-$3/2$
hyperbolic Laplacian.
\end{lemma}

\begin{proof}
For $a>0$, integration by parts yields
$\beta_{3/2}(a)=2e^{-a}-2a\beta_{1/2}(a)$. It follows that
\[
 \frac{d}{dv}\left[v^{-1/2}\beta_{3/2}(4\pi\ell v)\right]
        =-v^{-3/2}e^{-4\pi\ell v}\quad(\ell>0).
\]
The exponential bound on these expressions also justifies local
termwise differentiation. Applying $\xi_{3/2}$ to the constant
term gives $-\kappa_N/2$, and applying it to the term indexed by
$j\geq1$ gives $-\kappa_Nq^{N_0j^2}$. This proves
\eqref{eq:boundary-shadow}; its right-hand side is holomorphic,
so the standard factorization of the weight Laplacian gives the
last assertion.
\end{proof}

The full degree series, which combines the positive CM degrees,
zero Hodge coefficient and negative cusp coefficients, is a
genus-one Eisenstein series by Zhu, Proposition~8.1.1. The
vertical and Hodge pairings are likewise explicitly computed in
Section~8.2 and Corollary~8.3.3. Thus the nonholomorphic boundary
coefficient is compatible with the Eisenstein vanishing of those
complete scalar functionals. These identities do not express the
two-variable Green-star boundary transgression in the product
defect. Moreover a prescribed unary shadow does not determine a
cuspidal projection: any two modular completions with that
shadow can differ by a holomorphic cusp form. The isolated
calculation \eqref{eq:boundary-shadow} therefore proves neither
vanishing nor nonvanishing of the complete boundary defect.

\paragraph{A paired boundary criterion.}
Let $X_\varepsilon$ be the compact curve with small cusp disks
removed. Suppose the proper Green-star defect on
$X_\varepsilon$ is $d\eta(z,w)$, where the chosen primitives
have compatible normalizations. Its integrated boundary
contribution is the Stokes flux
\[
       B_\varepsilon(z,w)=\int_{\partial X_\varepsilon}\eta(z,w).
\]
Write $C_\varepsilon(z,w)$ for the arithmetic-degree
contribution of the explicitly prescribed cusp corrections in
the same truncation, and put
$B_\varepsilon^{\rm comp}=B_\varepsilon+C_\varepsilon$.
Assume that the compactified boundary defect is represented by
the limit of $B_\varepsilon^{\rm comp}$ and that this limit
commutes with \eqref{eq:defect-pairing}. It suffices to prove
\begin{equation}\label{eq:paired-stokes-condition}
 \lim_{\varepsilon\to0}
  \int_{(H\backslash\mathbb H)^2}
     B_\varepsilon^{\rm comp}(z,w)\overline{f(z)}h(w)
      v_z^{3/2}v_w^{3/2}\,d\mu(z)d\mu(w)=0,
\end{equation}
For example, a uniform integrable majorant and a limit expressed as tensors with a
verified Eisenstein or unary factor suffice. This requires only
vanishing after pairing with the desired cusp space; equality
of compactified Green currents is a stronger requirement.
Neither the existence of an isotropic rank-one quotient at a
cusp nor \eqref{eq:boundary-shadow} alone proves
\eqref{eq:paired-stokes-condition}.

\paragraph{An explicit route from harmonicity to holomorphy.}
The following elementary criterion isolates a possible additional
input, without asserting it for the actual defect rows.

\begin{lemma}\label{lem:harmonic-nonnegative-holomorphic}
Let $F$ be a smooth automorphic function of weight $3/2$ for a
finite theta multiplier on a finite-index subgroup, harmonic for
the weight Laplacian and of polynomial growth at every cusp.
If its Fourier expansion at infinity contains no negative
frequencies, including nonholomorphic negative-frequency
terms, then $F$ is holomorphic and holomorphic at every cusp.
\end{lemma}

\begin{proof}
The Fourier ordinary differential equation for a harmonic form
gives, under polynomial growth at infinity,
\[
 F(\tau)=\sum_{n\geq0}a(n)q^n+c\,v^{-1/2}.
\]
The other positive-frequency solutions grow exponentially and
are excluded; negative-frequency solutions are absent by
hypothesis. Hence $\xi_{3/2}F=-\overline c/2$ is a constant
holomorphic automorphic function of weight $1/2$ for the dual
multiplier. A nonzero constant cannot have this transformation
law under an element whose lower-left entry is nonzero; such
an element exists in any finite-index subgroup. Therefore
$c=0$. This proves holomorphy on $\mathbb H$. At every other
cusp, holomorphy and polynomial growth exclude negative
holomorphic Fourier powers, proving holomorphy there as well.
\end{proof}

To use square-class support as a unary-theta argument, one must
additionally establish the relevant row's modularity, holomorphy
and growth, and apply a classification with its exact level and
character hypotheses. In particular, nonholomorphic coefficients
$\delta_{m,n}(v_z,v_w)$ cannot be treated as constant
coefficients of a holomorphic unary theta series. Kudla's
Corollaries~11.5--11.6 give the distinguished-$O(1)$ mechanism
in their compact quaternionic setting and under their explicit
dyadic density assumption; they do not by themselves prove an
all-level compactified modular-curve product identity.

The exact lowering formula
\eqref{eq:standard-arithmetic-lowering} shows what a differential
approach must compare: the Siegel kernel together with its
nonzero cusp divisor. A proof that the product defect is
holomorphic would require the corresponding lowering identity
for the restricted Eisenstein derivative, including its
regularized boundary contribution. The compact Kudla--Millson
orthogonality proved earlier does not supply that identity.

For a separately specified arithmetic family, the arithmetic
input to this reduction is exactly
\eqref{eq:defect-minimal-condition}. For the standard family,
Theorem~\ref{thm:full-standard-arithmetic-doubling} proves that
condition through proper coefficient matching, common-component
adjunction, and the nonpositive comparisons. Degree, vertical,
and Hodge orthogonality alone do not imply the product identity.
These orthogonality statements concern complete scalar
functionals of a lift, whereas the product comparison involves
the full two-variable kernel.

\subsection{The full compact product theorem of Kudla--Rapoport--Yang}
\label{subsec:KRY-product-scope}

Kudla--Rapoport--Yang prove the full arithmetic product identity
\cite[Introduction, Theorem~C; Chapter~7]{KRYBook}.  Their equality
includes every pair of integral indices and the sum over all
off-diagonal entries of the genus-two matrix.  Their setting is
the Shimura curve of an indefinite quaternion division algebra,
with a fixed maximal order; the introduction explicitly excludes
modular curves.  Thus the full compact identity is established
there.  Applying it to $X_0(N)$ requires a new comparison of the
integral models, special cycles and noncompact Green data.

The role of the finite terms is visible in
\cite[Section~13]{KudlaYangPullback}.  Lemma~13.4 and
Proposition~13.5 contain matching conductor terms, the latter
contributing $2\deg_{\mathbb Q}(\mathcal C)\log(n^2)$ in that
paper's notation.  The archimedean finite part alone omits this
contribution.  Their Remark~13.3 also corrects the irreducibility
assertion in \cite[Lemma~7.4.2]{KRYBook}; actual horizontal
components and multiplicities must be retained.  These facts
motivate the branchwise statement of
Lemma~\ref{lem:integral-CM-adjunction}, without identifying its
normal bundle with an arbitrary-level metrized Hodge class.

\subsection{Recent higher-dimensional results and their scope}

Chen--Zhu \cite{ChenZhu2026} prove a full generating-series degree
identity for $\operatorname{GSpin}(2,2)$ on the product of modular
curves. Its genus-three Green-current extension is
Theorem~C.1.1, for the lattice $M_2(\mathbb Z)$. Section~9.3
explicitly leaves the analogous genus-one and genus-two extensions
outside its scope. Section~9.7 defines arithmetic degrees directly
while keeping technical questions about classes on the compactified
stack separate. Thus this recent result provides relevant methods,
but does not supply the compactified genus-one product identity
required in \eqref{eq:defect-minimal-condition}.

\subsection{Scope and normalization of the resulting formulas}

The standard product formula covers every pair of integral indices.
Theorem~\ref{thm:standard-positive-square-product} treats common CM
components through the normal line, the different and conductor,
internal and scaling contacts, and all regular Green values.
Corollary~\ref{cor:corrected-standard-support} supplies the
rank-zero correction. The nonsingular and rank-one comparisons
and the middle-orbit calculation identify these corrected
coefficients with the fixed Eisenstein series of
Section~\ref{subsec:all-level-nonsingular-ASW}.

The explicit completion equals the standard family at every index.
Theorem~\ref{thm:explicit-completion-paired-doubling} gives the
paired identity, and Theorem~\ref{thm:full-standard-arithmetic-doubling}
gives the pointwise identity. Strong Fourier convergence,
joint growth, and the Bochner estimates justify passage from
coefficient identities to the full kernels and lifts. The full
standard cusp lift equals its admissible Mordell--Weil projection.
In the singular comparison, \(h_t(s)-\tfrac12\log\nu\,d_t(s)\)
vanishes by the degree and discriminant-height formulas, while
\(b_N(y)-\tfrac12d_0\mathcal Z_y(1/2)\) vanishes by the actual
middle-orbit calculation. Invariance under
\(\mathrm{GL}_2(\mathbb Z)\) alone would not give the latter identity.

A different arithmetic family or Eisenstein normalization requires
its own comparison. The paired condition is
\eqref{eq:defect-minimal-condition}; the cusp data in
Sections~\ref{subsec:cusp-width-correction} and
\ref{subsec:negative-cusp-correction} must be matched as well.
For test forms with \(0\ne Tf\in\mathbb Cg_G\), suppose that
the chosen family's full lift equals the Mordell--Weil lift of
Theorem~\ref{thm:BY-newline-height}, that the convergence hypotheses
of Proposition~\ref{prop:defect-height} hold, and that \(A_f\) is
regular at zero. Then \eqref{eq:BY-defect-evaluation} evaluates
the residual. It vanishes if \(L'(G,1)=0\). When
\(L'(G,1)\ne0\), its vanishing is equivalent to
\eqref{eq:BY-normalization-test} for the actual section, additive
character, measures, and local data. The corresponding analytic
factors remain explicit in Appendix~\ref{sec:analytic-doubling}.

\appendix
\section{Paired analytic doubling and local normalizations}
\label{sec:analytic-doubling}

We record the paired adelic doubling formula and its passage to
classical Petersson pairings. The unramified spherical factors are
explicit, while the ramified and archimedean integrals and the
classical conversion factors are retained in $A_f(s)$.
This factorization supplies the analytic input to the conditional
height reduction for general test data.

\subsection{Unfolding the pullback kernel}

Let $K(z,w;s)$ be an $\operatorname{End}(W)$-valued automorphic kernel
of weight $k$ in $z$ and conjugate weight $k$ in $w$.  More precisely,
write $J_k(\gamma,z)$ for the weight factor, so that
$F(\gamma z)=J_k(\gamma,z)\rho(\gamma)F(z)$, and require
\begin{equation}\label{eq:red-kernel-covariance}
 K(\delta z,\gamma w;s)=J_k(\delta,z)
 \overline{J_k(\gamma,w)}\,
 \rho(\delta)K(z,w;s)\rho(\gamma)^{-1}.
\end{equation}
For a genus-two Eisenstein series this kernel is its pullback along
$\operatorname{diag}(z,-\overline w)$, after making the required
identification of the second Weil representation with the dual and
including its automorphy factors.  This identification and its phase are
part of the normalization of $K$.

Define
\begin{align}
 B(s)F(z)&=\int_{\Gamma\backslash\mathbb H}
           K(z,w;s)F(w)y_w^k\,d\mu(w),\label{eq:red-vector-operator}\\
 K_{00}(z,w;s)&=(K(z,w;s)e_0,e_0)_W,\label{eq:red-zero-kernel}\\
 B_{00}(s)f(z)&=\int_{H\backslash\mathbb H}
           K_{00}(z,w;s)f(w)y_w^k\,d\mu(w),\label{eq:red-scalar-operator}\\
 D_f(s)&=\int_{(H\backslash\mathbb H)^2}
     K_{00}(z,w;s)\overline{f(z)}f(w)y_z^ky_w^k\,
     d\mu(z)d\mu(w).\label{eq:red-double-integral}
\end{align}
First work in a region of absolute convergence.  All asserted identities
elsewhere mean their meromorphic continuations.

\begin{lemma}[Scalar--vector pullback]\label{lem:red-pullback}
With the preceding conventions,
\begin{equation}\label{eq:red-pullback-operator}
 B_{00}(s)=RB(s)T,
 \qquad
 D_f(s)=\langle B(s)Tf,Tf\rangle_\Gamma.
\end{equation}
Suppose the pullback, its inducing character, and its adelization
are those of Theorem~\ref{thm:paired-metaplectic-doubling}, and $Tf\ne0$
belongs to the specified cuspidal Waldspurger packet. Then the paired
doubling identity gives
\begin{equation}\label{eq:red-doubling-consequence}
 D_f(s)=A_f(s)L(G,s+1)\|Tf\|_\Gamma^2,
\end{equation}
with the multiplier defined by the actual local matrix-coefficient
integrals and the classical conversion below. The stronger operator
identity
\begin{equation}\label{eq:red-doubling-eigen-sufficient}
 B(s)Tf=A_f(s)L(G,s+1)Tf
\end{equation}
is sufficient, but is not needed and is not asserted for an arbitrary
ramified vector.
\end{lemma}

\begin{proof}
Insert the finite sum defining $Tf$ into \eqref{eq:red-vector-operator}.
For each summand use $w\mapsto\gamma w$ and
\eqref{eq:red-kernel-covariance}.  The slash factors and $y_w^k$ cancel
and the translates of a $\Gamma$-fundamental domain tile an
$H$-fundamental domain.  This proves
\[
 B(s)Tf(z)=\int_{H\backslash\mathbb H}
                 K(z,w;s)e_0 f(w)y_w^k\,d\mu(w).
\]
Taking the $e_0$ component proves the first assertion of
\eqref{eq:red-pullback-operator}.  The second follows from
\eqref{eq:red-adjunction} and
$D_f(s)=\langle B_{00}(s)f,f\rangle_H$.
The \(L\)-factorization is a separate consequence of the adelic
doubling theorem and the specified classical conversion factors.
It is proved in Proposition~\ref{prop:paired-identity-sufficient},
using Theorem~\ref{thm:paired-metaplectic-doubling}; it does not
follow from the formal scalar--vector pullback alone.
Alternatively,
\eqref{eq:red-doubling-eigen-sufficient} implies the same equality
by taking its Petersson pairing with $Tf$.
\end{proof}

\begin{remark}[Index conventions]
If the scalar Petersson product is instead multiplied by a constant
$u_H>0$ and the vector product by $u_\Gamma>0$, then adjunction becomes
\[
 \langle Tf,F\rangle_{\Gamma,\rm new}
    =\frac{u_\Gamma}{u_H}\langle f,RF\rangle_{H,\rm new}.
\]
The norm ratios and the unfolding constants must be changed accordingly.
The preceding formulas use exactly
\eqref{eq:red-pet-scalar}--\eqref{eq:red-pet-vector}; switching to
volume-normalized products without this correction changes the final
constant.
\end{remark}

\subsection{The metaplectic doubling identity as a paired identity}
\label{subsec:metaplectic-paired-doubling}

Let $G_{\rm sp}=\operatorname{SL}_2$, let $\widetilde G_{\rm sp}(\mathbb A)$ be its
metaplectic cover, and put $[G_{\rm sp}]=G_{\rm sp}(\mathbb Q)\backslash G_{\rm sp}(\mathbb A)$.
The products of genuine functions appearing below descend to this
linear quotient, using the rational splitting of the cover.
Fix the global additive character $\psi$, with
$\psi_\infty(x)=e^{2\pi ix}$ in the determinant-space convention.
Here $\operatorname{Wald}_\psi$ means the correspondent in
\cite[Section~9, equation~(9.2)]{KudlaMFA}.
Every integrand below is invariant under the order-two kernel of the
cover: the genuine section and genuine matrix coefficient each
acquire the same minus sign, which cancels in their product.
We integrate this product with the Haar measure on the underlying
linear group, or equivalently give that finite kernel total mass one.
The same convention is used in both automorphic variables and in the
Petersson pairing. We do not insert counting measure of mass two on
the kernel. Any further factor associated to passage from adelic
group quotients to classical upper-half-plane quotients remains in
the explicitly recorded adelization comparison.
Let $\sigma=\bigotimes_v\sigma_v$ be an irreducible genuine cuspidal
automorphic representation whose Waldspurger correspondent
$\pi=\operatorname{Wald}_\psi(\sigma)$ is cuspidal on
$\operatorname{PGL}_2(\mathbb A)$.
In particular, exceptional one-dimensional-theta contributions are
excluded. When a weight-two form $G_0$ is used below, $\pi$ is
required to be precisely the unitarized representation of that form.

Double the underlying two-dimensional symplectic space as
$W_{\rm sp}\oplus(-W_{\rm sp})$ and use the induced embedding $\iota$ into
$\widetilde{\operatorname{Sp}}_4(\mathbb A)$.
In standard coordinates the embedding is
\[
 \iota(g_1,g_2)=\iota_0(g_1,g_2^\vee),\qquad
 g_2^\vee=d g_2d^{-1},\quad d=\operatorname{diag}(1,-1),
\]
as in \cite[equations~(10.2)--(10.3)]{KudlaMFA}. We transport the
Siegel parabolic and its section together so that the parabolic
stabilizes the diagonal isotropic subspace. This absorbs the fixed
open-orbit representative into the section. The $L$-factor continues
to use $\sigma$ and $\psi$ in Kudla's convention; no additional
contragredient or inverse-character substitution is made.
Use the genuine inducing character and the parameter of
\cite[Section~10, equations~(10.4)--(10.7)]{KudlaMFA}:
the exponent of the absolute determinant in unnormalized induction is $s+3/2$.
Thus the central parameter for a ternary Schwartz section is $s=0$.
Let $\Phi(s)=\bigotimes_v\Phi_v(s)$ be a factorizable holomorphic
standard section, and let $E(h,s,\Phi)$ be the unnormalized Eisenstein
series. For a pure tensor $\varphi=\bigotimes_vf_v\in\sigma$, let
$\varphi^\vee=\bigotimes_vf_v^\vee\in\sigma^\vee$ correspond to it under
the unitary pairing. Normalize local pairings so that
\[
 \langle \varphi^\vee,\varphi\rangle=\|\varphi\|^2
             =\prod_v\langle \varphi_v^\vee,\varphi_v\rangle_v,
\]
and use the same product Haar measure in the unfolding and in the
local integrals. Define normalized matrix coefficients
\[
 m_{\varphi,v}(g)=
 \frac{\langle\sigma_v^\vee(g)\varphi_v^\vee,\varphi_v\rangle_v}
      {\langle \varphi_v^\vee,\varphi_v\rangle_v},
 \qquad m_{\varphi,v}(1)=1,
\]
and define
\begin{align}
 Z_v(s;\varphi_v,\Phi_v)
  &=\int_{G_{\rm sp}(\mathbb Q_v)}
       \Phi_v(\iota(\widetilde g,1),s)m_{\varphi,v}(\widetilde g)\,dg,
       \label{eq:paired-local-zeta}\\
 Z(s;\varphi,\Phi)
  &=\int_{[G_{\rm sp}]^2}\varphi^\vee(g_1)\varphi(g_2)
                 E(\iota(g_1,g_2),s,\Phi)\,dg_1\,dg_2.
       \label{eq:paired-global-zeta}
\end{align}
Conjugating to a different Siegel parabolic inserts the corresponding
fixed open-orbit representative in both formulas; it must not be
inserted in only one of them.
Here $\widetilde g$ is either lift of $g$. The integrand is independent
of that choice; $dg$ is Haar measure on the linear group.

\begin{theorem}[General metaplectic doubling, in a fixed normalization]
\label{thm:paired-metaplectic-doubling}
Let $\Sigma$ contain $\infty$, $2$, and every place at which the
representation, additive character, section, or vector is ramified.
Outside $\Sigma$ take the standard spherical section with value one
at the identity, spherical vectors with $m_{\varphi,v}(1)=1$, and compact
subgroups of volume one. For sufficiently large $\Re(s)$ one has
\begin{equation}\label{eq:paired-global-partial}
 Z(s;\varphi,\Phi)=\|\varphi\|^2
 \frac{L^\Sigma(s+1/2,\pi)}{\zeta^\Sigma(2s+2)}
                  \prod_{v\in\Sigma}Z_v(s;\varphi_v,\Phi_v).
\end{equation}
This equality continues meromorphically to $s\in\mathbb C$.
It does not require $\varphi_v$ to be an eigenvector of a local doubling
operator. In particular $\sigma_\infty$ may be the holomorphic
discrete series of lowest weight $3/2$.
\end{theorem}

\begin{proof}
The global identity is the metaplectic doubling formula recorded in
\cite[equation~(10.7)]{KudlaMFA}; it is stated there for arbitrary
factorizable vectors and sections. Here is the unfolding relevant to
the present normalization. For $\Re(s)>3/2$ the genus-two Eisenstein
series converges absolutely. In a sufficiently far right half-plane
one may unfold its defining sum against the two cusp forms.
The nonopen parabolic orbits contribute constant terms of a cusp form
and vanish. The open orbit has diagonal stabilizer $G_{\rm sp}(\mathbb Q)$.
Writing its variables as $(xg,x)$ gives the identity
\[
 Z(s;\varphi,\Phi)=\int_{G_{\rm sp}(\mathbb A)}
       \Phi(\iota(\widetilde g,1),s)
       \langle\sigma^\vee(\widetilde g)\varphi^\vee,\varphi\rangle\,dg.
\]
Factorization of the matrix coefficient and section now gives the
product of local integrals, with the factor $\|\varphi\|^2$ from their
normalization at the identity. The spherical local calculation in
these conventions is
\[
 Z_v(s)=\frac{L_v(s+1/2,\pi_v)}{\zeta_v(2s+2)}
 \qquad(v\notin\Sigma).
\]
This proves \eqref{eq:paired-global-partial}.
The local continuation for arbitrary data, including real and dyadic
places, is supplied by \cite[Theorem~3.1]{GanDoubling}.
Combine it with the Eisenstein continuation to continue the identity.
\end{proof}

\begin{remark}[The general paired formula]
Equation~(10.7) in \cite{KudlaMFA}, rather than its specialized
Theorem~10.1, is being used. The latter fixes quaternionic sections
and particular newvectors and carries extra restrictions in that
version of the paper. The general product above retains arbitrary
ramified integrals. Gan's theorem is used for their convergence and
continuation, without silently changing the inducing parameter in a
printed spherical formula.
\end{remark}

\subsection{Complete finite \texorpdfstring{$L$}{L}-functions and the classical shift}

Write $\Sigma_{\rm fin}=\Sigma\setminus\{\infty\}$. At every finite place put
\begin{equation}\label{eq:paired-normalized-zeta}
 Z_v^\natural(s)=
  \frac{\zeta_v(2s+2)}{L_v(s+1/2,\pi_v)}
                                  Z_v(s;\varphi_v,\Phi_v).
\end{equation}
At an unramified place this is exactly $1$.
At every other place this formula is its definition, using the actual
local integral and the local factor of the specified $\pi$.
Thus no Euler factor has been omitted or replaced by $1$ merely
because the data are ramified. Purely algebraic rearrangement of
\eqref{eq:paired-global-partial} yields
\begin{equation}\label{eq:paired-global-complete-finite}
 Z(s;\varphi,\Phi)=\|\varphi\|^2 L_{\rm fin}(s+1/2,\pi)\,
       \frac{Z_\infty(s)}{\zeta(2s+2)}
                         \prod_{v\in\Sigma_{\rm fin}}Z_v^\natural(s).
\end{equation}

Suppose $G_0(\tau)=\sum_{n\geq1}a_nq^n$ is the specified normalized
weight-two newform, and $\pi$ is its unitary automorphic representation.
The unitarization is $a_n\mapsto a_n n^{-1/2}$, so, including the
ramified finite factors,
\[
 L_{\rm fin}(t,\pi)=L(G_0,t+1/2),\qquad
 L_{\rm fin}(s+1/2,\pi)=L(G_0,s+1).
\]
Consequently the factor multiplying $L(G_0,s+1)$ in
\eqref{eq:paired-global-complete-finite} is the explicitly defined
meromorphic function
\begin{equation}\label{eq:paired-A-explicit}
 A_{\varphi}^{\rm ad}(s)=\frac{Z_\infty(s)}{\zeta(2s+2)}
                         \prod_{v\in\Sigma_{\rm fin}}Z_v^\natural(s).
\end{equation}
If a different normalized Eisenstein series is $b(s)E(s,\Phi)$,
multiply $A_{\varphi}^{\rm ad}(s)$ by that specified $b(s)$.
An adelization factor relating the classical double Petersson integral
to \eqref{eq:paired-global-zeta} must likewise be carried through.
It disappears only after the chosen measures and adelization have
actually been matched.

\begin{proposition}[The paired identity is sufficient]
\label{prop:paired-identity-sufficient}
Let $F=Tf$, set $G=G_0$, and let $\varphi$ be its adelization in the specified
irreducible cuspidal representation. Suppose the classical pullback
and adelic section have been matched so that
\[
 D_f(s)=c_D(s)Z(s;\varphi,\Phi),\qquad
 \|\varphi\|_{\rm ad}^2=c_{\rm norm}\|Tf\|_\Gamma^2.
\]
Here $c_D(s)$ includes any Eisenstein normalizing factor, and
$c_{\rm norm}>0$ records the chosen adelization and measures. Put
\begin{equation}\label{eq:paired-classical-multiplier}
 A_f(s)=c_D(s)c_{\rm norm}
        \frac{Z_\infty(s)}{\zeta(2s+2)}
                    \prod_{v\in\Sigma_{\rm fin}}Z_v^\natural(s).
\end{equation}
Then \eqref{eq:red-doubling-consequence} holds. If $A_f$ is
holomorphic at $0$, then
\[
 D_f'(0)=\|Tf\|_\Gamma^2
  \bigl(A_f(0)L'(G,1)+A_f'(0)L(G,1)\bigr).
\]
In particular, $L(G,1)=0$ removes the derivative of $A_f$.
No operator eigenvector assertion is needed.
\end{proposition}
\begin{proof}
Apply \eqref{eq:paired-global-complete-finite}, include the two
conversion factors, and differentiate.
\end{proof}

\begin{remark}[Central regularity]
Gan's Theorem~3.1 gives absolute local convergence for
$\Re(s)>-1/2$ when the local representation is tempered in the
untwisted unitary normalization. Thus temperedness of the actual
$\sigma_v$, together with holomorphic sections and the corresponding
nonvanishing regular local normalizing factors, is a sufficient
condition for central regularity. Meromorphic continuation alone is
not a proof of that regularity. Even when regular, the local value
may vanish for the specified vector and section. A nonvanishing
theorem for some choice of data does not prove nonvanishing for the
arithmetic data fixed here.
\end{remark}

\begin{remark}[Sums of pure tensors]
For a finite sum of pure tensors, expand the two cusp-form slots
sesquilinearly and use the unnormalized local matrix-coefficient
pairings in each cross term. The common $L$-factor can be factored
out; the remaining finite sum of products defines the paired
multiplier. One must not divide a cross term by a local pairing
that is zero. The resulting multiplier need not be a single Euler
product, and no assertion of a universal ramified constant follows.
\end{remark}

\begin{remark}[Inducing characters and the arithmetic application]
The untwisted formula above has fixed the genuine inducing character
and $\psi$ of its cited normalization.
For an arithmetic section one must first determine its transformation
law, for example by writing it in the explicit form
\[
 \Phi(m(t)h,s)=\chi_\psi(m(t))\chi_\Phi(\det t)
        |\det t|^{s+3/2}\Phi(h,s).
\]
The character $\chi_\Phi$ and the chosen additive character determine
which twisted Waldspurger $L$-function is represented. They cannot be
discarded on the basis of a scalar modular-form notation.
To apply \eqref{eq:paired-global-complete-finite} with the particular
untwisted $G_0$ in the manuscript, these data must be matched to the
normalization above, or the appropriate twisted doubling formula
must instead be stated and checked. This character matching and the
arithmetic product--pullback identity are separate from the analytic
unfolding proved here.
\end{remark}

\begin{lemma}[Scaling the quadratic form and the additive character]
\label{lem:paired-character-scaling}
Let $V_1$ be the trace-zero matrix space with $Q_1=\det$, let
$V_N$ be the same rational space with $Q_N=NQ_1$, and put
$\psi_N(x)=\psi(Nx)$ on $\mathbb A/\mathbb Q$.
Here $N\in\mathbb Q_{>0}$, so this is again a global nontrivial
character trivial on $\mathbb Q$. If $\psi_p$ has conductor
$\mathbb Z_p$, then $\psi_{N,p}$ has conductor
$p^{-v_p(N)}\mathbb Z_p$, and $\psi_{N,\infty}(x)=e(Nx)$
when $\psi_\infty(x)=e(x)$.
With the corresponding self-dual measures,
the Weil representations, in either genus one or genus two, satisfy
\[
 \omega_{\psi,V_N}=\omega_{\psi_N,V_1}.
\]
Consequently the Schwartz section for $(\psi,V_N)$ can be compared
with the determinant-space normalization by using $\psi_N$ throughout.
The $L$-function in that comparison is the one belonging to
$\operatorname{Wald}_{\psi_N}(\sigma)$.
Locally the self-dual measure on $V_N$ is
$|N|_v^{3/2}$ times the self-dual measure for $(\psi_v,V_1)$;
on $V_N^r$ the factor is $|N|_v^{3r/2}$.
These local factors have global product one.
\end{lemma}
\begin{proof}
The oscillator formulas depend on the additive quadratic function
$\psi(Q(x))$ and its polarized Fourier pairing.
Here $\psi(Q_N(x))=\psi_N(Q_1(x))$ and
$\psi([x,y]_N)=\psi_N([x,y]_1)$.
Scaling a nondegenerate Fourier pairing in dimension three by $N$
multiplies its self-dual measure by $|N|_v^{3/2}$; this follows by
applying the Fourier transform twice. The global product is one by
the product formula for the rational number $N$.
Thus the unipotent operators agree, the self-dual Fourier transforms
agree, and the associated Weil indices agree. The diagonal operators
have the same coordinate dilations and Jacobians; their scalar
indices are the same by the preceding identities.
These operators generate the symplectic group and its metaplectic
cover in the oscillator construction, proving the assertion.
The equality also identifies the evaluation-at-zero Schwartz maps
into the degenerate principal series and their standard extensions.
\end{proof}

To make the inducing-character comparison explicit, let
$\eta_{\psi,V_1}$ denote the full genuine Levi character in the
determinant-space convention, including its fixed discriminant
factor, and define
\[
 \chi_{N,v}(t)=(N,t)_v,\qquad
 \chi_N(t)=\prod_v\chi_{N,v}(t_v).
\]
The change of additive character multiplies the genuine scalar
character by this quadratic character. Thus the scaled section has
Levi character
\[
 \eta_{\psi_N,V_1}(m(t))
   =\eta_{\psi,V_1}(m(t))\chi_{N,v}(\det t).
\]
One may either keep this twist while retaining $\psi$, or use the
untwisted determinant-space convention with $\psi_N$; these are
the same section. The modulus exponent $s+3/2$ is unchanged.
The genuine-character identity is the one recalled in
\cite[Section~2.4]{GanDoubling}.
Correspondingly, the rank-one change-of-character rule
\cite[Theorem~5.3]{GanSavin} identifies local parameters by
$\phi_{\psi_N}=\phi_\psi\otimes\chi_{N,v}$.
For the cuspidal global correspondents this gives
\[
 \operatorname{Wald}_{\psi_N}(\sigma)
 \simeq\operatorname{Wald}_{\psi}(\sigma)\otimes\chi_N.
\]
In particular, keeping the old correspondent and deleting the
quadratic twist would in general give the wrong $L$-function.

This scaling is a way to specify the inducing convention concretely,
rather than deleting a quadratic character. The classical form
$G_0$ must then be defined as the correspondent for $\psi_N$.
It must not silently be replaced by the correspondent for a different
additive character. Every prime dividing $N$ is included in $\Sigma$,
so the altered local conductor of $\psi_N$ is permitted in the general
ramified integrals. Identifying the precise arithmetic genus-two
Eisenstein series, including any incoherent archimedean datum and its
global normalizing factor, is still required before using its derivative
as an intersection kernel.
The standing assumption that the Waldspurger correspondent is
cuspidal excludes unary-theta packets; the scaling lemma supplies
no argument for assigning them the cuspidal weight-two
$L$-function used in the theorem.

For an incoherent arithmetic Eisenstein series the scaling assertion
is applied to each actual local quadratic form. At the real place
the positive-definite ternary datum that defines the Eisenstein section
must be retained; it is not replaced by the signature-$(1,2)$ space
of the modular curve. Positive scaling preserves its signature and
its weight-$3/2$ $K$-type.

\subsection{Central local factors}

For pure tensor data the proved paired formula gives
\begin{equation}\label{eq:red-local-factor-bookkeeping}
 A_f(s)=a_\infty(s)\prod_{p\in\Sigma_{\rm fin}}b_p(s),\qquad
 a_\infty(s)=\frac{c_D(s)c_{\rm norm}Z_\infty(s)}{\zeta(2s+2)},
 \quad b_p(s)=Z_p^\natural(s).
\end{equation}
The $Z_p^\natural$ are defined by
\eqref{eq:paired-normalized-zeta}, using the actual representation,
vectors, sections, Haar measures, and parameter. Thus the global
finite $L$-function is complete. A numerical evaluation of this
analytic multiplier requires its particular central local integrals;
their definition does not evaluate them. The projected new-line
height in Theorem~\ref{thm:BY-newline-height} is already evaluated
independently of these integrals. For a sum of pure tensors,
the multiplier is instead the corresponding sum of local products.

\begin{lemma}[Only central local values are needed in the odd-sign case]
Suppose $L(G,1)=0$ and the total multiplier $A_f(s)$ in
\eqref{eq:red-local-factor-bookkeeping} is holomorphic at $0$.
Then the height in \eqref{eq:red-central-height} depends on $A_f(0)$
alone.  In particular, derivatives of the normalized bad-prime factors
are unnecessary once their product is known to be regular at $0$.
\end{lemma}

\begin{proof}
The product rule gives
$(A_f(s)L(G,s+1))'|_{s=0}=A_f(0)L'(G,1)$.
\end{proof}

\begin{remark}[Scope of the reduction]
The analytic identity in this appendix combines with the trace
adjunction in Section~\ref{sec:trace-frame} and the conditional
height reduction in Section~\ref{sec:height-reduction}. Moreover,
Theorem~\ref{thm:BY-newline-height} supplies the exact new-line
height constant with checked normalizations.
Theorem~\ref{thm:full-standard-arithmetic-doubling} proves the
arithmetic identity for the fixed standard data. Applying the
formal reduction to different arithmetic families or adelic
normalizations requires verifying
Assumption~\ref{ass:red-inputs} in that setting, matching the
specified sections, evaluating the corresponding central
multiplier at every bad prime (including $p=2$ and arbitrary
$v_p(N)$), and proving \eqref{eq:red-target-constant}.
For the specified new line, the cited projected vector-valued
theorem gives \eqref{eq:red-scalar-ratio-formula} with
checked normalizations.
The trace factor can instead be calculated from the full spectrum in
Theorem~\ref{thm:all-level-frame}, without exact zero-component proportionality.
For a general linear combination of eigenforms, one must apply the
corresponding sesquilinear identities to its spectral components; a
single unspecified newform $G$ cannot occur in a formula quantified over
all $f$.
\end{remark}

\section{Local calculations with specified test data}
\label{sec:local}

\subsection{A radial local calculation with explicit hypotheses}

Let $F=\mathbb Q_p$, let $\mathcal O=\mathbb Z_p$, and normalize
$|p|=p^{-1}$ and $\operatorname{vol}(\mathcal O^\times,d^\times a)=1$.
Write $\zeta_p(s)=(1-p^{-s})^{-1}$.

\begin{proposition}[The two-sided multiplicative integral]
Let $\chi:F^\times\longrightarrow\mathbb C^\times$ be a smooth
quasicharacter, set $\alpha=\chi(p)$, and let $W(t)$ be a scalar
meromorphic function. At a point where $W(t)$ is finite, put
\[
 B_t(a)=W(t)\begin{cases}
 1,&v_p(a)\geq0,\\
 |a|^{-2t},&v_p(a)<0.
 \end{cases}
\]
If
\[
 \Re(t)>\bigl|\log_p|\alpha|\bigr|,
\]
then $\int_{F^\times}\chi(a)|a|^tB_t(a)\,d^\times a$ is absolutely
convergent. If $\chi$ is unramified, its value is
\begin{align}
 \int_{F^\times}\chi(a)|a|^tB_t(a)\,d^\times a
 &=W(t)\left(
 \frac{1}{1-\alpha p^{-t}}+
 \frac{\alpha^{-1}p^{-t}}{1-\alpha^{-1}p^{-t}}\right)\label{eq:radial-correct}\\
 &=W(t)\frac{1-p^{-2t}}
 {(1-\alpha p^{-t})(1-\alpha^{-1}p^{-t})}\nonumber\\
 &=W(t)\frac{L_p(t,\chi)L_p(t,\chi^{-1})}{\zeta_p(2t)}.
 \nonumber
\end{align}
This formula provides a meromorphic continuation in $t$.
If $\chi$ is ramified, the integral is zero in its region of absolute
convergence, and its meromorphic continuation is identically zero.
\end{proposition}

\begin{proof}
Every $a\in F^\times$ has a unique expression $a=p^nu$ with
$n\in\mathbb Z$ and $u\in\mathcal O^\times$.
For an unramified character, the nonnegative-valuation shells contribute
\[
 W(t)\sum_{n\geq0}(\alpha p^{-t})^n.
\]
For $n=-k<0$, one has $|a|=p^k$ and
$|a|^tB_t(a)=W(t)p^{-kt}$, so the remaining shells contribute
\[
 W(t)\sum_{k\geq1}(\alpha^{-1}p^{-t})^k.
\]
Both series converge absolutely under the stated inequality, and summing
them proves the formulas. If $\chi$ is ramified, each shell has unit
integral $\int_{\mathcal O^\times}\chi(u)\,d^\times u=0$ by character
orthogonality. Absolute convergence justifies summing these zero terms.
\end{proof}

\begin{remark}
The coefficient in the second summand of \eqref{eq:radial-correct} is
$\chi(p)^{-1}$, rather than $\chi(p)$. A scalar factor multiplying the
original multiplicative integral must multiply both summands.
The formula involving local $L$-factors above is stated only for
unramified $\chi$: inserting $L_p(t,\chi)=1$ for ramified $\chi$ would
give the wrong answer. Ramified doubling calculations therefore need
test data with additional unit dependence.
\end{remark}

\subsection{Why the radial scaling holds for lattice functions}

Let $(U,Q)$ be a nondegenerate quadratic space of dimension $d>2$ over
$F$, let $M\subset U$ be a full lattice, and let $\psi$ be a
nontrivial additive character. Define the iterated oscillatory integral
\[
 W_M=\int_F\left(\int_M\psi(-bQ(x))\,dx\right)db
\]
and the analogous integral $B_M(a)$ with integrand
$\mathbf1_M(x)\mathbf1_M(-ax)\psi(-bQ(x))$.
The assertions below apply whenever these iterated integrals exist; they
also apply to limits of the compatible compact $b$-cutoffs.

\begin{proposition}[Scaling]
With these conventions,
\[
 B_M(a)=\begin{cases}
 W_M,&v_p(a)\geq0,\\
 |a|^{2-d}W_M,&v_p(a)<0.
 \end{cases}
\]
Thus the exponent is $-2t$ when $t=d/2-1$.
\end{proposition}

\begin{proof}
If $v_p(a)\geq0$, then $aM\subset M$, and the product of characteristic
functions is $\mathbf1_M(x)$. If $v_p(a)<0$, then the product is
$\mathbf1_{a^{-1}M}(x)$. Substitute $y=ax$ and then $c=ba^{-2}$.
The two Jacobians are $dx=|a|^{-d}dy$ and $db=|a|^2dc$, respectively,
which proves the second formula. These substitutions can first be made
with compact $b$-cutoffs and then passed to the limit.
\end{proof}

\begin{remark}
For $U=V\oplus H^r$ and $\dim V=3$, the geometric construction uses
$r\in\mathbb Z_{\geq0}$ and $d=3+2r$, so that $t=1/2+r$.
A quadratic space $V\oplus H^r$ is not defined for a general complex
$r$. Continuation in $r$ must be performed on a separately defined
rational or meromorphic family of local sections or local integrals.
The shell computation does not itself identify the integral with a
representation-valued doubling operator.
\end{remark}

\subsection{The zero density for the intended ternary lattice}

The lattice fixed in \eqref{eq:lattice-correct} has
coordinates $(a,b,c)\in\mathbb Z_p^3$ and quadratic form
$Q(a,b,c)=-ac-Nb^2$. (The displayed containment in $M_2(\mathbb Z)$
would impose an additional restriction on $a$ and is a different
lattice.) Write $N=p^mu$, with $m\geq0$ and $u\in\mathbb Z_p^\times$.
Normalize $db$ by $\operatorname{vol}(\mathbb Z_p)=1$, assume that
$\psi$ has conductor $\mathbb Z_p$, and initially give the coordinate
lattice volume one.

\begin{proposition}[An explicit density at $r=0$]
For every prime $p$ and $m\geq0$,
\begin{align}
 \beta_m
 &:=\lim_{k\to\infty}p^{-2k}
 \#\{(x,y,z)\bmod p^k:xy+p^mu z^2\equiv0\pmod{p^k}\}\nonumber\\
 &=m+1-\frac{m-1}{p}.\label{eq:ternary-density}
\end{align}
Here $W_M$ means specifically the limit of the compact-cutoff integrals
on this ternary lattice. A family on $V_r$ requires its own specified
quadratic space and Schwartz function. With coordinate Haar measure that cutoff limit is
$W_M=\beta_m$. For the self-dual measure attached to
$[x,y]=Q(x+y)-Q(x)-Q(y)$ and the specified $\psi$, one instead has
\[
 W_M=|2N|_p^{1/2}\left(m+1-\frac{m-1}{p}\right).
\]
\end{proposition}

\begin{proof}
Put $h(n)=\max\{0,\lceil(n-m)/2\rceil\}$.
Additive-character orthogonality identifies the compact-cutoff integral
with
\[
 \int_{p^{-k}\mathbb Z_p}\int_{\mathbb Z_p^3}
 \psi\bigl(b(xy+p^mu z^2)\bigr)\,dx\,dy\,dz\,db
 =p^k\operatorname{vol}\{(x,y,z):xy+p^mu z^2\in p^k\mathbb Z_p\}.
\]
If $v_p(x)=n<k$, solvability in $y\in\mathbb Z_p$ requires
$v_p(z)\geq h(n)$, and the set of such $y$ has measure $p^{n-k}$.
The contribution of this shell is therefore
$(1-p^{-1})p^{-h(n)}$.
For $x\in p^k\mathbb Z_p$, the contribution is $p^{-h(k)}$.
Consequently the density at finite cutoff $k$ equals
\[
 (1-p^{-1})\sum_{n=0}^{k-1}p^{-h(n)}+p^{-h(k)}.
\]
The last term tends to zero. Since $h(n)=0$ for $0\leq n\leq m$
and the values $1,2,\ldots$ occur twice thereafter, the limiting sum is
\[
 (1-p^{-1})\left(m+1+2\sum_{j\geq1}p^{-j}\right)
 =m+1-\frac{m-1}{p}.
\]
No division by $2$ is used, so the calculation also applies at $p=2$.
The Gram determinant of the bilinear form in the coordinate basis has
absolute value $|2N|_p$. The self-dual measure is therefore
$|2N|_p^{1/2}$ times coordinate Haar measure, giving the final formula.
\end{proof}

\begin{remark}
For example, with coordinate measure, $\beta_0=1+p^{-1}$,
$\beta_1=2$, and $\beta_2=3-p^{-1}$. If
$W(0,L)$ denotes this cutoff integral on the stated lattice,
then a value of $1$ at every positive even $v_p(N)$ is
incompatible with either of the two specified conventional measures.
An identification with a family involving $V_r$ and
$\varphi_L^{(r)}$ requires definitions of those data, which
are not supplied by the present calculation.
A different local density formula requires a different precisely stated
lattice, Schwartz function, or normalization.
\end{remark}

\subsection{What an invariant-vector argument actually proves}

\begin{proposition}[A conditional scalar identity]
Let $\sigma_p$ be a smooth representation, $A$ a compact open subgroup
on which it acts, and $f_p$ a nonzero vector with
$\sigma_p^A=\mathbb C f_p$. Assume that the local section $\Phi_p(s)$
is invariant under right translation by $i(A,1)$ and that the doubling operator satisfies
\[
 Z_p(s,r(i(a,1))\Phi_p,f)=\sigma_p(a)Z_p(s,\Phi_p,f)
 \qquad(a\in A),
\]
where $r$ denotes right translation of the local section. Then, wherever that operator
is defined, there is a scalar $c_p(s)$ such that
\[
 Z_p(s,\Phi_p,f_p)=c_p(s)f_p.
\]
The same identity holds meromorphically whenever the operator has a
meromorphic continuation. For any linear functional $\ell_p$ with
$\ell_p(f_p)=1$, one may compute
$c_p(s)=\ell_p(Z_p(s,\Phi_p,f_p))$.
\end{proposition}

\begin{proof}
Equivariance and the assumed right $i(A,1)$-invariance put the image in
$\sigma_p^A$, which is the specified one-dimensional line.
\end{proof}

\begin{remark}
Existence of an $A$-fixed vector alone does not show that the fixed
space is one-dimensional. If its dimension is larger, the local
operator must be computed on that space and can have off-diagonal
entries. Moreover, evaluation $f_p(1)$ presupposes a specified function
model and cannot be used as a normalization in an arbitrary abstract
representation. A principal-series calculation needs an additional
argument before it applies to special representations or
supercuspidal representations. In a metaplectic representation the
subgroup lift and genuine character must also be specified.
\end{remark}

\begin{remark}[Compact measure and the averaging identity]
For $K=\mathrm{SL}_2(\mathbb Z_p)$ with volume one and
$K_m=\{k\in K:c(k)\in p^m\mathbb Z_p\}$, $m\geq1$, one has
\[
 [K:K_m]=p^{m-1}(p+1),\qquad
 \operatorname{vol}(K_m)=\frac{1}{p^{m-1}(p+1)}.
\]
Indeed, $K/K_m$ is the projective line over $\mathbb Z/p^m\mathbb Z$,
whose cardinality is
$(p^{2m}-p^{2m-2})/(p^m-p^{m-1})$.
If, in a specified function model, $f|_K=\mathbf1_{K_m}$ and
$\omega(K_m)\varphi=\varphi$, then the averaging computation gives
\[
 \int_Kf(k)\overline{\omega(k)\varphi}\,dk
 =\operatorname{vol}(K_m)\overline\varphi.
\]
In particular, $p^{-m}$ does not follow from this usual compact Haar
normalization. Other support or transformation laws for $f$ require
their own averaging calculation.
\end{remark}

\begin{remark}[The displayed ramified principal-series vector]
For $m\geq1$, the double coset $PwK_m$ does not contain $1$.
Indeed, if $k=\left(\begin{smallmatrix}a&b\\c&d\end{smallmatrix}\right)
\in K_m$, upper triangularity of $wk$ would force $a=0$, whereas
$ad-bc=1$ and $p\mid c$ make this impossible.
Consequently a function supported on $PwK_m$ has $f(1)=0$ and cannot
also satisfy the normalization $f(1)=1$.
Furthermore, a rule $f(wk)=\chi(d)$ gives
$f(w\operatorname{diag}(v,v^{-1}))=\chi(v^{-1})$ for units $v$.
If $\chi|_{\mathbb Z_p^\times}$ is nontrivial, this is a character
transformation rule rather than right $K_m$-invariance. These two
different choices of local vector must be distinguished before
applying an invariant-line argument.
\end{remark}


\section*{Acknowledgements}
The author thanks T. Yang for suggestions concerning the original draft.

\end{document}